# A. A. Cournot

## Exposition of the Theory of Chances and Probabilities

> The more difficult it seems to establish
> logically what is variable and
> obeys chance, the more delightful is
> the science that determines the results
>
> Huygens[1]

## Contents

**Note by Translator**

Antoine Augustin Cournot (1801 – 1877) was a mathematician, philosopher, economist and educator (Feller 1961; *Etudes* 1978). Here, I discuss his contribution of 1843 reprinted in 1984 (Paris, Libraire J. Vrin) as vol. 1 of his *Oeuvres Complètes* complete with an Introduction and Commentary by B. Bru.

In this Introduction, Bru remarked that in 1828 and 1829 Cournot had published two notes on the *calculus of chances and combinations*. In 1834, Cournot translated J. Herschel's astronomical treatise and appended a discussion of cometary orbits which constitutes here a large part of Chapter 12. A long paper on the application of the theory of probability to judicial statistics followed in 1838 and was largely reprinted here, in Chapters 15 and 16.

Bru's commentary certainly demanded great efforts; by his permission, I quoted some of them adding his initials (B. B.). Regrettably, many of his comments are too short and for the same reason some of his references are not readily understandable.

He repeatedly mentions Condorcet, Lacroix, and certainly Laplace and Poisson as the main authors from whom Cournot had issued and he also notes (see p. 318, comment to p. 108, line 12) that Cournot had seldom indicated them.

Cournot reprinted Kramp's table of the normal distribution appended to his book (1799); it is not included in the translation.

There exist translations of Cournot's book into German (1849) by C. H. Schnuse, who also translated, in 1841, Poisson (1837), and into Russian (Moscow, 1970) by N. S. Chetverikov, the closest student of Chuprov.

A few of Cournot's terms ought to be explained. *Element* is parameter; *philosophical criticism* apparently means philosophical discussion; a *commensurable* number is rational. And the *same randomness* is the same random variable. See explanation of Bernoulli *theorems* in § 30. And, finally, I have introduced the then still unknown notation $C_m^n$ and *n*!

The Preface begins by an explanation: Cournot wished to *make assessable* his subject *to those unacquainted with the higher chapters of mathematics*. In § 123 he even repeated the formulation of the Pythagorean proposition (excluded in the translation) and in § 69 he says that he adduced a table (Kramp's table) of a certain function (of the exponential function of a negative square) but he did not provide its analytical expression. I doubt that that was good enough, but then, his Chapter 12 was certainly beyond the reach of ordinary readers.

As Bru noted in his Introduction, the book *was not understood in its entirety either by mathematicians or other scientists* [!], *frequently quoted* […] *but rarely read*.

Cournot's sentences are long-winded, up to 12 and 13 lines (§§ 117 and 240/4). In many cases, perhaps copying Poisson (1837), he connected the parts of complex sentences by semicolons rather than words which is not easy to understand. Demonstrative pronouns are often lacking, but unnecessary repetitions are plentiful.

**[1]** Cournot was obviously ignorant of precise observations (measurements) and he ignored Gauss. His Chapter 11 is therefore barely useful.

**[2]** Cournot (§ 145) considered himself a pioneer in applying statistics to astronomy, but he forgot to mention William Herschel, and he certainly had no means for studying the starry heaven.

**[3]** He did not study the application of statistics to meteorology although Humboldt, in 1817, had introduced *isotherms*, cf. Cournot's definition of the aims of statistics in § 103!

**[4]** While discussing statistics of population (Chapter 13), Cournot had not mentioned Daniel Bernoulli's classical study of prevention of smallpox published in 1766, and Gavarret's contribution (1840) escaped his notice. For many decades, in spite of the work of Graunt, Süssmilch and Daniel B., later statisticians had avoided medical statistics.

**[5]** Many elementary calculations (in §§ 13. 70, 165 – 167, 170, 182, 203, 204) are wrong which had not, however, affected Cournot's general conclusions.

**[6]** Cournot's description of the Bayes rule (§ 88) is superficial: he did not notice that Bayes had treated an unknown constant as a random variable. True, he (§ 89) remarked that without prior information that rule leads to a subjective result, and he (§ 95) attempted to prove that with a large number of observations that result becomes objective. Cf. Note 9 to Chapter 8.

His treatment of the Petersburg game is interesting, but he failed to refer to Condorcet (1784, p. 714) who had remarked that the possibly infinite game nevertheless only provided one trial so that many such games should have been discussed. Freudenthal (1951) expressed the same opinion and provided pertinent recommendations.

Poisson's law of large numbers is ignored; during Poisson's lifetime, Cournot (1838) had, however, at least twice mentioned it. It is generally known that the reason of that about-face was Bienaymé's attitude.

Cournot's description of tontines (§ 52) was completely wrong; see also Note 17 to Chapter 14.

Poisson (see Preface) indicated that Cournot *discern*[ed] the difference between chances and probabilities. However, in § 12, see also § 240/3, the latter stated that probability was the ratio of the pertinent chances. And he almost indifferently applied the terms *theory*, or *doctrine*, *of probability*, and *of chances*.

Now, however, I turn to other points whose positive aspect much prevail over the negative sude.

**[1] Probability.** Cournot's subjective philosophical probabilities (§§ 43, 233 and 240/8) can be related to expert opinions whose study undoubtedly belongs to mathematical statistics. Laplace (1812, Chapter 2) had introduced them, noted their possible application to decisions of tribunals and elections, but did not introduce any special term.

In § 18 Cournot offered a definition of probability covering both the discrete and continuous cases (i. e., and geometric probability). He appropriately introduced the ratio of the extents (étendue) of the pertinent chances; nowadays, we would have said of the measures. In § 45 he objected to Laplace's belief that probability is relative in part to our ignorance and in part to our knowledge.

Poisson's letter to Cournot (see Preface) indicates that the latter *discern*[ed] the […] *difference between the words <u>chance</u> and <u>probability</u>* […]. However, in § 12, see also § 240/3, Cournot called probability the ratio of the pertinent chances and thus contradicted Poisson's inference. And, see above, he almost indifferently applied the terms *theory* or *doctrine of probability*; or *of chances*.

**[2] Probability density.** In 1709 Niklaus Bernoulli introduced a continuous distribution (and its density) and many later scholars, including Laplace and Gauss, applied such distributions and densities. Cournot (§§ 64 – 65, see also § 31) followed suit. His term was curve of probability; it *appropriately represented* […] *the law of probabilities of different values of a variable magnitude*. Indeed, in § 73, although after Poisson (1837, § 53), Cournot introduced a *grandeur* fortuitously taking a series of distinct values. True, in 1756 and 1757 Simpson, in an error-theoretic context, had effectively applied such variables.

Cournot (§ 73) also described the determination of the density of a function *u of a magnitude x which takes a series of various fortuitous values*. Actually, he (§ 74) considered the case of a function of two independent variables and a linear function of many variables. Supposing that $u = |x - y|$, he concluded that for $0 \leq a \leq 1$ $P(u \geq a) = (1 - a)^2$. The probability of the contrary event would have described the once-popular *encounter problem* (Whitworth 1886 and possibly 1867; Laurent P. H. (1873, pp. 67 – 69): two persons are to meet during a specified time interval but their arrivals are independent and occur at random and the first to arrive only waits a specified time.

Before Cournot Bessel (1838, §§ 1 and 2) determined the densities of two functions of a continuous and uniformly distributed variable and Laplace solved such problems even earlier.

In § 81 Cournot studied a mixture of densities. Let $n_1$, $n_2$, … observations have densities $f_1(x)$, $f_2(x)$, …, then the mixture of those observations will have a density equal to the weighted arithmetic mean of $f_1(x)$, $f_2(x)$, …

**[3] Median.** It was Cournot (§ 34) who introduced this important parameter.

**[4] Randomness.** Cournot (§ 40) defined a random event as an intersection of (two) independent chains of other events and thus revived an ancient idea (Aristotle). In § 45 he mentioned the *mathematical theory of randomness* acting in the *proper field of science* and declared that randomness has *a notable role* […] *in governing the world*. Lacking was the dialectical link between randomness and necessity which Kant (1781/1911, p. 508) had clearly indicated: *Randomness in a single case nevertheless obeys a rule in a totality*.

Cournot (§ 43) connected randomness with physical impossibility, a very important notion, as he stated. It is physically impossible for a right circular cone to stand on its apex, as he remarked, but regrettably did not mention randomness. Indeed, here (as also when two chains of events are intersecting) a small cause leads to a considerable effect which is Poincaré's main and generally known explanation of randomness.

Physical impossibility is contrary to moral certainty which Cournot did not mention but which Descartes introduced in 1644 and Huygens mentioned in a letter of 1673 (Sheynin 1977, pp. 204 and 251), and Jakob Bernoulli recommended for application in law courts. And in 1693 Leibniz (Couturat 1901, p. 232) stated that there existed three degrees *of security in judgements: logical certainty, physical certainty* (*which is only logical probability*) *and physical probability*.

In § 42 Cournot indicated another aspect of randomness: *For properly understanding randomness, we should only attach to it the [...] idea of independence or absence of solidarity between different series of facts or causes*. This is an interesting idea. According to one of the modern approaches to identifying a random numerical sequence, it should only have a small number of regularities. Note also that solidarity (a notion which Cournot repeatedly applies) had been known to astronomers from the antiquity. Thus, refraction (see § 230) is a common cause altering the zenith distances of all the stars, and horizontal refraction caused by meteorological factors became known in geodetic operations. In general, any observations or measurements are fraught with solidarity, − with systematic errors.

Cournot returned to the notion of randomness in his later contributions. Thus, he (1851/1975, § 33, Note 38) recalled Lambert's forgotten attempt at formalizing *randomness* of the digits of irrational numbers by an intuitive notion about normal numbers.

Cournot appropriately mentioned Poisson (and had previously dedicated a contribution (1841) to his memory). Still, he had not hesitated to criticize Poisson (§§ 61, 93, 149 Note, 225 and 237). In § 93 his criticism seems unwarranted: Poisson, whom he had not directly mentioned, did not provide the non-existing statistical data necessary for solving a problem about the births of boys and girls. Cournot thought that this circumstance was *unbecoming of eminent geometers*. In § 225 he criticized useless mathematical considerations, for example concerning the study of facts testified by a chain of witnesses, cf. Poisson (1837, § 39).

I take this opportunity to note Mises' exaggerated opinion (1928/1930, p. 243) about Poisson (1837): it is *one of the most remarkable books in the history of the development of mathematical theories*.

**[5] Jurisprudence.** Following Poisson, but several times disagreeing with him, Cournot applied stochastic reasoning to verdicts and decisions of judges and jurymen in law courts. Unlike Poisson, he did not need to introduce a preliminary probability of guilt of an accused and he attempted to study the dependence between the voters. He concluded that cases should be separated into categories so as to ensure a useful analysis of data, to choose a certain category for further study, but at the same time (§§ 111 − 114) warned about the possible ensuing pitfalls: another classification could have suggested another category for additional investigation. When a large number of observations was available, the number of categories should be

increased, Cournot remarked (§ 115). Quetelet (1846, p. 278), however, believed that too many subdivisions of the data is a *luxe de chiffres*, a kind of *charlatanisme scientifique*. Finally, in connection with his study of classifying the data, Cournot introduced a pattern of stratified sampling (Stigler 1986, pp. 196 – 197).

[6] **Statistics.** I repeat that Cournot omitted applications of statistics to meteorology and medical statistics although emphasized its possible use in astronomy (and statistically studied planets and comets) and, apparently, chemistry (Note 9 to Chapter 11). According to his definition of statistics (§ 103), it collects and coordinates facts, it should *appreciably* exclude anomalies of chance and discover regular causes acting together with randomness. It (§ 105) should have *its theory, its rules and principles*, should *penetrate* into *the essence of things* (§ 106).

Cournot himself had not formulated any statistical rules or principles and did not mention statistics at all in his **Summary** (§ 240). *Penetration* etc had not then been generally accepted. Fourier (1821, pp. iv – v) stated that *the spirit of dissertations and conjectures is in general opposed to the veritable progress of statistics* and the just established London (now, Royal) Statistical Society declared that statistics did not discuss causes or effects (Anonymous 1839). This was the viewpoint of Staatswissenschaft, for many decades the rival of statistics as understood today.

On the other hand, opposite opinions had also been formulated: *Statistics should investigate not only why, but even the why of the why, to explain the present state of a nation by its past* (Gatterer 1775, p. 15). Quetelet, in spite of his carelessness and even happy-go-lucky attitude (Sheynin 1986), had been advocating *penetration* into the essence of phenomena. Thus, he (1869, t. 1, p. 419) recommended to study (no doubt, statistically) the changes brought about by the construction of telegraph lines and railroads. And Cauchy (1845/1896, p. 242) thought that statistics offered a means for judging doctrines and institutions.

Bru (Introduction) stated that Cournot had left *an incomparable testimony about the European* [statistical] *thought of the first half of the 19th century*. However, those thoughts had not been united in a single school (see above) and, once more, Quetelet comes to mind.

Chuprov several times mentioned Cournot, *one of the most profound thinkers of the 19th century* (1909/1959, p. 30), *the real founder of the modern philosophy of statistics* (1925/1926, p. 227). The first statement (although likely representing all other achievements of his hero as well) is an obvious exaggeration. The latter opinion is difficult to evaluate. Kruskal (1978, p. 1082) called statistics a neighbour of philosophy, *a part of philosophy of science*. For a philosopher, statistics is a method of stochastic reasoning, partly inductive, and partly deductive. In any case, Cournot's contribution was indeed philosophical. And still, Kries (1886) had barely noticed him and, much worse, Lexis (1879), who originated the Continental studies of the stability of statistical series, did not mention him (or Poisson).

Chuprov's viewpoint is understandable since almost all treatises on statistics at least until the beginning of the 20[th] century had been completely unphilosophical. Those contributions included a course of lectures written by Chuprov's own father, an eminent nonmathematical statistician A. A. Chuprov, and first published in 1886.

I am citing my website www.sheynin.de which is being diligently copied by Google: Oscar Sheynin, Home.

## Bibliography

**Anonymous** (1839), Introduction. *J. Stat. Soc. London*, vol. 1, pp. 1 – 5.
**Bessel F. W.** (1838), Untersuchung über die Wahrscheinlichkeit der Beobachtungsfehler. *Abhandlungen*, Bd. 2. Leipzig, 1876, pp. 372 – 391.
**Cauchy A. L.** (1845), Sur les secours que les sciences du calcul peuvent fournir aux sciences physiques ou même aux sciences morales. *Œuvr. Compl.*, sér. 1, t. 9. Paris, 1896, pp. 240 – 252.
**Chuprov A. A.** (1909), *Ocherki po Teorii Statistiki* (Essays on the Theory of Statistics). Moscow, 1959.
--- (1925), *Grundbegriffe und Grundprobleme der Korrelationstheorie*. Leipzig – Berlin. Russian version: Moscow, 1926 and 1960.
**Condorcet M. J. A. N.** (1784), Sur le calcul des probabilités. *Hist. Acad. Roy. Sci. Paris 1781 avec Mém. Math. et Phys. pour la même année*, pp. 707 – 728.
**Cournot A. A.** (1838), Sur l'applications du calcul des chances à la statistique judiciaire. *J. Math. Pures et Appl.*, t. 3, pp. 257 – 334.
--- (1841), *Traité élémentaire de la théorie des fonctions et du calcul infinitésimal*. *Œuvr. Compl.*, t. 6/1. Paris, 1984.
--- (1849), *Grundlehren der Wahrscheinlichkeitsrechnung*. Braunschweig.
--- (1851), *Essai sur les fondements de nos connaissances … Œuvr. Compl.*, t. 2. Paris, 1975. English translation: New York, 1956.
**Couturat L.** (1901), *La logique de Leibniz*. Paris.
*Etudes* (1978), *A. A. Cournot. Etudes pour le centenaire de sa mort*. Paris.
**Feller J.** (1961). Cournot. *Dict. Scient. Biogr.*, vol. 9, p. 983.
**Fourier J. B. J.** (1821 – 1829, 1821), *Recherches statistiques sur la ville de Paris* … tt. 1 – 4. Paris.
**Freudenthal H.** (1951), Das Peterburger Problem in Hinblick auf Grenzwertsätze der Wahrscheinlichkeitsrechnung. *Math. Nachr.*, Bd. 4, pp. 184 – 192.
**Gatterer J. C.** (1775), *Ideal einer allgemeinen Weltstatistik*. Göttingen.
**Gavarret J.** (1840), *Principes généraux de statistique médicale*. Paris.
**Kant I.** (1781), *Kritik der reinen Vernunft*. Werke, Bd. 3. Berlin, 1911. *Critique of Pure Reason*. Kindl publishers, 2010.
**Kries J.** (1886), *Die Prinzipien der Wahrscheinlichkeitsrechnung*. Tübingen, 1927.
**Kramp Chr.** (1799), *Analyse des réfractions astronomiques*. Leipzig – Strasbourg.
**Kruskal W. H.** (1978), Statistics: the field. In Kruskal W. H., Tanur J. M., Editors, *Intern. Enc. of Statistics*, vols 1 – 2. New York, pp. 1071 – 1093.
**Laplace P. S.** (1812), *Théorie analytique des probabilités. Œuvr. Compl.*, t. 7. Paris, 1886.
**Laurent P. H.** (1873), *Traité du calcul des probabilités*. Paris.
**Lexis W.** (1879), Über die Theorie der Stabilität statistischer Reihen. *Jahrbücher f. Nationalökonomie u. Statistik*, Bd. 32, pp. 60 – 98. Reprinted in author's *Abhandlungen* … Jena, 1903, pp. 170 – 212.
**Mises R.** (1928), *Wahrscheinlichkeit, Statistik und Wahrheit*. Wien. Russian translation: Moscow, 1930. *Probability, Statistics and Truth*. New York, 1981.
**Poisson S.-D.** (1837), *Recherches sur la probabilité des jugements* … Paris, 2003. English translation: www.sheynin.de downloadable file 53.
**Quetelet A.** (1846), *Lettres sur la théorie des probabilités*. Bruxelles.
--- (1869), *Physique sociale*, tt. 1 – 2. Bruxelles.

**Sheynin O.** (1977), Early history of the theory of probability. *Arch. Hist. Ex. Sci.,* vol. 17, pp. 201 – 259.

--- (1986), Quetelet as a statistician. Ibidem, vol. 36, pp. 281 – 325.

--- (2001), Social statistics and probability theory in the 19$^{th}$ century. *Hist. Scientiarum*, vol. 11, pp. 86 – 111.

**Stigler S. M.** (1986), *History of Statistics*. Cambridge (Mass.) − London.

**Whitworth W. A.** (1959), *Choice and Chance*. New York, reprint of the edition of 1901. First edition, 1867, one of the later editions, 1886.

# Preface

Here, I am setting myself two goals. First, I aim to make accessible the rules of the calculus of probability to those, unacquainted with the higher chapters of mathematics. Without that, it is impossible to conceive clearly either the precision of measurements obtained in the sciences of observation, or the values of the numbers provided by statistics, or the conditions of success of many commercial enterprises. And, second, I wished to correct the mistakes, to eliminate the ambiguities and dissipate obscurities from which, as it seems to me, the works of the most able geometers studying that delicate subject are not at all free. Since mistakes and obscurities concern the principles of the calculations rather than purely mathematical deductions, I thought that both these goals are compatible, so that instead of writing a book only for geometers I will seize the opportunity of inserting remarks useful for those attracted by that theory even if the exposition will be purely speculative.

And so, I attempted to ensure that the reading of my book will not require any other knowledge except elementary algebra, or even, strictly speaking, algebraic notation. Otherwise, I would have been compelled to replace it by verbiage at the expense of conciseness and clarity. I also wished to indicate the results of calculations and, if possible, to elucidate their meaning without entering into technical details of the pertinent proofs. Explanations containing necessary symbols of the infinitesimal calculus are placed in notes, but even in such cases I had often indicated rather than demonstrated the results.

The calculus of probability is only really important if applied to sufficiently large numbers, and for ensuring practicable results we therefore have to use approximations. I invariably had to make use of such formulas and, consequently, appended a table for applying those formulas to all the provided numerical examples without needing to know anything except ordinary arithmetic. To determine exactly the approximation furnished by those formulas and the conditions under which they can be safely used is an extremely difficult analytical problem. It is not yet solved in any complete way and I did not allow myself to touch it at all.

In the theory of probability, there occurs something similar to the mathematical theory of heat[2]. If a body is somehow heated and then subjected to the action of regular and constant sources of heat or cold, the temperature at each of its points gradually approaches a level called *final* so that all traces of the initial irregularities disappear. However, before reaching that final condition (which in the strict mathematical sense would have required infinite time) the temperature at each point passes through a certain state called *penultimate*. While it lasts, the law of the variation of temperature without sensible error obeys a regular and simple mathematical expression.

Just the same, the objects of the theory of probability are certain numerical relations which take constant and quite determined values when the number of trials made on the same randomness indefinitely increases. And when the number of trials is yet finite, those values are the closer to their final stage, and they oscillate between limits which contract the more, the larger is that number. The mathematical relation between the amplitude of the oscillations and the number of the trials is greatly simplified when that latter becomes appropriately large. By analogy, it is then possible to call the situation the *penultimate stage* obeying the laws which are the most important in the mathematical theory of chances and represent the main object of this work. In general, it is not necessary to carry out a very large number of trials for achieving the transition to the penultimate stage. I insert many examples proving that, but I did not aim at establishing with mathematical precision the conditions for the appearance of that stage. Moreover, that would have been more interesting speculatively than practically.

To tell the truth, I attach most importance to that part of my work which aims at explaining the philosophical value of the notions of chance, randomness, probability as well as the real sense of the results of the calculations to which we are led by the development of those fundamental concepts. I provided repeated explanations of independence and solidarity of causes, of the double sense of the word *probability* as a certain measure of our knowledge and of the possibility of things irrespective of our knowledge about them. These explanations seem proper for resolving the difficulties which until now led eminent authors[3] to suspect the entire theory of mathematical probability.

Here, the reader will find definitions and ideas which I believe to be new or at least had not been properly formulated. They brought me to consider the doctrine of posterior probabilities and most of the applications connected with them in way quite different from those adopted by really celebrated authors[4]. I will hardly yield to the authors' illusion by suggesting that my ideas at least deserve to be discussed; that they will be able to interest philosophers as well as geometers; that in attracting the attention of those engaged in this field of human knowledge they can foster its later perfection.

The term *probability* had been the source of so many ambiguities[5] that I had at first intended to abandon it completely and to apply in appropriate cases either *chance* or *possibility*[6]. Then, however, I found it more inconvenient to reject a term so rooted in geometry. I also believe, and attempted to prove that the word *probability* has other meanings differing from, but sufficiently close to those applied in calculations of the geometers so that they should not be completely isolated from the latter in a philosophical exposition.

Circumstances permitting, some day I will try to develop the ideas which are only indicated here in the last chapter[7]. If I insist on my ideas still more, I fear that reproach for admixing too much metaphysics with geometry will follow. *Est modus in rebus* [There is measure in everything]. However, I have necessarily borrowed without scruple two epithets, *objective* and *subjective*, for radically separating the two meanings of the term *probability* as applied in the calculus. Here, however, I am following the example of Jakob Bernoulli[8].

I hope that the readers will find here a selection of sufficiently differing applications for ensuring them a fair idea about the usefulness of the theory of chances and that those who search statistics for something in addition to raw results[9] will be able to try out their own possibilities on a path to new applications. I have considered in detail two curious problems which I had already treated before. One of them considers the distribution of the cometary orbits in space and the other has to do with the theory of probability of judgements and its application to statistical documents published in France by the Administration [Ministry?] of justice. Competent readers will pronounce their opinion about the value of my solution of this theoretical problem and compare it with those provided by other authors, notably by Poisson in his great work of 1837.

As to the application of judicial statistics, I have for the first time made use of the more precise data appearing in the *Comptes généraux* of the Administration of criminal justice as a result of changes in jury panels introduced by the law of 9 Sept. 1835.

I am not concluding this Preface, perhaps already protracted, without expressing gratitude to my excellent friend Mr. Bienaymé, inspector general of finances, whose work in statistics and probability is well known to geometers. For a long time, not knowing each other, we have been occupied with the same objects of study, then becoming closer because of a singular conformity of ideas and inclinations. During the printing of my book, he gladly helped me and was obliging to such an extent that even re-read the proofs and calculated anew a part of the numerical computations. He was all the more unselfish since long ago he came to the theory of posterior probabilities by issuing from considerations quite different from those that guided me. As it seems to me, this theory largely reappears here in Chapter 8. It will become possible to access our similarities and dissimilarities if he decides to publish his own researches[10]. At present, however, I hasten to acknowledge that the originality is all his and that his ability deftly to apply analysis undoubtedly enabled him to discover much of what had escaped my attention.

I ought to be excused for publishing here a letter which Poisson wrote me in response to my own in which I had sketched the subject of this book. I do it much less for maintaining my priority over some ideas than as a testimony of friendship the memory of which will invariably be precious to me, and for testifying about some opinions of the celebrated man.

<div align="right">Paris, 26 Jan 1836</div>

*Sir, With great pleasure will I read the work on the <u>Doctrine of chances</u> which you propose to publish. What I am now completing will not hinder you at all, and I am leaving enough space for a more comprehensive book. I discern the same difference between the words <u>chance</u> and <u>probability</u> as you, and strongly insist on it. As to your approach to the main problem, the probability of judgements, I will compare it with my own after reviewing and definitively accomplishing that part of my work. I only have to finish that and to copy the entire text before being able to begin printing. There are some problems whose solutions I will include in one of the last chapters provided that I complete them at least to my own satisfaction.*

*Finally, you will find in that work some metaphysical considerations and see that I am not at all denying that branch of human knowledge. I am writing this letter during a sitting of the Council[11], and my turn to speak is about to begin, so I am unable to continue.*

*Be assured, Sir, of my attachment and undivided devotion.*

*Poisson*

## Notes

**1.** Following a nasty tradition, Cournot did not indicate the source of that statement. Bru found it: Huygens' covering letter to F. van Schooten accompanying his contribution (1757) inserted in van Schooten's book. Cournot quoted the Latin text published in 1760 whereas my translation is from the French version.

**2.** Cournot's analogy between probability theory and the mathematical theory of heat is based on Fourier (1819; 1822). [B. B.]

**3.** D'Alembert, Poinsot (1836), Comte. [B. B.] Poinsot denied the application of the theory of probability beyond natural sciences. O. S.

**4.** Bayes, Condorcet, Laplace. [B. B.]

**5.** Pascal. [B. B.]

**6.** De Moivre applied the term *chance* in the title of his contribution (1718).

**7.** See Cournot (1851). N. S. Chetverikov, translator of Cournot into Russian.

**8.** See the first lines of the *Ars Conjectandi*. [B. B.]

**9.** Guerry (1864), Moreau de Jonnès (1847). [B. B.]
  The former described in detail the early history of statistics. O. S.

**10.** I. J. Bienaymé (1796 – 1878), a greatest statistician of the 19th century, apparently never decided to publish his researches largely only known by short extracts published by the Société Philomatique and the Paris Academy of Sciences. All traces of his course of probability (Sorbonne 1848) are lost. See Heyde & Seneta (1977). [B. B.]

**11.** Conseil Royal de l'Instruction Publique. [B. B.]

## Bibliography

**Cournot A. A.** (1851), *Essai sur les fondements de nos connaissances* etc. Paris, 1975.

**De Moivre A.** (1718), *Doctrine of Chances*. London, 1738, 1756 ; New York, 1967.

**Fourier J. B. J.** (1819), Sur la théorie analytique des assurances. *Annales de chimie et de physique*, t. 10, pp. 177 – 189.

--- (1822), *Théorie analytique de la chaleur*. Paris.

**Guerry A. M.** (1864), *Statistique morale de l'Angleterre comparée à la statistique morale de la France*. Paris.

**Heyde C. C., Seneta E.** (1977), *I. J. Bienaymé*. New York.

**Huygens C.** (1757), De calcul dans les jeux de hazard. *Oeuvr. Compl.*, t. 14. La Haye, 1920, pp. 49 – 91. In French and Dutch.

**Moreau de Jonnès A.** (1847), *Eléments de statistique*. Paris.

**Poinsot L.** (1836), Discussion of Poisson's note in *C. r. Acad. Sci. Paris*, t. 2, pp. 377 – 382.

## Chapter 1. Combinations and Order

**1.** Among the abstract ideas not arbitrarily created by the human mind, but suggested to it by the very nature of things, that of *combination* is one of the most general and simple. After considering individual objects in isolation, we are led to understand that, according to their nature, those objects combine, or unite, one to another one, two to two other ones, three to three, etc, and form certain systems or complex objects which in turn can combine with each other and form other groups or more complicated systems, etc.

The theory of combinations which the Germans called *syntactics*[1] is an abstract and purely logical science like the science of numbers and geometry. It is intimately connected with all branches of mathematics, notably algebra so that perfection, or, as it is called, elegance of algebraic formulas, achieved by aptly chosen notation, provides the greatest obviousness to the laws of combination.

In essence, each scientific synthesis successfully combines certain principles or primordial facts. From that viewpoint, logic, general grammar, chemistry, just as algebra, issue from combinatorics. Notably logic, in the theory of syllogisms, offers us a curious example of combinatorial synthesis. Logicians[2] invariably require [of their followers] to form carefully all the possible combinations. Omitting only one of them is sufficient for the reasoning to become illegitimate. Actually, little did the logicians think about tracing the rules for securely compiling a complete record of the combinations involved. Not unreasonably, they thought that in simple cases such rules will be useless and impracticable otherwise.

A systematic and regular approach is unknown to most people or its application had been too slow for satisfying practice. Therefore, the art of forming combinations embracing all at once a more or less large number [of elements] depends on the aptitude and education of individuals. It is usual to understand as calculative minds those who more eminently possess that power of combination and apply it to various objects although most often calculations proper do not involve numerical reckoning. The mechanician, the geometer, the tactician, the chess-player is distinguished, each in his field, by his ability to form and classify combinations, and the common opinion admits an affinity between these aptitudes however different are the pertinent objects.

**2.** It is easy to understand that the rules ensuring the formation of all the possible combinations should implicitly include rules for calculating their number without having to form them one by one or to pay attention to each in particular. To provide a simplest example[3]: Suppose that a chemist is trying out all combinations of *m* acids with *n* alkalis. He is assured to omit no combination when denoting each acid and each alkali by signs or numbers in succession and combining acid No. 1 with each alkali, then repeating this procedure with acid No. 2, etc. However, it is seen at once that the total number of all the necessary trials or of all the abstractly possible combinations is *mn*. And if *m* and *n* are considerable numbers, 100 or 200, say, the product *mn* will be so large that the formation or even the enumeration of all the possible combinations becomes a long-winded and tiresome operation. However, nothing is simpler and quicker than multiplying *m* by *n*.

We will see that the solution of very important and very curious problems about which this book is destined to provide an idea, essentially depends on the possibility of indicating the number of combinations of a certain type of objects or at least the ratio of such numbers for two different types of objects. At the same time, it would have been barely important or often impossible to consider these combinations one by one.

Therefore, in a more special sense, the theory of combinations is understood as a science aiming to assign the number of combinations of a given kind. Synthesis, thus reduced to determining their numbers, is naturally included in those branches of mathematics with which, as I said, it is intimately connected.

**3.** Suppose that *m* objects *a, b, c, …, k, l* of some kind should be combined. For exhausting all the possible combinations of two from two (*binary* combinations) we can combine the object, or element *a* with all the ($m-1$) other elements *b, c, …, k, l*; then combine the element *b* with all the ($m-1$) other elements *a, c, …, k, l* etc. We will thus have $m(m-1)$ combinations; however, each combination, for example *ab*, can evidently be obtained twice, as *ab* and *ba* so that the number of different binary combinations will be $m(m-1)/2$.

Combinations of three from three, or *tertiary* combinations, will be exhausted if each binary combination *ab* is successively joined with each of the ($m-2$) elements *c, …, k, l*. We will thus obtain the same combination *abc* three times: $ab + c$, $ac + b$ and $bc + a$, and the number of the tertiary combinations will be $m(m-1)(m-2)/2·3$. It is not necessary to continue; by an evident induction we conclude that the number of different combinations of *m* elements from *n* is $C_m^n$. […]

Connect the pure idea of combination with certain relations of order or situation so that combination *ba* is now regarded as differing from *ab*. Then *m* elements will provide $m(m-1)$ binary combinations, $m(m-1)(m-2)$ tertiary combinations, and, in general,

$$m(m-1)(m-2) \ldots (m-n+1) \qquad (3.1)$$

combinations (arrangements, selections) of order *n* from *m* elements. Various authors have not at all agreed about how to name those combinations in which not only the included elements, but also their order is taken into account. Most convenient, as it seems to me, is to call them *ordered combinations* without involving very particular considerations. On the contrary, we will call those, in which the order of the elements is not essential, *absolute combinations*, or simply *combinations*.

**4.** The number of ordered combinations is expressed by formula (3.1), and that of absolute combinations is $C_m^n$. It follows that the denominator of $C_m^n$, or $n!$, expresses the number of all changes of the order, or *permutations* of *n* elements. This is easy to prove by direct reasoning which will also show that the theory of permutations and of order in general is actually the same as the theory of combinations presented from a different point of view.

Suppose that in some order in space or time or even (if it can be understood) independent from space and time there are *m* determined places denoted by numbers 1, 2, …, *m*, say. Two elements, *a* and *b*, are put there. We will exhaust all the possible *arrangements* (selections)[4] if *a* is put on place 1, and *b*, on each of the other places, 2, 3, …, *m* in succession; then *a* is put on place 2, and *b*, on each of the other places 1, 3, …, *m* in succession, etc. Therefore, the number of arrangements of order 2 from *m* elements is $m(m-1)$. If the system includes a third element *c*, each of the previous arrangements, for example *ab*, can be supplemented by *c* taking each of the $(m-2)$ places 2, 3, …, *m*. The number of the arrangements will now be $m(m-1)(m-2)$. And in general when $m = n$, the product (3.1) will become $n!$.

We can suppose that the ordered elements are letters written one after the other and form a series of the kind which geometers call *linear*[5]. This assumption, however, is only intended to ease the idea of applying a sign (? - O.S.)[i]. Indeed, there is nothing essential here except the notions of elements and order understood in the most general sense. We can imagine *n* points somehow distributed in space or *n* spheres of various radiuses whose centres successively coincide with each of those points. The product $n!$ indicates also the number of the arrangements or the various configurations offered by the system of those spheres. We can also imagine *n* people occupying a social hierarchy having different functions and able to exchange their positions. Then the number $n!$ will, as previously, indicate how many ways there are for modifying that hierarchic system.

The idea of *absolute order* can be restricted by introducing particular conditions and thus obtaining a smaller number of different arrangements. For example, if *n* elements are arranged in a *circular* series, or a periodic series indefinitely increasing or decreasing, without accounting for the definite places they occupy and only taking into consideration the order in which they follow, the number of arrangements will be $(n-1)!$. And, without distinguishing either left or right or increasing or decreasing order, that number should be halved.

**5.** When taking *n* elements out of *m* so as to form a combination of *n* elements from *n*, the total group will be separated in two, one of these consisting of *n* elements, the other, of $(m-n)$ elements. Then $C_m^n$ will indicate in how many ways this separation can be done. It follows that that quotient $C_m^n$ ought also to express the number of different combinations characterizing the second partial group so that its value will not change if *n* is replaced by $(m-n)$ […] and $C_{m+n}^m = C_{m+n}^n$.

Expression $C_{m+n}^m$ can be directly obtained by considerations applied in the theory of order. Suppose that a system of $(m + n)$ places is separated in two, A and B, of *m* and *n* places. The number of different arrangements for the initial system is $(m + n)!$. But then, if only considering the distribution of the elements among the two partial groups, the arrangements will only differ by permutations of order in A or in B. Therefore, the product $(m + n)!$ should be divided first by $m!$ then by $n!$.

The same problem of order can be presented otherwise. Suppose that A has $m$ letters and B, $n$ letters. Then $C_{m+n}^{m}$ will express the number of different ways for them to follow each other. Indeed, when considering that letters are individually different, the number of permutations will be $(m + n)!$. But in the contrary case the arrangements which only differ by their order in the partial groups should be regarded identical. If letters A and B denote respectively $m$ and $n$ events of the same nature following each other successively, the quotient $C_{m+n}^{m}$ will express the number of ways of different successions.

An analogy will indicate well enough that, when separating a total group of $m + n + p$ elements in three parts consisting of $m, n$ and $p$ elements respectively, the quotient

$$\frac{(m+n+p)!}{m!n!p!}$$

will express the number of different ways for achieving that purpose. It will also show how many different series can be formed with $m$ letters A, $n$ letters B and $p$ letters C etc. It is needless to continue generalizing this reasoning.

Products $n!$ invariably appear in the theory of order and combinations and therefore in its applications to other branches of mathematics, and analysts3 had called them *factorials* and attempted to study their properties. The properties of the numbers $C_{m+n}^{m}$ which combine three factorials had also been thoroughly studied. They are connected with those of certain *functions* appearing in higher analysis and are at present known as the *Euler functions* since Euler essentially advanced their theory. However, the nature of this book only allows me to indicate for some readers the connection of our present subject with other abstract speculations having less immediate applications.

**6.** If the same element can be repeated in the combinations (like letters in alphabetic, or digits in numerical combinations) factorials are replaced by powers of numbers. Thus, with $m$ letters *a, b, c, ..., k, l* we can form $m^2$ binary combinations *aa, ab, ac, ..., ak, al; ba, bb, bc, ..., bk, bl* etc differing from each other either by composition or order of letters. There will be $m^3$ tertiary combinations, and generally $m^n$ combinations of $n$ letters or

$$\frac{m(m+1)...(m+n-1)}{n!} \tag{6.1}$$

if only taking into account absolute combinations, i. e. regarding combinations formed by the same differently ordered elements as identical. Conforming to the remark in § 5, instead of the expression above we can write $C_{m+n-1}^{m-1}$.

Thus, when throwing two dice with six faces each, each face of the first die combines with each face, coinciding or not, of the second die. The number of combinations is $6^2 = 36$. However, if we only consider the outcome, for example, 2 and 3, without examining on which die had each of those points appeared, the number of different throws is reduced to 6·7/1·2 = 21. For 3 dice the number of ordered combinations will increase to $6^3 = 216$ and that of absolute combinations or differing combinations will be 6·7·8/1·2·3 = 56.

**7.** The preceding clearly enough shows the connection between the theory of combinations and one of the four fundamental arithmetical operations, *multiplication*. That operation is actually characterized in that if each factor is considered to be complex, formed by adding up some *terms*, the total product will be the sum of the partial products obtained by multiplying each term of one of the factors by each term of another. Thus, when considering two natural numbers *m* and *n* as a sum of *m* and *n* unities, the product *mn* will contain the same number of partial products, i. e. as many unities as binary combinations can be formed by taking one element from the first series of *m* of them and the other element from the second series of *n* of them.

The proof of the first theorem of the number theory immediately follows: the product of two numbers does not change whichever of them is considered the multiplicand or multiplier. And it is not difficult to conclude that the product of some number of factors does not change if their order is somehow permuted. It also follows that the algebraic multiplication of $(a + b + c + …)$ by $(a' + b' + c' + …)$ contains all the binary combinations which can be formed when combining a letter from the first polynomial by another from the second.

When ordering the expansion of the product of *m* binomials $(x + a)$, $(x + b)$, $(x + c)$, …, $(x + k)$, $(x + l)$ by the decreasing powers of *x* the coefficient of $x^{m-1}$ will be the sum of *a, b, c, …, k* and *l*; of $x^{m-2}$ and $x^{m-3}$, the sums of the different products formed by the binomial and tertiary combinations of those letters, etc with the free term being the product of those *m* magnitudes. Suppose that all the magnitudes *a, b, c, …* are the same so that their sum is *ma*; the sums of their binary and tertiary combinations will be $C_m^2 a^2$ and $C_m^3 a^3$, etc, and $C_{m-n+1}^n$ will be the numerical coefficient of the term $a^n x^{m-n}$ in the expansion of $(x + a)^m$.

That expansion, fundamental for algebra, is called the *Newton binomial* after the great [physicist and] geometer who discussed it. Somewhat earlier Pascal had provided its equivalent by constructing his *arithmetical triangle*. He had not written it in an algebraic manner which deprived it of the immense advantage attached to algebraic formulation[6].

The coefficients of the binomial formula occur in many other formulas playing a considerable role in the higher parts of analysis. The cause of these analogies is easy to see: they obviously follow from the law of those coefficients being justified by a rule of combinatorial synthesis quite independently from the nature of the represented calculative procedure or the secondary idea of multiplication which can be associated in the elements of algebra with the abstract notion of combination.

**8.** If $x = a = 1$ the binomial formula

$$(x + a)^m = x^m + (m/1)ax^{m-1} + C_m^2 a^2 x^{m-2} + \ldots + a^m, \qquad (8.1)$$

after the first term on the right side is transferred to the left, provides

$$2^m - 1 = m/1 + C_m^2 + C_m^3 + \ldots + 1. \qquad (8.2)$$

Thus, $2^m - 1$ is the number of all possible combinations out of $m$ elements taken one from one, two from two, etc, and, finally, all of them together. Now let $x = 1$ and $a = -1$, transfer the same term to the left and change all the signs:

$$1 = m/1 - C_m^2 + C_m^3 - \ldots \qquad (8.3)$$

The positive part of the right side is the sum of the combinations of *odd orders*; and the negative part is the sum of the combinations of the *even orders*. When adding up the equations (8.2) and (8.3) the combinations of the even orders will disappear and the sum of the combinations of the odd powers will be $2^{m-1}$. Therefore, the sum of the combinations of the even powers will be $2^m - 1 - 2^{m-1} = 2^{m-1} - 1$. The first sum exceeds the second by a unity whether $m$ is even or odd. This result was thought to be strange and attempts were made to justify it by prior considerations, but, while treating a more general problem, I (1829) had shown that that opinion was unfounded.

**9.** I am concluding this extremely concise explanation of the most general principles of the theory of combinations by a few numerical examples. It is common knowledge that the previous Lottery of France[7] consisted of 90 numbers 5 of which were extracted at each drawing. These 90 numbers provided 90 simple extractions, 90·89/1·2, 90·89·88/3!, 90·89·88·87/4! and 90·89·88·87·86/5! combinations of 2, 3, 4 and of all 5 numbers.

In a drawing of the 5 numbers 5 simple extractions were possible as well as 5·4/1·2, 5·4·3/3!, 5·4·3·2/4! and 5!/5! = 1 combinations of 2, 3, 4 and of all 5 numbers. A gambler playing on 3 numbers had 10 winning combinations out of 117,480 etc.

In the game of piquet the operation called *dealing* means distributing 32 cards in 4 groups, 2 of them of 12 cards given respectively to each gambler and 2 other of 5 and 3 cards which constitute the *talon*. The number of combinations taking place is 32!/12!12!5!3!. Since this number is enormous, when recalling the date assigned to the invention of card games, we are assured by a very simple calculation how much time is required for the cards in that game to be distributed in all possible ways. However, in piquet the cards only differing by suit are considered identical so that the number of various distributions is considerably less.

If among all the combinations we are only interested in those in which the 4 aces are in one of the 2 groups of 12 cards, − in that, for example, which is dealt out to the first gambler, − we can find the number of such combinations by imagining that the aces are removed from the pack and that the remaining 28 cards are distributed in all possible ways in 4 groups of 8 (for the first gambler), 12, 5 and 3 cards. The required number is 28!/8!12!5!3! and its ratio to the number calculated above is 9·10·11·12/29·30·31·32 = 99/7192 = 0.01137653 which is easy to calculate without determining those two enormous numbers.

**10.** The examples provided show how rapidly increases the number of combinations even when the number of the combined things increases inconsiderably. Soon it becomes impossible not only to form or examine these combinations one after another, but even to accomplish the operations necessary for calculating the number of combinations. Imagine an assembly consisting, as our Chamber of Deputies now is, of 459 members separated by chance into 9 sections with 51 members each.

The number of possible distributions will be $x = 459!/(51!)^9$, but the calculation of this number by ordinary arithmetic is either impractical or excessively long. Tables of logarithms or, better, of sums of logarithms were specifically compiled for such aims (Degen 1824). We can therefore determine, at least with 12 decimal points, that $\lg x = 428.445 \ldots$ In our decimal number system $x$ will have 429 digits before the decimal point.

### Notes

**1.** Bru named the main workers of the German combinatorial school beginning with K. F. Hindenburg.
**2.** Arnauld & Nicole (1662, pt. 3, Chapter 19, § 4) and Jakob Bernoulli. B. B.
**3.** Cournot (1847, pp. 13 – 14) repeated this example. B. B.
**4.** In elementary algebraic treatises in which the theory of combinations is only taught from the viewpoint of its applications to algebra, it is now usual to call *arrangements* what we denoted as *ordered combinations*. This term is improper. To *arrange*, in ordinary language signifies putting things in a certain order rather than choosing or combining them. *Permutation* is an operation replacing an arrangement by another one with the arranged things remaining the same. A. A. C.
**5.** Unusual term. [B. B.]
**6.** The first traces of the theory of combinations [in the new times] are found in the correspondence of Pascal and Fermat and in Pascal's treatise on the arithmetical triangle. We should not forget that Leibniz' glorious career began in 1666 with his thesis on combinations. Traces of his first youthful speculations [on combinations] reappear in all parts of the great man's philosophical system and notably in his views about the *Caractéristique universelle*. A. A. C.

The expansion of $(a + b)^n$ for natural numbers $n$ was known before Newton, but he extended it on the case of negative and fractional values of $n$ (with the number of terms becoming infinite). On the arithmetic triangle see Edwards (1987). O. S.
**7.** The Lottery of France existed from 1758 (initially under a different name) to 1836, but was prohibited from 1793 to 1797. [B. B.]

### Bibliography

**Arnauld A., Nicole P.** (1662 French), *The Art of Thinking*. Indianapolis, 1964.
**Cournot A. A.** (1829), *Bull. sci. math*. de Ferrusac, t. 11, p. 93.
--- (1847), *De l'origine et des limites de la correspondance entre l'algèbre et la géométrie*. Paris, 1989.
**Dogen C. (K.) F.** (1824), *Tabularum ad faciliorem et breviorem probabilitatis computationem utilium enneas*. 2012.

**Edwards A. W. F.** (1987), *Pascal's Arithmetic Triangle*. Baltimore – London, 2002.

## Chapter 2. Chances and Mathematical Probability

**11.** I (§ 2) have indicated that the theory of combinations is mainly applied in cases in which a great number of combinations can be distributed according to some viewpoint in a small number of categories, so that it will only be interesting to know how many combinations out of all of them are included in one and another category.

For example, suppose that a gambler intends to play 30, say, sets of a game. The number of hypotheses or combinations taking place in the uncertain succession of losses and gains (§ 8) is $2^{30}$. It is clear, however, that after all he is only interested in the number of won or lost sets. It follows that he could include in the same category all combinations that only differ in the order in which the won and lost sets succeed each other. Evidently, in accord with formula (6.1) there will only be 31 different hypotheses. It will be otherwise if the sum at the gambler's disposal can be exhausted by consecutive losses making it impossible for him to play a fixed number of sets. The order of the succession of gains and losses will not be anymore indifferent so that the possible combinations should be distributed in a larger number of different categories.

In general, when having no cause for supposing in advance that one of the realizable combinations and hypotheses will appear rather than another one, the rules of the so-called games of chance determine a certain number of them. Those hypotheses or combinations are called *chances*, and they are naturally distributed in two categories, favourable for the gambler and leading to his winning, and unfavourable, causing losses. It is clear that he and those who share his hopes and fears are interested not in enumerating or examining all the possible chances one after another which is almost never practicable, but in being able to calculate directly how many chances are favourable and contrary.

**12.** It is easy to see that the gambler will know everything important for him if only being able to calculate the ratio of the favourable and contrary chances or, which comes to the same, of the former to the total number of chances. That proportion means that it is indifferent to the gambler if the numbers of favourable and contrary chances increase or decrease proportionally. This proposition can be regarded as one of those fundamental notions which, when being developed, are subjected to possible obscurity. Here is a physical model to which we may resort following Laplace[1] for rendering that proposition more obvious.

Suppose that an urn contains 20 white and 15 black tickets and that another urn has 40 and 30 of those. An extraction of a white and black ticket means that a gambler, respectively, wins or looses. We say that for the gambler it is indifferent whether the ticket was drawn from the first or the second urn. […] We conclude that when considering random events[2] the interest lies not at all in determining the total number of chances or the absolute numbers of favourable and contrary chances for the occurrence of such events, but only in finding out *the ratio of the number of chances favourable for that event to the total number of chances*, which does not change when those two magnitudes vary proportionally. That ratio should be named so as to dispense with the need for incessantly repeating the stated definition, and it is called *mathematical probability*, or simply *probability* of the event.

**13.** The described substitution of calculating the two magnitudes constituting a ratio by determining that ratio transforms the theory of combinations. It becomes extended which is difficult to overestimate. In Chapter 1 we saw how rapidly increase the numbers determined by the combinatorial synthesis with the increase in the number of the combined elements. These formulas soon lead to impossible calculations although the ratio of two incalculable numbers only largely depends on the first significant digits of those numbers and can still be estimated.

For example (§ 9), the total number of combinations taking place when dealing out the cards in the game of piquet and of the combinations for the 4 aces to be dealt out to a gambler, are enormous although can be calculated without too much work. However, their exact ratio, 99/192 = 0.0137653 (or the *probability* that the 4 aces will be dealt out to a gambler) can be determined without calculating either of these two numbers. Suppose that [by ordinary arithmetic] this is impossible and that we calculate these numbers by logarithms arriving at 15,928 with 11 zeros and 219, again with 11 zeros. Then the probability sought will approximately be 219/15,928 = 0.0137169 which coincides with its exact value at least to about 1/10,000.

The calculus of probability exclusively consists of largely similar methods by including mathematical probability into the theory of combinations and thus singularly multiplying its applications. However, the more complicated become the problems the greater insight and knowledge of mathematical analysis they demand. In general, the obstacles encountered in the applications of the exact sciences are engendered either by the very nature of things which make them inaccessible to our calculations, or by the duration and complication of calculations which become impractical even when understanding their theory.

It is mostly the obstacles of the second kind that occur in the calculus of probability, and for overcoming the difficulties of work analysts have to apply all the opportunities offered by perfecting the analysis. The nature of my contribution does not allow me to describe their methods, but it is possible to provide at least an idea of their goals[3].

**14.** Not only the chances or combinations monstrously multiply with the number of the combined elements, but the combined things themselves and all the more the combinations can become infinite in number and indefinite. However, when distributing them in two categories, the number of combinations in each is [still] infinite and indefinite, but maintains a finite and assignable ratio between them. Since a fundamental notion is here discussed, we would like to elucidate it by an example and make it quite clear to all our readers.

Suppose that a billiard board one meter long is divided into 10 equal parts. A billiard ball thrown randomly hits the board at a certain point. If that point belongs to part 7, we say that point 7 was hit. Throw the ball twice, then the difference between the hit points can vary from 0 to 9 inclusive. It is required to determine the probability that it will not be less than 3. The number of chances or combinations of possible points, when having no reason to suppose one of them to occur rather than another, is $10^2 = 100$ with 44 of them providing a difference less than 3. The required probability is 0.56.

Divide now the same board into 100 parts, and the same probability will become 0.497 and 0.4907 when the board is divided into 1000 parts. Following in the same way, we will find that probability equal to 0.49007, 0.490007, … We see that it incessantly tends to 0.49 and will soon only differ from that value by an extremely small fraction.

Now formulate that problem otherwise, without imagining that the board is divided into parts. Suppose that we measure the distance of each hit point from an end of the board. It is required to determine the probability that the difference between two measures thus obtained will not be less than 0.3 of the board's length. It is clear that he number of the combined *points* and, all the more, of the combinations or chances, is infinite. Indeed, each throw of the ball can lead to infinity of values of those measures. And it is also clear that, when successively dividing the board in 10, 100, 1000, … we ever closer approach the present case since we neglect centimetres, then millimetres, then …

Mathematicians know general calculating procedures for determining the limits to which certain ratios ever closer approach when the terms of those ratios vary by ever smaller degrees[4]. In the presented case, the application of those rules provided 0.49 exactly for the required probability

**15.** And so, the calculus of mathematical probability which at first presented itself as a branch of synthesis or the theory of combinations, became wider than the entire synthesis in the sense that it is applicable to cases in which either the formation of the combinations one after another, or the calculation of their unbounded number is impracticable.

The discussed transformation is of special interest in the real world where the number of combinations or chances is usually infinite since in nature almost everything varies continuously rather than by leaps. *Natura non facit saltus* (Nature does not leap), as the scholastics of old used to state. This saying should not be understood literally, but it is true in the sense that for natural phenomena continuity is the rule and leaps or discontinuities are exceptions. For combinations engendered by human hands the inverse is true.

One of the best discoveries made in the exact sciences was the determination of methods for passing from discontinuity to continuity. This is how geometers pass from considering polygons whose sides sharply change their directions to curves. And, not to abandon our subject, the ideas of chances and probabilities, formed when considering games only offering a finite number of combinations, were extended by applying them to cases of nature in which the ratios and combinations can vary infinitely.

**16.** For example, if a married couple paid a sum for assuring a pension to the surviving spouse, the problem for the insurer was to determine the probability that the *difference* between the lifetimes of husband and wife counted from the beginning of the insurance contract did not exceed a determined period. These lifetimes take an infinity of various values, so that this problem was similar to that which we provided in a purely geometric form but differing from it in one essential aspect. We were able to assume there that at each throw of the ball the distance of the hit point from the chosen end of the board took all values from 0 to 1 *m* without having a cause to suppose that one of those values was preferable to another. Here, we are not allowed to adopt the same viewpoint, to regard the realization of all the hypotheses about the lifetimes of the spouses as indifferently possible.

Each of these hypotheses has its own proper probability, its own ratio between the chances or combinations leading to its realization and the total number of chances. Those ratios vary from one hypothesis to another, and the law of that variation should be known in advance for solving by calculation the problem about the probability interesting for the insurer. The geometric problem of § 14 can be modified so that it will be more similar to that of the insurance. Suppose that the ball is twice thrown at random on a circular surface with radius of 1 *m* [cf. Poisson (1837, § 102) − B. B.] The distance of the hit point to the circumference can take all values from 0 to 1 *m*, but they correspond to the same number of unequally probable hypotheses. Indeed, the ball can hit by chance any point of the circle, and when considering two portions of the circle having the same area, there will be no reason for the ball to hit one of them rather than the other. [Cournot then derived the probability of the ball hitting a point having a given distance from the circumference.]

**17.** We provide one more example of very simple geometric conditions of the game of *franc-carreau* (Buffon 1777, § 23). A floor is paved with regular hexagonal tiles and a coin randomly falls on it. One gambler bets on the coin to rest completely within a tile, the other one, on it to fall on a joint.

Choose one of the hexagons and construct another such figure inside it with its sides parallel to the corresponding sides of the initial hexagon and the distances between such sides equal to the coin's radius. The first gambler obviously wins if the coin's centre falls inside the internal hexagon and loses if it falls between the two figures. It is seen that his gain is measured by the ratio of the areas of those hexagons.

**18.** We have defined mathematical probability as the ratio of the number of chances favourable for an event to the total number of chances. That definition presumed that all the chances can be enumerated and constitute the same number of *discrete* unities. For modifying that definition and thus rendering it applicable for an infinite number of chances with a passage from one chance to another performed without discontinuity, the numbers should be replaced by continuous magnitudes. Among the concepts of such magnitudes we most immediately imagine that of extent. Thus, we can also define mathematical probability[5] as *the ratio of the extent of chances favourable for an event to the total extent of chances*. However, the word *extent* is only generally used by assimilation although it can also be applied in its proper sense so that probability will be immediately determined by the ratio of geometric magnitudes just like in the examples above.

**19.** Following the natural inclinations of the human mind, we will now pass over from the case in which the chances can be enumerated as so many differing hypotheses to that in which the chances constitute a continuous whole whereas geometers often move in the opposite direction (? - O.S.); when an enumeration of chances, although theoretically possible, leads to impracticable calculations, they introduce a fictitious continuity and this is indeed one of the most fruitful methods of approximation for evaluating ratios of large numbers (§ 13). In essence, this artifice is the same as invariably practised in most usual circumstances. For example, instead of counting grains, they are measured by considering them a continuous mass, and the ratio of volumes of grain of the same kind should not remarkably differ from the ratio of the numbers of grain contained there.

**20.** Here, we provide some frequently applied general principles immediately following from the notion of mathematical probability. For simplifying the exposition we suppose that there is a finite number of chances so that the probabilities will be expressed by commensurable fractions and considered as continuous magnitudes.

**I.** *The probability of an event which can occur according to various unequally probable hypotheses is the sum of the probabilities of each hypothesis favourable for that event*.

To fix this idea by a simplest example, suppose that an urn contains $N$ balls, $n$ of them white, $n'$ red, $n''$ yellow, and some of other colours. The random event is the gain of a gambler occurring if a ball of any of those three colours is extracted. The probability of his gain is evidently

$(n + n' + n'')/N = n/N + n'/N + n''/N$.

**21.** It is possible to enquire not about the absolute probability of the gambler's gain but about its relative probability that the extracted ball was white rather than red or yellow. If a person bets on it, and another person bets against it, they will disregard all drawings and chances leading to the former's loss, and the required probability will be $n/(n + n' + n'')$. The absolute probabilities of extracting a white, and a red or yellow ball will be

$n/N$ and $(n' + n'')/N$

and therefore

**II.** *The relative probability of an event is the quotient of its absolute probability divided by the sum of the absolute probabilities of the events which are compared with it.*

**22.** In the case just above, a compound bet could have been put on the gambler's gain if the colour of the extracted ball was not only one of the three mentioned, but white. The probability of that event is

$n/N = [(n + n' + n'')/N][n/(n + n' + n'')]$

and the following rule can be therefore stated:

**III.** *The absolute probability of an event composed of two other events the second of which can not occur before the first one is the product of the absolute probability of the first event multiplied by the probability that the second will appear if the first one did, or by the relative probability of the second event.*

The application of this principle often facilitates the calculation of probability as seen in this very simple example borrowed from Lacroix.

Suppose that we collected in a random order the 13 cards of the same suit from a pack of 52 cards. It is required to determine the probability that the first two of them are an ace and a deuce. The probability of the ace occurring in the first place is 1/13 and 1/12 is the probability of the deuce occurring in the second place with the second event being subordinated to the first one. The required probability is $(1/13)(1/12) = 1/156$.

For solving this problem directly by enumerating the chances we remark that the number of the possible arrangements of the 13 cards is (§ 4) 13!; after assigning the first two places to the ace and deuce, the number of permutations of the 11 other cards is 11! and the required probability is $11!/13! = 1/156$, just as derived previously without needing to know the formula for the number of permutations.

**23.** We have considered a compound event resulting from the concurrence of two other events of which the second was subordinated to the first. However, a compound event frequently results from the concurrence of two or many events independent from one another, each having its own proper probability. It is required to determine the probability of that compound event.

Suppose that we have two urns, one of them with $m$ and $m'$ white and black balls, the other, with $n$ and $n'$ balls of those colours. It is required to determine the probability of extracting two white balls, one from each urn. There is evidently as many combinations or equal chances as unities in the product obtained by multiplying together the total numbers of balls in both urns (§ 2). And there are as many combinations or chances favouring the compound event as unities in the product obtained by multiplying together the numbers of white balls in both urns. The required probability is thus

$$mn/[(m + m')(n + n')] = [m/(m + m')][n/(n + n')].$$

By generalizing this reasoning we may formulate the following rule:

**IV.** *The product of the probabilities of many events independent from each other is the probability of the compound event resulting from their concurrence.* Or, more briefly: *The probability of a compound event is the product of the simple probabilities*.

Suppose that $p$ and $q$ are the probabilities of two contrary events A and B such that one of them necessarily occurs [at each trial]; let also A′ and B′ be two other contrary events with probabilities $p'$ and $q'$, etc. Then

1 = p + q = p′ + q′ = …

The product

(p + q) (p′ + q′)(p″ + q″) …

being expanded in a series of terms like pp′q″… corresponding to compound events AA′B″ … resulting from the concurrence of simple events A, A′, B″, … has as many such terms and as many compound events as possible combinations from a set of random events (§ 7). The sum of all these terms is unity as it should be since one of the possible compound events formed from all possible combinations of simple events should necessarily occur.

If *p, q, r* are the probabilities of events A, B, C one of which should necessarily occur in a random trial, $p + q + r = 1$. Then the binomial factor (p + q) in the product above should be replaced by a trinomial factor (p + q + r) etc.

**24.** Here is one more principle, an evident corollary of the precedent.

**V.** *The absolute probability of an event having differing probabilities according to different hypotheses is the sum of the compound probabilities obtained by multiplying the probability of that event according to each hypothesis by the probability of that hypothesis.*

And so, suppose that two urns contain *m* and *n* white, and *m′* and *n′* black balls. It is required to determine the probability of extracting a white ball out of an urn selected by chance. The probability of selecting each urn is 1/2, and that of drawing a white ball is $m/(m + m′)$ or $n/(n + n′)$, and the required probability is

1/2[m/( m + m′)] + 1/2[n/( n + n′)].

The probability of extracting a black ball is

1/2[m′/( m + m′)] + 1/2[n′/( n + n′)].

The sum of both probabilities is unity since a ball of either colour will necessarily be drawn.

It will be a grave mistake to decide that the probability of extracting a white ball is the ratio of the total number of white balls in both urns to the total number of all the balls there without considering the combinations resulting from the distribution of the balls among the urns. Suppose for example that the first urn contains 1 white ball and 2 black balls, and the second urn, respectively, 5 and 3. The number of white balls exceeds the number of black balls, and it is possible to believe that there are more chances to draw a white ball, or that the probability of that event exceeds 1/2. However, according to the method of extracting the balls, that probability is by the preceding formula

(1/2)(1/3) + (1/2)(5/8) = 23/48 < 1/2.

**25.** One and the same problem about combinations or probabilities can be presented in different ways and we should select that which leads to the most elegant or simplest solution. Let us take an example which occurs at the yearly military recruitment. Denote by $N$ the number of young men entered in a canton's list, by $N'$, those who have a legal cause for being exempted, and by $c$, the canton's contingent. It is required to determine the probability that the drawings will reach number $n$, $c < n < c + N'$.

We may suppose that there are $(N - N')$ white balls and $N'$ black balls and a box with $N$ pigeonholes numbered 1, 2, …, $N$. The total number of permutations obtained by successively occupying each pigeonhole with each ball is $N!$. All the arrangements ($S$) in which the number of black balls contained in the first $n$ pigeonholes is not less than $(n - c)$ correspond to chances of the drawings reaching number $n$. Then $S/N!$ is the required probability and the problem is reduced to solving a problem in permutations for determining $S$.

But we can also imagine an urn containing $(N - N')$ white and $N'$ black balls and the required probability will be the same as extracting from that urn at least $(n - c)$ white balls in $n$ consecutive drawings without replacement. That method of drawing which occurs in many other problems has no physical resemblance with the one applied in military recruitment. However, the problem stated in that form can be easily solved by applying the principle of compound probabilities, as we will indicate in the next chapter.

**26.** Sometimes very particular considerations following from the physical conditions of a problem can dispense with all enumeration of chances and calculations in general. The game of *passe-dix* offers such an example. When tossing 3 dice on a table a gambler bets against his adversary on the sum of the appeared points to exceed 10. Among the 216 possible combinations (§ 6) we should enumerate those providing a number larger than 10. It is easy to find formulas which will dispense with the need to enumerate the combinations one by one, but it is even simpler to benefit from the following remark.

On ordinary dice the points are disposed so that their sum on each of the two opposite faces is 7: 1 is opposite to 6 etc. And even if the manufacturer did not adopt that usage, it is always possible, without changing the conditions of randomness, to admit that the points are arranged as usual. According to those suppositions, the sum of the appeared points together with the opposite points on the faces lying on the table is 21[6]. Therefore, to each combination favouring the gambler that bet on *passe-dix*, there corresponds another leading to his loss, − the one which can be obtained by turning over the three dice. It is thus seen that each gambler has the same number of chances favouring and contrary to him so that they can bet even money.

### Notes

1. See Laplace (1814/1995, p. 4). B. B.
2. Concerning that term, Bru refers to Jakob Bernoulli's *Ars Conjectandi*, p. 212 (Chapter 1 of Pt. 4). Bernoulli's discussion was, however, more philosophical than mathematical.

**3.** Suppose that it is required to calculate the product *x*! by logarithms. If *x* is a large number, and a certain precision is necessary, we should use logarithms calculated with many more digits than provided by ordinary tables and their addition becomes laborious. Fortunately, this work can be avoided by means of the very remarkable Stirling formula which is a particular case of a much more general formula discovered by Euler for transforming sums into integrals and vice versa. The Stirling formula is

$$\ln x! = \ln \sqrt{2\pi} + (x + 1/2)\ln x - x + \frac{1}{12x} - \frac{1}{360x^2} + \ldots \qquad (13.1)$$

It was by means of this formula that the tables mentioned in the Note to § 10 were compiled. The Stirling formula can be considered typical of those which are applied in the theory of probability for numerical calculations. All of them have the singular and characteristic property of engendering series whose consecutive terms at first decrease very rapidly when the natural number which serves as the variable is only a few dozen, but then always take very slowly increasing values. This suffices for including them in the class of divergent. The series of that class can nevertheless be safely applied for numerical calculations of functions whose expansions they are if only we are able to assign a superior limit for the error made when stopping at some term. These limiting values are then of the order of magnitudes which we are allowed to neglect. Such limits were established for the Stirling formula; in general, however, they essentially exceed the involved error and do not therefore provide a proper idea of the obtained approximation.

They apparently require to take more observations than really needed or to assign to *x* larger values than really sufficient (? - O.S.). That grave imperfection is more or less characteristic of all similar formulas applied in the theory of probability. It can be said that they provide more than they promise in the sense that they ensure a sufficient approximation in the cases in which it is not yet possible to prove rigorously that the approximation is sufficient. Suppose that $x = 10$, then $x! = 3,628,800$ and $\lg x = 6.5597630$. When only taking into account the term $1/12x$ that logarithm will be 6.5597642 and its error barely exceeds $1/10^6$.

Another formula can be derived from (13.1):

$$x! = \sqrt{2\pi} x^{x+1/2} e^{-x}[1 + \frac{1}{12x} + \frac{1}{288x^2} + \ldots].$$

If applying it for calculating 28!12!/[32!8!], the number that expresses the probability that in the game of piquet after dealing out the cards a gambler will have all four aces, when neglecting the terms $1/12x + \ldots$ of the series, we will find 0.013807 with an error of only 1/335 of the true value. When the term $1/12x$ is included, we will have that value 0.0137653 with precision $1/10^7$. A. A. C.

De Moivre is known to have also derived the Stirling formula although he did not notice that the constant included there was the square root of $2\pi$, see De Moivre (1718/1756, p. 244). In that source, his note of 1733 was reprinted on pp. 243 – 254. De Moivre also published a table of lg *n*! for *n* = 10(10)900 with 14 digits correct to 11 – 12 digits with a misprint in the fifth digit of lg 380!. Bayes (1764) denied the possibility of applying divergent series (of the Stirling formula in the first place) but his criticism remained unheeded. Fichtenholz (1947/1951, § 501, p. 820) rigorously proved that the Stirling series diverged. O. S.

**4.** The rapprochement of mathematical principles of the calculus of probability with those of infinitesimal calculus clearly shows the similarity proper for being indicated to a reader familiar with both. Mathematical probability is a ratio between two terms which can increase to infinity with that ratio converging to a finite and assignable limit. A fluxion, or a derivative, or a differential coefficient (since all these terms are identical) is a ratio between two terms decreasing indefinitely with that ratio converging to a finite and assignable limit.

Considering directly, by reasoning and calculation, mathematical probability instead of combinations (that is, the ratio instead of its terms), is the same as instead of dealing with infinitesimal magnitudes in Leibniz' theory, operating directly with fluxions or derivatives in the theories of Newton and Lagrange. In both cases we introduce an auxiliary symbol and substitute a direct procedure corresponding to the nature of things by an artificial method adapted to our intellectual organisation. A. A. C.

**5.** Bru noted that Jakob Bernoulli (*Ars Conjectandi*, p. 227, this being the end of Chapter 4 in Pt. 4) actually discussed the possibility of introducing geometric probability and that Buffon (1777, e. g., § 23) forcefully introduced that concept. A thought experiment introducing geometric probability was due to Newton, but his manuscript of 1664 – 1666 was only published in 1967, see Sheynin (2003, p. 42).

**6.** This is difficult to understand.

## Bibliography

**Bayes T.** (1764), A letter … to J. Canton. *Phil. Trans. Roy. Soc.*, vol. 53, pp. 269 – 271.

**Buffon G. L. L.** (1777), Essai d'arithmétique morale. *Oeuvr. Phil.* Paris, 1954, pp. 456 – 488. English translation: Isf.lu/eng/Research/Working-Papers/2010

**De Moivre A.** (1718), *Doctrine of Chances*. London, 1738, 1756. Reprint of 3$^{rd}$ edition: New York, 1967.

**Fichtenholz G. M.** (1947), *Kurs Differenzialnogo i Integralnogo Ischislenia* (Course in Differential and Integral Calculus), vol. 2. Moscow – Leningrad, 1951. Single numbering of chapters and pages for all 3 volumes.

**Laplace P. S.** (1814 French), *Philosophical Essay on Probabilities*. New York, 1995. Translator A. I. Dale.

**Sheynin O.** (2003), Geometric probability and the Bertrand paradox. *Hist. Scientiarum*, vol. 13, pp. 42 – 53.

### Chapter 3. Laws of Mathematical Probability of Repeated Events

**27.** The mathematical theory of chances would have only been attractive for speculations if restricted to finding out how many favourable and contrary chances has an isolated event, repeated either not at all, or only under very rare circumstances.

We will show, however, that it becomes very important even for practical purposes when trials of the same randomness are repeated many times under similar circumstances. We can liken all such repeated trials to repeated drawings with replacement from an urn containing balls of different colours, so that the chances at each extraction remain without change. The solution of all possible problems concerning repeated trials is implicitly contained in the rule derived in § 23 about the principle of compound probabilities. Suppose that events A′, A″, … are repetitions of event A, events B′, B″, … repetitions of event B, … Then

$$p = p' = p'', \ldots q = q' = q'', \ldots$$

and if there are *m* trials, the product

$$(p + q)(p' + q')(p'' + q'') \ldots$$

becomes $(p + q)^m$. Its general term $C_m^{m-n} p^n q^{m-n}$ expresses the probability that in *m* trials the events A and B occur *n* and $(m - n)$ times[1].

The sum of the terms of that binomial until and including its general term expresses the probability that event A arrives at least *n* times or that the contrary event B will not occur more than $(m - n)$ times. It is required to determine the probability of obtaining one point at least twice in four tosses of a die. We have $p = 1/6$, $q = 5/6$, $m = 4$, and that probability is

$$p^4 + 4p^3q + 6p^2q^2 = 171/1296,$$

a fraction between 1/7 and 1/8. Then, how many trials are necessary for the event A to occur at least once with probability 1/2? Or, which is the same, for the event B to occur invariably with the same probability? The unknown *m* should be derived from $q^m = 1/2$. […] Suppose that the event A is the arrival of a double-six when tossing two dice at once. Then (§ 6)

$$p = 1/36, q = 35/36, m = \lg 2/(\lg 36 - \lg 35) = 24.6 \ldots$$

It is therefore beneficial to bet on the occurrence of a double-six in 25 tosses, but disadvantageous to agree on 24 tosses. This is the only possibility of interpreting, in this case, the derived incommensurability of number *m* which by its essence should be natural[2].

**28.** Each term of the binomial $(p + q)^m$ corresponds to a possible hypothesis about the ratio of the numbers of events A and B for the case of *m* trials. The sum of all the terms or of all the probabilities corresponding to these different hypotheses is unity. And since the number of the terms or the hypotheses is (*m* + 1), it will not be difficult to understand that the absolute (? - O.S.) values of the different terms should become ever smaller as the number of the trials increases. However, while decreasing, they maintain certain ratios between them. The law of those ratios is now our most important subject.

The general term of the binomial, $C_m^{m-n} p^n q^{m-n}$ is preceded by term $C_m^{m-n-1} p^{n+1} q^{m-n-1}$ and the ratio of the former to the latter is

$$\frac{n+1}{m-n} \cdot \frac{q}{p}. \qquad (28.1)$$

If this ratio is larger than unity, or *n* + 1 > *p*(*m* + 1), the former term is larger than the latter and vice versa.

If *p*(*m* + 1) is a natural number, call it *k*, then there will be number *n* = *k* in the number sequence 0, 1, …, *m*, and the term *k* of the binomial corresponding to the event A occurring *k* times, and B, (*m* – *k*) times, will be followed by another term of the same value. Otherwise, denote by *k* the largest natural number contained in *p*(*m* + 1), and the term *k* will be larger than both its preceding and following terms, i. e., will be the largest term of the expansion.

If *pm* is a natural number, it will indeed be the largest *k* contained in *p*(*m* + 1), and (*m* – *k*) will be another natural number equal to *qm*. And the largest term of the expansion will correspond to the combination for which the ratio of the numbers of the occurrences of events A and B is the same as *p*/*q*. In any case, the largest natural number *k* contained in *p*(*m* + 1) will differ less than by unity from *pm*. And if neglecting that fraction of unity as compared with numbers *pm* and *qm* (which is justified when the number *m* is very large, at least when *p* or *q* are not extremely small fractions) we may say that in general the most probable combination is that for which the number of events A is to the number of events B as *p*/*q*.

It is in addition evident that, if *n* = *pm*, *m* – *n* = *qm*, the ratio (28.1) becomes

$$\frac{pm+1}{qm} \cdot \frac{q}{p}$$

and approaches unity the nearer the larger becomes *m*. And comparing the largest term not only with the immediately preceding and following terms, but with those situated two, three, four, … places further, we become assured in that in both directions from that term the decrease becomes ever less rapid as the number *m* increases.

Those terms whose values are the next largest and whose sum constitute the greater part of the total sum of all the terms of the expansion accumulate in the vicinity of the maximal term.

**29.** For elucidating these notions by an example, I suppose that a ball should be extracted at random from an urn containing 2 white balls and 1 black ball. Event A will be the arrival of a white ball whose probability is 2/3. The contrary event B with probability 1/3 is the appearance of the black ball. I denote by (*u, v*) a compound event, an extraction of *u* white and *v* black balls. The following numbers are the nominators of the fractions that denote the corresponding probability; for 3, 6 and 9 drawings the denominators are 27, 729 and 19,683. According to the binomial formula, we will have for those drawings

(3, 0, 8); (2,1, 12); (1, 2, 6); (0, 3, 1)

(6, 0, 64); (5, 1, 192); (4, 2, 240); (3,3, 160);
(2, 4, 60); (1,5, 12); (0, 6, 1)

(9, 0, 512); (8, 1, 2304); (7, 2, 4608); (6, 3, 5376); (5, 4, 4032);
(4, 5, 2016); (3, 6, 672); (2, 7, 144); (1, 8, 18); (0, 9, 1)

In these three series of probabilities the largest term corresponds to the combination in which the number of white balls exactly twice exceeds that of black balls, and the terms decrease in both directions from those maximal terms. The ratios of the largest terms to those immediately preceding or following them diminish and approach unity as the number of terms increases. On the contrary, the ratios of the largest to the extreme terms invariably increase since the latter decrease very rapidly whereas the largest term decreases as well, but much slower.

For the last series of 10 different terms or combinations the sum of the largest term and the immediately preceding and following terms amounts to more than 0.7 of the sum of all the terms taken together. And, taking a series of 90 extractions, we conclude by aid of logarithmic tables (§ 10) that the sums of the largest term and the two terms neighbouring it from both sides are

(62, 28, 0.081817); (61, 29, 0.087460); (60, 30, 0.088918);
(59, 31, 0.086049); (58, 32, 0.079327)

The sum of those 5 terms is 0.423571, more than 2/5 of the sum of all the terms. On the contrary, the numerical values of the extreme terms are excessively small; for the term (90, 0) the probability is a fraction with 1 in the numerator and a number with 16 digits in the denominator. For the term (0, 90) that fraction is even incomparably smaller: with the same nominator its denominator is a number with 43 digits[3].

**30.** Bearing in mind that we are dealing with fundamental propositions, we now summarize and complete what was said.

**I.** *When the arrivals of events A or B depend on a random trial and when these trials are repeated many times, the most probable distribution is that for which the ratio of the numbers of these events* [the former ratio] *is equal to the ratio of their probabilities* [the latter ratio] *or deviate from it as little as possible. The probabilities of the other distributions decrease as the former ratio ever more deviates from the latter.*

**II.** *As the trials multiply, the number of possible distributions increases, and the probability of each former ratio lowers the more*

*rapidly the more that ratio deviates from the latter ratio, and vice versa*.

**III.** *Therefore, there exists an ever heightening probability that the former ratio will not deviate from the latter ratio beyond certain given limits. And, however narrow these limits are chosen, that probability can arbitrarily close approach unity, if only the number of the trials is sufficiently increased*.

For these theorems, we are obliged to Jakob Bernoulli who provided them in Part 4 of his posthumous *Ars Conjectandi* of 1713.

We should not forget that in these various statements the term *probability* is only understood in its mathematical sense (§ 12). Proposition I therefore signifies that, among all the combinations or hypotheses which do not take into account the order of the succession of the events A and B the number of those indicating that the former ratio is equal to the ratio of their chances exceeds the number of the others.

**31.** For rendering obvious the law obeyed by the numerical values of the different terms of the binomial expansion we can draw a segment AB and divide it into *m* equal intervals. Erect then perpendiculars to AB from the (*m* + 1) obtained points including points A and B and mark off distances proportional to the value of the first term (*m*, 0) on perpendicular Aa, the distance proportional to the value of the second term (*m* – 1, 1) on perpendicular $A_1a_1$ etc. If *m* is a considerable number, the points of the division of AB will be very close to each other, and each perpendicular or *ordinate* will little differ from the neighbouring ordinates. We can join the ends of the ordinates by a curve whose course will represent the law which we wished to show. That curve has a maximal ordinate Kk. By a generally known rule the straight line touching the curve at point k is parallel to AB. Finally, according to another rule of geometry the area [under a portion of the curve] is approximately equal to the product of one interval of AB by the sum of the included ordinates plus half the sum of the extreme ordinates.

Suppose for the sake of greater simplicity that *p* is commensurable and that the given value of *m* renders *pm* a natural number. If *m* increases but invariably obeys that condition, and if the same construction is repeated, we will obtain another curve whose maximal ordinate Kk will begin at the same point K, but in virtue of law II, the ordinates which delimited a portion of the initial curve and have their bases in their former places will shorten more rapidly than Kk and the partial area [under the curve] will constitute a larger portion of the whole area. And, by law III, we can sufficiently increase *m* for that portion to differ arbitrarily little from the whole area.

**32.** When *m* becomes very large, a direct calculation of a sum of a large number of terms of the binomial expansion will be impractical. We turn to formulas of approximation whose use exactly corresponds to constructing a curve as above through the ends of the ordinates and substituting by the mentioned rule the calculation of a portion of the curve's area by summing the corresponding ordinates. Among formulas or *functions* analysts select for this goal those which can be written algebraically whose expressions are not too complicated and whose courses best accord with that of the curve described above. The same method is applied in a similar case for passing from a real discontinuity to a fictitious continuity (§ 19).

**33.** Denote by *P* the probability that the number of events A in *m* trials is contained within $m(p - l)$ and $m(p + l)$ or that the ratio *w* of that number to the total number of the events will be contained within $(p - l)$ and $(p + l)$. According to formulas of approximation, for large values of *m* the value of *P* only depends on

$$t = l\sqrt{m / 2p(1-p)} \qquad (33.1)$$

so that, if *t* does not change (but *l, m* and *p* can change) the probability *P* will not vary either.

As the algebraists say, *P* is a *function* of abstract number *t*. Now, *t* varies proportionally to *l* which is the limit of the difference between *p* and *w*, proportionally to the square root of the number *m* of trials and inversely proportional to the square root of $p(1 - p)$. Therefore,

**[1]** After assigning certain values to *l* and *m* and calculating *P*, we will successively try out values of *l* equal to 1/2, 1/3, 1/4, 1/10 of its initial value. The number of trials should be increased 4, 9, 16, 100 times for obtaining the same probability that the random deviation $\pm (p - w)$ is contained within these new limits. In other words, to obtain the same probability that the anomalies of randomness concerning the determination of *w* are contained in ever narrower limits it is necessary to increase the number of trials inversely proportionally to the squares of these limits.

**[2]** For values of *p* very little differing from 0 and 1 the product $p(1 - p)$ is very small. It attains its largest value at $p = 1/2$. Therefore, the more different are the probabilities of A and the contrary event, the less is the need to multiply the trials for obtaining the same probability *P* that the random deviation $\pm (p - w)$ will be contained within the same limits; or, the narrower they become if the probability and the number of trials are the same.

These rules are derived by approximate calculations and are themselves only approximately exact. However, the approximation that they provide is quite sufficient when *m* is of the order of hundreds; still better, of thousands, tens of thousands, … Such numbers rarely occur in random trials taking place in agreements between individuals, but they are usual in physical and social phenomena for which the theory mainly ought to be established[4].

**34.** After assigning the numbers *p, m, l* the magnitude *t* is determined and it is necessary to calculate the corresponding value of *P* by approximate formulas whose origin we have only indicated. Or, which is much better, it is necessary to calculate once and for all a table of the values of *P* corresponding to a series of values of *t* sufficiently close to each other. Those alien to higher mathematics can apply that table without knowing the theory of its compilation just the same as we invariably use the tables of logarithms and sines without knowing either the theory of these magnitudes or the methods of compiling these tables.

At the end of this book I provide a table calculated for values of $t$ from 0 to 3 increasing by hundredths. It will be barely useful to continue that table further since for $t = 3$ $P$ is already 0.999978 so that the probability $(1 - P)$ of a deviation larger than the $l$ corresponding to $t = 3$ becomes an extremely small fraction 0.000022 which is lower than the probability of randomly extracting the single black ball from an urn containing 45,000 balls of other colours.

The Table indicated that the value $P = 1/2$ corresponds to the value of $t$ between 0.47 and 0.48 and calculation provides 0.476937. We will call the value of $l$ which, for given magnitudes $m$ and $p$, leads to that value of $t$ and to $P = 1/2$, the *median* value. Then there will be as many chances for $\pm (p - w)$, or to the numerical value of the deviation, to be contained within or beyond that value. In that case and in similar cases authors have applied the expression *probable value* which is not at all proper. On the one hand, all possible values of a deviation have their own chances or probabilities, and, on the other hand, the law of their probabilities (§ 16) is such that as the numerical value assigned to the deviation decreases its probability heightens, and the value of the deviation which we called *median* is actually *less probable* than any other of its smaller values.

Already at $t = 2$ we have such a value of $l$ that the probability of a larger deviation is 0.00468 or lower than the probability of extracting by chance a black ball from an urn having only one such ball out of 212. Finally, if $t = 2.87$, the probability $(1 - P)$ becomes equal to that of a drawing by chance of the only black ball out of 20,000. This value of $t$ is remarkable in that it is approximately 6 times larger than the value providing $P = 1/2$. And, since for constant values of $m$ and $p$, $t$ and $l$ increase in the same proportion, it is easy to remember that the median value of a deviation is 1/6 of the value that can be selected as that limit beyond which a deviation will have probability not higher than 1/20,000.

**35.** Suppose, just like in § 29, that an event A consists in drawing a white ball from an urn containing 2 white balls and 1 black ball. For a series of 9000 trials the median value of the deviation is 0.003348. There is probability 1/2 that the number of events A will be between 5970 and 6030; probability 211/212 for these limits to be 5874 and 6126; and an extremely high probability of 19,999/20,000 that they will be 5819 and 6181.

Suppose now that there are 9 million trials. Then the median value of the deviation will be about 32 times smaller. The limits corresponding to the same probabilities as above will be 5,999,047 and 6,000,953; 5,996,000 and 6,004,000; 5,994,260 and 6,005,740.

**36.** We (§ 28) have considered all the repetitions of the same random trials consisting of drawings with replacement so that the random conditions of the successive extractions did not change. Drawings without replacement are also important, − less important but, because of similarity, proper to be treated here briefly.

Suppose we have $a$ white and $b$ black balls in an urn from which balls are extracted one by one without replacement. It is required to determine the probability of drawing $n$ white and $(m - n)$ black balls in $m$ extractions, We continue to denote by A a simple event consisting of drawing a white ball, and by B, a contrary event of extracting a black ball. Then AB will indicate a compound event consisting of an arrival of a white ball followed by a black ball etc. It is obvious that, according to the principle of compound probabilities (§ 23), such a compound event as AAB will have probability

$$\frac{a}{a+b} \cdot \frac{a-1}{a+b-1} \cdot \frac{b}{a+b-2} = \frac{a(a-1)b}{(a+b)(a+b-1)(a+b-2)},$$

whereas the compound event ABA, which only differs from the previous by the order of the succession of the simple events, has probability only differing from the previous expression by the order of the numerator's factors.

The generality of this remark is evident, so that the required probability is

$$\frac{a(a-1)...(a-n+1)b(b-1)...[b-(m-n)+1]}{(a+b)(a+b-1)...(a+b-m+1)}$$

taken as many times as there are permutations in the order of the events (§ 5). So that probability is

$$C_m^{m-n} \frac{a(a-1)...(a-n+1)b(b-1)...[b-(m-n)+1]}{(a+b)(a+b-1)...(a+b-m+1)}. \qquad (36.1)$$

If consecutively $n = 1, 2, \ldots (m-1)$, we will have the probabilities of the arrival of $1, 2, \ldots, (m-1)$ white, and $(m-1), (m-2), \ldots, 1$ black balls in $m$ drawings. The probability of only extracting white balls is

$$\frac{a(a-1)...(a-m+1)}{(a+b)(a+b-1)...(a+b-m+1)}.$$

The sum of all the probabilities is unity since it corresponds to the same number of hypotheses one of which ought to occur. Therefore,

$$1 = \frac{a(a-1)...(a-m+1)}{(a+b)(a+b-1)...(a+b-m+1)} +$$

$$C_m^1 \frac{a(a-1)(a-2)...(a-m+2)b}{(a+b)(a+b-1)...(a+b-m+1)} + ... +$$

$$C_m^{m-n} \frac{a(a-1)...(a-n+1)b(b-1)...[b-(m-n)+1]}{(a+b)(a+b-1)...(a+b-m+1)} + ... +$$

$$\frac{b(b-1)...(b-m+1]}{(a+b)(a+b-1)...(a+b-m+1)}. \qquad (36.2)$$

The denominator of the last fraction is

$$a(a-1) \ldots (a-m+1) + C_m^1 \, a(a-1) \ldots (a-m+2)b + \ldots +$$
$$C_m^{m-n} a(a-1) \ldots (a-n+1)b(b-1) \ldots [b-(m-n)+1] + \ldots +$$
$$b(b-1) \ldots (b-m+1). \qquad (36.3)$$

This last formula is remarkable for its similarity with the binomial expansion with multipliers replacing powers of the latter. This analogy already sensed (§ 6) is justified by the very principles of combinatorial synthesis. That formula also expresses an algebraic relation which should persist whichever are the numerical values of *a, b, m* with *a + b > m* and therefore even if those letters do not anymore denote natural numbers. It can also be proved by pure algebra, but we saw how easy it followed when considering compound probabilities. Thus, for example, by artificially introducing an alien element it is sometimes possible to simplify a successive course of certain abstract truths and, for example, by considerations borrowed from mechanics simpler prove certain propositions of pure geometry.

When summing (*m – n* + 1) first terms of (36.2) until term (36.1) inclusive, we will obtain the probability that not less than *n* white balls are extracted in *m* drawings. The term (36.1) is the largest in the expansion (36.2) since *n*/(*m – n*) and therefore (*a – n*)/(*b – m + n*) are equal to *a/b* or differ from it as little as possible. The terms of the expansion (36.2) decrease on each side of the largest term following laws similar to those considered in § 28 and the next sections.

Instead of extracting *m* balls one by one, it is evidently possible to draw all of them at once without changing the probability of obtaining *n* white and (*m – n*) black balls.

**37.** We (§ 25) have indicated an application of a problem treated in § 36. It can be most immediately made use of in many problems about trials of sorts taking place in political assemblies whose members are usually separated in two parties, or thought to belong to two large factions. For example, an assembly has 459 members (§ 10) of which 240 belong to the majority, and 219, to the minority factions. Suppose that a deputation or commission of 20 is elected by chance. What will be the probability that the majority or the minority of the assembly will compose the majority or the minority of the commission?

If because of random causes such as illnesses acting independently from the parties, 30 members are absent at the voting, we can require the probability that the majority of the assembly will constitute a minority in the commission. That question is the same as asking for the probability that in 30 drawings by chance at least 26 white balls will be drawn from an urn containing 240 white and 219 black balls. Our formulas applied with the aid of logarithmic tables provide for this probability the value 0.000049547 ≈ 1/20,000.

The right to challenge peremptorily a certain number of judges or jurymen suggests similar problems. For example, out of 36 jurymen the accused and the public prosecutor have the right to challenge 12 of them respectively. The challenge is made by extracting the names of the jurymen from an urn in which 12 names should remain. For considering the simplest case, I suppose that the prosecutor has no reason to use his right of challenging, and that the accused wishes to reject 6 jurymen. It is required to determine the probability that the accused will not benefit from his right. The problem is tantamount to determining the probability that in 12 drawings only white balls will be extracted from an urn containing 30 white and 6 black balls. That probability is $0.069102 \approx 7/100^5$.

**Notes**

**1.** The meaning of $p$ and $q$ is obvious.

**2.** The problem concerning the abovementioned tricktrack is famous since it became the occasion (? - O.S.) of Pascal's first researches [in probability] and thus originated the calculus of probability. The correspondence of this great man tells us that the pertinent question was posed by a man about town Chevalier De Méré, who was remote from mathematics. Because of this fortunate circumstance, his name is since belonging to the history of science. A. A. C.

**3.** Suppose that $p = q = 1/2$ and $m = 100$. The terms equally remote from the middle term are equal to each other, and calculations provide [Cournot lists their values with 7 decimals for (50, 50) − (75, 25) and continues:]

The 50 terms situated furthest from the same middle term can evidently be considered absolutely negligible. Each of the values of the two extreme terms is $2^{-300}$ which is a fraction with numerator 1 and denominator, a number with 31 digits. The sum of the seven middle terms (53, 47) − (47, 53) is 0.5158814, larger than a half of the sum of all the 101 terms of the expansion. A. A. C.

**4.** In mathematical treatises it is shown that, neglecting magnitudes of the order of $1/m$,

$$P = \frac{2}{\sqrt{\pi}} \int_0^t \exp(-t^2)dt + \frac{\exp(-t^2)}{\sqrt{2\pi p(1-p)m}}. \qquad (33.2)$$

For the sake of greater simplicity we suppose that, throughout this book, the value of $P$ is only represented by the first term of (33.2). The same holds concerning similar expressions which we invariably apply. That simplification is all the more allowed since ordinarily it is much less required to calculate with a close approximation the value of $P$ rather than to assign to that fraction an inferior limit. The second term of (33.2) is positive so that, when neglecting it and finding out that the probability $(1 – P)$ is less than some number, the consequences of that fact will all the more persist if $P$ is corrected by that term.

Function $\exp(-x^2)$ should be considered as the algebraic type of functions [!] which very rapidly decrease symmetrically on both sides of the origin of the variable $t$. Its numerical value never exactly disappearing is excessively small for quite small values of $t$. It is for this reason (? - O.S.) that that function is included in all formulas constructed by analysts for applying them in the theory of chances. Because of that property of $\exp(-x^2)$ the curve passing through the ends of a large number of ordinates representing the consecutive terms of the expansion of the binomial $(p + q)^m$ can almost exactly coincide, especially near point $k$, with the curve whose ordinates are those of $(p + q)^m$ and abscissas measured along AB with origin at $K$. […] Denoting by $a$ and $b$ the abscissas KA and KB we will have approximately by the known theorem of geometry and integral calculus

$$P = \int_{-t}^{t} \exp(-t^2)dt \div \int_a^b \exp(-t^2)dt.$$

However, owing to the smallness of the ordinates Aa and Bb and the extreme rapidity with which the function $\exp(-x^2)$ decreases at larger numerical values of $t$, the integral in the denominator can be extended over $[-\infty, \infty]$ and therefore becoming equal to $\sqrt{\pi}$, so that this $P$ becomes equal to the first term of the formula (33.2).

Denote by *k* the largest natural number contained in *p*(*m* + 1) and let numbers *l* and *t* increase by steps with $l\sqrt{m(m+1)}$ remaining a natural number, call it λ. The second term in (33.2) expresses the probability that the number *n* of events A is exactly equal to (*k* − λ) or (*k* + λ). When the second term is subtracted from the first instead of being added to it we will get the probability that *n* is contained between (*k* − λ) and (*k* + λ) but dose not reach these limits. The first term of (33.2) expresses the probability that *n* is contained within (*k* − λ) and (*k* + λ + 1) or (*k* − λ − 1) and (*k* + λ). Finally, the complete value of *P* is the probability that *n* is contained within (*k* − λ − 1) and (*k* + λ + 1). The role of the second term and of the error made when neglecting it is thus better understood.

Although formula (33.2) is only considered exact to within magnitudes of the order of 1/*m*, it usually provides a much better approximation. Suppose that *p* = 1/2 and *m* = 100. According to the formula, the probability that the number *n* of events exceeds 39 but is less than 61 is 0.9653 whereas our calculation (Note 3) provides 0.9648, so that the error of the formula is only 0.0005 instead of one or a few hundredths as could have been feared when keeping to the exact terms of the ordinary demonstration. A. A. C.

Formula (33.2) is due to Laplace (1812/1886, p. 284) who had thus developed de Moivre's formula of 1733. It also occurs in Poisson (1837, § 79). On the role of the second term of formula (33.2), see Poisson (1837, § 79). B. B. See also Note 6 to Chapter 1.

**5.** This reasoning is incomplete: the accused would have probably been glad if even less than 6 undesirable jurymen were challenged.

## Bibliography

**Laplace P. S.** (1812), *Théorie analytique des probabilités. Oeuvr. Compl.*, t. 7. Paris, 1886.

**Poisson S.-D.** (1837), *Recherches sur la probabilité des jugements* … Paris, 1837, 2003, 2012. English translation: www.sheynin.de   downloadable file 53.

**Chapter 4. Randomness. Physical Possibility and Impossibility**

**38.** Until now, we had in a certain sense discussed pure arithmetic. We enumerated combinations; asked for the ratios between the numbers that express how many combinations lead to some result and how many are contrary to it; assigned the limits of those ratios when the pertinent numbers increased to infinity because of passing from discontinuity to the continuous; examined how the values of those ratios for compound events depended on the values calculated for the simple events.

Now, however, we aim to find out whether all that theory is only a jeu d'esprit, a curious speculation, or, on the contrary, it strives to discover very important and very general laws which govern the real world. For passing from the idea of an abstract ratio to the notion of an efficient law of realities and phenomena, mathematical reasoning based on a series of identities is obviously insufficient. We should turn to other notions, to other principles of knowledge. In a word, we should apply philosophical criticism. We direct complete attention of readers to this delicate point. Although this subject had been, in our opinion, imperfectly understood or described by philosophers both being geometers[1] or not, we do not despair of making that subject sufficiently clear for preventing any ambiguity or mistakes in its future applications.

**39.** No phenomenon or event is produced without a cause. This is a supreme principle and regulator of the human mind when facts of reality are investigated. The cause of a phenomenon often eludes us or we assume as a cause something which it is not. However, neither our helplessness in applying the principle of causality nor the mistakes made while applying it can shake our attachment to this principle understood as an absolute and necessary rule.

We move from an effect to its immediate cause. In turn, that cause is considered as an effect, etc, without at all understanding any limits of this law of regression. In turn, the actual effect becomes or can become a cause of a subsequent effect, and thus it can continue to infinity. That infinite chain of causes and effects exists in time[2]. Actual phenomena form its links and constitute a linear series (§ 4). An infinity of such series can coexist in time, they can intersect when a phenomenon, engendered by many conducive phenomena, is an effect of many different series of generating causes or in turn engenders many series of effects which remain different and perfectly separated beyond the initial term.

A simple idea of such intersections and such isolation is formed by comparing that with human generations. Through his parents a man is connected with two series of ancestors. And, while ascending, both these series fork with each generation. In turn, the man can originate and be the common ancestor of many descending lines, which, after issuing from him, do not intersect or only intersect accidentally by marriages within the family. After some time, each family or each genealogic branch contracts marriages with many others, but these other branches much oftener propagate collaterally and remain perfectly distinct and isolated from each other. Or, if they have a common origin, its genuineness is not justified by science or historical research.

In the ascending order, each human generation only provides a division by pairs. However, we can imagine a possible existence of a much larger number of divisions in both temporal directions when some causes and effects are considered. Some phenomenon can be imagined as being conditioned by a multitude of different causes. Passing in most cases from discontinuity to continuity and thus leading us to believe in the infinity the number of conducive causes to infinity even appears to conform to the general plan of nature (§ 15). Therefore, the branches of intersecting lines which in our imagination represent the causal connection of phenomena similar to a pencil of rays, widen, and concentrate without tearing their tissues.

**40.** However, whether we consider the number of generating causes of some phenomenon finite or infinite, according to the principle of common sense there will be series of *solidary* phenomena depending on each other and other series which develop on parallel or successively without any dependence, any solidarity. Actually, some philosophers[3] imagined that everything in the whole world is connected and proved it by subtle arguments or ingenious trifles. However, neither their subtleties, nor their nothings can prevail over beliefs of common sense. No one can seriously think that kicking the Earth will hamper a navigator sailing somewhere at the antipodes or stir the system of Jupiter's satellites. When wishing to admit theoretically the existence of disturbances of that kind caused by such causes, it will be necessary to recognize that these perturbations are imperceptible so that we have no means for finding their traces in phenomena. In other words, the alleged solidarity does not manifest itself by any sensible sign as though among observable facts it does not exist.

*Events causally produced by combinations or encounters of phenomena belonging to independent series[4], are those which we call fortuitous or resulting from randomness.*

**41.** We elucidate that statement by examples. I suppose that brothers serving in the same corps perished in the same battle[5]. There is something startling us in the common misfortune, but, after deliberating we see that these two circumstances can be not independent one from another and that randomness alone did not lead to that pernicious rapprochement. Perhaps the younger brother decided on the military career following his elder brother; and it is therefore natural for them to try serving in the same corps. They are thus exposed to the same danger, feel it necessary to help each other and, if in great danger, both could have succumbed which would not be surprising.

Causes independent from their ties of relationship could have played a role in that event, but the coincidence of their being brothers and their common fate was not due to pure chance. Now I suppose that they served in different armies, one on the northern frontier and the other, at the foot of the Alps. Battles occurred on the same day in both places and the brothers perished. We are justified to consider this coincidence as a result of chance.

At a great distance apart the operations of the two armies consisted of two series of facts. General orders could have been issued from a common centre but later developments were perfectly independent from each other and accommodated local circumstances. Those leading to the battle on some day and at some place in the first army had no connection with the circumstances which, on the same day, brought about the battle in the second. And if those corps to which the brothers belonged participated in both battles respectively, if the fights were desperate and both brothers succumbed, there will be nothing in their blood relation which justifies that coincidence.

A man surprised by a storm took refuge under an isolated tree and was killed by a lightning. That accident is not purely random. Physics tells us that electric fluids tend to discharge on treetops and any spikes. And it stands to reason that a man ignorant of physical principles chose a tree as a refuge and it is also justifiable that a lightning struck exactly that place. On the contrary, if a man is killed in the steppe or forest, that accident would have been fortuitous since there is no connection between the causes that brought the man to that point and the reason for the lightning to strike him at that moment.

A man who can not read extracts letters one after another from disordered type of a printing house, and they form the word *Alexander*. That is an accidental fact or the result of randomness since there is no connection at all between the causes that direct the man's hand and those that impose the name *Alexander* of the famous conqueror which was later attributed to other historical personalities so that that name became popular and one of the best known in the language.

**42.** Such events are rare and surprising but that is not at all the reason for choosing them as examples of the result of randomness. On the contrary (as we will explain and as it was possible to foresee by what was said in the preceding chapter), it is because randomness engendered them rather than many others caused by other combinations that they are rare. And since they are rare, they amaze us. There is nothing rare or surprising for a blindfolded man to extract a white or a black ball from an urn containing the same number of balls of both colours. And nevertheless either event will be justifiably considered a result of randomness since there is manifestly no connection between the causes for some ball to be taken by the man's hand and its colour.

It is quite true that in ordinary language we voluntarily apply the expression *randomness* when discussing rare and surprising combinations. Four uninterrupted extractions of a black ball from an urn containing the same number of white and black balls is thought to be an effect of a rare chance. The same will perhaps not be said if 2 white balls were followed by 2 black balls, and all the less said if the 4 balls followed each other less regularly. Nevertheless, the causes directing the operator's hand and those determining the colour of the balls are perfectly independent.

It was remarked that the death of both brothers on the same day was due to chance, but it would not have been indicated, or remarked less had they died a month, or three or six months apart. Still, in any of these cases there is no solidarity between the causes leading to the death of the elder brother on a certain day and those that brought about the death of the younger brother on some other, or those that made them brothers. When a blindfolded man extracts disordered letters, no attention would have been paid to their sets which do not represent any articulate sounds or words used in a known language. Nevertheless, there is never any connection between the causes that successively direct the operator's hand on a certain piece of metal and those which stamped certain letters on these pieces.

However, that nuance of the expression attached to the word *randomness* in usual conversation and the language of the town is vague and barely defined and it ought to be removed from the exact language of science and philosophy. To understand the notion of randomness well enough, we should only attach to it the fundamental and categorical; that is, the idea of independence or absence of solidarity between different series of facts or causes.

**43.** Another notion connected with randomness with very important consequences for theory and practice, is *physical impossibility*. It is now appropriate to turn to examples for facilitating the understanding of these abstract generalities.

We consider it physically impossible that a material cone remains in equilibrium on its apex; that an impulsion communicated to a sphere is exactly directed along a straight line passing through its centre and therefore does not lead to any rotation; that the centre of a disc falling on a floor covered with square tiles lands on the intersection of the diagonals of a tile; that an angle-measuring instrument is exactly centred; that a balance is rigorously exact; that a certain measure rigorously conforms to the standard, etc.

All these physical impossibilities are of the same nature and are evidently connected with the notion of chance encounters or independence of causes as described above. Suppose that it is required to find the centre of a circle. The ability of the performer and the precision of his instruments assign the limits of the possible error or the distance between the veritable centre and the determined point he indicates as the centre. On the other hand, within certain limits differing from the former and separated by a shorter interval the performer is not anymore guided by his senses or instruments.

The central point in a more or less small area is undoubtedly determined by some causes, but *blind* causes independent from geometric conditions which determine the veritable centre. There exists an infinity of points on which the blind causes can fix the performer's instrument with no reason occasioned by the nature of the task for choosing one point rather than another. The coincidence of the chosen point and the veritable centre is an event completely similar to blindly drawing a white ball from an urn containing a single white ball and infinitely many black balls.

*A physically impossible event is therefore such whose mathematical probability is infinitely low.* That sole remark provides thoroughness and an exceptional objective value for the theory of mathematical probability.

Just the same, if a sphere collides with a body moving in space, because of causes independent from the presence of that sphere in a certain place it is physically impossible, and it never happens, that among the infinitely many possible directions of that body the causes of its motion lead to its exactly passing through the centre of the sphere. We therefore admit the physical impossibility of the sphere not to begin rotating in addition to translating. If the impulsion is communicated by an intelligent being with a restricted sense of perfection who aimed at that result it would have still been physically impossible to achieve it. […] In the same way we can explain the physical impossibility admitted by the whole world of retaining a cone in equilibrium on its apex. Similar reasoning holds for all the described cases.

The notion of physical impossibility doubtless essentially differs from mathematical or metaphysical impossibility, but we are unable to establish the transition from one to another. A physically impossible thing is understood as mathematically or metaphysically possible although never happening. There is no reason for the only combination of facts or independent causes which can lead to it, to be present in preference to infinitely many others. That general and abstract notion of independence of causes and the infinite multitude of possible combinations provides a foundation of physical impossibility without turning to empirical notions about the material world conveyed to us by our senses. For this reason it is perhaps better to call physical impossibility *actual impossibility* opposing it to mathematical or metaphysical impossibility, apparently more properly called *rational* or *absolute impossibility*.

**44.** An actually or physically impossible event is that whose mathematical probability is infinitely low. It can be likened to a blind extraction of a single white ball from an urn also containing infinitely many black balls. However, having any finite ratio of white to black balls, when repeating ever more trial drawings we will obtain, according to the Bernoulli theorems (? - O.S.), a heightening probability that the ratio of the extracted white to black balls ever less deviates from the ratio of their mathematical probabilities. For an infinity of extractions we will get an infinitely low probability or a physical impossibility of those ratios to differ by a given and arbitrarily small fraction. And, assigning a sufficiently large number of trials and convenient limits of that deviation, the entire doctrine of mathematical probabilities will be attached to the notion of physical impossibility.

Mathematical probability will not be anymore a simple abstract ratio caused by our mind but the expression of the ratio maintained by the nature itself of the things established by observation when, under the influence of independent accidentally combining causes, the trials of the same randomness indefinitely multiply as it occurs all the time in natural phenomena and social facts.

In the strict language proper for abstract and absolute truths in mathematics and metaphysics a thing is possible or not; there are no degrees of possibility or impossibility. However, in the world of facts and realities, when two contrary phenomena can be and are realized according to random combinations of certain variable causes and other causes or constants, it is natural to regard a phenomenon as endowed with the greater ability to occur, or to occur with the greater physical or absolute possibility the oftener they are reproduced in a large number of trials. Mathematical probability then becomes the measure of *physical possibility* so that these expressions are interchangeable. This, after all, is only a definition of terms. The advantage of the term *possibility* (usage has already recognized the truth of what we describe) is that it clearly denotes the experience of a ratio which exists between the things themselves and does not depend on our manner of judging or sensing varying from one individual to another according to their circumstances and the degree of their knowledge.

Finally to apply the technical language of learning, the term *possibility* expresses an *objective* sense whereas *probability* ordinarily implies a *subjective* sense[6] and therefore deceived excellent minds, caused so many misunderstandings and corrupted ideas which should be formed about the theory of chances and mathematical probabilities.

**45.** Thus, often repeated is Hume's conception that *Properly speaking, there is nothing random but there is something equivalent to it: our ignorance of the real causes of events*. Laplace himself (1814/1995, p. 3) formulated the following principle: *Probability is relative in part to our ignorance and in part to our knowledge*. It follows that for a superior intelligence which would have discerned all the causes and the resulting effects[7] the science of probabilities disappears owing to the absence of a subject.

All these reflections are however wrong. Undoubtedly, the word *randomness* indicates an idea rather than a substantial cause, an idea of combinations of systems of causes or facts, each of them developing in its own series, independently one from another. A superior mind only differs from the human mind in that it is less often wrong, or never wrong. It will not risk considering series causally influencing each other as independent entities, or, inversely, imagine actually independent causes as depending on each other. A superior mind would have more surely or even quite exactly separated the part of randomness in the development of successive phenomena. It would have assigned in advance the results of the coincidence of independent causes which we most often are unable to do.

For example, an irregular die should be tossed a large number of times. At each toss causes, independent from those acting in the following trials, determine the intensity, direction and point of application of the impulsive forces. Then that mind, unlike ours, would have almost exactly known the ratio of the number of tosses leading to a determined outcome to the total number of them. And that knowledge would have been certain, whether or not the superior mind knew the acting forces, and was able to calculate the effects of each toss.

In a word, it advances further than we do and better applies the knowledge of those mathematical ratios which are connected with the notion of randomness and became the laws of nature in the world of phenomena. In this sense it is correctly stated (and very often repeated) that randomness governs the world or rather that it has its role, and a notable role at that, in governing the world.

This does not at all contradict the idea which we should form about the supreme and providential direction[8]. It either takes care of only mean and general results[9] ensured by the laws themselves of randomness, or the supreme cause disposes the details and particular facts for coordinating them with aims surpassing our sciences and theories. When remaining in the world of secondary causes and natural facts, which are the proper field of science, the mathematical theory of randomness appears as a widest application of the quantitative science and justifies in the best way the saying *Mundum regunt numeri* (The world is governed by numbers). Actually, in spite of the thoughts of some philosophers[10], nothing authorises us to believe that the foundation of all phenomena is found in the notions of extension, time, movement, and, in a word, in the notions of continuous and measurable magnitudes, the object of geometry.

At the level of our knowledge, the acts of intelligent and spiritual living beings can not be explained at all and we can fearlessly state that they will never be explained by mechanics (? - O.S.) of the geometers. They do not at all find themselves on the side of geometry or mechanics, in the numerical domain. However, they are led there since the notions of combinations and chance, cause and randomness are higher in the abstract world than geometry and mechanics and applicable to the facts of living nature[11], to the intellectual and moral field just as the phenomena produced by the movement of inert matter.

**46.** In essence, the theory of chances and mathematical probabilities is applicable to two clearly distinct fields: to problems of *possibility* which objectively exist, as we have explained, and of *probability* which are actually relative in part to our knowledge and in part to our ignorance[12]. When we say that the probability of achieving a double-six in tricktrack is 1/36 (§ 6), we can think about possibilities so that it means that, had the dice been perfectly regular and homogeneous cubes, there would have been no reason caused by their physical structure for one face to appear rather than another one. It then follows that a double-six will arrive in approximately 1/36 of the total number of tosses because the directions of the impulsive forces are absolutely independent from the numbers stamped on the faces of the dice.

However, we can also think about simple probability. Without enquiring whether that regularity of structure exists or not, it suffices that we do not know how its irregularities act if they exist. We will then have no reason to suppose that one face arrives rather than another one and the occurrence of a double-six for which there is only one combination out of 36 is for us less *probable* than the arrival of 2 and 1 favoured by two combinations. This conclusion takes place in spite of the latter outcome being perhaps less physically possible or even impossible (? - O.S.).

If a gambler bets on a double-six, and another one, on 2 and 1, and if they agree to disregard all other outcomes, there will be no other means (as we explain in the next chapter) to regulate their stakes than assuming the ratio of 1:2. And this will be as fair as when we are certain that the structure of the dice is perfectly regular. At the same time, if an umpire knew that the dice were fraudulent his agreeing with that ratio will be unfair and favourable to one of the gamblers.

This almost restricts [exhausts] the applications of the theory if it aims at simple judgements of probabilities varying according to the knowledge and ignorance of men. However, if the theory is transferred to discussions of natural phenomena and social facts, the consequences derived from such judgements can lead to mistakes which we will illustrate by examples and which can undermine the confidence in legitimate applications.

**47.** The calculus of chances was born in connection with such regulation of stakes (*compositio sortis*, casting lots[13]). A remarkable passage from Pascal's letter shows that he did not at all think about applications of the *geometry of chance* in the field of judging possibilities to the economy of natural facts. The great geniuses of the 17th century, Fermat, Leibniz, Huygens who were occupying themselves with the calculus of combinations and chances simultaneously with, or a few years later than Pascal, only thought about the problem of points[14]. Jakob Bernoulli, in his *Ars Conjectandi*, formally determined the essential aim, the objective value of the theory of chances. At the same time, however, the continuous application of the terms *probability, conjecture*, etc paved the way for misunderstandings leading to confusing expositions and uncertain applications.

The title of the first edition of De Moivre's *Doctrine of Chances* which appeared in 1718, five years after the publication of the *Ars Conjectandi*, lacked that inconvenience, but even now authors sometimes discuss *chances* of an event in the sense of *possibility* and *probability*. It is regrettably inconvenient to apply the same term in two meanings, for denoting either each random combination leading to a determined event, or the ratio of their number to the number of all random combinations with both numbers being either finite or infinite but their ratio converging to an assignable limit when those two numbers increase unboundedly.

To comply with the most ordinary usage, we continue to apply the word *probability* as a synonym of physical possibility with the exception of the case in which the discourse indicates a subjective meaning. At the end of this book we examine whether there are other probabilistic judgements in addition to those connected with the mathematical theory of chances and randomness and try thus to complete the exposition of our subject.

**48.** The term *probability* taken in the subjective sense corresponds to *certitude*, and it is often said that if probabilities are measured by fractions, unity is the measure of certainty[15]. Actually, if all the chances or possible random combinations favour an event, it certainly occurs and the probability of the contrary event is exactly zero. However, on the other hand it is recognized that there exists an essential difference[16], not only with respect to magnitudes, but between probability and absolute certainty.

It is absurd to say that absolute certainty is composed of the sum of two or more probabilities. In the entire doctrine of mathematical probabilities or possibilities the term for comparison is not a rationally, metaphysically or absolutely certain event, but a *physically certain* event, whose probability only differs from unity by an infinitely small magnitude, or the event whose contrary is physically impossible as explained above both by general reasoning and examples. Thus is homogeneity re-established. It should be invariably present in things subjected to measurement and calculation and only a semblance of difficulty remains for those familiar with the sense of expressions in mathematics.

## Notes

**1.** D'Alembert, Condorcet. [B. B].
**2.** This idea goes back to Cicero. [B. B.]
**3.** D'Holbach and others. [B. B.]
**4.** The stated definition of randomness goes back to Aristotle. [B. B.]
**5.** That example is due to Cicero. [B. B.]
**6.** Bru named Jakob Bernoulli but had not justified his statement.
**7.** See Laplace (1814/1995, p. 2) who followed d'Holbach and Buffon. [B. B.] I add: see also Maupertuis (1756, p. 300) and Boscovich (1758, § 385). O. S.
**8.** See *Ars Conjectandi,* beginning of Chapter 1 in pt. 4. [B. B.]
**9.** See De Moivre (1733/1756, p. 253). [B. B.]
**10.** Bru named Newton and Clarke. In 1715 – 1716 the latter exchanged letters with Leibniz in which he defended Newton. In 1717 he published this correspondence, see its edition of 1998.
**11.** Bru named Laplace (1814/1995, p. 62) and Quetelet.
**12.** In the beginning of § 45 Cournot quoted Laplace who had formulated that idea but added that *All these reflections are ... wrong*!
**13.** *And then an absolutely new scientific work whose subject is not until now studied: distribution of chances in games which obey them; in French, it is called "faire les partis des jeux"* [divide the stakes]. *The uncertain fortune so well submits to the fairness of calculation that each gambler always gets exactly <u>what belongs to him by right</u>. And this is certainly what should be determined by reasoning the more <u>the less possible it can be found by experience</u>.*

*Actually, the results of ambiguous lots are justifiably attributed to fortuitous contingency rather than to natural necessity. This is why the problem remained uncertain to this day. However, now, what <u>had been resisting experience</u>, can not avoid the dominion of reason. Owing to geometry, we have reduced the problem quite surely in an exact manner, and its certainty has partly daringly advanced. And so, combining rigour of scientific demonstration and incertitude of randomness and conciliating things apparently contrary to each other, it is possible to elicit its name from both and rightfully appropriate for it the stupefying title Geometry of randomness, aleae geometrie.*

This is a passage from Pascal (1654). The Academy of Sciences was only founded in 1666. A. A. C.

Cournot quoted this passage in its original Latin. Bru provided its French text, now also published, and my translation is from French. O. S.

**14.** This is only true with regard to Fermat (who only corresponded with Pascal about games of chance). Leibniz left not less than five manuscripts first published in 1866 on *Staatswissenschaft* and political arithmetic and Huygens is known to have stated that the new theory was only in the making. In 1669 he studied problems of mortality but that work was only published in 1895. And Cournot should have certainly mentioned Halley.
**15.** Bru mentioned the beginning of Chapter 1 of pt. 4 of the *Ars Conjectandi* and Lacroix.
**16.** Bru mentioned Buffon and Euler's *Letters to a German Princess*.

## Bibliography

**Bernoulli J.** (1713), *Ars Conjectandi*. *Werke*, Bd. 3. Basel, 1975, pp. 107 – 259. German translation of 1899 (Frankfurt/Main, 1999) is modernized. My English translation of pt. 4 is entitled *On the Law of Large Numbers*. Berlin, 2005. Also www.sheynin.de  downloadable file 8. See my review of an illiterate English translation of the whole book (2006) in *Hist. Scientiarum*, vol. 16, 2006, pp. 212 – 214.

**Boscovich R. J.** (1758), *Philosophiae naturalis theorie*. Chicago – London, 1922. Latin – English edition.

**De Moivre A.** (1733, Latin), A method for approximating the sum … Translated by author. In author's *Doctrine of Chances*. London, 1756 and New York, 1967, pp. 243 – 256, an extended version.

**Laplace P. S.** (1814, French), *Philosophical Essay on Probabilities*. New York, 1995. Translated by A. I. Dale.

**Maupertuis P. L. M.** (1756), La divination. *Œuvres*, t. 2. Lyon, 1756, pp. 298 – 306.

**Pascal B.** (1654, Latin), A la très illustre Académie Parisienne de Mathématique. *Oeuvr. Compl*. Paris, 1963, pp. 102 – 103, Latin and French.

## Chapter 5. Sale Prices of Chances and Probabilities. The Market of Chances and Games in General

**49.** If a lottery[1] offers a commercial object, each of its tickets representing an eventual right to own it can be in turn offered for sale, and its sale price is the price of the chance of the eventual right as ensured by the ticket. There is absolutely no reason to value one ticket more than another so that two people having *m* and *n* tickets respectively possess values in the ratio of *m/n*.

What we say about chances in a lottery usually represented by tickets can be equally applied to any kind of chances and it follows that when many people have eventual rights to a commercial object which can in turn be commercialized, the sell price of such an object is necessarily proportional to the respective probabilities of obtaining it. That consideration is not yet sufficient for establishing the absolute value of each chance, and it is actually clear that each person can determine its sale price as well as that of any other merchandize according to its particular convenience.

And just as the course of things usually dealt with in commerce, is established, so also a course is developed for chances which can become objects of everyday speculations. Therefore, the price of each chance is to the price of the thing to which the chance ensures a risky right as unity is to the total number of chances. And if the course assigns a lower price for each chance, the owner of the thing will not get the usual price and will not raffle it. On the contrary, if the course assigns a higher price for each chance, speculators will benefit by buying that or a similar thing for distributing its value among negotiable chances. The emerging competition will lead to the lowering of the course until it returns to the level which they exceeded for some time.

In this purely theoretical reasoning we have for the sake of greater simplicity abstracted ourselves from the invariably involved overheads and salaries. Instead of the actual state we substituted a fictitious state which is the closer to the real processes of commerce the freer are those processes, and this is how the laws of commercial equilibrium should be researched.

**50.** Because of a rather bizarre association of words the product obtained by multiplying the value of a thing in monetary units by the fraction expressing the mathematical probability of gaining it is called *mathematical expectation*. According to the above, it is the limit to which by the laws governing free trade invariably converges the sale price of chances owned by each who claims a thing or the sale price of the probability of its gain. If the risky right is obtained for moneys variations of the course are excluded and the mathematical expectation of each claimant becomes fixed as soon as the probabilities of gain of each of them is known just as it is in games of pure chances when the combinations can be calculated.

If the gamblers agree to stop playing they ought therefore to share the stakes proportionally to their probabilities of gaining since each chance represents an equal right to get them. The rule of mathematical expectation is therefore reduced to the *division rule* (§ 47) which became the occasion for the first researches of mathematical probability. If Peter bets on event A and Paul, on event B, their money whose sum forms the stakes should be respectively proportional to the mathematical probabilities of A and B. Indeed, if all the chances or equally possible random combinations which can produce A and B are enumerated, there will be no reason for betting on one of them rather than on another. The sum bet on A can be considered as a total of the sums bet on each of the chances bringing about that event, and the same is true for event B. When examined from this other point of view, the rule of mathematical expectation is confounded with the *division rule*[2].

**51.** Selling a thing for a fair price means selling it for a price established by a free competition of buyers and sellers, and just as well the mathematical expectation is therefore the *fair price* of chances or the limit to which that price approaches when the overheads of the transaction diminishes. If the demanded price is different, equitable conditions do not anymore govern the market of chances. The same happens to all other markets if one of the contractors benefits from the advantage of his position, − from the needs, passions or ignorance of the other contractor, and gives him in exchange less than would have been determined in the absence of any illusions by free competition.

It does not follow that the same price of one and the same thing suits everyone or that a certain expense is reasonable only because its price was not higher than the course. The *acceptable value* as opposed to the commercial price is evidently subordinated to the buyer's particular situation and fortune. It is impossible to measure and subject to calculation that price of things, − of chances and of all other commodities.

It is undoubtedly well understood that when a person buys a chance which is something uncertain he risks the more the larger is its given certain price relative to his fortune. Common sense also tells us that the importance of an amount of money decreases for the person whose fortune increases. Thus, for a worker who saved a thousand francs and risks a half of it in the game of *passe-dix* (§ 25), the 500 francs of possible gain cost less than the 500 that he risks[3]. That relative value of chances is called *moral expectation*, and various rules for evaluating it were proposed, all of them arbitrary and lacking real applications. Calculus should not be abused if desired that it preserves its authority over things situated in its domain. And in general we run the risk of discrediting logical argumentation (with calculus only being one of its branches), when it is transferred beyond the sphere of logical combinations.

**52.** Let us return to the lottery with moneys proposed as the prize and its tickets representing the corresponding chances of gaining them. If the government does not exercise a monopoly on such an enterprise for itself or for the concessionaire, the price of a ticket does not exceed the attached mathematical expectation more than justifiably needed for covering the expenses of management, of selling [the tickets] and the commercial interest on the necessarily involved capital.

Here, however, we evidently have a bad application of a part of the capitals and productive forces of a nation. First of all, that application is unproductive for the country or for the society of its citizens since some of them only gain what others lose. A prolonged conduct of such an enterprise tends to impoverish many people for enriching a few favoured by chance and thus to widen the inequality of fortunes beyond the difference that can be demanded by the natural laws of society. And, after all, these instant gains not at all compensating work incline to prodigality, luxury and unproductive expenses; they harm the society in a purely economic way and corrupt morals which is much more pernicious.

It would have been otherwise had the lottery been the only possible means for producing a useful thing, as publishing, for example, a luxurious book whose copies are too expensive to be sold but are the lottery's prizes. We see that under some circumstances speculating on chances for supplementing a productive speculation can become an advantageous and laudable means for making use of a part of the capitals and productive forces of a nation.

Suppose that [a group of] workers each aged thirty draw the same moneys out of their savings and deposit them in a joint account under the condition that after thirty years those moneys complete with interest will be shared by those still living. Such associations are called *tontines*,[4] and the chances there, similar to those in lotteries, can be useful if not directly productive if each worker sacrifices a certain sum without which he can painlessly manage while being able-bodied and thus ensures his subsistence in case he becomes old and unable to live by working. On the contrary, if he could have been able to ensure a decent living in old age just by economizing but preserving the fruits of his labour for his children, then his deposit would have only aimed at easily enriching himself or his family by trusting chances, − then the tontine should be condemned like lotteries, and for the same reason[5].

Economists have justly remarked that considerable salaries, attached to certain lofty functions or the large possible profits in certain professions for those who excel in them and managed to become eminent, act like a large random bonus, like a considerable prize offered for many but obtained by few. These random bonuses allow the maintenance of money payments for many public and private services on a more modest level and excite the activities which would have numbed. And in these various aspects they can in proper limits foster progress and the well-being of the social corps.

**53.** In general, randomness intervenes in all things of our world. In the economic life, each speculation has more or less the nature of a market of the fortuitous. All kinds of commercial affairs invariably include the buying and selling of chances. When possible and desirable to save a commercial speculation or a single private affair from its inherent fortuitous condition, the contract is called *assurance*, and we treat it in a special chapter. The contract of assurance is always favourable since it dissipates uncertainty restraining productive activity and generates free development. It extends the power that man had acquired by his free intelligence, his foresight over physical nature which only obeys the laws of fatality.

The addition of the fortuitous bonuses to speculation to which it does not necessarily belong is an operation inverse to insurance. We can consider it favourably or unfavourably depending on whether it is a useful auxiliary of a productive operation or an absorbent acting contrary to productive operations, becomes a fancy dress of sorts and degenerates into hullabaloo, into a pure game. It is easy to understand that between these extremes there can be innumerable nuances which can not become an object of precise determination.

**54.** We (§§ 49, 52) supposed that the lottery runs under conditions of free competition. However, if a similar enterprise becomes an object of legal monopoly ran by the government or tax-farmers, it can benefit the operator, the state revenue after covering all the expenses of the exploitation.

Pursuing fiscal aims, the government could have supported passions and disorders, moral or economic, occasioned by fortuitous speculations and at least would have benefited from the opinion that these passions and disorders are indestructible for attempting to gain a certain monetary advantage for profiting the political corps instead of conceding that advantage to private speculation.

It is useless to return to a problem treated many times for promoting order and morals and happily legally solved in France. In cases of embarrassment governments often resort to fortuitous bonuses for favouring its loans. And if the accumulation of capitals in a country on the way to prosperity renders this expedient useless, we see nevertheless that commercial deals concerning public funds, useful in themselves, serve as an excuse for organizing a vast market of chances whose misuse the governments neither can nor desire to suppress. It is not our intention to treat these political or financial problems which are only in a roundabout way connected with the mathematical theory of chances. The indicated suffices.

**55.** The entrepreneur of the fictitious lottery which we considered as typical, or the *banker*, does not gamble, they only distribute the chances. Usually, however, in public lotteries or similar markets of chances the banker plays against the *punters* or those whom he sells the chances. The punters do not gamble against each other through the banker, all play against him. The same random event leads to the gain of some punters and to the loss of others depending on their preference for some chances based on fortuitous circumstances, whims or *systems* of gambling.

For simplifying calculations, we assume that in each of a long number $m$ of drawings the punter stakes the same sum $a$ with probability $p$ of gaining $b$. If the game is fair, the stake is exactly equal to the mathematical expectation [of gain] $pb$ or only exceeds it for covering the banker's expenses. However, once the bank has a financial aim, or the monopolist's interest comes into play, the difference ($a - pb$) or the banker's *advantage* will amount to a notable portion of $a$. In the previous *Lottery of France* it was 1/6 of the stake in cases of gambling on one number out of 5, about 1/3 and almost 22/25 for gambling on 2 and 4 numbers.

Denote by $P$ the probability that the number of sets won by the punter is contained between $m(p - l)$ and $m(p + l)$. For large values of $m$ the value of $P$ only depends on the ratio (33.1)

$$t = l\sqrt{m/2p(1-p)}$$

and we adduce an appropriate table sufficient for all practical purposes.

Now, *P* is also the probability that the total sum gained by the punter is contained between $mb(p - l)$ and $mb(p + l)$. His stakes totalled *ma*, and if $mb(p + l) < ma$, or if $l < a/b - p$, *P* will be the probability that the punter's loss is contained between

$ma - mb(p - l)$ and $ma - mb(p + l)$.

If, on the contrary, $l > a/b - p$, *P* will be the probability that the punter's loss or gain do not exceed $[ma - mb(p - l)]$ and $[mb(p + l) - ma]$ respectively. Suppose for example[6] that $m = 3000$, $p = 1/18$ and the banker's advantage, as in the example above, 1/6. Then even money can be bet on the final loss of the punter to be contained between 373 and 627 times his stake, and 20,000 against 1, that his gain does not exceed 265 times, and his loss does not exceed 1265 times that stake. When the stake is unity, these limits are expressed by large numbers, but they become 15 times smaller if we assume that $b = 1$.

In lotteries understood in their proper sense, drawings follow each other very slowly so that the same punter can not repeat trials of the same randomness many thousand times, but in public games, on the contrary, the sets end so rapidly that such numbers do not present anything extraordinary. And in such games the banker only reserves for himself a very small advantage not to discourage the punters and at the same time the prompt repetition of the sets multiplies his benefits and assures them.

**56.** A person who usually plays games of pure chance against the first comer is like a punter for whom the public is the banker. In such a case he is however playing a fair game, i. e. his stake is equal to his mathematical expectation, and most often he plays when the chances are equal so that the probability of gain is 1/2. For a series of 3000 sets it is even money that his loss or gain will be less than 19 times his stake, and a bet of 20,000 to 1 on the loss or gain not to exceed 111 times his stake[7]. If the number of sets becomes 4, 9, 16, … times larger, those limits should be multiplied by 2, 3, 4, … The probable loss or gain always increases with the number of sets although in a much less rapid and invariably slowing progression.

If the same two gamblers are playing all the time against each other, and the game lasts indefinitely, that progression, however slow it becomes, will certainly ruin one of them. Depending on the ratio of the stakes to their fortunes, one of them will be likely ruined after a more or less large number of sets. It is, however, proper to remark that if that ratio is not much larger than usual, the number of sets leading with a considerable probability to the ruin of one of the gamblers is larger than practically possible.

**57.** Anyway, these calculations are based on the hypothesis that the capital of each gambler allows him to complete the number of sets which we denoted by *m*. This can always be supposed with respect to those who only play for amusement. However, the contrary can regrettably happen when the game is passionate, and then we ought to find out whether the chances of a final loss will not disappear as compared with those of the anticipated ruin. Calculation of the latter chances becomes here the more necessary since we usually have a vague idea about their influence on the gambler's fate and since it is believed necessary, with very laudable intentions, to exaggerate that influence[8]. We think that it is still better to keep to arithmetical rigour.

To offer a simple example, suppose that A and B play staking equal sums and having equal chances, and that the stake is 1/50 of A's capital. Suppose at first that A and B are equally rich or disposed of exactly the same capitals. Then there will be probability 0.8859 or a bet of almost 9 against 1 on A not being ruined at the 1000-th set and 0.4954 or almost 1:1 that he will be ruined not later than at the 10,000-th set.

Suppose now that B's capital is twice as large as A's. Then the former probability will not appreciably change but the latter will become 0.604; about 3 can be bet against 2 on A's ruin not later than at the 10,000-th set. The influence of the superiority of B's fortune becomes appreciable, but much less than usually thought. Finally, suppose that B's capital is infinite or inexhaustible. The probability of A's ruin not later than in the 1000-th set still remains appreciably the same, but the probability of his ruin not later than in the 10,000-th set will become 0.617 which is only a bit higher than before.

If A's capital increases twice, three or four times, the number of sets necessary for his anticipated ruin with the same probabilities should be increased 4, 9, 16 times[9].

**58.** As a summary, the calculations confirm that, as indicated by common sense, when unequally rich gamblers play high, the richer of them, other things being equal, is at advantage since he can longer sustain misfortune. However, at the same time calculation shows that that advantage is much feebler than people tend to believe and insensible if the stake of a gambler in each set is not a considerable fraction of his capital and the number of sets is not very large.

The usual superiority of the richer gambler, if shown by well ascertained observations can also be explained otherwise. The richer gambler less feels losses and freer exercises his mind, whereas his opponent despaired by the reverses of fortune usually begins to increase ever more his stakes which suffices for ruining him in a much less number of sets. It is proper to remark here that all the calculations above supposed that the stakes were invariably constant.

The disadvantage of the gambler playing a large number of sets with a richer adversary becomes more justified when he plays the same number of sets against the first comers. Indeed, this is tantamount to playing with an infinitely rich opponent who can ruin him without risking to be ruined himself.

The capital of the tax-farmer of a gambling house is enormous as compared with those possessed by the punters. Therefore, each punter taken alone is at a disadvantage, but it does not follow that that inequality of position suffices for assuring the banker's benefit. If all the other conditions of the game are equal, the banker can equally lose or gain. It is wrong to see here a contradiction to what was said about the disadvantage of punters. Although it is more probable that each individual punter will be ruined, it does not follow that it is more probable that all of them will be ruined or that they will not gain.

**59.** The profits of an entrepreneur of public games[11] rest on a much more solid base than the inequality of the conditions of playing between the punters and him. We saw (§ 55) that in an unequal game the mean loss of a gambler being at a disadvantage increases proportionally to the number of sets whereas the interval between the limits in which the loss oscillates increases proportionally to the square root of that number. Therefore, a moment will arrive when the mean loss imposed by the rules of the game and the gambler's disadvantage becomes incomparably larger than the variations of that loss occasioned by the anomalies of chance. When for example the mean loss is counted by the million those variations are counted by the thousand, since a thousand is the square root of a million.

In the public games sets rapidly succeed one another and when many punters play at the same time betting different chances and varying their stakes capriciously or systematically, it is impossible, at least by experience, to calculate the number of sets after which there will be a certain probability that the accidental variations of the banker's gain remain contained within some limits. This precision is however barely needed. Suffice it that the theory and experience agree in ascertaining the banker's mean benefit resulting from the rules of the game.

The same remark all the more applies to lotteries. The number of yearly drawings is not considerable, and when considering them as sets, it seems that a great number of years should pass for the lottery to rely certainly on benefits. On the other hand, if all the stakes are equal one to another; if the chances were of the same kind; and if each punter selected the numbers accidentally, – then each drawing can be considered as providing as many sets as there were stakes, and this returns us to the very simple example of § 55.

However, all these various assumptions are inexact. Certain prejudices prevail among regulars of lotteries[12] and act in the same way on many stakes preventing a sensibly equal distribution of the chosen numbers or their combinations which takes place when the selection only depends on irregular and accidental causes. Only experience can measure the effect of those prejudices and only it can therefore indicate the number of drawings and stakes sufficient for assuring the lottery's benefit.

The official table of the results of the *Lottery of France* shows that from year VI [1793] to 1835 inclusive they varied with the increase and decrease of the [number of] stakes and depended on the unequal distribution of the different chances. However, that inequality much greater than it would have been in case of purely random causes never threatened the Treasury by loss. Experience also convincingly showed the administration of the Lottery the possible influence of those prejudices at each drawing, and it never exercised its right to *withdraw* the overloaded numbers.

**60.** The feebler was the probability of gain corresponding to the different chances of the previous lottery (when gamblers selected 1, 2, … numbers) the more advantageous it was for its administration since it raised the price of the tickets ever higher than the mathematical expectation [of gain]. It thus heavier taxed greater greed, but it also understood the need to ensure greater security in cases of less frequently repeated chances which can, although very unlikely, notably reduce the lottery's reserve. The gain on the *quine*, i. e. on a gamble on 5 numbers, was only a million times more than the stake whereas its probability was 1/43,949,268 so that the administration's benefit amounted to 42/43 [43/44] of the stake.

Nevertheless, the administration finally suppressed that version of gambling either to save itself from troubles or because the *quine* was played too rarely for the benefit of speculating on that chance to justify the complications of reckoning. It is well known that there should be a limit of the smallness of negotiable chances or probabilities. A random extraction of the one white ball contained in an urn among $10^8$ or $10^9$ balls is so unlikely that no one will wish to speculate on it. And if someone accidentally makes such an attempt, it will be an exception too rare for that chance to acquire a current commercial value or to be included in the price list of a company enjoying a monopoly on chances.

**61.** These considerations will provide a most natural solution of an entertaining problem called the *Petersburg game*[13]. By its captivating form it resembles the celebrated sophisms of ancient Greece. It can be formulated thus. Peter and Paul play *passe-dix* (§ 25) under the condition that Peter pays Paul 1 franc if the number of points exceeds 10 at the first throw, 2 francs, if that only occurs on the second throw, 4 francs, if only on the third throw, etc. The game does not end until that event happens and it is required to determine the expectation of Paul or how much he ought to stake in advance.

According to the fundamental notions of the calculus, Paul has probabilities 1/2, 1/4, 1/8, … of gaining 1, 2, 4, … francs depending on Peter throwing more than 10 points in the first, the second, the third, … toss. The value of his expectation is the sum of these random gains multiplied respectively by the corresponding probabilities. However, each of these products is 0.5 francs so Paul should stake 50 francs if the game will necessarily end at the 100th throw; 500, if at the 1000th; and more than any assignable sum, or an infinite amount if the gamblers agree to continue playing as long as necessary for Peter to throw more than 10 points. And still, as is remarked, no sensible person will risk here not only an infinite sum (which no one has) but even any considerable part of his fortune.

For solving this paradox, most geometers insert their hypotheses about *moral expectation* (§ 51) according to which the *useful value* of a sum of money increases less rapidly than its nominal value or even, as some of them hold, does not increase above a certain limit. Those explanations seem to us too arbitrary for adopting them.

Poisson very simply remarked that Peter can not pay more than he has and if he possessed 50 *mln* francs, an exorbitant sum for an individual, he will not be able loyally to commit himself to continue playing beyond the 26th throw. Indeed, at the 27th throw his debt to Paul in case of loss, $2^{26}$ = 67,108,864 francs, will exceed his fortune. On the other hand, Paul, if knowing Peter's fortune, will not commit himself to the same and will not risk more than 13 francs. Supposing that the number of throws is not restricted, we find that, although he can not at all receive from Peter more than 50 *mln* francs, the mathematical value of his expectation is not more than 13.5 francs.

    This remark does not however get to the bottom of the difficulty. The proper value of a thing should not be confused with the relative value which results from the degree of the debtor's solvency.

    Suppose that a public lottery is organized under the same conditions. Peter is substituted by a blind instrument throwing dice; the lottery's administration issues tickets numbered 1 ensuring their holders 1 franc if 10 points are exceeded on the first throw; tickets numbered 2 ensure 2 francs, if the same is only achieved on the second throw, etc. Because of the monopoly enjoyed by the lottery, it can price the tickets numbered 1 at higher than 0.5 francs and sell them easily. Similar pricing of tickets numbered 2, 3, … will be successful, but finally a number arrives for which there is no buyer, or they are found so rarely that the administration suppresses that chance. And still the game will be fair because the administration's solvency guaranteed by the state, can not be doubted. Peter's lot as described here initially, is the same as the punter's who bought a ticket of each kind.

    **62.** Systematic gamblers, i. e. those who know a system of gambling ensuring them, as they believe, a benefit or at least preventing their loss, are often met in games played in societies as well as in public games. Their stakes follow some progression or they prescribe rules for entering a game and leaving it. The bounds of our work do not allow a detailed description of these systems which can vary most widely. Suffice it to express as a mathematical truth which immediately follows from definitions and which it is easy to show by example of some system. A gambler, whichever system he applies, can not ensure a probability of 100 to 1 of gaining a franc in any fair game without running the risk measured by probability of 1 to 100 of losing 100 francs.

    The products obtained by multiplying the possible gain and loss by their probabilities should always remain strictly equal. If the game is not fair, no method of playing can abolish the inequality of the conditions of the opposing parties. In each case there is a product of two factors whose value invariably results from the conditions of the game which the gambler can not change by his method. However, he certainly can increase one of these factors by his method to the detriment of the other. He can choose to decrease or increase the eventual gain or loss by proportionally increasing or decreasing the probability of either. Just the same, the vis viva expended by various machines can serve to double the distance and halve the [carried] mass or halve the distance and double the mass […]**14**. To search for a machine which creates rather than absorbs the vis viva means to share the chimera of the inventors of perpetual motion. […]

The overheads of the game represent the wasteful absorption of the vis viva by the machine, and any system of gambling can be considered as a machine by whose aid the gambler can vary the two elements according to his views and motives, but their product never varies.

**63.** It is proper to indicate here that people otherwise sensible are prone to become victims of an illusion. Each of us has a confused feeling that in a long succession of events the anomalies of chance should likely be almost compensated. Therefore, it is imagined that if an event lacking many chances is reproduced more often during a certain period it will occur less often in the next one, as though the independence of successive events does not exclude any influence of the passed chances on the future chances. However, imagination has difficulties in seeing that the laws of chance are consequences of mathematical laws governing combinations; it is always tempted to provide chance with a substantial and productive property with its own energy and an aim of sorts. Suffice it to indicate such illusions for safeguarding any reasonable person against them.

If, in a long run of random trials, the ratio of the number of events A to that of the contrary event B appreciably differs from the ratio of their probabilities, it reveals some irregularity in the pertinent chance mechanism or, more generally, the existence of a cause under whose influence the allegedly equally possible combinations serving for calculating the probabilities of those events were not really such.

If, for example, the double-six occurred a thousand times in ten thousand tosses of two dice (§ 6), the structure of the dice was certainly irregular. It can also be possible that the gambler, who puts the dice in their box, throws them either habitually or deftly in such a manner that the probability of the arrival of the double-six is heightened. Instead of being able to derive the value of that probability in advance, it should be determined by experience, as we will describe below.

If a game depends on skill and chance with the same chances for both gamblers, and one of them wins much more than a half of a long series of sets, it was certainly due to the superiority of his skill. Conversely, in a long run superiority of skill should prevail over irregularity of chance. What we attribute to that superiority of the gambler's skill and sangfroid, others, if they desire, impute to a mysterious fatality which pursues some people and seems to be pleased to favour the rest.

This belief is one of those rooted in the human heart. It results from a confused feeling of an inconceivable supernatural system not yielding to reasonable discussion. It contributes at least as much as greed and ambition to maintaining an inclination to venturesome enterprises and fortuitous speculations. The history of that belief, its origin and effects belong to the field of moralists and psychologists and have no place in an exposition of a mathematical subject.

**Notes**
**1.** A law of 1836 suppressed all but charitable lotteries. [B. B.]

**2.** Concerning that *other point of view*, Bru refers to De Moivre (1718/1756, p. 3). There, De Moivre first applied the word *expectation* irrespective of its definition which he added a few lines afterwards.

**3.** Bru refers here to Buffon (1777, § 13). See also Note 13.

**4.** Tontines were named after an Italian banker Lorenzo Tonti (1630 – 1695). Contrary to Cournot's explanation, they were groups of annuitants of about the same age considered by the entrepreneurs (usually, the appropriate state) as single entities. A tontine distributed yearly payments among its still living members and those living long enough came to enjoy considerable moneys. Bru noted that tontines were *ruinous* for the state (apparently only because their financial conditions had not been properly considered) and suppressed in 1740 but reappeared later and existed even in the mid-19[th] century.

**5.** Bru noted that in 1821 a commission of the Paris Academy of Sciences (Fourier, Poisson, Lacroix) negatively reported on the establishment of a tontine and that it had been felt that tontines harmed the development of insurance.

**6.** This contradicts the pertinent result of § 34. [B, B,]

**7.** Bru stated that Cournot's estimates should have been 30 and 222.

**8.** Bru referred here to Ampère and Lacroix.

**9.** Suppose that gambler A plays a fair game with opponent B whose fortune is reputed infinite and each time stakes $1/\alpha$ of his capital. Denote by $\Pi$ the probability that he will be ruined not later than at the *n*-th set and let

$$t = \frac{\alpha}{\sqrt{2n}}.$$

We will have almost exactly[10]

$$1 - \Pi = \frac{2}{\sqrt{\pi}} \int_0^t \exp(-t^2) dt + \frac{1}{\sqrt{\pi}} \frac{t \exp(-t^2)\alpha}{2n} [1 - \frac{2}{3n}].$$

In most cases the second term can be neglected (§ 33, Note 3). The value of the first term is provided in the table adduced at the end of this work. A. A. C.

**10.** Bru referred to Laplace and Fieller (1931).

**11.** Public games were allowed in Paris and spa towns from 1818 to 1838. [B. B.]

**12.** Many authors beginning with Montmort (1708/1713, pp. vi − vii) mentioned these prejudices.

**13.** It was Nikolaus Bernoulli (Montmort 17108/1713, pp. 280 – 285) who invented the Petersburg game. Below, in the same section, Cournot stated that *many geometers* studied it. The main author among them was Daniel Bernoulli (1738) who published his memoir in Petersburg. At the end of the 19[th] century economists, issuing from Bernoulli's ideas, began to develop the theory of marginal utility.

**14.** Bru noted that Mises applied the same comparison with the vis viva for justifying one of his axioms.

## Bibliography

**Bernoulli D.** (1738, Latin), Exposition of a new theory of the measurement of risk. *Econometrica*, vol. 22, 1954, pp. 23 – 36.

**Buffon G. L. L.** (1777), Essai d'arithmétique morale. *Oeuvr. Phil*. Paris, 1954, pp. 456 – 488. English translation: Isf.lu/eng/Research/Working-Papers/2010

**De Moivre A.** (1718), *Doctrine of Chances*. London, 1738, 1756. New York, 1967.

**Fieller E. C.** (1931), The duration of play. *Biometrika*, vol. 22, pp. 377 – 404.

**Mises R.** (1928), *Probability, Statistics and Truth*. New York, 1981.

**Montmort P. R.** (1708), *Essay d'analyse sur les jeux de hazard*. Paris, 1713, New York, 1980.

## Chapter 6. The Laws of Probability. Mean Values and Medians

**64.** We know that in natural phenomena the number of combinations or chances which can occur by a fortuitous conjunction of independent causes is usually infinite (§ 15). Therefore, ordinarily a magnitude whose determination depends on these conjunctions can without discontinuities take all values contained between certain limits or even be unbounded either above or below. Since these values are infinite in number, the probability of each is infinitely low. It is physically impossible for a bet on some precise value not to be lost. Generally, however, those infinitely low probabilities are not the same at all. Between themselves, they preserve certain finite and assignable ratios which only in particular cases are reduced to unity.

This was indeed borne in mind when the problem of § 16 was treated. Its formulation can be modified for generalizing it as desired. Suppose that we have a flat figure limited by segment *AB* of the *x*-axis, perpendiculars *Aa* and *Bb* and some curve *ab*. A sphere is randomly thrown on that figure so that the abscissa of the point of contact can take all values between *OA = a,* and *OB = b* where *O* is the origin of coordinates. It is requested to determine the probability that that abscissa takes a certain intermediate value *OI = x*. In other words, the probability that the point of contact is situated on the perpendicular *Ii* to *AB* between *Aa* and *Bb*. That probability is infinitely low because the point can just as well be on infinitely many other perpendiculars to *AB* situated between *Aa* and *Bb*. However, when comparing the probabilities of the fall of the sphere on *Ii* and another perpendicular *Hh*, it follows from the formulation itself of the problem that the ratio of these probabilities equals to that of the lengths of *Ii* and *Hh* so that these infinitely low probabilities are not at all equal. […]

The same fact can be presented otherwise. Choose points $I_1$ and $I_2$ close to each other and situated on different sides from *I*. The probability of the point of contact of the sphere and the figure being within the curvilinear trapezoid $I_1I_2i_2i_1$, or that the probability of the abscissa of that point being between $OI_1$ and $OI_2$ is evidently expressed by the ratio of the areas of that trapezium and of the entire figure. […]

**65.** In general, whichever are the conditions of randomness that assign one of the infinitely many values contained within limits *a* and *b* to a certain magnitude *x*, that case can be likened to the abovementioned. The curve *ab* will represent the *law of probabilities*[1] of the different values of *x* situated between those limits. The ratio of *Ii* to *Hh* is that of the infinitely low probabilities of *x* taking exactly the values *OI* and *OH*. [Here, Cournot introduces the terms abscissa, ordinate (its *function, as the geometers say*, and axes of coordinates.]

For the sake of brevity, we call the *curve of probability*[2] that curve which appropriately represents, as we have explained above, just like in § 31, the law of probabilities of different values of a variable magnitude.

**66.** The curve of probability extends to infinity, if, for example, all the values of magnitude *x* from 0 to ∞ are really possible. However, for the notion of probability to make sense it is necessary for the total area of the figure between the coordinate axes and the curve *ab* (Fig. 1) extended to infinity in both directions to preserve a finite value.

This condition supposes that a current ordinate will decrease and become less than any limit [any arbitrarily small magnitude]. But this is not sufficient for ensuring the existence of the first condition, at least when considering it abstractly and purely mathematically.

In the physical reality it always happens that for a certain value of the abscissa *OB* the corresponding ordinate becomes so small that the portion of all the area beyond it can be neglected without an appreciable error. The values beyond *OB*, although strictly speaking are possible, occur so rarely that they can be disregarded as physically impossible which they always are. Thus, the probability of living until 110, 120, 130 years is not exactly zero since some people exceeded them.

Probably there even are no mathematical conditions or other of mathematical rigour which determine the limit of an absolutely unattainable age. Nevertheless, in all problems belonging to the domain of the calculus of probabilities it is quite admissible to treat probabilities of living until 110 years or more as zeros, and it can be assumed that no one had lived 200 years. If the appropriate magnitude randomly takes negative just as positive values, the curve of probabilities can extend to infinity in both directions. […]

**67.** Imagine now that after tracing the curve of probabilities we separate the interval *AB* between the extreme values of the abscissa in *i* equal parts $AA_1, A_1A_2, \ldots$ by equidistant ordinates $A_1a_1, A_2a_2, \ldots$ Denote by $\Omega$ the total area *Abba* and by $w_1, w_2, \ldots$ the partial areas $AA_1a_1a, A_1A_2a_2a_1, \ldots$ Then, as *i* becomes ever larger, the quotient

$$\frac{w_1 OA_1 + w_2 OA_2 + \ldots}{\Omega}$$

ever closer approaches a certain value *M* represented by segment *OG* situated between *OA* and *OB*. Thus defined, *M* is the *mean* of all the values that *x* can randomly take within *OA* and *OB*. Because of its proper probabilities, each of these particular values contributes to the formation of the mean *M* which should appreciably coincide with the arithmetic mean of the particular values provided by a very large number *N* of fortuitous trials.

Separate the total series of these particular values in *i* partial series. The first of them, $n_1$ in number, being between *OA* and $OA_1$; the second, $n_2$ in number, being between $OA_1$ and $OA_2$; … Since intervals $AA_1, A_1 A_2, \ldots$ are short, the mean μ of the total series is appreciably equal to

$(1/N)[n_1 OA_1 + n_2 OA_2 + \ldots]$.

On the other hand, since each partial series is supposed to contain a very large number of particular values, we have approximately

$n_1/N = w_1/\Omega, n_2/N = w_2/\Omega, \ldots$

Therefore, the fixed value which the mean μ ever closer approaches when the numbers $i$, $n_1$, $n_2$, ... and the more so, $N$ ever increase, will indeed be $M$ as defined above and represented by the segment $OG$.

According to the elementary notions of statics, if the figure *Abba* is one of the faces of a ponderable plate of uniform thickness and density, its centre of gravity will be on ordinate $Gg$. And if segment $AB$ represents a ponderable bar of uniform thickness and density changing as the ordinates of the curve *ab*, point $G$ will be its centre of gravity.

**68.** Suppose that $OI$ is the abscissa whose corresponding ordinate $Ii$ divides the total area Ω in two equal parts. Then $OI$ is what we call the *median* value of $x$. Gamblers betting on $x$ being smaller or larger than $OI$ have equal chances of gaining. For a very large number of randomly determined values of $x$, the quotient of those, larger (or smaller) than $OI$, to the total number of trials very little differs from 1/2.

As I have already remarked (§ 34), until now, authors usually but very improperly call this median value *probable*. Generally, this value does not coincide with the one to which corresponds the *maximal* ordinate of the curve *ab* and can not therefore be considered more probable than the others. Nothing prevents it from corresponding to the *minimal* ordinate of that curve and even to a zero ordinate in which case the median will not anymore be one of the infinitely many values randomly taken by $x$.

If the ordinates of the curve invariably increase from $A$ to $B$, the median exceeds the mean value and vice versa. If the curve is symmetric with respect to some ordinate $Gg$, the mean value and the median coincide with the abscissa $OG$ which is also the half-sum of the extreme values, and to which the largest or the smallest ordinate ordinarily correspond.

**69.** It is not difficult to understand that, depending on the interval between the limits within which oscillates the random magnitude, and on the form of the curve representing the law of probabilities of [its] different values, the mean value μ determined by a large number of trials should more or less rapidly tend to the *absolute* mean $M$ as defined in § 67. Then, depending on the considered case, the number of the trials should be more or less large for securing a given probability that the anomalies of chance will only result in deviations contained within the assigned limits.

*Modulus of convergence*[3] or simply *modulus,* as I call it, is the number which measures in each particular case the rapidity with which the means provided by the trials converge to the absolute mean. The value of that modulus is obtained in advance by the rules of the integral calculus after assigning the type of the function or the form of the curve representing the law of probabilities.

Denote by $g$ the modulus, by $m$ the number of trials and by $P$, the probability that the mean μ determined by these trials does not deviate in either direction from the abscissa of the absolute mean $M$ by more than $l$. For large values of $m$ $P$ only depends on number[4]

$$t = lg\sqrt{m} \qquad (69.1)$$

so that for a constant *t* (numbers *l, m* and *g* can vary) probability *P* does not vary either. Note that *P* is the same function of *t* as in § 33 for which we provide a table. For determined values of *g* and *m* the value of *l* corresponding to $t = 0.476937$ and $P = 1/2$ is that which we call the median deviation. The value approximately equal to 6 such deviations corresponds to $t = 2.87$ and $P = 19,999/20,000$ and can be considered as the extreme limit of deviations.

If the modulus of convergence is constant, the limit of the deviations for the same values of *P* varies inversely proportional to the square root of the number of trials. And if that latter number is constant, the limit of the deviations varies inversely proportional to the modulus[5].

**70.** *First example*. Points are randomly distributed on a segment 1 *m* long. In § 14 we similarly supposed that they are the points of contact of a ball thrown at random with a billiard board in such a way that there is no reason to suppose that the ball hits one point rather than another. The distance of the points of contact to one of the extremities of the segment is a magnitude that with the same probability can randomly take all the values from 0 to 1 *m*. The curve of probabilities becomes a segment parallel to the *x*-axis. The mean value coinciding in this case with the median value is 1/2 *m* and the modulus of convergence $\sqrt{6} = 2.4495$. It follows that in a series of 1000 trials the median value of the deviations is 0.006519 *m*, or somewhat larger than 6 *mm*. 20,000 to 1 can be bet on the deviation not to exceed 36 *mm*. These limits of deviations can be halved if there will be 4000 trials.

*Second example*. Points are randomly distributed on a unit circle as was indicated in § 16. The distance of a point to the centre of the circle is a magnitude that can randomly take all values from 0 to 1 *m*, but they are unequally probable. The probability heightens from one value to another proportionally to the distance of the point from the centre of the circle. The curve of probabilities becomes a segment not anymore parallel, but at an angle, whose tangent is 2, to the *x*-axis and passes through the origin *O*.

The mean value is 2/3, the median value that should exceed it (§ 68) is $1/\sqrt{2} = 0.7071$ *m*, and the value of the modulus is the natural number 3. The limits of the deviations as compared with those in the preceding example are reduced in the ratio 300:245.

*Third example*. Points are randomly distributed in space within a unit sphere. The distance of a point to its centre is a magnitude *x* which once more randomly takes all the values within 0 and 1 *m* and their probabilities are proportional to the surface of a sphere of radius *x*, or to $x^2$. The curve of probabilities becomes a parabola with its vertex in point *O*, *OY* as its axis and the focus at distance 3/4 *m* from the vertex. The mean value is 3/4 *m*, the median value, $1/\sqrt{2}$ [$1/\sqrt[3]{2}$] $= 0.7937$ *m*, and the value of the modulus, $2\sqrt{30}/3 = 3.5683$ [3.6533].

It follows that the limits of the deviations as compared with the first example are reduced approximately in the ratio 357/245, a little less than 3/2.

**71.** *Fourth example*. Imagine a globe with the poles, an equator and meridians and parallels just like on terrestrial or celestial globes. It is randomly thrown and lands on the floor and its point of contact with the floor is carefully marked. Each of these points has a longitude and latitude, the former varying from 0 to 360°. The latter, if the essence of the problem allows us, as below, to disregard its sign, is regarded as a magnitude susceptible to vary fortuitously from 0 to 90°. Supposing that the globe is really spherical and homogeneous, there will be no reason for it to stop at one rather than another region of its surface so that each longitude will be equally probable with its mean being 180°. If assuming the length of a meridian as unity, the value of the modulus will be the same as in the first example. Therefore, after a series of 1000 trials the median value of the deviation will be

360°×0.006159 = 2°13′2″.064.

For reducing that value to 1° we should have more than 4000 trials.

With latitudes, it is otherwise. Each of its values is the less probable the nearer it is to 90°, or the closer is the appropriate point to one of the poles. Indeed, two latitudinal circles very close to each other, with the latitudinal difference being 1′, say, circumscribe a zone on the spherical surface whose area is proportional to the cosine of the latitude.

The mean value of the latitude is equal to the complement of the arc of the same length as its radius, or 32°42′15″.2. The median value, 30°, is in this case smaller than the mean value (§ 68). The value of the modulus is 2.9518 if a quarter of the circumference, or the distance between the limits within which the latitudes can oscillate, is unity. And then, having a series of 1000 trials, the median value of the deviation will be

90°×0.005111 = 0°27′35″.964.

The limits of the deviations as compared with those for the mean longitude, are reduced not only because each particular value only varies in a four times shorter interval, but also since the modulus of convergence is increased.

**72.** It is easy to justify these results of calculations by a very simple geometric consideration. Actually, when a point on a plane or in space experiences a displacement *z* measured along a certain straight line, its distance from a fixed point varies less than *z* except when the displacement occurs along a radius of a circle or a sphere with that fixed point as its centre. It follows that in general the influence of chance inequalities in the distribution of points on the mean distance from a fixed point should lessen when the distribution on a straight line is replaced by a distribution on a plane, and then in space. The same occurs in the example of § 71: the nearer is a point to the poles of a sphere, the easier its light displacement alters the longitude without noticeably influencing the latitude. However, calculations are indispensable for precisely measuring the effects which are only vaguely discerned by geometric considerations.

**73.** A certain magnitude *u* can be in a known connection with magnitude *x* which takes a series of various fortuitous values in a succession of trials. Therefore, *u* can also be considered as indirectly obtaining the same number of fortuitous determinations. And, owing to the known connection between $u$ and $x$, the law of probabilities of $u$ can be derived from that of $x$. The same holds for the limits of that function's oscillations, its mean and median values, and the modulus of convergence.

Depending on that modulus being larger or smaller than the modulus of $x$, the influence of the anomalies of chance on the difference between the absolute mean and that provided by a series of fortuitous trials will be weakened or strengthened when passing from $x$ to its function $u$. If the variations of $u$ and $x$ are proportional, or if $u = b + cx$, the mean value of $u$ will correspond to the mean value of $x$, but in general this is not so. Let, for example, $u = x^2$. Then the mean value of $u$ or of $x^2$ will always exceed the square of the mean value of $x$. It is the difference between these two magnitudes on which depends the value of the modulus of convergence for the variable $x$ (§ 69, Note 5). On the contrary, since $u$ invariably increases with $x$, the median value of $u$ will necessarily correspond with the median value of $x$.

**74.** If (Fig. 2) function $u$ depends, according to a known law, on many magnitudes $x, y, z, \ldots$, each independently one from another taking a fortuitous value at the same trial, we derive the law of probabilities of magnitude $u$ from the laws of probabilities of those independent magnitudes. Also determined, although not as simple as before, are the limits within which it oscillates and its modulus of convergence.

In the problem considered in § 14, the function $u$ is the difference (without taking into account its sign) between two magnitudes, $x$ and $y$, each taking fortuitous values at each trial, or at each pair of joint trials. Both $x$ and $y$ can indifferently take all values between 0 and 1 and the function $u$ can also take them although they are not equally probable. Draw a unit square $OACB$; each point such as $m$ within the square or on its edge will have coordinates $x, y$ and to each of such infinitely many points there will correspond equally probable hypotheses about the system of fortuitous values of those coordinates. Let

$OP = OQ = a, AQ' = AP, BP' = BQ$.

All the points of the square for which the function $u$ takes the particular value $a$, will be situated on one of those equal and parallel segments $PQ'$ and $QP'$. Therefore, the probability of function $u$ taking a particular value $a$ between 0 and 1 is proportional to the lengths of those segments or to $(1 - a)$, and the median value of $u$ is 0.2928 …, mean value 1/3 and its modulus will be 3 (§ 70). The probability that the fortuitous value of $u$ will not be less than $a$ is equal to the ratio of the sum of the areas of the right triangles $APQ'$ and $BQP'$ to the area of the square and therefore equal to $(1 - a)^2$. If $a = 0.3$, that probability will be 0.49, as we indicated in § 14.

Suppose now that the magnitudes $x$ and $y$ can take, as previously, all the values from 0 to 1 but that they are not equally probable. For example, their probabilities, just like in the geometric problem of § 16, are respectively proportional to $(1 - x)$ and $(1 - y)$. Imagine that the square *OACB* is a ponderable plate whose density at each point is proportional to the product $(1 - x)(1 - y)$. The ratio of the sum of weights of triangles *APQ′* and *BQP′* to the weight of the square is the probability that the fortuitous value of *u* will not be less than *a*. That ratio, by the rules of the integral calculus, is $(1 - a)^3[1 + a/3]$ or 0.3773 if $a = 0.3$.

Let us also consider the case of *u* being the sum of two magnitudes indifferently taking all the values between 0 and 1. Then *u* can take all unequally probable values from 0 to 2. Trace, as in the previous example, square *OACB* and let $OP = OQ = a$. All the points of the square for which the function *u* takes the particular value *a* are situated on segment *PQ* whose length is proportional to the probability that *u* will take that value at least if *a* remains less than unity, or less than the side of *OACB*. If $a > 1$, or if points *P* and *Q* are situated beyond *A* and *B*, the probability of *a* becomes proportional not to the length of *PQ*, but of *P′Q′* which is part of *PQ* intercepted by the two other sides of the square. It follows that the function measuring the probability of each value of *u* will be equal to *u* when its values are contained within 0 and 1, and to $(2 - u)$ when within 1 and 2. At the value $u = 1$ that function experiences what we will call discontinuity of the second order represented by an ordinate of line *OGB* broken at *g* and formed by two equal sides of a right isosceles triangle whose height *Gg* is unity and base is twice longer.

And now we find that the modulus of convergence for the function $u = x + y$ is $\sqrt{3}$. For each magnitude *x* and *y* it is $\sqrt{6}$ (§ 70). The limits of the deviations proper for *x* and *y* should be increased in the ratio $\sqrt{2}:1$ for obtaining the limits of the deviations corresponding to *u*. The values of *u* oscillate in a twice wider interval as compared with that between the extreme values of *x* or *y*, but the limits of its deviations do not at all reach twice those limits for *x* and *y*.

Constructions in space similar to those discussed in a plane serve to solve the possible problems about the law of probabilities of a function *u* of three variables *x, y, z* with known laws of probabilities and fortuitously taking values independently one from another. In general, if *u* is a linear function of a certain number of variables,

$u = b + c_1 x + c_2 y + c_3 z + \ldots,$

where $b, c_1, c_2, c_3, \ldots$ are positive or negative constants, the mean value of the function *u* coincides with what is obtained if *x, y, z,* … are replaced by their mean values $M_1, M_2, M_3, \ldots$:

$M = b + c_1 M_1 + c_2 M_2 + c_3 M_3 + \ldots$

If *u* ceases to be a linear function, its mean value will not in general be obtainable in that way. This is what we had remarked in § 73 about functions of one single variable. The case of linear functions merits, however, special attention because any function can be artificially linearized, as proved in pure mathematics, if the magnitudes on which they depend only experience very small variations.

## Notes

**1.** Bru notes that Cournot's Chapter 6 had been repeated (repris) by all authors of the 19th, and some of the 20th century. Now, *repeated* is certainly too strong, and why was not this remark inserted at the beginning of the chapter?

**2.** The term *curve of probability* is due to Laplace [B. B.]

**3.** This modulus occurred in Gauss (1809, § 177) as parameter *h* of the normal law:

$$\varphi(\Delta) = \frac{h}{\sqrt{\pi}} \exp(-h^2 \Delta^2) .$$

In § 178 Gauss called it the measure of precision. In 1823, Gauss, having rejected the universality of the normal law, introduced the variance as $1/2h^2$. Cournot ignored Gauss, see Note 8 to Chapter 7 and Note 7 to Chapter 11.

**4.** Formula (69.1) which Cournot denoted by (L), is due to Laplace [B. B.]

**5.** Suppose that *a* and *b* are the inferior and superior limits of the possible values of *x* and *fx* expresses the law of possibility [!] of those different values. Then the modulus of convergence will be[6]

$$g = 1 \div \sqrt{2[\int_a^b x^2 fx dx - (\int_a^b xfx dx)^2 ]}.$$

Function *fx* is necessarily subordinated to the condition

$$\int_a^b fx dx = 1.$$

Three integrals

$$\int_a^b fx dx, \ \int_a^b xfx dx, \ \int_a^b x^2 fx dx \qquad (69.2)$$

express respectively: the area of the curve whose ordinates represent the law of probabilities; the mean value *M* of the variable *x*; and the mean value of $x^2$. And $1/g^2$ is equal to twice the difference between the mean value of the square and the square of the mean value. In addition,

$$\frac{1}{g^2} = \int_a^b \int_a^b (x - x')^2 fx fx' dx dx'.$$

This means that $1/g^2$ is the mean of all the infinite number of values which the square of the difference between two values fortuitously assigned to variable *x* can take according to its law of probabilities.

If *x*, instead of taking all the infinite values between *a* and *b* only took a finite number of different values $x_1, x_2, \ldots, x_n$ with probabilities $p_1, p_2, \ldots, p_n$, the integrals included in the expression of *g* will be replaced by sums:

$$g = 1 \div \sqrt{2[p_1 x_1^2 + p_2 x_2^2 + \ldots + p_n x_n^2 - (p_1 x_1 + p_2 x_2 + \ldots + p_n x_n)^2]}$$

and if $p_1 x_1 + p_2 x_2 + \ldots + p_n x_n = M$,

$$g = 1 \div \sqrt{2[p_1(x_1 - M)^2 + p_2(x_2 - M)^2 + \ldots + p_n(x_n - M)^2]} =$$

$$1 \div \sqrt{2[p_1 p_2 (x_1 - x_2)^2 + p_1 p_3 (x_1 - x_3)^2 + \ldots + p_2 p_3 (x_2 - x_3)^2 + \ldots]}.$$

It is seen that that value of *g* becomes minimal and equal to $\sqrt{2}/(b-a) = \sqrt{2}$ if $(b-a) = 1$ and

$x_1 = a, x_n = b, p_1 = 1/2, p_2 = 0, \ldots, p_{n-1} = 1/2$.

The minimal value of $\sqrt{1/2p(1-p)}$ in formula (33.1) is also $\sqrt{2}$ corresponding to $p = 1/2$. Therefore, when fortuitously determining a large number *m* of particular values; calculating the mean value μ; separating that series in *n* terms smaller, and (*m* – *n*) larger than the median value, the deviations $(1/2 - n/m)$ and $(M - \mu)/(b - a)$ will oscillate fortuitously with the same probability *P* between unequally spaced limits. The interval between the limits of the first deviation will always be larger than the second one.

It is also evident that *g* has no maximal value. If one of the limits *a* and *b* or they both is/are infinite, the second integral (69.2) can become infinite but the first one will still be unity. Then, strictly speaking, the absolute mean *M* to which the mean μ can converge as *m* ever increases does not exist. If the first and the second integral (69.2) remain finite, the third can still become infinite, and then, strictly speaking, the modulus of converging will not exist. We do not consider here the various singular cases possibly occurring in the applications of those formulas.

Generally

$$P = \frac{2}{\sqrt{\pi}} \int_0^t \exp(-t^2) dt$$

is the probability that the deviation $(M - \mu)$ is contained within the limits $\pm l$ determined by equation (69.1). This formula proved in mathematical treatises is reputedly exact to within magnitudes $1/m$ but actually much more precise.

For illustrating this, suppose that all the values of *x* are equally probable. Then[6]

$$P = \frac{1}{m!}[(\alpha m)^m - C_m^1(\alpha m - 1)^m + C_m^2(\alpha m - 2)^m - \ldots \\ - (\beta m)^m + C_m^1(\beta m - 1)^m - C_m^2(\beta m - 2)^m \pm \ldots].$$

For the sake of brevity denote $(M + l)/(b - a) = \alpha$, $(M - l)/(b - a) = \beta$ and terminate each series when the appropriate term is not positive anymore. Let $(b - a) = 1$, $M = 0.5$ and $l = 0.1$, then even for $m = 10$ that expression becomes $2,585,698/3,268,800 = 0.71255$. For the same numerical values $t = \sqrt{3/5} = 0.7746$ and our table provides the corresponding $P = 0.7266 \ldots$ The difference is less than 0.015 although we could have feared, if strictly keeping to the usual derivation, that the difference will amount to one or many tenths which would have rendered the approximate formula illusory. It can be safely thought that for $m = 100$ the error will be quite negligible. A. A. C.

**6.** Bru referred to Poisson (1829; 1837, § 102) and remarked that Cournot had actually introduced an example of the Schwarz – Buniakovsky inequality which goes back [indirectly] to Laplace (1812/1886, p. 316).

**7.** This formula is due to Lagrange, 1776, and occurred in Laplace (1812/1886, p. 260) and Poisson (1837, § 110) who attributed it to Laplace. [B. B.]

## Bibliography

**Gauss C. F.** (1809, Latin), *Theory of Motion* … New York, 1963. Translator C. H. Davis. First edition, 1857.

**Laplace P. S.** (1812), *Théorie analytique des probabilités*. *Oeuvr. Compl.*, t. 7. Paris, 1886.

**Poisson S.-D.** (1824 – 1829), Sur la probabilité des résultats moyens des observations. *Conn. des temps* pour 1827, pp. 273 – 302; pour 1832, pp. 3 – 32.

--- (1837), *Recherches sur la probabilité des jugements* … Paris, 2003, 2012. English translation: www.sheynin.de   downloadable file 53.

**Chapter 7. On the Variability of Chances**

**75.** Until now, we have supposed that during the repetition of the trials the chances of the same event did not change. This hypothesis is underlying the theorem of Jakob Bernoulli, see Chapter 3, and the rules of the convergence of the mean values (Chapter 6). However, in general the chances of the same event change according to their nature from trial to trial or from one of their series to another if accomplished under other circumstances and by other instruments[1]. For example, when tossing a coin, the probability of the appearance of *heads* is not strictly equal to 1/2 because of the coin's irregular structure which is always necessary to suppose. It does not change during successive trials if always tossing the same coin when other circumstances, such as the density and the velocity (? - O.S.) of the air, remain the same.

However, when the coin is changed from one trial to another, the probability of the appearance of *heads* will also change. Assuming that all the available coins are perfectly identical the probability will [nevertheless] change from one series to another if the applied coins were differently minted. When always applying the same coins they will wear so that by the end of a series *heads* can experience progressive variations and acquire an essentially different value as compared with the beginning of the series.

What we say about an event insignificant in itself which can only become useful as a result of a conditional consent, but it is applicable to fortuitous natural and social economic phenomena of great importance. It is therefore essential to examine how the laws of probability are modified owing to the variability of chances.

**76.** Suppose that there are *n* urns with white and black balls in various proportions. In $n_1$ of them the probability of drawing a white ball is $p_1$, in $n_2$ of them, $p_2$ etc. At first we suppose that at each trial an urn is chosen by chance and a ball randomly extracted from it. The probability of obtaining a white ball is calculated by the rules of compound probabilities (§ 23) and it evidently does not change from trial to trial. The probability of choosing an urn of the first kind is $n_1/n$ and $n_1 p_1/n$ is the probability of the compound event. Therefore, the probability of extracting a white ball from some urn is

$$p = \frac{n_1 p_1 + n_2 p_2 + ...}{n} = \frac{n_1 p_1 + n_2 p_2 + ...}{n_1 + n_2 + ...}.$$

In other words, it is the arithmetic mean (§ 67) of the probabilities of that event. If the extracted ball is each time replaced in the appropriate urn which is then inserted by chance among the other urns, it does not change from one trial to another[2]. Nothing should be therefore changed in the law of probability as it was described above, suffice it to understand that the fraction *p* instead of denoting a magnitude constant for a given urn, is now a mean of the probabilities varying from one urn to another.

**77.** I will, however, adduce a useful remark. In § 33 we saw that for the same number of trials and the same probability *P* the fortuitous deviation (*p* – *w*) will be contained within the limits ± *l*, and that the value of *l* varied proportionally to the square root of $p(1 - p)$ at least when the number *m* of the trials was not less than several hundred. We have applied this hypothesis owing to the simplicity that it introduces in the calculations and the exposition of the theory.

[For the sake of brevity we denote] $n_1/n = k_1$, $n_2/n = k_2$, ... Then

$$p(1-p) = (k_1p_1 + k_2p_2 + ...)(1 - k_1p_1 - k_2p_2 - ...).$$

Carry out the multiplication on the right side, replace $k_1^2$ by

$$k_1(1 - k_2 - k_3 - ...),$$

substitute similar expressions instead of $k_2^2$, $k_3^2$ ..., then

$$p(1-p) = k_1p_1(1-p_1) + k_2p_2(1-p_2) + ... + k_1k_2(p_1-p_2)^2 + k_1k_3(p_1-p_3)^2 + ... \quad (77.1)$$

This proves that always

$$p(1-p) > k_1p_1(1-p_1) + k_2p_2(1-p_2) + ..., \quad (77.2)$$

i. e. that the value of the product $p(1-p)$ where $p$ is a mean value, always exceeds the mean of the values of that product for each urn in particular. In addition, by virtue of the principle that a mean of squares always exceeds the square of the mean value (§ 73), we have another inequality which is also easy to verify,

$$k_1p_1(1-p_1) + k_2p_2(1-p_2) + ... > [k_1\sqrt{p_1(1-p_1)} + k_2\sqrt{p_2(1-p_2)} + ...]^2.$$

Then, all the more

$$p(1-p) > [k_1\sqrt{p_1(1-p_1)} + k_2\sqrt{p_2(1-p_2)} + ...]^2,$$

$$\sqrt{p(1-p)} > k_1\sqrt{p_1(1-p_1)} + k_2\sqrt{p_2(1-p_2)} + ...$$

Denote by $l_1$, $l_2$, ... the new values of *l* corresponding to the same values of *m* and *P* and *p* replaced successively by $p_1$, $p_2$, ... The previous inequality will be equivalent to

$$l > k_1l_1 + k_2l_2 + ...$$

and we can therefore maintain that the value of the limit *l* that measures the influence of the anomalies of chance exceeds the mean of the values $l_1$, $l_2$, ... in case the same number of trials is made with each urn. It is even possible that *l* exceeds the largest of those values.

**78.** Suppose now that the urns are not anymore selected by chance but that for $m_1$ trials we choose an urn from those having probability $p_1$ of drawing a white ball, for $m_2$ trials, an urn from those having probability $p_2$ of the same event etc[3]. As before, the total number of trials is

$$m = m_1 + m_2 + \ldots$$

Calculations show that the ratio $w$ of the number of the extracted white balls to the total number $m$ of them converges to

$$p = m_1 p_1/m + m_1 p_1/m + \ldots \qquad (78.1)$$

or to the mean of the probabilities of drawing a white ball from each of the selected urns. And we have probability $P$ that the difference $(p - w)$ is contained within the limits $\pm l$ such that $l$ is connected with $t$ (§ 33) and therefore with $P$ by the formula

$$t = l\sqrt{m \div 2[m_1 p_1(1-p_1)/m + m_2 p_2(1-p_2)/m + \ldots]}.$$

These consequences take place if $m$ is a sufficiently large number and (what is very important to note) when each of the numbers $m_1$, $m_2$, … making up the large number $m$ is small or even reduced to unity[4]. The anomalies of chance are then compensated not within each particular series, but in the total series.

Suppose for a moment that the numbers $m_1$, $m_2$, … are proportional to $n_1$, $n_2$, … considered in the assumption made in §§ 76 and 77. Then

$$m_1 p_1(1-p_1)/m + m_2 p_2(1-p_2)/m + \ldots = k_1 p_1(1-p_1) + k_2 p_2(1-p_2) + \ldots$$

so that the inequality (77.2) will indicate that for the same value of $m$ the anomalies of chance are contained in limits narrower than before. This can be known in advance without calculations. Actually, if the selection of the urns is fortuitous, with an increasing number of trials the numbers $\mu_1$, $\mu_2$, … of the chosen urns of the first, second, … kind tend to become proportional, respectively, to the numbers $n_1$, $n_2$, … and, in addition, the ratio $w$ converges to the value of $p$ defined by equation (78.1). Therefore, the influence of the anomalies of chance will decrease if for some reason governing the selection of the urns the numbers $\mu_1$, $\mu_2$, … will become equal to $m_1$, $m_2$, … and consequently, by the assumption made, will be strictly proportional to the numbers $n_1$, $n_2$, …

According to the hypothesis of a fortuitous selection of urns, the ratio $w$ which converges to the value $(n_1 p_1 + n_2 p_2 + \ldots)/n$ or to the mean of the probabilities of extracting a white ball from all the urns, at the same time converges to $(\mu_1 p_1 + \mu_2 p_2 + \ldots)/m$ or to the mean of the values of probabilities of extracting a white ball from the $m$ actually applied urns. If $m$ becomes ever larger, these two means tend to, but never strictly coincide. Common sense tells us that the random difference between $w$ and the second mean should oscillate (with the same probability) within narrower limits than the random difference between the same ratio and the first mean. This statement assumes that

$$p(1-p) > \mu_1 p_1(1-p_1)/m + \mu_2 p_2(1-p_2)/m + \ldots$$

or that (§ 77)

$$k_1p_1(1-p_1) + k_2p_2(1-p_2) + \ldots >$$
$$\mu_1p_1(1-p_1)/m + \mu_2p_2(1-p_2)/m + \ldots \qquad (78.2)$$

This inequality can be verified[5].

**79.** The third hypothesis which we should consider consists in that the numbers $m_1, m_2, \ldots$ are fixed in advance but that the urns to be applied in the first, the second, … series are selected by chance[6]. Assume also, for avoiding excessive complications of formulas, that $m_1 = m_2 = \ldots$ so that $m = im_1$ where $i$ is the number of the partial series. The ratio $w$ will converge to

$$p = (n_1p_1 + n_2p_2 + \ldots)/m$$

just as under the first hypothesis but less rapidly.

The limit of the deviation and the probability $P$ will be connected by means of a supplementary magnitude $t$ according to a remarkable equation (Bienaymé 1839) which we will provide in the form

$$t = l\sqrt{m \div 2\{p(1-p) + (m_1 - 1)[k_1(p-p_1)^2 + k_2(p-p_2)^2 + \ldots]\}}. \quad (79.1)$$

If $m_1 = 1$ we return, as it should have been, to the formula taking place under the first hypothesis. The factor in the square brackets

$$k_1(p-p_1)^2 + k_2(p-p_2)^2 + \ldots \qquad (79.2)$$

is certainly positive. It expresses the mean of the squares of the differences between each of the probabilities $p_1, p_2, \ldots$ and their mean $p$. If that factor is a very small fraction, for example 1/10,000, $m_1$ can be of the order of tens or hundreds without $l$ (for the same $P$) being appreciably altered when passing from the first to the present hypothesis. However, in general that factor (79.2) is a fraction comparable with $p(1-p)$ so that the value of $l$ will considerably increase due to that passage. For small values of $m_1$ such as 10 or 12 and with $m$ remaining a very large number, it easily trebles. This occurs all the more if $m_1$ is a large number; the fraction $p(1-p)$ will generally become negligible as compared with the product of the factor (79.2) by $(m_1 - 1)$. The radical in formula (79.1) instead of being of the order of $\sqrt{m}$, in general has order $\sqrt{m/m_1}$ or $\sqrt{i}$.

For narrowing the interval between the limits of the anomalies of chance it is necessary that both the total number $m$ of trials and the quotient $i$, or the number of the partial series, are large. Instead of supposing that $m_1 = m_2 = \ldots$, or even that they are unequal but assigned in advance by causes not at all accidental, we may assume that the number of trials constituting the first, the second, … series is determined fortuitously with the total number $m$ of trials remaining the same. If chance enters the choice of the urn or the group of urns for each partial series, it will be its second intervention. Then, the chance inherent in the extraction of the balls is its third appearance. Common sense indicates, without calculations whose complications necessarily ever increase, that the more systems of fortuitous trials are interlocked, the more considerable are the anomalies of chance that corrupt the final result. We should increase the total number of trials, or in our example the total number of drawings from which the final system is constituted, for including these oscillations in the same limits with the same probability.

**80.** Generally, let us admit that there are $n$ urns subdivided into groups of $n_1, n_2, \ldots$ urns with probabilities of extracting a white ball being $p_1, p_2, \ldots$ A first long series of $m^{(1)}$ trials is made, but the system of causes, whether random or not, acts in such a way that the numbers $m_1^{(1)}, m_2^{(1)}, \ldots$ are not proportional to $n_1, n_2, \ldots$ denoting respectively the number of drawings made from an urn belonging to the first, the second, … group. Therefore, the ratio $w^{(1)}$ for the series of $m^{(1)}$ trials very little differs from the mean

$$p^{(1)} = (m_1^{(1)}p_1 + m_2^{(1)}p_2 + \ldots)/m^{(1)}$$

but can appreciably differ from the mean $p$ found with a good approximation according to the hypothesis of § 76.

When making a second series of $m^{(2)}$ trials, with the system of causes influencing the selection of urns fortuitously or not, we will find that some ratio $w^{(2)}$ will very little differ from the mean

$$p^{(1)} = (m_1^{(2)}p_1 + m_2^{(2)}p_2 + \ldots)/m^{(2)}.$$

Here, $m_1^{(2)}, m_2^{(2)}, \ldots$ are similar to $m_1^{(1)}, m_2^{(1)}, \ldots$ This mean can appreciably differ from both means, $p^{(1)}$ and $p$.

However, when having a very large number of similar series all the accidental and irregular in the causes that jointly acted in selecting the urns will be compensated and therefore lacking and the mean

$$\frac{m^{(1)}w^{(1)} + m^{(2)}w^{(2)} + \ldots}{m^{(1)} + m^{(2)} + \ldots}$$

will converge to a certain fixed limit, and to $p$ if all the causes determining the selection were purely fortuitous. Suppose that the chances which, taken together, determined the selection, are rigorously defined, then we can in advance (without the complications of calculation) assign the number of partial series $m^{(1)}, m^{(2)}, \ldots$ each consisting of a large number of trials, which should be obtained for arriving at an appropriately fixed mean approximately equal to $p$. But if the conditions of the selection or of the extractions of the balls progressively change in time, the derivation of an appropriately fixed mean will not be always possible.

**81.** Similar remarks are applicable to the laws of probability defined in Chapter 6 and to the mean values derived from them. Suppose, just like in § 71, that an equator, the poles, meridians and parallels are traced on an appreciably spherical globe. Throw it fortuitously and mark the longitude and latitude of the points of its contact with the earth. Each of these coordinates, the latitude for example, is a magnitude taking infinity of values between certain limits. For a strictly spherical and homogeneous globe the law of probabilities of these values is easily assigned, but in general it depends on its form and internal structure and varies from one globe to another.

Suppose that we have a large number of randomly collected globes and trials are successively made with a globe each time fortuitously chosen rather than with the same one. Let $n_1, n_2, \ldots$ be the number of globes for which the probability of value $x$ (? - O.S.) is proportional to function $f_1x, f_2x, \ldots$ and $n = n_1 + n_2 + \ldots$ is the total number of the globes. In each trial the probability of $x$ resulting from the selection of the globe and its fortuitous throw is proportional to

$$fx = \frac{n_1 f_1 x + n_2 f_2 x + \ldots}{n} = \frac{n_1 f_1 x + n_2 f_2 x + \ldots}{n_1 + n_2 + \ldots}.$$

In other words, it is the arithmetic mean of the probabilities of that value for each globe. By tracing the curve of probabilities for each globe we will obtain the ordinate $fx$ of the mean curve of probabilities by matching each abscissa to the mean of all the corresponding ordinates.

A median value and a mean value (§§ 67 and 68) correspond to each of the functions $fx, f_1x, f_2x, \ldots$, but there is no simple relation or a general formulation connecting the median value of $fx$ with those proper for functions $f_1x, f_2x, \ldots$ However, denoting by $M, M_1, M_2, \ldots$ the mean values for $fx, f_1x, f_2x, \ldots$, we have

$$M = \frac{n_1 M_1 + n_2 M_2 + \ldots}{n_1 + n_2 + \ldots}$$

so that $M$ is the arithmetic mean of the values $M_1, M_2, \ldots$ repeated as many times as there are globes to which they correspond.

**82.** By calculations similar to those of § 77, we find that the value of the modulus of convergence for the mean function $fx$ is smaller than the mean of those moduli for each applied globe[7]. It can even happen that the value of that modulus is smaller than the minimal modulus of the different particular forms of that function.

Denote by $l$ the limit of the difference between the veritable mean and the mean resulting from the $m$ trials made with fortuitously chosen globes and by $l_1, l_2, \ldots$ the limits of equally probable differences when the trials are successively made with a globe from the first, the second, … series. Then it follows that

$$l > \frac{n_1 l_1 + n_2 l_2 + \ldots}{n_1 + n_2 + \ldots}.$$

**83.** A second hypothesis corresponding to that introduced in § 78, consists in supposing that the selection of the globes is not anymore fortuitous, but made by non-random causes or by their combination with fortuitous causes. The series of $m$ trials consists of $m_1, m_2, \ldots$ trials made with globes for which the law of probabilities is $f_1x, f_2x, \ldots$ Then the mean resulting from the $m$ trials will converge to the value

$$(m_1 M_1 + m_2 M_2 + \ldots)/m$$

very little differing from it if $m$ is a large number, even if $m_1, m_2, \ldots$ decrease until unity. The modulus of convergence becomes[8]

$$1 \div \sqrt{\frac{m_1}{mg_1^2} + \frac{m_2}{mg_2^2} + \ldots}$$

and, assuming that the numbers $n_1, n_2, \ldots$ are respectively proportional to $m_1, m_2, \ldots$, exceeds the modulus obtained under the previous hypothesis. Nevertheless, we always have

$$l > \frac{m_1 l_1 + m_2 l_2 + \ldots}{m_1 + m_2 + \ldots}.$$

The values of $l, l_1, l_2, \ldots$ and the meaning of this inequality were indicated above clearly enough[9].

**84.** After the statements of § 79 it is understandable that if some cause has fixed in advance the numbers $m_1, m_2, \ldots$ and the globes to be applied in the series $(m_1), (m_2), \ldots$ are selected fortuitously, the mean resulting from the trials converges to the value $M$ just like under the hypothesis of a fortuitous selection for each isolated throw, even when $m_1, m_2, \ldots$ are very small numbers, but notably slower. Otherwise a large $m$ is not sufficient for the anomalies of chance to be contained within narrowly spaced limits; it is then necessary for each of the partial series $(m)$ to include a large number of throws.

**85.** So as to formulate now the most general hypothesis, imagine a combination of many series composed of very large numbers of trials $m^{(1)}, m^{(2)}, \ldots$ The causes, whether fortuitous or not, which determine the selection of globes vary from one series of trials to another so that the numbers $m_1^{(1)}, m_2^{(1)}, \ldots$ express how many globes of the first, the second, … kind were applied in the first series, whereas the numbers $m_1^{(2)}, m_2^{(2)}, \ldots$ express the same numbers pertaining to the second series etc. None of these numbers supplied with super- and subscripts ought to be of a certain order of magnitude, some can be reduced to 1 or even to 0. The means derived from the first, the second, … series of trials very little differ from, respectively,

$$(m_1^{(1)} M_1 + m_2^{(1)} M_2 + \ldots)/m^{(1)}, \quad (m_1^{(2)} M_1 + m_2^{(2)} M_2 + \ldots)/m^{(2)}, \ldots$$

but they can appreciably differ one from another, and all of them can appreciably differ from the mean $M$. However, compensation in a large number of similar series will remove all the accidental and irregular in the selection, so that, denoting by $\mu^{(1)}, \mu^{(2)}, \ldots$ the means resulting from the first, the second, … series, the general mean

$$\frac{m^{(1)} \mu^{(1)} + m^{(2)} \mu^{(2)} + \ldots}{m^{(1)} + m^{(2)} + \ldots}$$

converges to a fixed limit coinciding with *M* if the selection is purely fortuitous. It should be borne in mind that if the conditions of selection or extraction (? - O.S.) experience progressive variation in time, it can happen that, whichever is the number of the partial series $m^{(1)}$, $m^{(2)}$, … composing a total series, we will never arrive at an appreciably fixed mean.

**Notes**

**1.** In his correspondence with Jakob Bernoulli, Leibniz referred to changes in various circumstances and thus objected to stochastic reasoning. That problem was studied by Poisson and Bienaymé and originated later developments. [B. B.]
 Bernoulli answered Leibniz without naming him at the end of Chapter 4 of pt. 4 of his *Ars Conjectandi*: [*R*]*esume observations with the pebbles if it is assumed that their number in the urn is variable*. O. S.

**2.** Poisson (1837, § 55) did not agree. [B. B.]

**3.** That pattern was studied by Laplace (1812/1886, pp. 430 – 431) and Poisson (1837, § 94). [B. B.]

**4.** Poisson (1837, § 95) indicated an exception. [B. B.]

**5.** The magnitude

$$k_1 p_1(1 - p_1) + k_2 p_2(1 - p_2) + \ldots \qquad (78.3)$$

is the absolute mean of the function $p(1 - p)$. For large values of *m* it should very little differ from

$$\mu_1 \, p_1(1 - p_1)/m + \mu_2 \, p_2(1 - p_2)/m + \ldots \qquad (78.4)$$

which is provided by the *m* fortuitous trials. At least if the magnitude

$$k_1 k_2 (p_1 - p_2)^2 + k_1 k_3 (p_1 - p_3)^2 + \ldots$$

is not very small, this should lead to inequality (78.2). However, in that other case the modulus of convergence of the function *p* will be very large, see my Note 5 in Chapter 6, as obviously that modulus of $p(1 - p)$. This circumstance will even more decrease the difference between the magnitudes (78.3) and (78.4) but the inequality (78.2) will still be preserved. A. A. C.

**6.** The third hypothesis is due to Bienaymé (1839) whom Cournot followed. However, the latter considered a similar pattern in 1834. [B. B.]

**7.** Suppose, just like in § 77, for shortening the calculations, that

$$n_1/n = k_1, \; n_2/n = k_2, \; \ldots$$

so that

$$fx = k_1 f_1 x + k_2 f_2 x + \ldots$$

Then, by the formulas in my Note 5 to Chapter 6,

$$\frac{1}{2g^2} = k_1 \int_a^b x^2 f_1 x dx + k_2 \int_a^b x^2 f_2 x dx + \ldots -$$

$$- [k_1 \int_a^b x f_1 x dx + k_2 \int_a^b x f_2 x dx + \ldots]^2.$$

Expand the right side of this equation and replace $k_1^2$, $k_2^2$, … respectively by

$$k_1(1 - k_2 - k_3 - \ldots), \; k_2(1 - k_1 - k_3 - \ldots), \; \ldots$$

Then, denoting by $g_1$, $M_1$, $g_2$, $M_2$, the analogues of $g$, $M$, we find for each function $f_1$, $f_2$, …

$$\frac{1}{2g^2} = \frac{k_1}{2g_1^2} + \frac{k_2}{2g_2^2} + ... + k_1 k_2 (M_1 - M_2)^2 + k_1 k_3 (M_1 - M_3)^2 + ...$$

Therefore,

$$\frac{1}{g^2} > \frac{k_1}{g_1^2} + \frac{k_2}{g_2^2} + ...$$

However, according to the principle that the mean of squares always exceeds the square of the mean value (§ 73),

$$\frac{k_1}{g_1^2} + \frac{k_2}{g_2^2} + ... > \left[\frac{k_1}{g_1} + \frac{k_2}{g_2} + ...\right]^2$$

and the more so

$$\frac{1}{g} > \frac{k_1}{g_1} + \frac{k_2}{g_2} + ..., \quad l > k_1 l_1 + k_2 l_2 + ...$$

From the last but one inequality it follows that

$$g < 1 \div [k_1/g_1 + k_2/g_2 + ...],$$

but on the other hand

$$1 \div [k_1/g_1 + k_2/g_2 + ...] < k_1 g_1 + k_2 g_2 + ...$$

since that last inequality can be written as

$$\frac{k_1 k_2 (g_1 - g_2)^2}{g_1 g_2} + \frac{k_1 k_3 (g_1 - g_3)^2}{g_1 g_3} + ... > 0$$

and the more so

$$g < k_1 g_1 + k_2 g_2 + ... \qquad\qquad \text{A. A. C.}$$

**8.** That the variance of a sum of independent random variables is equal to the sum of the variances is due to Gauss (1823, § 15). The modulus is evidently that of the mean globe. [B. B.]

**9.** Letter $l$ denotes the limit of the difference between the veritable mean and the mean resulting from the $m$ trials and $l^2 = \sum(m_i/m)l_i^2$. [B. B.]

### Bibliography
**Bienaymé I. J.** (1839), Théorème sur la probabilité des résultats moyens des observations. […] *L'Institut*, 284, t. 7, pp. 187 – 189.
  **Gauss C. F.** (1823 Latin), Theory of combinations of observations […] Philadelphia, 1995. Latin and English. Translated by G. W. Stewart.
  **Laplace P. S.** (1812), *Théorie analytique des probabilités. Oeuvr. Compl.*, t. 7. Paris, 1886.
  **Poisson S.-D.** (1837), Recherches sur la probabilité des jugements […] Paris, 2003, 2012. English translation: www.sheynin.de downloadable file 53.

# Chapter 8. Posterior Probabilities

**86.** For random events whose conditions are not determined by man, the causes providing them certain chances or determining the law of probabilities of the various values of a variable magnitude, are almost always either of an unknown nature and manner of acting, or so complicated that it is impossible to analyze them rigorously or subject them to calculation. Even in games in which everything is of human convention and invention, the construction of the instruments of chance is prone to irregularities and the influence of the modification which they impart on chances is impossible to evaluate in advance.

When throwing a homogeneous die whose rectangular (? - O.S.) faces are not strictly equal [congruent], the required chances of the arrival of each of them, although seemingly a very simple problem of mechanics with all its conditions being strictly defined, can not be determined given the present state of mathematical analysis. And in games in which probabilities only depend on a purely arithmetical enumeration of combinations without any influence of mechanical or physical conditions, the solution of the arithmetic problem can still exceed the power of analysis. Suppose for example that it is required to determine the advantages of the *main* in piquet or of *dé* in tricktrack[1], or establishing the probability that the gambler having them wins the set […]. This demands such an inextricable enumeration of chances which no existing procedure of calculation can ensure.

It is therefore really necessary for the application of the theory of chances to be able to ascertain posteriorly by experience those chances whose direct determination actually and evidently exceeds the powers of calculation. From all said until now it follows that the principle of Jakob Bernoulli leads to such an experimental determination. Indeed, denote by $x$ the unknown chance of the arrival of an event, and by $n$, the number of its appearances in $m$ trials. It will then be always (§ 33) possible to obtain a probability $P$, such that the fortuitous difference $(x - n/m)$ is contained within limits $\pm l$; if only numbers $m$ and $n$ are sufficiently large, the number $l$ and the difference $(1 - P)$ will be less than any assignable magnitudes. It is clear that, if nothing restricts the number of trials, the probability $x$ can be determined with an arbitrary precision. It is possible, for example, to become sure that the difference between the ratio $n/m$ given by experience and the unknown $x$ will be less than 1/100,000. The existence of a larger difference, although possible in a strict sense, will be an event of the kind reasonably reputed as physically impossible, so that we may disregard it when describing various phenomena (§ 43).

By issuing from the Jakob Bernoulli theorems (? - O.S.) whose sense and importance their inventor had perfectly well understood, we are now able to pass immediately to the applications which they enjoy in the sciences of facts and observation. However, a rule first announced by the Englishman Bayes on which Condorcet, Laplace and their followers wished to construct a doctrine of posterior probabilities became a source of many ambiguities which should first be cleared up, and grave errors which should be rectified. They are rectified now that we recognize the spirit of the fundamental distinction between probabilities existing objectively and measuring the possibility of things, and subjective probabilities relative partly to our knowledge and partly to our ignorance and varying from one mind to another depending on their capacities and the information provided them (§ 46).

**87.** Imagine three groups of urns containing balls, 3 white balls in each urn of the first group, 2 white balls and 1 black ball, and 1 white ball and 2 black balls in urns of the other two groups[2]. All groups consist of the same number of urns. An urn is selected by chance and a white ball fortuitously extracted from it. It is required to determine the probabilities that that urn belonged to one of those groups.

The probability of selecting an urn from the first group is 1/3, and then the drawing of a white ball is certain. For the other groups these probabilities are 1/3 and 2/3, and 1/3 and 1/3. Therefore, the prior (§§ 20 – 23) probability of the extraction of a white ball is

$1/3·1 + 1/3·2/3 + 1/3·1/3 = 2/3$.

Denote by $A_1$, $A_2$, $A_3$ the events consisting of extracting a white ball from an urn of those groups. After repeating the described procedure a large number $m$ of times, the number of those events will appreciably be $m/3$, $2m/9$ and $m/9$.

Three gamblers could have bet that the white ball will be extracted respectively from an urn of those groups and agreed to disregard arrivals of a black ball. Their stakes should be fixed in the proportion of 3:2:1 not at all because, being ignorant of the particular causes that determine the random event, there is no reason to fix them otherwise, but since the chances of their gain are in that proportion. This, as we will explain, becomes manifest after a long series of trials.

**88.** On the contrary, if a white ball is extracted under the same conditions but it is not known from which urn, the three gamblers should regulate their stakes the same way. Their probabilities of gaining (and probability is here applied in the objective sense as equivalent to possibility) are respectively proportional to those numbers. Therefore, if the same bets were repeated many times under the same circumstances the number of sets gained by the gamblers will be approximately in the same proportion of 3:2:1.

It is in this sense that the rule attributed to Bayes (1764)[3] should be understood. It can be thus formulated:

*Probabilities of causes or hypotheses are proportional to the probabilities that these causes provide to the observed events. The probability of one of these causes or hypotheses is a fraction whose numerator is the probability of the event being due to that cause and the denominator, the sum of such probabilities relative to all the causes or hypotheses.*

In our example, the cause or the previous event was the selection of an urn from the three groups. The subsequent and observed event, i. e. the extraction of a white ball previously had different probabilities depending on the selected urn to belong to one or another group. Thus understood, the Bayes rule is a theorem which can not lead to any ambiguities and whose truth can not be contested.

**89.** Suppose now that an urn is taken from many of them, each containing 3 balls, white or black. The ratios of the balls of the two colours is unknown; also unknown is the ratio of urns containing only white balls to the total number of urns, etc and neither is it known whether the selection of urns is fortuitous or influenced by non-random causes.

Suppose also that a white ball was fortuitously extracted from an urn and that gambler $A_1$ bets that that urn only contained white balls, gambler $A_2$, that it contained 2 white balls and 1 black ball, and gambler $A_3$, 1 white ball and 2 black balls. It is required to determine how their stakes should be regulated.

Before selecting the urn 4 hypotheses could have been formulated about its contents: 1) 3 white balls; 2) 2 white balls and 1 black ball; 3) 1 white ball and 2 black balls; and 4) 3 black balls. According to the conditions of the problem there is no reason to prefer one of them rather than another. Three cases correspond to each of them, and, when issuing from the formulated question, there was no reason to prefer one of them rather than another. The observed event excluded the fourth hypothesis with its 3 cases; it also excluded 2 out of 3 cases of the third hypothesis and 1 out of 3 of the second. Only 6 cases remained with no reason to prefer one of them rather than another, three of them corresponding to the first hypothesis favour $A_1$, two corresponding to the second favour $A_2$, and one corresponding to the third favours $A_3$. The gamblers' stakes should therefore be regulated by the ratio 3:2:1.

However, when declaring that the probabilities of gaining are in that ratio and when the Bayes rule is understood in that sense, the term *probability* will be purely subjective, variable from one individual to another depending partly on their knowledge and partly on their ignorance. This meaning does not at all signify that in a large number of similar betting the number of sets gained by the three gamblers will approximately be in the ratios of 3:2:1. It can happen that gambler $A_1$, for example, will invariably lose which invariably takes place if the urn at each trial is selected from a group with no urns having only white balls […].

The Bayes rule thus applied for determining subjective probabilities can only be useful for fixing the stakes when some hypothesis is about things both known and unknown to the arbiter. It leads to unjust fixing if the arbiter knows more than supposed about the real conditions of the random trial[4].

**90.** Return to the hypothesis of § 87 and suppose that after selecting an urn and extracting a white ball it is required to determine the probability of extracting another such ball. The first ball is returned to the urn so that the conditions of both trials are the same. We saw that after the first drawing the probabilities of the three possible hypotheses about the content of the urn are 3/6, 2/6 and 1/6. The respective probabilities of extracting a second white ball under these hypotheses are 1, 2/3 and 1/3. Therefore[5],

$$3/6 + 2/6 \cdot 2/3 + 1/6 \cdot 1/3 = 7/9$$

is the probability of extracting a second white ball and probability is here understood in its objective sense. Randomly select a very large number of urns and each time consecutively extract two balls from the chosen urn. Count the number *m* of first trials resulting in the drawing of a white ball and the number *n* of those among these *m* trials in which the second trial also resulted in an extraction of a white ball. Then the ratio *n/m* will little differ from 7/9.

However, suppose that under the hypothesis of § 89, after the extraction of the first white ball, two gamblers bet on the appearance of a white and a black ball respectively at the second trial. Their stakes should be justly regulated at a ratio of 7:2 but only because of the hypothesis about known and unknown things in the conditions of the random trial. Indeed, those conditions can be such that either the gambler betting on a black ball will invariably lose or in a large number of similar bets the ratio of the numbers of gained bets to those gained by his adversary can very much differ from 2:7.

**91.** In ordinary applications of the Bayes rule absolutely nothing is known about the contents of the urn and it is admitted that the chances possibly vary continuously, or, in other words, that the urn contains an infinity of balls and that the ratio of the number of white balls to the total number of balls can take all the values, infinitely many of them, contained between 0 and 1. In advance, all these values have equal and infinitely low probabilities. The extraction of a white ball assigns another law of probabilities to these values (§ 65) but the probability of each particular value remains infinitely low.

Just like in the second example of § 70, the curve of probabilities is a segment passing through the origin of coordinates at the angle to the *x*-axis having tangent 2. The mean value resulting from that law of probabilities is 2/3, and it precisely expresses the probability of an appearance of a white ball in the next drawing from the same urn. The median value is $1/\sqrt{2} = 0.7071$. […] A bet of 3:1 can be taken on the chance of that event to exceed 1/2 or on the urn containing more white balls than black ones.

**92.** More generally, if in *m* drawings with replacement made from the selected urn there appeared *n* white and (*m* – *n*) black balls the ordinate of the curve of probabilities will be expressed by the fraction

$$\frac{x^n(1-x)^{m-n}}{C_m^n}$$

with *x* being the value of the chance of extracting a white ball.

The curve of probabilities passes through the origin and touches the *x*-axis once more at point *x* = 1. The maximal ordinate OK = *n/m*, the mean value, also expressing the probability of extracting a white ball in the next drawing is[6]

OG = (*n* + 1)/(*m* + 2)

which is smaller than OK if *n* > *m* – *n*, and the median value is situated between OG and OK.

**93.** All these results should be understood in the objective sense applicable for measuring the possibility of events if the urn is indeed fortuitously selected from infinity of other urns so that the probability of choosing an urn for which the chance $x$ of extracting a white ball does not change with $x$. However, otherwise and in ordinary applications these results can only lead to arbitrary regulating the conditions of a bet without knowing the real circumstances of a fortuitous trial. But still in most cases in which these circumstances are unknown we sufficiently understand their nature for avoiding the temptation to regulate the stakes of a bet by the abovementioned rules based on the Bayesian principle.

Suppose that we have a pile of recently minted coins. Select randomly one of them, and toss it. *Heads* appeared but we will not bet 2:1 on the appearance of *heads* in a toss of another randomly chosen coin. However, if someone takes such a bet and repeats it many times under the same circumstances, he will not at all gain 2/3 of his bets. Indeed, although certainly varying a bit from coin to coin because of the irregularities of their physical structures, for any coin the probability of throwing *heads* can not differ much in either direction from the fraction 1/2. In that case, it will therefore be contrary to our notion about the conditions of chance to attribute indifferently all the values from 0 to 1 to the prior chance of throwing *heads*.

We know absolutely nothing about the chances of each woman to conceive a child of one or another sex; these chances certainly vary from one of them to another. We only know the mean values of these chances as derived from the statistical results of a very large number of births. If a woman gave birth to a boy, we should in our ignorance perhaps regulate the bets of two gamblers on the birth of a boy or a girl after the second pregnancy by the ratio 2:1. However, that regulation only motivated by our ignorance has no relation to the real chances of the two events.

If the same bet is repeated many times under similar circumstances, the ratio of the bets gained by the gamblers will, to all appearances, much differ from 2:1. To find out a value with good approximation, it is necessary to study the registers of the certificates of birth and to establish how many times a birth of a firstborn boy was followed by births of a second boy or a girl. As far as we know, that interesting research had not yet been done[7]; and, for the time being, the application of the Bayes rule as we saw only leads to futile and illusory consequences. Nevertheless, such poorly based applications are fearlessly made to problems essentially interesting for the society and morals like those relating to judicial decisions and testimonies and thus leading to aberrances unbecoming of eminent geometers[8].

**94.** If successive drawings are made not from the same urn but from urns each time fortuitously selected from the same group, the problem does not change, but the letter $x$ will denote (§ 76) the mean chance of extracting a white ball from that group of urns. If the conditions of selecting change from trial to trial according to an unknown law, nothing can be concluded about future events given the observed events. Condorcet provided and discussed at length absolutely illusory formulas pertaining to that case.

**95.** When the numbers *m* and *n* from § 92 are very large, the points K and G will appreciably coincide and the result obtained by the Bayes rule will not anymore essentially differ from that following from the Jakob Bernoulli theorem. And this should certainly happen because the verity of that theorem is independent from any hypothesis about the previous selection of the urn. It is not as though (as many authors apparently imagined) the Bernoulli rule becomes exact when it approaches the Bayes rule; on the contrary, the latter becomes exact and acquires its lacking objective value[9] by coinciding with the former.

Let us delve in pertinent details seemingly delicate or even trivial but demanded by the importance of the subject. Given an extraction from an urn of *n* white balls in *m* drawings the Bayes rule provides a certain probability *P* that the chance *x* of the appearance of a white ball is contained within the interval

$$[n/m - l, n/m + l]. \qquad (95.1)$$

Here, *l* is a magnitude, ever decreasing for the same value of *P* when *m* and *n* increase and can finally become smaller than any given magnitude.

Suppose that Kk is the maximal ordinate of the curve of probabilities and KI and KL, ordinates on both sides of Kk such that KI = KL = *l*. Then the probability *P* is represented by the ratio of the area of the figure between KI, KL, the *x*-axis and the curve of probabilities to the area of the figure *below* that curve[10].

If the chance of selecting an urn for which the chance *x* of extracting a white ball does not change with *x*, probability *P* will have an objective value. In other words, if after fortuitously selecting an urn and extracting from it *n* white balls in *m* drawings, I *decide* that the chance *x* of drawing a white ball from that urn is contained within the interval (95.1) and if I repeat the same judgement after a large number *N* of similar results obtained with the same number of different urns, the ratio of the number of correct judgements to that of mistaken judgements, will be approximately equal to $P/(1 - P)$. And, once more in other words, $(1 - P)$ precisely measures the chance or the possibility of an error inherent in the first judgement.

However, in general this will change if the chance of selecting an urn varies with the value of *x* for that urn. Represent the law of probabilities of the values of *x* when the urn is selected by curve *o'k'b'* (Fig. 3), the law according to which actually varies the chance of selecting the urn for which *x* denotes the ratio of the number of white balls to the total number of balls. Let the lengths of *OI'*, *OK'* and *OL'* be equal respectively to $n/m - l$, $n/m$, $n/m + l$. Suppose also for the sake of definiteness that the ordinate *K'k'* is minimal. The chance to draw exactly *n* white balls in *m* extractions is the feebler the more the value of *x* for the selected urn deviates from *OK'*.

But on the other hand the chance to select the urn for which the chance of extracting a [white] ball is *x*, becomes the higher, the more *x* differs from *OK'*. Therefore, and also because the area *I'L'l'i'* is only a small fraction of the total area *OB'b'k'o'*, it can happen that for a large number of events consisting in fortuitously selecting an urn and extracting *n* white balls in *m* drawings, the number of cases in which *x* is beyond the interval (95.1) certainly exceeds by the Bernoulli principle the number of cases in which *x* is contained within it although probability *P* defined above continues to exceed 1/2 or is even very close to unity.

It follows that if, after extracting *n* white balls in *m* drawings, it is stated that the value of *x* for the urn fortuitously selected for the drawings is contained within the interval (95.1) the chance of an error of that judgement will not in general be $(1 - P)$ but $(1 - P')$ which can essentially differ from $(1 - P)$ and remain unknown until the law of probabilities represented by curve *o'k'b'* remains unknown. In a large number *N* of identical judgements made under the same circumstances the ratio of the number of the correct to that of the mistaken judgements will not anymore be approximately in the ratio $P/(1 - P)$. The probability *P* resulting from the Bayes rule can only be understood in the subjective sense serving for regulating conditions of betting since we do not know anything about the form of the curve *o'k'b'*.

But suppose now that *m* and *n* denote large numbers. Then, owing to the Bernoulli principle for the values of *x* such as *OE'* which is much shorter than *OI'*, the event consisting of extracting *n* white balls in *m* drawings will be very rare and actually impossible. Therefore, although the ordinate *E'e'* is much larger than *K'k'*, in approximate calculations the appearance of that event when drawing from an urn in which *x* has value *OE'* or a smaller value can be disregarded. Such an event will occur extremely rarely even when *N* becomes very large. The same remark is applicable to values of *x* much larger than *OL'*. Therefore, when calculating *P'* as denoted above, if the curve *o'k'b'* is given, it is only necessary to consider its portion close to *k'* in which ordinates little differ from *K'k'* or *n/m*. Indeed, the other portions of that curve whichever form they have do not essentially influence the value of *P'*.

Since the variations of the ordinates in the portion of the curve close to *k'* are small, the error made when implicitly supposing that by the Bayes rule that ordinate [*K'k'*] is constant, is very small. Therefore, the value of *P* can be assumed with a sufficient approximation as the value of *P* [of *P'*]. For sufficiently large *m* and *n* probability *P* thus acquires an objective value independent from the form of the unknown function. In other words, the contrary probability $(1 - P)$ provides an approximate but sufficiently precise measure of the chance of an error which actually affects our judgement in pronouncing, after the extraction of *n* white balls in *m* drawings, that, for the urn which fortuitous causes or causes independent from those governing the drawings have selected among many others for our experiment, the chance *x* is contained in the interval (95.1).

**96.** When *m* and *n* are large numbers, calculations based on the abovementioned remarks prove that there exists a connection between the probability[11] *P* and the limit *l* expressed through an auxiliary variable *t* (§33) by means of equation

$$t = lm\sqrt{m/[2n(m - n)]}. \qquad (96.1)$$

This result can be understood by a very simple reasoning[12]. Actually $x$ always denotes the unknown chance and (§ 33) for sufficiently large values of $m$ and $n$ there exists an equation

$$t = l\sqrt{m/[2x(1-x)]}.$$

This formula will not essentially change if the unknown $x$ is replaced by a very little differing ratio $n/m$, differing the less the larger are $m$ and $n$. Formula (96.1) is indeed the result of that substitution.

**97.** It can be required to determine the probability that in another series of a large number $m'$ of drawings from the same urn the ratio $n'/m'$, where $n'$ is the new value of $n$, is contained within the interval $n/m \pm l'$. And if $l$ is exactly zero, the probability $P$ that $(n/m - n'/m')$ is contained within $\pm l'$ will be provided by formula[13]

$$t = l'm\sqrt{m'/[2n(m-n)]}. \qquad (97.1)$$

However, we know that the value of $l$ corresponding to the same probability $P$ should necessarily increase because of the fortuitous and very small, but differing from zero values as admitted by the deviation $(x - n/m)$. Calculations prove that the preceding formula should therefore be replaced by

$$t = l'm\sqrt{\frac{mm'}{2n(m-n)(m+m')}}. \qquad (97.2)$$

If $m'$, although considerable, is very small as compared with $m$ (for example, if it is of an order of a thousand, and $m$, of millions), formulas (97.1) and (97.2) appreciably coincide, and if $m' = m$, the values of $l'$ in (97.2) and (97.1) are in the ratio $\sqrt{2}:1$[14]. That ratio ever increases with the increase of $m'$ and with $m$ remaining constant although the absolute value of $l'$ in (97.2) incessantly decreases. Finally, when $m'$ becomes very large as compared with $m$, that value $l'$ appreciably coincides with the value of $l$ in formula (96.1).

This result was easy to foresee because the ratio $n'/m'$ can not essentially differ from $x$. It is also possible to arrive at formula (97.2) when determining by the Bayes rule the probability that after $n$ white balls had been extracted in $m$ drawings, $n'$ were extracted in $m'$ drawings, and then, by formulas of ordinary approximation to pass on to the case in which $m, n, m'$ and $n'$ are large numbers. However, the main point is that the obtained result is independent from the conditions of selecting the urn and therefore from the Bayes rule. They are only consequences of the Bernoulli principle and the hypothesis which assigns large values to those four numbers.

**98.** The limit $l'$ corresponding to probability $P$ is provided by formula (97.2) as a function of $m, n$ and $m'$ or of $m, n, m'$ and $n'$ according to formula[15]

$$t = \frac{l'mm'\sqrt{mm'}}{\sqrt{2[m^3n'(m'-n') + m'^3n(m-n)]}} \qquad (98.1)$$

in which these four numbers enter symmetrically.

It is understandable that the two values of *t* in formulas (97.2) and (98.1) can not in general be identical because of the fortuitous deviations still inherent in *n'* after assigning a determinate value to *m'*. However, these two values become identical and equal to

$$t = \frac{l'm\sqrt{m\alpha}}{\sqrt{2n(m-n)(1+\alpha)}} \qquad (98.2)$$

if the two ratios, *m/n* and *m'/n'*, are strictly equal so that we can state that, at the same time, *m'* = α*m* and *n'* = α*n*. It follows that the two values of *t* in formulas (97.2) and (98.2) very little differ from each other if the difference between those two ratios remains very small. And this should be admitted if *m, n, m'* and *n'* are large numbers and both series of extractions are made, as supposed, from the same urn.

For applying formulas (97.2) and (98.2) it is evidently not necessary for the first series of drawings to be already made, so that, before all the trials the deviation (*m/n* − *m'/n'*) between the two series of future extractions from the same urn or from two urns giving the same chance *x* of drawing a white ball remains within the limits ± *l'* with the same probability *P*.

The value of *l'* is connected with those of *t* and *P* by formula (98.1) which we prefer because of its symmetric composition.

**99.** Suppose that the two series of extractions are made from two urns or from the same urn when its contents can change between those series. Denote by $x_1$ and $x_2$ the unknown chances of extracting a white ball. Experience provides

*n/m* − *n'/m'* = δ,

in which for the sake of definiteness δ is a positive fraction. And, having the observed deviation δ, it is required to determine the probability Π that $x_1 > x_2$ at least by α which can be smaller or larger than δ.

As before, denote by *P* the probability that existed before the trials when admitting that $x_1 = x_2$ and that (*n/m* − *n'/m'*) will be contained within ± (δ − α). The value of *P* will be provided by the auxiliary magnitude *t*:

$$t = \frac{\pm(\delta-\alpha)mm'\sqrt{mm'}}{\sqrt{2[m^3n'(m'-n') + m'^3n(m-n)]}}.$$

The double sign indicates that (δ − α) should always be positive for *t* to remain positive. Then

$$\Pi = (1 \pm P)/2$$

with the signs corresponding to α being smaller or larger than δ.

If α = 0, Π will be the probability that $x_1 - x_2 > 0$ however small this difference is supposed to be. In other words, it will be the probability that the deviation δ can not be only attributed to the anomalies of chance but that it indicates a variation of chances from one series of trials to another.

Suppose that $n'/m' = w$ and let $m'$ increase to infinity. Then Π will denote the probability that $x_1 > w$ at least by α. In this case the value of $t$ which determines $P$ and therefore Π, becomes

$$t = \frac{\pm(n/m) - w - m\sqrt{m}}{\sqrt{2n(m-n)}}.$$

For example, if $n > m/2$, the probability that $x_1 > 1/2$ or that the urn contains more white balls than black ones will be equal to the value of Π corresponding to

$$t = \frac{(2n - m)\sqrt{m}}{2\sqrt{2n(m-n)}}.$$

**100.** Imagine that each ball contained in the urn, whether white or black, is marked in red ink by letters $a$ or $b$ but that it is not known whether the marks were applied by a blind agent without distinguishing the colours of the balls, or, on the contrary, by some cause so that one of these letters is preferred for balls of a certain colour.

At each drawing the mark on the extracted ball is registered and the result of analysing the extractions is the separation of the total series of the $m$ drawn balls in two partial series, $m_1$ balls marked $a$ of which $n_1$ are white, and $m_2$ balls marked $b$ of which $n_2$ are white. It is required to find out whether the deviation $(n_1/m_1 - n_2/m_2) = \delta$ should be attributed to anomalies of chance or, on the contrary, does it indicate with a sufficient probability that the chances of extracting a white ball are not identical in both series, and that the letters $a$ and $b$ were not distributed completely independently from the colours of the balls.

Before the analysis there was probability $P$ corresponding to

$$t = \frac{\delta m_1 m_2 \sqrt{m_1 m_2}}{\sqrt{2[m_2^3 n_2(m_2 - n_2) + m_1^3 n_1(m_1 - n_1)]}}. \qquad (100.1)$$

that the deviation will be contained within the limits ± δ provided that the chance $x$ remained identical in both series. After the analysis there appeared probability

$$\Pi = (1 + P)/2$$

that the deviation $\delta$ indicates that $x_1 > x_2$. These letters denote the values of the chance $x$ in the two series. It follows that if $\Pi$ only differs from unity by a very small fraction, it should be considered almost certain that the chance of extracting a white ball varied from the first series to the second.

**101.** When comparing a partial series not with the other one, but with the total series, the probability $P$ that the deviation $(n/m - n_1/m_1)$ is contained within the limits $\pm \delta$ is provided as a function of $t$

$$t = \frac{\delta m \sqrt{mm_1}}{\sqrt{2n(m-n)(m-m_1)]}}$$

due to Bienaymé (1840)[16]. It is a reasonable analogy with that of § 97.

**102.** It is quite possible to imagine that the balls are also distinguished by inscribing on them letters $a'$ and $b'$. We suppose that the distributions of those new letters and of $a$ and $b$ are independent. However, we also suppose that we do not know in advance whether the former distribution is independent from the colours of the balls. An analysis of the $m$ drawings separates the total series in two new partial series; one of them is composed of balls marked $a'$, $n_1'$ of them white, and $(m_1' - n_1')$ black, and the other, of balls marked $b'$, $n_2'$ of them white, and $(m_2' - n_2')$ black. To the observed deviation $(n_1'/m_1' - n_2'/m_2') = \delta$ correspond probabilities $P'$ and $\Pi'$ sufficiently defined by the above.

Nothing restricts the number $s$ of the binary systems of contrary letters followed by the same number of separations of the total series. Denote by

$m_1^{(i)}, n_1^{(i)}, m_2^{(i)}, n_2^{(i)}, \delta^{(i)}, x_1^{(i)}, x_2^{(i)}, P^{(i)}, \Pi^{(i)}$,

the numbers in the system $(a^{(i)}, b^{(i)})$ being the analogues to

$m_1, n_1, m_2, n_2, \delta, x_1, x_2, P, \Pi$,

the numbers of the system $(a, b)$.

Suppose also that all the terms of the sequence $\Pi, \Pi', \ldots, \Pi^{(s-1)}$ except $\Pi^{(i)}$ are sufficiently low for believing at once that all the deviations $\delta, \delta', \ldots, \delta^{(s-1)}$ except $\delta^{(i)}$ should be attributed to anomalies of chance. It is required to find out what follows from the determined values of $\delta^{(i)}$ and $\Pi^{(i)}$. In this connection it is important to make an important distinction. If the experimenter wishes to separate at once, without any analysis, the total series and prefers the system $(a^{(i)}, b^{(i)})$ whether chance suggested him the idea to study it in preference to all the other ones or he had some prior motive to believe that the chances $x_1^{(i)}, x_2^{(i)}$ are unequal then, as explained above, the number $\Pi^{(i)}$ undoubtedly measures the probability that after the drawings $x_1^{(i)}$ really exceeds $x_2^{(i)}$. And if $\Pi^{(i)}$ very little differs from unity that excess should be regarded as almost certain.

However, if the experimenter is only led to consider the system ($a^{(i)}$, $b^{(i)}$) in preference to the other ones by studying the results of the drawings, his conclusion can not be the same. Actually, if the chances $x_1$ and $x_2$ are identical, even having a high probability $P$ that for each system ($a, b$) the deviation $\delta$ is contained within certain limits $l$, if the number $s$ is very large there can be a high probability that at least for one system that deviation will exceed these limits. In that case the observed anomaly concerning the system ($a^{(i)}$, $b^{(i)}$) considered among many others can indeed be fortuitous. It is even likely that, when multiplying the numbers of systems and drawings we will finally fortuitously encounter such an anomaly.

And there is nothing surprising that essentially unequal chances of error effect the same judgement about the same fact depending on how was the judgement made. Thus, it is quite understandable that the experimenter is only mistaken once in a thousand cases if only, according to a preconceived idea, he at once separates the total series according to the letters $a^{(i)}$, $b^{(i)}$ and finding that $\Pi^{(i)} = 0.999$ announces that $x_1^{(i)} > x_2^{(i)}$. On the contrary, he can be mistaken 999 times out of a thousand in formulating the same judgement when only encountering a deviation $\delta^{(i)}$ after a large number of attempts by trial and error and having this result as the only reason for concentrating his attention on the system ($a^{(i)}$, $b^{(i)}$) in preference to many other ones.

There follows a singular consequence[17]. A person not knowing how the data were analysed and whom the experimenter told the result of that analysis concerning the system ($a^{(i)}$, $b^{(i)}$), but not how many attempts he made to achieve that result, is unable to judge with *a determined chance of error* whether the chances $x_1^{(i)}$ and $x_2^{(i)}$ are equal or not. Actually that person could have had prior reasons to believe in their inequality and by similar reasons the system ($a^{(i)}$, $b^{(i)}$) rather than many other equally possible systems interested the experimenter independently from the results of the analysis. However, to appreciate these motives is not equivalent to [revealing] a measurable probability having an objective value and representing the veracity or error really affecting a judgement when the conditions of randomness are strictly defined. In the next chapter devoted to the applications of statistics we return to this singular consequence and explain it more sensibly by less abstract examples.

**Notes**

**1.** The term *advantage* is here connected with equalizing an unfair game by differing the gamblers' stakes (Montmort 1708/1713, p. 74). The outcome of a game of piquet, see Laplace, 1774 and Poisson (1837, § 65) depends on chance and the gamblers' abilities. Jakob Bernoulli, in his *Lettre à un amy …* appended to his *Ars Conjectandi* studied such a game. [B. B.]

**2.** Cf. § 89 and Poisson (1837, § 28). [B. B.]

**3.** Bayes did not introduce that formula which is due to Laplace.

**4.** Poisson (1837) repeatedly stated that two persons having different knowledge about a disputed thing can derive differing conclusions about it.

**5.** Poisson (1837, § 32) offered a general formula for such probabilities. [B. B.]

**6.** This is the famous *rule of succession* (Zabell 1989).

**7.** Poisson (1837, § 59) investigated that problem, but, as Cournot stated, no one had studied the appropriate statistical data.

**8.** Poisson (1837, p. 2) noted that it was Condorcet who had remarked on the possibility of applying the Bayesian principle to the probability of judicial judgements and Laplace (1812/1886, pp. 528 – 530) took up that recommendation later criticized by Bienaymé (1838, p. 208). [B. B.]

**9.** Concerning the objective value of the Bayes rule see Mises (1964, pp. 342 – 343) who noted that the ignorance of prior probabilities becomes less important with the increase of the number of observations. [B. B.]

**10.** I did not reproduce Cournot's pertinent figure (several others either), hence my careless word *below*.

**11.** Probability $P$ is the probability that chance $x$ is contained within the interval $(n/m - l, n/m + l)$. Formula (96.1) a few lines below is due to Laplace (1812/1886, p. 287) and was applied by Poisson (1837, § 83). [B. B.]

**12.** See that simple reasoning in Poisson (1837, § 83). [B. B.]

**13.** Formula (97.1) is due to Poisson (1837, § 83) and formula (97.2), to Laplace (1812/1886, pp. 393 – 394). [B. B.]

**14.** The ratio $\sqrt{2}/1$ was noted by Laplace in 1774 and 1780 and Poisson (1837, § 87). [B. B.]

**15.** Formula (98.1) is due to Poisson (1837, § 87). [B. B.]

**16.** Heyde & Seneta (1977, pp. 108 – 111) discussed Bienaymé's note (1840) who had introduced the formula provided later by Cournot. [B. B.]

**17.** Cournot was one of the first to turn attention to that phenomenon. [B. B.]

## Bibliography

**Bayes T.** (1764), An essay towards solving a problem in the doctrine of chances. *Biometrika*, vol. 45, 1958, pp. 293 – 315.

**Bienaymé I. J.** (1838), Probabilité du jugements et des témoignages.[…] *L'Institut*, 235, t. 6, pp. 207 – 208.

--- (1840), Principe nouveau du calcul des probabilités avec ses applications aux sciences d'observation. […] *L'Institut*, 333, t. 8, pp. 167 – 169.

**Heyde C. C., Seneta E.** (1977), *Bienaymé*. New York.

**Mises R.** (1964), *Mathematical Theory of Probability and Statistics*. Edited and complemented by Hilda Geisinger. New York – London.

**Montmort P. R.** (1708), *Essay d'analyse sur les jeux de hazard*. Paris, 1713, New York, 1980.

**Poisson S.-D.** (1837), *Recherches sur la probabilité des jugements* … Paris, 2003. English translation: www.sheynin.de   downloadable file 53.

**Zabell S. L.** (1989), The rule of succession. *Erkentniss*, Bd. 31, No. 2 – 3, pp. 283 – 321.

## Chapter 9. Statistics in General
## and Experimental Determination of Chances

**103.** Statistics is a quite modern science[1]. The geniuses of the ancient times had not been voluntarily busying themselves with work demanding precision. Means of research and communication did not exist and finally (what is most of all surprising) in spite of the variety of their philosophical speculations they apparently had not suspected the existence of the principle of compensation[2] which in the long run always manifests the influence of regular and permanent causes and ever more weakens that of irregular and fortuitous causes.

In our day, on the contrary, statistics developed, so to say, abundantly and we should be guarded against its premature and excessive applications which can discredit it for a while and delay the so desirable epoch when the materials of the experience will serve as a certain basis for all the theories aimed at the diverse parts of the social organization. Actually, statistics (as indicated by its name) is understood as a collection of facts which take place because of the agglomeration of men into political societies. However, for us that word takes a wider acceptance. We believe that *statistics is a science aimed at collecting and coordinating numerous facts of each kind for obtaining numerical ratios appreciably independent from anomalies of chance and indicating the existence of regular causes whose action is joined with that of fortuitous causes*.

**104.** That distinction between regular or permanent causes and accidental or fortuitous causes frequently occurs in this work and it is proper to attach quite an exact sense to it and to see how is it connected with our notion of randomness and physical possibility of events. When a die of irregular structure is tossed many times in succession, the appearance of its certain face at each throw is the effect depending on the direction and intensity of the impulsive force as well as on the form of the die and the manner of the distribution of its mass. However, in general we may admit that the causes determining the intensity and direction of the impulsive force as well as its point of application at one throw are perfectly independent from those determining them at the next throws.

But the irregularity of structure, for example the [existence of a] distance of the die's centre of gravity from its geometric centre each time acts[3] in the same sense and favours the appearance of a certain face oftener than some other. These invariably present causes whose influence extends over a whole series of trials are those which we call regular or permanent whereas those which dominate each trial separately without leaving any trace of solidarity between their action from one case to another, are called accidental or fortuitous. The effects of their irregular variations compensate themselves and disappear in the mean result of a large number of trials. They do not influence the measure of possibility of an event A or of [a contrary] event B although in each particular case they efficiently unite and determine the appearance of A or B.

The permanent causes are those whose influence determines the possibility of each event susceptible of appearing or not according to the combination taking place between permanent and fortuitous causes. And the elimination of accidental causes and the investigation of permanent causes is the essential aim of statistical researches.

**105.** For statistics to deserve the name of science it should not only consist of a compilation of facts and numbers; it should have its theory, its rules and principles and be applicable to physical and natural as well as to social and political facts. In this sense, the phenomena occurring in the celestial spaces can be subjected to statistical rules and investigations just like agitations of the atmosphere, perturbations in the animal economy and still more complex facts generated by the state of the society, by the friction between individuals and nations[5].

**106.** The study and critical analysis of documents should lead in each branch of statistics to difficulties and special rules with which we will not occupy ourselves. When admitting that we have collected necessary materials or documents having the required exactness and authenticity, they should be processed, ordered, their elements restored to proper form and connected with the original data whose values implicitly determine everything else but can sometimes be inaccessible to direct observation. For example, when treating judiciary statistics, it is seen that statistics can not directly provide the chance of an error corrupting sentences returned by a judge or a tribunal which can nevertheless be indirectly derived from other statistical numbers like effect is derived from its cause.

Statistics is a science of observation. Numbers are the instruments applied by statisticians and their precision is made comparable by formulas of the theory of chances. But the essential goal of the statistician, just like of any other observer, is to penetrate as deep as possible into the knowledge of the essence of things. To achieve this, he should by a rational discussion separate as distinctly as possible the immediate data of observation and their modifications introduced solely by the observer's point of view and his means of observation.

**107.** The ordinary immediate aim of statistical enumerations and tables is to find out the chances of the occurrence of an event that, depending on fortuitous combinations, can happen or not under given circumstances, or to determine the mean value of a variable magnitude susceptible of fortuitous oscillations within certain limits, or, finally, to assign a law of probabilities of an infinite number of values which a variable magnitude can take under the influence of fortuitous causes. It is natural to treat at first the statistical determination of the chance of an event or to provide a measure of its possibility.

It was seen in § 96 that if an event A whose unknown probability is $p$ has appeared $n$ times in $m$ observations or trials, there is probability $P$ that the error made when assuming that $p$ is equal to $n/m$, is contained within the limits $\pm l$ where $l$ is connected with an auxiliary magnitude $t$ and therefore with probability $P$ by equation (96.1)

$$t = lm\sqrt{m/[2n(m-n)]}.$$

We have explained that probability $P$ has an objective value. It actually measures the possibility of an erroneous judgement when assuming the above. Even when among the indefinite multitude of facts to which statistical observations are applicable, owing to the unknown reasons certain values of *p* can be produced more frequently than the others, the number of proper judgements as stated above will be to the number of those erroneous as $P/(1 − P)$ if only we collect a sufficiently numerous series of judgements for the fortuitous anomalies to be appreciably compensated.

For the same value of *P* the interval 2*l* of the errors is inversely proportional to

$$\frac{m\sqrt{m}}{\sqrt{2n(m-n)}} \qquad (107.1)$$

so that this number can be adopted as the measure of *precision* of the determined unknown *p* when its value is believed to be *n/m*. In the words of other authors, (107.1) measures the *weight* of the result obtained by statistical observations. For rendering the degree of precision of statistical results easily comparable it would have been proper to supplement the ratios *n/m* by their appropriate weights (107.1).

**108.** In a new series of *m'* observations event A will occur *n'* times and it should be regarded impossible that *n'/m'* strictly coincides with *n/m*. If the chance *p* did not change between the series of trials, there will be probability *P* (98.1) that the deviation (*n/m* − *n'/m'*) is contained between the limits ± *l* determined by the equation

$$t = \frac{l'mm'\sqrt{mm'}}{\sqrt{2[m^3 n'(m'-n') + m'^3 n(m-n)]}}.$$

Suppose that experience provided that

*n/m* − *n'/m'* = δ,

then, in the preceding equation *l'* = δ and after calculating the corresponding values of *t* and *P*, the fraction $\Pi = (1 + P)/2$ will provide the probability that the chance *p* of the event A in the first series of trials exceeds its chance *p'* in the second series or that the probability of the deviation δ can not be attributed to anomalies of chance.

**109.** In general, *p* and *p'* should be understood as the means of a multitude of the distinct possible values of the chance of event A when passing from one category to another[6] or even from one individual case to another (§ 75). These means can differ because either the chances of the event in the different categories vary between the series; or the same categories enter in these series in differing proportions; or the categories in the first series do not enter the second or vice versa; or, finally, because of the coincidence of all these circumstances.

All these causes of the variation of the mean *p* between the series can be fortuitous and irregular. Thus, an event comparable to a throw of dice can determine the proportions in which categories (a), (b), (c),

… combine to form the series (m) and (m′) [the series consisting of *m* and *m*′ trials]. Another event also comparable to a throw of dice can determine the values of the chance of event A for (m) and then for (m′) for each category (a), (b), (c), … In this sense it will not be exactly true to say that the deviation δ indicating a likely variation of the mean *p* from one series to another can not be attributed to irregularity of chance. Nevertheless, we will continue to apply the same expressions in the following sense.

It is clear that the event which determines the value of the chance *p* for each category (a), (b), (c), …of the series (m), although being as though a fortuitous and irregular cause for a system formed by a large number of series (m), (m′), (m″), … is a constant cause (§ 104) of sorts with respect to individual observations of the series (m) because, in the same time and the same manner, it is solidary[7] and affects all its observations belonging to each category separately. Just the same, the event which determines the proportions of the categories forming series (m), being as though a fortuitous cause with respect to the system of series (m), (m′), …, affects all at once the individual cases of each series.

Let us mentally group all the influences or causes affecting in the same way all the observations of a series or a part of them. Group also all the influences which affect each observation independently from all the others. The first group determines for each series the means *p, p′*, … which either change fortuitously and irregularly or not from one series to another. What remains from the action of the second group in the observed results and does not disappear due to compensation is regarded as the part of randomness, of fortuitous anomalies in each series. And when we say that the deviation δ indicates a variation of the mean *p* not to be attributed to anomalies of chance, we mean that δ is not, or at least is not wholly due to the influences of the second group. It is engendered at least partly by the modifications introduced by the influences of the first group although they can also result from fortuitous and irregular causes whose effects are compensated when collecting a sufficient number of observations.

**110.** Suppose that *m* observations of the first series are distributed among two very numerous categories (a) and (b); that in the category (a) are $m_1$ observations and $n_1$ events [occurrences of event] A, and in category (b), $m_2$ and $n_2$ respectively. Let also

$n_1/m_1 - n_2/m_2 = \delta$,

then there will be probability $\Pi = (1 + P)/2$ that the mean chance of event A for category (a) will exceed that for category (b). The value of *P* will correspond to the value of *t* as provided by equation (100.1).

**111.** We should now return to the remark made in § 102. It is clear that nothing restricts either the number of viewpoints from which statistically researched natural events or social facts can be considered, or, as it follows, the number of indications according to which they can be distributed into many groups or distinct categories. Suppose for example that it is required to determine, by issuing from a large number of observations collected in a country such as France, the chance of a male birth known to exceed 1/2. We can at first distinguish those born in and out of wedlock, and find out, when having a large number of observations, that there is a high probability that that chance is much higher in the former case. We can also distinguish births in the countryside and in towns and arrive at a similar conclusion. These two classifications so naturally come to mind that they became the object of study for all statisticians.

It is clear that the births can also be classified according to primogenitures (? - O.S.), age, profession, fortune, religion of the parents. We can distinguish first and second marriages and births in different seasons of the year. In a word, we can study many accessory circumstances and an indefinite number of indications as a basis for the same number of distributions among categories. It is also evident that with that number increasing without limit it is ever more probable in advance that solely by the effect of randomness at least one of them will provide essentially differing rates of the number of male births for two contrary categories.

Therefore, as we have already explained, for the statistician occupied by grouping and comparisons, the probability that a given difference is not attributable to anomalies of chance will take very different values depending on the number of groups tested before encountering that difference. We invariably suppose that large numbers are available so that by virtue of the indicated principles (§ 95) in each system of tests that probability will have an objective value as being proportional to the number of bets which the experimenter will actually win if the same bet is repeated a large number of times always after perfectly similar tests and if he has a sure *criterion* for distinguishing the cases in which he was in the right or not.

However, unsuccessful tests usually leave no traces[8]; the public only knows the results which the experimenter thought to be deserving notice. It follows that a person alien to the testing is absolutely unable to regulate bets on whether the result is, or is not attributable to anomalies of chance. Even approximately assigning the rate of erroneous judgements, when having a very large number of similar judgements made under identical circumstances, will be impossible. For that person the probability which we denoted by $\Pi$ corresponding to deviation $\delta$ loses all objective consistency. He will differently appraise the same magnitude of the deviation depending on the idea about the *intrinsic value* of the indication selected as the basis of the corresponding division into categories.

**112.** For elucidating that remark without leaving our example suppose that births are divided into two categories, those occurring between summer and winter solstices or, on the contrary, between winter and summer solstices. It is plausible in advance that the chance of a male birth in those periods is not strictly invariable or even that it appreciably varies from one of them to another.

The time of birth is connected with the time of conception[9] so that births of the first category mostly correspond to conceptions in winter, between the autumn and spring equinoxes. And it is indeed normal to presume that the differences of temperatures, diet, work and habits occasioned in our climate by the passage from summer to winter can appreciably influence the chance of conceiving a male[10]. If the magnitude of the observed difference corresponds with that conjecture, its attribution to the anomalies of chance will be extremely unlikely. Only in a very large number of divisions into categories chance all by itself can produce appreciable differences in some of them, and in advance it is very unlikely that exactly in these latter cases we would have thought that the difference was occasioned by a real variation of chances.

Suppose now that we distribute the births in two categories depending on their occurrence during even or odd days and that the resulting difference is very appreciable. We will still believe that that difference is likely attributable to anomalies of chance. Indeed, from the very beginning the distribution made should seem to be extremely arbitrary and artificial. It does not correspond to any natural phenomenon or to habits of social life. Seasons, weeks, even phases of the Moon circulate indifferently to even and odd days. We therefore have every reason to believe that a notable difference corresponds to a distribution of categories that can occur fortuitously because of the multiplicity of the distributions and only discovered due to the patience of the calculator.

If, nevertheless, the difference persists in a new series of observations we will be obliged to admit, however bizarre it seems in advance, that the chance of a male birth is not the same for even and odd days of the month. The result of the first experience signals that the pertinent distribution should be checked by a new series of trials and, since there is an infinite multitude of possible distributions according to arbitrarily selected indications, it would have been extraordinary for randomness all by itself to lead twice in succession to a notable difference.

**113.** It follows that the probable judgement, pronouncing that an observed deviation is not attributable to anomalies of chance, results from two elements. One of them can be precisely and mathematically determined; it is the rate denoted until now by *P* of the fortuitous combinations which provide a smaller deviation for a randomly selected distribution. The other element is the preliminary judgement according to which we consider the distribution leading to the observed deviation as one of those among their infinite possible multitude which it is natural to study, but not only because it is one of those to which the observed deviation turns our attention.

That preliminary judgement, which in our opinion should direct statistical studies to some distribution rather than to another, is based on motives which can not be rigorously appreciated and can be differently appraised by different minds. It is a conjectural judgement also founded on probabilities which can not be reduced to enumeration of chances and whose discussion does not properly belong to the doctrine of mathematical probabilities. We will return to them in our last chapter.

Although a judgement made after inspecting statistical tables includes a variable element unyielding to precise measure, we should guard ourselves against concluding that the mathematical theory of chances is useless for statisticians. It is evident that the importance [the significance] of the deviation δ as an observational fact depends on its magnitude and the number of applied observations. But what is the law of that dependence? Only the theory of chances can tell us and show us how to calculate the ratio *P* corresponding to δ. As to probability denoted above by $\Pi$, in its application to statistics it really has no objective consistency. It does not at all measure the chance of verity or error inherent in a definite judgement.

**114.** We do not conceal the delicate in all the discussions and wish to multiply the examples throwing light on it. Suppose that we have the distribution of births by sex for France in its entirety and for a department in particular. There will be a certain deviation δ and a prior probability $P$[11] that the deviation is smaller than the observed provided that the chance of a male birth is the same in both cases and a posterior probability $\Pi$ that that chance differs. However, for that latter probability to have an objective value the department should be selected randomly […], all the deviations should be purely random and the real chance of a male birth, should not really change from one department to another.

In many respects naturally interesting for the statistician exceptional conditions in the department of Seine can essentially influence the chance of a male birth and the ratio $\Pi$ will once more assume its original sense. Indeed, had it been probable that a game of chance all by itself produces for one of the departments such a large deviation, as provided by observations in the department of Seine, it would be extraordinary that chance directs that deviation to that department […]

However, it is unknown whether what surprises all minds when discussing the department of Seine, will surprise everyone, or surprise to the same extent, when discussing the department of Corsica, or du Nord[12], or many others which also seem in advance to be placed in exceptional conditions. How to estimate the value of statistical experience and the derived posterior probability for each department separately? Here evidently enters a variable element resisting mathematical determination.

A similar remark is applicable to the study of yearly births. If I choose a year by chance and derive its notable deviation from the mean of many successive years, there will be a certain probability that the mean chance for that year is not the same as for all the series of collected observations. However, if we only turned our attention to that year after studying each year separately, since its deviation was the largest, the probability of a variation of the chance will become very different. It can retain its primitive value if that year was distinguished by a large climatic perturbation, high prices, [unusual] morbidity and therefore deserved preferential statistical study.

**115.** When observations are chronologically classified, as in the example above, it is impossible to extend the same statistical investigations and the obtained results should be accepted with their inherent indeterminateness. However, when having reason to believe, or when the observations themselves indicate that in another system of classification the chances do not vary in time, we will dispense with all the subtle distinctions. When indefinitely continuing the observations, we will determine with any desired precision the pertinent proper chances for each system of categorical division.

With the increase in the number of observations, the number of categories can be multiplied while being invariably guided by the notions acquired in accordance with the conditions of the studied natural or social phenomenon. And the complicated fact, the object of the first numerical determination (§ 106), will be gradually decomposed into its elements. The Bernoulli principle to which we at last always return as to the only solid basis of all the applications of the theory of probability[13] will be sufficient for the statistician whereas mathematical formulas incessantly provide the measure of the attained degree of precision.

**116.** Those formulas, as provided throughout our book, are only approximate but their precision usually suffices if observations are counted by the hundred. The most eminent authors, notably Laplace, did not hesitate to apply them in such cases (§13, Note 3; § 33, Note 2; and § 69, Note 5). And even when the number of observations is too small for allowing this practice, those formulas can be suitable for revealing variations of chances from one category to another. However, in such cases numerical calculations of probabilities should be done by formulas whose complication often renders their usage very tiresome.

The applied number of observations can also be too small for providing very high probabilities to the mentioned variation of chances between the categories even when making use of the total series of observations. However, the manner of decomposing such series can render the existence of those variations very likely. Thus, we will obtain a certain posterior probability $\Pi$ of an inequality of chances when considering the difference $\delta$ between two series of 50 observations of each category, and another probability $\Pi'$ when the difference between two series of 150 observations each is $\delta'$ which can be either larger or smaller than $\delta$.

In general, however, $\Pi'$ will exceed $\Pi$ although both can be of such an order that the difference will be explained by a fortuitous anomaly. On the contrary, when decomposing each series of 150 observations into 3 consecutive series of 50 observations each[14] the persistence of a difference of the same order for the pairs of partial series can render the hypothesis of a succession of fortuitous anomalies so unlikely that quite reasonably all doubts about the existence of an inequality of the chances of the same event in the two categories disappears[15].

**117.** It is often insisted on the need to combine a very large number of observations for arriving at appreciably fixed results ridden from all the irregularities of chance. And, as it is very often borne in mind, we should distinguish fortuitous influences affecting each observation independently from all the others of the same series, and other influences dominating all the observations of a series or its part but nevertheless being fortuitous since irregularly varying from one series to another (§ 104). The effects of their variation are compensated when collecting a large number of series and therefore a very large number of individual cases.

Among the causes of solidarity of causes or influences which dominate diverse trials of the same randomness and thus demanding an accumulation of a very large number of trials for the mean results to be stable we shall first of all mention the closeness of trials in space or time. We have seen in an example of § 79 borrowed from Bienaymé that if the value of a chance varies not fortuitously from one trial to another in a total series of *m* trials but is fortuitously determined in a partial series of $m_1, m_2, \ldots$ trials, it will not in general be sufficient to have a large number *m* of observations for the fortuitous anomalies to be considerably compensated. This aim will only be achieved when the numbers $m_1, m_2, \ldots$ are not so large, but the total series (m) is composed of a large number of partial series $(m_1), (m_2), \ldots$

Actually, the example of § 79 is not directly applicable to statistical facts. We can not in general suppose that the chances remain rigorously constant for a whole series of neighbouring trials and change sharply and fortuitously from one series to another. On the contrary, they are subjected to progressive modification and like all the magnitudes in nature even in their irregular variations usually obey the law of continuity. However, it is not always that in general chances differing but little dominate a large number of neighbouring trials, so that we can not at all regard as independent two trials immediately following each other or even as absolutely independent those not separated by such a large number of intermediates that the trace of the initial state of the chance is erased.

This is how the unevenness of the terrestrial surface or of the surface of an agitated sea is considerably independent from point to point at large distances whereas close points in spite of their irregular course necessarily remain on very little differing heights. We can not in advance reduce to formulas the influence of this connection between chances of neighbouring trials[16], but the fictitious examples as those in § 79 sufficiently prove that that influence can be very large and experience confirms this presumption.

In many applications to social [political] economy and to such branches of natural sciences as meteorology[17], we should accumulate observations in a number much exceeding those assigned by the formulas based on the hypothesis of independence of chances from one trial to another for obtaining considerably stable means.

**118.** A fact seemingly singular at first sight which authors did not fail to note[18] is that things, engendered by the developing activities of man, which apparently result from a multitude of very complex causes such as the ratio of the number of the accused for committing crimes to the number of the inhabitants of the country, or the rate of conviction, experience less annual variations than things depending on the coincidence of blind forces of nature. However, after contemplation that result ceases to surprise.

It is easily understood that there only exists little or no solidarity between the causes whose coincidence determines the perpetration of different crimes, or between condemning different accused, whereas there evidently exists a very strong solidarity between causes whose fortuitous coincidence leads to rain at the same place today and tomorrow. It is therefore quite understandable that for things, depending on the individual activity of man, the mean chances appear more stable and actually less often experience irregular perturbations. On the contrary, there is every reason to believe that the succession of slow transformations of the social situation leads to their secular variations not in general observed in physical phenomena either because they do not exist or occur extremely slowly.

**119.** Even when there is no solidarity resulting from neighbourhood in space or time between chances dominating each trial, the differences that reveal variations in the mean values of chances from one series to another which actually dominate each trial do not necessarily testify to a variation in the general system of fortuitous causes on which depend the phenomena studied by statistics. We (§ 109) have remarked that the proportions in which categories (a), (b), (c), … unite in forming the series (m), (m′), … can vary from one series to another owing to causes acting irregularly and fortuitously in differing series or their parts although influencing at once all, or a part of observations of a series and therefore affecting the value of the mean chance. Thus, the category of the accused for illegally cutting timber can increase if the year or the winter was severe and firewood expensive. Also increased can be the number of those brought to court for brawls during years when the price of alcoholic beverages is low and pubs are much frequented, etc.

Since the rate of conviction is different for each category of crime, its mean value for the totality of the yearly accused certainly varies. However, the causes of that variation are reasonably thought to be accidental and fortuitous, their effects compensate each other if a few years are combined. On the contrary, no one reckons changes of criminal legislation among fortuitous causes since they suppress a class of crimes or submit it to a court of lower instance.

**120.** The most important goal of statistics is the investigation of causes dominating physical and social phenomena[19]. By their enormity, the accumulated numbers satisfy the conditions for stabilizing the means, and it is now more necessary to decompose the fortuities from each other (§§ 79, 106, 109) and, so to say, clean the conditions of randomness. Individual cases should be agglomerated into series for the sole aim of compensating the effects of causes acting quite independently on each individual case. It will be therefore not necessary anymore to work with very long series.

In any case, if such studies sometimes lead to error because of the caprice of chance, they ordinarily provide veritable consequences. We should not be deprived of the most fruitful means of investigation because it does not guarantee absolute certainty; for that matter, it is rarely attained by man. We should therefore multiply the number of categories and choose observations sufficiently close to each other in space or time so that in each series the variations of mean chances will be inconsiderable.

### Notes

**1.** Statistics originated in mid-17th century (Graunt, Petty). [B. B.]. Cournot's opinion about ancient times should really mean that ancient science was qualitative rather than quantitative. O. S.

**2.** From at least the time of Tycho Brahe and Kepler the arithmetic mean had been understood as an estimator preferable to a single observation because of that same compensation, and in 1756 – 1757 Simpson quantitatively confirmed that belief for two distributions. Cournot, however, repeatedly mentions compensation as though it invariably and unconditionally happens in a large number of observations. Again, sums of random variables, unlike their means, do not compensate.

**3.** The words *cause*[4]*, action*, etc are understood here in the widest possible sense just like in ordinary language but not as rigorously as it is sometimes necessary in metaphysical analysis. […] A. A. C.

**4.** In 1774, Laplace defined *cause* relative to event E as a thing providing the proper chance to the appearance of E. [B. B.]

**5.** The next section can be included without changes in a modern work. B. B.

**6.** How to define a category? Cournot's reference to § 75 is only partly helpful. See also § 119. A few lines below he denoted series and categories in the same way!

**7.** Such causes were generally known in geodesy and practical astronomy perhaps even in antiquity. In triangulation, morning and afternoon observations often lead to systematically differing results. Having the same number of both, the observer can expect the mean to be essentially free from the ensuing error. Cf. Note 2.

**8.** A usual practice for unscrupulous investigators. Categorically forbidden in geodesy (at least in Soviet Russia in mid-20$^{th}$ century).

**9.** Bru referred in this connection to Villerme, 1831.

**10.** Here, Bru referred to Aristotle's *De generatione animalium*. Cournot's opinion a few lines below is unfortunate in that meteorological differences strongly influence the other mentioned factors.

**11.** If the deviation is obtained by comparing a partial and a total series the expression for *P* can not remain as it was provided in § 108 when two partial series were compared. It should now be taken from § 101, but this particular remark does not change anything in general reasoning. A. A. C.

**12.** Cournot had apparently selected departments differing in climate, and Bru remarked that Aristotle (Ibidem), Montesquieu and Buffon had discussed the influence of climate on the sex ratio at birth.

**13.** Following Bienaymé, Cournot did not even mention Poisson's law of large numbers. [B. B.]

**14.** Concerning the decomposition of observational series, Bru referred to Fourier, Bienaymé and Quetelet.

**15.** Let *m* be the number of observations and *n*, the number [of the occurrences] of event A in a total series. If $m_1$ of the observations are collected at random to form a partial series, the probability that it contains $n_1$ events A is evidently the same as when extracting without replacement $n_1$ white balls in $m_1$ drawings from an urn containing *m* balls, *n* of them white. That probability is equal to (§ 36)

$$\frac{m_1(m_1-1)...(n_1+1)n(n-1)...(n-n_1+1)}{(m_1-n_1)!m(m-1)...(m-m_1+1)} \times$$
$$(m-n)(m-n-1) \ldots [m-n-(m_1-n_1+1)]$$

and according to that section, if all those four numbers are very large, it will become likely that $n_1/m_1$ very little differs from $n/m$.

If the chance *p* of event A does not vary during the observation of the series (m), it will be possible to form the series ($m_1$) from the $m_1$ first observations of series (m) which is the same as though done fortuitously. And again for sufficiently large numbers *m, n, $m_1$* and $n_1$ the ratio $n_1/m_1$ will very little differ from $n/m$ even if that latter ratio notably deviates from the chance *p* because of an anomaly extremely low probable in advance. Bienaymé (1840) showed it by an elegant calculation which nevertheless we believe unnecessary since his conclusion was evident according to the preceding reasoning.

It is not necessary either to include the condition that the chance *p* remains constant; suffice it that its variations, if they occur, take place not chronologically. To prove this, suppose that the total series (m) is composed of two series, (m′) and (m″) made, for example, in different places and dominated respectively by chances *p*′ and *p*″. We imagine as previously that the total series (m) is disposed chronologically so that $m_1$ first observations form a category in which $m_1'$ are dominated by chance $p'$, and $m_1''$ observations dominated by chance $p''$. If the numbers $m'$, $m''$, $m_1'$ and $m_1''$ are of a suitable order, the ratio $m_1'/m_1''$ will little differ from $m'/m''$.

If, for example, we chronologically order the yearly births in two departments, the ratio of those births can be thought to be essentially the same for the first six months and for the whole year. And extracting $m_1$ first observations from the total series (m) is tantamount to drawing at random $m_1'$ balls from an urn containing $m'$ balls, $n'$ of them white and $m_1''$ balls from another urn containing $m''$ balls, $n''$ of them white. If $m_1' + m_1'' = m_1$ and $m_1'/m_1''$ very little differs from $m'/m''$ and if, as previously, the numbers are sufficiently large, the ratio $n_1/m_1$ will with high probability very little differ from $n/m$ even if the ratios $n'/m'$ and $n''/m''$ considerably deviate from the chances $p'$ and $p''$ owing to an anomaly very unlikely in advance.

The remarks which were the object of this Note do not at all refute what was said in the main text; in our opinion, they are even quite independent. And it is not less permissible to replace posterior probabilities calculated from the total series by the product of probabilities obtained when separately considering the partial series. However, it can and even will happen when there are very large numbers that the posterior probabilities calculated for two different series will be very close to each other since the results obtained for the total series will with high probability lead to similar results for each partial series. A. A. C.

**16.** The study of dependent variables properly began with Markov.

**17.** This statement is borrowed from Bienaymé (1839, p. 188). [B. B.]

**18.** Bru referred here to Bienaymé (1839, p. 187) and Poisson (1837, p. 11).

**19.** Bru called this *passage* remarkable and stated that that appraisal of the aims of statistics was only *definitively adopted* in the second half of the 20[th] century (Neyman 1960). The same problem much interested Chuprov from the very beginning of his scientific career (Chuprov 1903).

## Bibliography

**Bienaymé I. J.** (1839), Théorème sur la probabilité des résultats moyen des observations. *L'Institut*, 284, t. 7, pp. 187 – 189.

--- (1849), Principe nouveau … sur la constance des causes … *L'Institut*, 333, t. 8, pp. 167 – 169.

**Chuprov A. A.** (1903, Russian), Statistics and the statistical method. Their vital importance and scientific aims. In author's *Voprosy Statistiki* (Issues in Statistics). Coll. Articles. Moscow, 1960, pp. 6 – 42.

**Neyman J.** (1960), Indeterminism in science and new demands on statisticians. *J. Amer. Stat. Assoc.*, vol. 55, pp. 625 – 639.

**Poisson S.-D.** (1837), *Recherches sur la probabilité des jugements*. … Paris, 2003. English translation: www.sheynin.de downloadable file 53.

## Chapter 10. Experimental Determination of Mean Values by Observations and Formation of Tables of Probability

**121.** In statistics, mean values are determined for two different reasons which are preferably distinguished[1]. The means are often magnitudes whose knowledge is immediately interesting since their values directly influence physical and social phenomena. For example, the mean quantity of grain produced in a country directly influences its population and its entire economic system. The same can be said about the mean value provided by a tax, or import or export of some foodstuffs.

However, even more frequently means are only considered as results considerably independent from the oscillations of chance so that their variations can indicate more or less surely and rapidly the existence of changes in the intensity or manner of action of regular causes. Suppose for example that a general census of the population of a country is made and we calculate the mean age of all the inhabitants. In itself, this mean is not significant because there is no social fact directly depending on its value. Nevertheless, that value depends on the law according to which the population is distributed by age and the chances of longevity provided by the climate and habits of the people. Changes in all these circumstances or in one of them will be revealed by a change in the mean.

Suppose also that we extract the mean value of the stature of conscripts from the pertinent table[2]. If some new circumstance tends to improve or worsen the physical state of the country's population, we will discover it by the change in that mean. For the mind research is facilitated when it considers a simple expression, even if only leading to incomplete and indirect consequences, rather than a set of complicated facts.

**122.** It is natural to ask whether values [estimators] differing from the ordinary mean will not in some cases fulfil the aim raised by statistics better, more promptly get rid of the oscillations of chance, more properly indicate the influence of constant and perturbative causes, or even more properly ensure comparisons in certain problems in law and political economy not to mention all statistical investigations[3].

Thus, in France the law of 15 May 1818 prescribed that for estimating the tax on transfer of property the mean of market price-lists for the previous 14 years should be calculated after excluding four extreme values, two maximal and two minimal[4]. The mean of the rest 10 values then constituted the mean price of a *usual* year. The legislature perhaps wished to act in the interest of both the state and the taxpayers by protecting the collection of the tax as much as possible from the influence of chance. Whether that [probable] goal will be better attained in that way or by keeping to the ordinary mean, can only be known by comparing a large number of market price-lists calculated according to both methods.

Suppose that instead of 10 or 14 particular values we have 1000 or 1400. Let *aib* (Fig. 4) be the curve whose ordinates measure the pertinent probabilities and *OA* and *OB*, the extreme [particular] values. The usual mean will be close to *OG* with the ordinate *Gg* passing through the centre of gravity of ABbga (§ 67). After cutting off portions *ACca* and *BDdb* with the area of each of them being 1/7 of the total area the mean taken in accordance with the law of 1818 will be close to *OG′* to which corresponds ordinate *G′g′* passing through the centre of gravity of the rest area *CDdc*. However, depending on the form of the curve it can happen that the modulus of convergence (§ 69) for the rest area will considerably exceed the modulus for the entire area. The fortuitous deviation in both directions from the mean *OG′* will then be contained in a narrower interval than the fortuitous deviation in both directions from *OG* in spite of the number of the particular values applied for an approximate determination of the latter exceeding that of the former in the ratio of 7/5.

Denote by *OI* the median value of the total area (§ 68), i. e., the value whose corresponding ordinate *Ii* divides the total area in two equal parts. The larger are the cut off parts *ACca* and *BDdb* (remaining equal to each other) as compared with the total area, the closer will ordinate *G′g′* approach the ordinate *Ii*. It follows that in general a system similar to that which the law of 1818 provides for the *usual* value, and when the number of the applied particular values becomes very large, it will result in an intermediate value between the mean in the proper sense and the median value[5].

During the years of very high prices the price of foodstuffs such as grain tends to rise above the mean much more than it lowers below that mean in the years of abundance. This is the same as saying that the mean *OG* exceeds *OG′* and *OI*, so that the law of 1818 is more favourable for the taxpayer than it would have been if the ordinary mean were applied. However, independently from this circumstance, when considering that the price of foodstuffs usually oscillates within sufficiently narrow limits and only appears beyond them under the influence of violent and transient perturbative causes, we understand that the legislator wished to eliminate completely that influence by rejecting the extreme values which can not reckon in transactions between citizens and are wholly beyond the economic system of the country.

**123.** When means are determined for diverse parts of a complicated system we should attentively note whether they are compatible. When the means are determined separately one from another, the system can become impossible (§ 74). If for example a triangle should remain right when its sides vary, each side will have a mean value, but, taken together, they will not be compatible with a right triangle. […] If the sides and the angles of a triangle take various values, the mean values of the angles will remain compatible in that their sum will be equal to two right angles (? - O.S.) and a triangle, or even infinity of similar triangles can be thus constructed whose angles will be those mean values.

The means of their sides, if only each mean is smaller than the sum of the other two, will apparently also belong to possible triangles. However, there will generally be no triangle with those mean values of angles and sides, and the mean area of each triangle will not coincide with the area of the triangle having mean sides, etc. Just the same, when measuring the dimensions of the diverse organs of many animals of the same species we can, and likely will arrive at mean values incompatible either one with another and with the conditions of the species' life.

We insist on that very simple remark because it seems to have been overlooked in a work otherwise really meritorious in which the *homme moyen* was defined and determined by a system of means derived from the measurements of stature, weight, physical strength, … of a large number of individuals[7]. Thus defined, the *homme moyen* is not at all the type of sorts of mankind but simply impossible; at least until now nothing authorizes us to think him possible.

**124.** Let *aib* (Fig. 1, § 67) be the curve representing the law of probabilities of magnitude $x$. Each ordinate such as *Ii* is proportional to the probability of the particular value measured by the corresponding abscissa *OI* (§ 65). For magnitudes applied in statistical research that law of probabilities is in general unknown in advance. If the number $N$ of the observed particular values is extremely large, the number $n_1$ of them situated between $OA = a$ and $OA_1 = a_1$ will approximately be to $N$ as the area $AA_1a_1a$ to the total area *Abba* although the interval $AA_1$ can be very small as compared with the interval $AB$ between the largest and the smallest possible values.

In other words, if the total area *Abba* is unity, $n_1/N$ will be the approximate measure of the partial area $AA_1a_1a$. If the interval $AA_1$ or the difference $(a_1 - a)$ are very small, the partial area can be considered as an area of a rectangle. […] Therefore, when dividing $AA_1a_1a$ or the number $n_1/N$ by the difference $(a_1 - a)$, the quotient can be considered without an appreciable error as a measure of the very little differing ordinates $Aa_1$ and $A_1a_1$. More precisely, that quotient numerically expresses the value of the ordinate corresponding to a point situated on $AA_1$ at equal distances from $A$ and $A_1$.

It is thus possible to determine empirically the law of probabilities of magnitude $x$ by a sufficient number of statistically provided values. To achieve this, we decompose the total number $N$ in numbers $n_1$, $n_2$, … of partial values situated between $a$ and $a_1$, $a_1$ and $a_2$, … The set of quotients

$$n_1/N(a_1 - a),\ n_2/N(a_2 - a_1),\ \ldots \tag{124.1}$$

corresponding to the values of $x$ equal to

$$(a + a_1)/2,\ (a_1 + a_2)/2,\ \ldots \tag{124.2}$$

constitutes a *table of probabilities* and determine as many points of the curve *aib* as there are terms in the series above. These points can be joined by a continuous curve thus providing a graphical representation of the law of probabilities. By calculation we can also get a curve with algebraically connected ordinates and abscissas passing through all the points determined by the table; indeed, mathematical analysis furnishes general pertinent formulas. Whether constructing it graphically or obtaining it by calculation, it is evident that, for the abscissas situated between two consecutive terms of the series (124.2), we should regard the values of the ordinates as almost precise, if, as we suppose, they are very close to each other and, moreover, if the differences between the corresponding terms of the series (124.1) are sufficiently small as compared with one or another of these terms.

Therefore, each of the intervals $(a_1 - a)$, $(a_2 - a_1)$, … should be a small fraction, for example, 1/100 of the range of possible values; if that range is not given in advance, 1/100 part of the interval between the largest and the smallest observed particular values. For constructing and applying the table, it is convenient to take equal small differences $(a_1 - a)$, $(a_2 - a_1)$, … However, when indicating that in certain parts of the series (124.1) the differences between consecutive terms become too large for the second of the abovementioned hypotheses to be admissible, the interval between the corresponding parts of series (124.2) should be shortened.

**125.** For the value of some term of series (124.1), for example, of the $i$-th term, to be appreciably independent from the anomalies of chance, the corresponding number $n_1$ should be considerable, at least of the order of hundreds. However, the probability of a notable error of each term of series (124.1) taken separately can be very low although because of their large number there will be a high probability that at least one of them is affected with a considerable error. Suppose that one of those terms essentially deviates from the law, apparently followed by its neighbouring terms situated both before and after it, without there being any reason to suspect that that law underwent a sharp change in the vicinity of the corresponding term of series (124.2). We should not therefore hesitate to reject the anomalous term from series (124.1) and to replace it temporarily (? - O.S.) by a suitable value answering the law followed by those neighbouring terms and derived by the known interpolation formulas.

It nevertheless follows that for sufficiently guaranteeing that the random errors affecting each term in the table are contained within very narrow limits, the number $N$ should be very large, much larger for example than would be sufficient for deriving with high precision the mean value of magnitude $x$ or, under ordinary circumstances, the possibility of a simple event. When indicating that in addition the laws of probabilities studied by the statistician can experience notable perturbations during the time needed for accumulating such a large number of observations; and when reckoning with the chances of possible errors corrupting the treatment itself of observations, it will be understood that the construction of a table of probabilities is the most difficult work, as though a masterpiece of statistics. The numbers that can not be thought trustworthy in one study, can lead to concordant results in a new series of observation. Tables of probabilities with narrow intervals are only available for the duration of human life and it is necessary that they, being in agreement with each other, provided that perfect sureness which we discussed.

**126.** Let us return to the determination of mean values. If we constructed a table of probabilities having $N$ particular values, sufficiently large for considering random errors affecting each term of the table insignificant, the mean value $M$ will all the more be free from any appreciable error. A new series with a large number of $m$ particular values though much smaller than $N$, will provide another mean value μ possibly affected by an appreciable error. Since *M* tangibly coincides with the rigorous mean, (*M* − μ) will be the random error affecting the determination of the mean μ.

Before making a new series of observations, the probability *P* that the difference (*M* − μ) will be contained within the limits ± *l* is (69.1)

$$t = lg\sqrt{m}$$

where *g* is the value of the modulus of convergence which is connected with the law of probabilities of the magnitude *x* and can be calculated by two different methods. The first of them presumes that the preliminary construction of the table of probabilities consists in algebraically expressing by interpolation (§ 124) the ordinates and abscissas of the curve of probabilities or the function denoted by *fx* in Note 3 to § 69 and apply to it the formula of integral calculus indicated there as though observations effectively determined the continuous sketch of that curve rather than a finite number of points through which or very close to them it should pass.

The second method consists in immediately applying the system of particular values provided by observation. No interpolation, never exempt from arbitrariness, is needed. And here is a very simple pertinent rule[8]:

*Calculate the sum of the differences between the mean value and each of the particular values; divide the number of those values by twice that sum and extract the square root from the quotient. The result will be the modulus of convergence.*

For writing this rule algebraically, denote by $x_1, x_2, \ldots$ the *N* particular values. Then

$$g = \frac{\sqrt{N}}{\sqrt{2[(x_1 - M)^2 + [(x_2 - M)^2 + \ldots]}}.$$

We can also write (Note 3 to § 69) that

$$g = \frac{N}{\sqrt{2[(x_1 - x_2)^2 + (x_1 - x_3)^2 + (x_2 - x_3)^2 + \ldots]}},$$

$$g = \frac{1}{\sqrt{2[(x_1^2 + x_2^2 + \ldots)/N - M^2]}}$$

or even express *g* in other forms which can be successfully applied according to circumstances.

When *N* is very large, the calculation of the squares of all the differences is impracticable. However, in applying the rule, all the particular values contained within interval ($a_1 - a$) supposed to be very short can be regarded as being equal one to another and to ($a_1 + a$)/2. The number of the particular values is $n_1$; multiply it by $[(a_1 + a)/2 - M]^2$ and the sum of all such products for all the intervals ($a_1 - a$), ($a_2 - a_1$), … can be taken as that sum of the squares of differences between *M* and each included particular value.

Still more exactly, multiply by $n_1$ the square $(\alpha_1 - M)^2$ where $\alpha_1$ is the mean of all particular values contained between $a$ and $a_1$ and repeat the same for each partial interval. For determining the modulus of convergence this second method does not suppose a preliminary construction of the table of probabilities. Actually, this modulus can, like the mean, and for the same reason, be determined with a sufficient precision by a number of particular values not at all sufficient for precisely constructing that table.

**127.** It follows from that remark that, when denoting by $\xi_1, \xi_2, \ldots$ the $m$ particular values in the new series of observations, the number

$$\gamma = \frac{\sqrt{m}}{\sqrt{2[(\xi_1 - \mu)^2 + (\xi_2 - \mu)^2 + \ldots]}} \qquad (127.1)$$

will little differ from $g$, and the two fractions, $(g - \gamma)/g$ and $(M - \mu)/M$, will be in general of the same order of magnitude. The mean $\mu$ is, according to the hypothesis, very close to the mean $M$ and $\gamma$ is also approximately equal to the modulus of convergence $g$.

Suppose now that we only have the series of $m$ particular values $\xi_1, \xi_2, \ldots$ insufficient for determining the mean $M$ without an appreciable error but quite enough for the random error $(M - \mu)$ to be numerically very small. The limits $\pm l$ within which the error is contained with probability $P$ is provided without an appreciable error by the formula $t = l\gamma\sqrt{m}$. Now, $l$ is a very small magnitude and $\gamma$ also differs from the modulus $g$ by a very small magnitude, and the error that affects the value of $l$ because of replacing $g$ by $\gamma$ is very small even as compared with $l$, of an order neglected in approximate calculations. Therefore the application of mathematical symbols allows to present the demonstration in a more rigorous form although the basis of reasoning does not change[9].

For the same value of $P$ the interval $2l$ of the limits of the error is inversely proportional to $\gamma\sqrt{m}$, see (127.1). That product can be therefore taken as the measure of precision with which the unknown $M$ is determined when assumed to be equal to $\mu$ as provided by the system of $m$ particular values $\xi_1, \xi_2, \ldots$ In other words, that product is the *weight* of $\mu$. It would be appropriate (§ 107) to accompany the means $\mu$ in statistical tables by the corresponding weights $\gamma\sqrt{m}$ determined by the system itself of observations.

**128.** A new series of $m'$ particular values will provide a mean $\mu'$ differing from $\mu$ although the law of probabilities is not supposed to change from one series to another. And we have, after accomplishing the series (m) and before the series (m′) is observed, the probability $P$ that the difference $(\mu - \mu')$ will be contained within the limits $\pm l'$ as [indirectly] provided by the formula:

$$t = \frac{l'\gamma\sqrt{mm'}}{\sqrt{m + m'}}.$$

If the total series (m) is decomposed in two partial series formed by $m_1$ and $m_2$ particular values with magnitudes $\mu$ and $\gamma$ becoming respectively ($\mu_1$, $\gamma_1$) and ($\mu_2$, $\gamma_2$), and supposing that the law of probabilities of magnitude *x* is the same in both partial series, we will have in advance, before determining by experience the numbers $\mu_1$, $\mu_2$, $\gamma_1$, $\gamma_2$, probability *P* that the difference ($\mu_1 - \mu_2$) is contained within the limits ± *l'* given by formula

$$t = \frac{l'\sqrt{m_1 m_2}}{\sqrt{m_1/\gamma_1^2 + m_2/\gamma_2^2}}.$$

Now we can calculate the posterior probability Π that ($\mu_1 - \mu_2$) = δ as given by experience indicates a change in the law of probabilities between the series ($m_1$) and ($m_2$).

Instead of comparing the two partial series with each other we can compare ($m_1$) and (m) and the analogy formulated in § 101 will sufficiently indicate that the limits of the difference corresponding to probability *P* will be determined from the formula[10]

$$t = \frac{l'\gamma\sqrt{mm_1}}{\sqrt{m - m_1}}.$$

**129.** All the discussion above concerning the interpretation of changes occurred or seemed to have occurred in the chances of an event during the passage from one series of observations to another is evidently applicable to changes in the mean values. We will not reproduce it and restrict our considerations by remarking that any existing traces of solidarity of trials close to each other in space or time among the diverse fortuitous determinations of the same magnitude can be revealed. When combining all the particular values taken two from two, the mean of the squares of the obtained differences will depend on modulus *g* or coefficient γ (§ 126). Traces of solidarity still remain if it notably exceeds the mean of the squares of the differences obtained when combining each particular value with its preceding or immediately following, or with those situated sufficiently close[11].

### Notes

1. On the two types of means (see below in text) see Sheynin (2007, § 5).
2. Two authors had discussed the stature of conscripts before Cournot and Quetelet (1846) studied it after him. [B. B.]
3. Many scholars of the 18[th] century studied this problem. [B. B.]
4. In the theory of errors, rejection of outliers is a most delicate operation, in particular because of the unavoidable systematic errors.
5. From Note 5 in Chapter 6 it follows that the empirical determination of a mean value always gets more rapidly rid of anomalies of chance than the median value. This, however, does not prevent the empirical determination of a value situated by its definition between the values of the mean and the median to get rid of those anomalies still more rapidly[6]. A. A. C.
6. For some distributions, and when the distribution is unknown, the median is more reliable than the mean.
7. Bru remarks that Cournot indirectly referred to Quetelet (1835) and that the latter (vainly) objected to criticisms.
8. Cournot apparently followed Fourier (1826/1890, p. 532). [B. B.]

**9.** See Poisson (1829). [B. B.]
**10.** Cf. § 101. [B. B.]
**11.** Bru remarked that Cournot had thus offered the first test of independence of a sequence of observations based on the coefficient of autocorrelation.

## Bibliography

**Fourier J. B. J.** (1826), Sur les résultats moyens déduits d'un grand nombre d'observations. *Œuvres*, t. 2. Paris, 1890, pp. 525 – 545.

**Poisson S.-D.** (1829), Sur la probabilité des résultats moyens des observations. *Conn. des temps* pour 1832, pp. 3 – 22.

**Quetelet A.** (1835), *Sur l'homme*. Bruxelles.

--- (1846), *Lettres sur la théorie des probabilités*. Bruxelles.

**Sheynin O.** (2007), The true value of a measured constant and the theory of errors. *Historia Scientiarum*, vol. 17, pp. 38 – 48.

# Chapter 11. Means of Measurements and Observations

**130.** The theory of convergence of the mean values [to their limits] is applicable to a problem of great importance for all the physical sciences: to determining the probable limits of the error of a numerical result when taking the mean of a large number of values corrupted by some errors.

Denote by $fx$ the function expressing the probability of error $x$ in the measure of magnitude $a$; by $\varepsilon$, the mean of the absolute values of the error $x$; by $g$, the modulus of convergence whose value implicitly results from the form of the function $fx$; by $m$, the number of the particular values from which the mean was taken and which we suppose to be sufficiently large for applying our approximate formulas (Note 5 in Chapter 6); by $\alpha$, the mean of the particular measures $a_1$, $a_2$, … It will converge to a fixed magnitude $a + \varepsilon$ and there will be probability $P$ that the fortuitous error $(a + \varepsilon - \alpha)$ will be contained within limits $\pm l$ determined by formula (69.1). It follows that if $\varepsilon = 0$, by virtue of the function $fx$ $P$ will be the probability that the error still corrupting the mean $\alpha$ is contained within the limits $\pm l$.

When presuming that the constant $\varepsilon = 0$, it is supposed that in general the curve with ordinates $fx$ is symmetric with regard to the $y$-axis […]. Indeed, this second hypothesis, more particular than the first one, necessarily demands that the error of the mean value and the median is zero (§ 68). On the other hand, if that symmetry is lacking, such a concurrence of circumstances is needed for the disappearance of $\varepsilon$ that this case can be considered quite unlikely.

If the curve of probabilities is not symmetric with respect to the $y$-axis, and the constant $\varepsilon$ takes an appreciable value, it is said that a cause of a constant error corrupts the series of measurements due to a defect in the construction of the applied instruments, or in the observer's organs of sense, or his manner of work. A series of measures thus affected should be rejected as being improper for determining the real value of magnitude $a$. The experimenter's sagacity principally consists in discovering the means for getting rid of the influence of such causes of errors, for analysing, studying and eliminating their effect.

It is doubtless infinitely unlikely that the constant $\varepsilon$ is strictly zero, however thoroughly are the instruments manufactured and the measurements themselves taken for avoiding all the causes which render errors of one sense more (or less) probable than those of the contrary sense. In general, mathematical absoluteness cannot be achieved in those pursuits which depend on organs of sense and the interaction of man and the material world. Had it been otherwise, when sufficiently increasing the number of observations, formula (69.1) would have led to an arbitrary precision of the measured magnitude $a$; to its expression, for example, by 20 correct decimals just as in the case of an incommensurable root or the ratio of a circumference to its diameter. Such a consequence is absurd, and we will return to the discussion of the causes that necessarily restrict the precision of measurement of each kind of magnitudes whichever is the number of their observations.

In the first place we suppose, like everyone who had treated this problem expressly or tacitly did, that the constant ε = 0 or at least negligible as compared with the error which can corrupt each particular measure or with the mean of these errors without considering their signs. It is thought possible to determine after observation, simply by inspecting the particular values, whether they are compatible with the hypothesis of symmetry corresponding to the condition ε = 0. If, for example, $m$ is at least of an order of hundreds, and the ratio of the number of the particular values exceeding α to $m$ appreciably differs from 1/2, we are notified that the median value does not coincide with the mean value as demanded by the hypothesis of symmetry. Other, infinitely many other methods of testing can be proposed.

When the measures are taken with sufficient care, the probability of error $x$ should decrease very rapidly as $x$ numerically increases. And the probability of an error not very small as compared with $a$, for example its 1/20 or 1/30 part, should be regarded as very small or as approximately zero. The curve of probabilities then becomes bell-shaped. We can not assign a limit at which the chances of error become exactly zero (§ 66), but the chances of the values of $x$ being beyond certain limits are so slim, that we are quite authorized to neglect them[1].

**131.** In general, the form of the function $f$ is unknown and can not be assigned in advance. However (§ 127), for a large number of measurements it is permissible to assume the number

$$\gamma = \frac{\sqrt{m}}{\sqrt{2[(a_1) - (\ ^2 + a_2) - \ ...^2] +}}$$

given by those observations themselves as a sufficiently close value to the modulus $g$. After thus calculating γ, the product $\gamma\sqrt{m}$ can accompany the mean α as measuring its *weight* (§ 127).

**132.** When diverse measures are taken by different observers by different instruments or under dissimilar circumstances, the probability of error $x$ will vary from one measure to another. In such cases the function $fx$ should be understood as a mean of all values provided by the different laws of probabilities proper for each measure as I explained in a general manner in § 81.

Suppose that a total series of $m$ measures taken by the same observer, or resulting from observations made by the same instruments, or under similar circumstances are grouped in partial series denoted by subscripts 1, 2, …, $i$, then the weights $\gamma\sqrt{m}$ will be expressed respectively by $\gamma_1\sqrt{m_1}$, $\gamma_2\sqrt{m_2}$, …, $\gamma_i\sqrt{m_i}$.

Assuming

$(α_1 + α_2 + … + α_i)/i$, or, less defective, $(m_1α_1 + m_2α_2 + … + m_iα_i)/m$

as the value of magnitude $a$, we do not properly take into account the weights of each particular result. If, for example, these values in series $(m_1)$ very little differ from their mean $α_1$, much less than the values in ($m_2$) differ from their mean $α_2$, $γ_1$ will be much larger than $γ_2$ and $γ_1\sqrt{m_1}$ can much exceed $γ_2\sqrt{m_2}$ even when $m_1$ is appreciably smaller than $m_2$.

Common sense also tells us that a less numerous series composed of better concordant observations should inspire more confidence than a more numerous series manifesting larger differences. The theory ought to specify the indications of common sense and provide a formal rule for combining the partial results to narrow as much as possible (without changing the probability) the limits of the error which should be thought to corrupt the final result[2]. This rule consists in taking the mean

$$\frac{m_1 γ_1^2 α_1 + m_2 γ_2^2 α_2 + \ldots + m_i γ_i^2 α_i}{m_1 γ_1^2 + m_2 γ_2^2 + \ldots + m_i γ_i^2}$$

as the value of magnitude $a$, i. e. in entering each partial mean proportionally to the square of its weight and thus providing the general mean. We will then have probability $P$ that that error is contained within the limits $± l$ given by the equation [by the formula]

$$t = l\sqrt{m_1 γ_1^2 + m_2 γ_2^2 + \ldots + m_i γ_i^2}.$$

**133.** The required magnitude often depends on many other directly measured magnitudes. As a simplest example, suppose that it is required to calculate the height of a tower by solving a right triangle whose base beginning at the tower's foot and the adjacent angle are measured. The height is a function of those elements and its error depends on the errors of both measurements. Or, when calculating the area of a triangle as a function of its sides […], its error depends on the errors of the latter. Finally, imagine a network of triangles like those constructed during large geodetic operations. The error of one of its calculated sides results from the errors of base and angle measurements.

Denote by $a, b, c, \ldots$ the magnitudes whose direct measurements provided mean values $α, β, γ, \ldots$; by $h$, a magnitude depending on them; and by $η$, the corresponding value of $h$ […]

$h = F(a, b, c, \ldots)$                         (133.1)

if $a, b, c, \ldots$ are replaced by their means $α, β, γ, \ldots$ If $h$ is a linear function (§ 74) of $a, b, c, \ldots$, the value $η$ will be the mean of $h_1 = F(a_1, b_1, c_1, \ldots)$, $h_2 = F(a_2, b_2, c_2, \ldots)$, … […] In general, however, $F$ is some non-linear function and this property will not persist.

Denote

$h - η = δ$, $a - α = δ_1$, $b - β = δ_2$, $c - γ = δ_3, \ldots$

so that $δ_1, δ_2, \ldots$ are the errors corrupting the means of the measures of magnitudes $a, b, c, \ldots$, and $δ$ is the error in the result of determining $h$.

In general, according to the degree of complication of the function *F*, that error will be connected in a more or less sophisticated way with the errors that corrupt each direct measure. However, the aim of researching the probability of the resulting error is much simplified and admits of general solution when each constituent error $\delta_1, \delta_2, \ldots$ is very small. This premise is allowed when the discussed operations are precise and merit a rigorous discussion and when the means $\alpha, \beta, \gamma, \ldots$ result from a large number of partial measures.

Imagine that we replace *h, a, b, c,* … in equation (133.1) by their exact values

$\eta + \delta, \alpha + \delta_1, \beta + \delta_2, \gamma + \delta_3, \ldots,$

and that, because of the small values of the errors $\delta_1, \delta_2, \delta_3, \ldots$ we can neglect in calculations their products and powers. The theory of functions [!] tells us that there exists a linear equation between them

$\delta = C_1\delta_1 + C_2\delta_2 + C_3\delta_3 + \ldots$ (133.2)

where $C_1, C_2, \ldots$ are constant numbers, positive or negative. The smaller is the numerical value of $C_1$, of $C_2, \ldots$ the less $\delta_1$ or $\delta_2$, or … influences the error[3] affecting the evaluation of magnitude *h*.

And so, suppose that P is the probability that the constituent errors $\delta_1, \delta_2, \ldots$ are contained within the limits $\pm l_1, \pm l_2, \ldots$ Then, there will be the same probability *P* that the resulting error $\delta$ is contained within the limits

$\pm l = \pm\sqrt{C_1^2 l_1^2 + C_2^2 l_2^2 + C_3^2 l_3^2 + \ldots}.$

If *P* = 1/2, the limits $l_1, l_2, \ldots$ will be the median values of the constituent errors, and *l*, the median value of the resulting error.

Among the different systems of the elements *a, b, c,* … which can be applied for calculating the magnitude *h*, there is one which renders the smallest possible value to the previous expression of *l*. Such systems are the most *advantageous* for determining the magnitude *h*. Thus, when measuring the height of a tower we easily discover that in the most advantageous triangle the angle adjacent to the base is closest to 45°. Indeed, the error in measuring it least influences the calculated value of the height. In much more complicated operations, notably in large geodetic networks the derivation of the most advantageous systems becomes at the same time most delicate and most important.

**134.** It frequently happens that the directly measured magnitude is not at all the same in each observation. For example, when precisely determining the height of a tower we can measure different bases for which the adjacent angles will also differ instead of measuring many times the same base (and the same angle). In trigonometric operations and physical trials the circumstances can be changed almost arbitrarily, but in such sciences as astronomy trials are replaced by observations properly understood, and meaning that the conditions are not at the observer's will. Thus, for determining the elements of a comet it is necessary to find out its astronomical places at different times and its altitudes and hour angles to be [therefore] observed under very different circumstances, possibly very unequally favourable for the exactness of the results.

When combining a large number of observations for deriving very precise results but rendering calculations practicable it is necessary that the unknowns sought were linear functions of directly measured magnitudes. That essential condition is fulfilled when approximately having in advance the values of certain elements which only need *corrections* and the values which the observations should assign to directly measured magnitudes. Suppose for example that we wish to determine quite certainly not the elements of the parabolic motion of a comet, but those of the elliptic motion of a planet already known with a good approximation. The angular magnitudes (the right ascension and declination) which fix the astronomical places of the planet at given epochs can approximately be assigned in advance. Observations provide *corrections* for each of these angular magnitudes or the differences between the calculated and observed magnitudes. […]

**135.** We will begin by treating the case of correcting one element when each observation only measures one magnitude. Denote by $a_1$ the small correction of measurement $\alpha_1$, corrupted by error $\delta_1$. Then

$$a_1 = \alpha_1 + \delta_1.$$

Let also $x$ be the correction of the element resulting from all the observations. Then

$$a_1 = C_1 x + c_1$$

where $C_1$ and $c_1$ are the known numerical coefficients. Introduce for the sake of brevity

$$\alpha_1 - c_1 = \Delta_1,$$

then

$$\delta_1 = C_1 x - \Delta_1. \qquad (135.1)$$

Neglecting the error $\delta_1$ we will have $x = \Delta_1/C_1$. When repeating the same procedure with the results of a large number of observations, we will derive the same number of different values of $x$ and can calculate their arithmetic mean as though the particular values of $x$ were provided by direct measurements. However, this manner of operating does not conform to the indication of reason. Actually, it is clear that the magnitude of the error $\delta_1$ influences the value of $x$ derived from formula (135.1) the less the larger is the value of the coefficient $C_1$. However, when simply adopting the mean of all the values of $x$ we consider all pertinent particular observations without distinguishing among them the more or the less advantageous.

For avoiding this inappropriateness the geometer Cotes had proposed a rule which astronomers have been applying owing to its simplicity[4]. It consists in adopting

$$x = \frac{\Delta_1 + \Delta_2 + \Delta_3 + \ldots}{C_1 + C_2 + C_3 + \ldots} \qquad (135.2)$$

where $\Delta_2$, $C_2$, $\Delta_3$, $C_3$, ... are similar to $\Delta_1$, $C_1$ for observations 2, 3, ... That rule is reduced to supposing that the sum of the errors disappears: $\delta_1 + \delta_2 + \delta_3 + \ldots = 0$ or to taking the mean of the values $\Delta_1/C_1$, $\Delta_2/C_2$, ... given by observations whose numbers are proportional to $C_1$, $C_2$, ... Thus, each observation influences the definitive value of the unknown $x$ the stronger the more advantageous it is by itself for deriving that value.

All this is nevertheless only a point of view based on an arbitrary hypothesis that the errors of all the observations oscillate within the same limits according to the same law of probabilities. Moreover, even under that hypothesis the described procedure is inexact. Laplace[5] had proved that we should adopt

$$x = \frac{C_1\Delta_1 + C_2\Delta_2 + C_3\Delta_3 + \ldots}{C_1^2 + C_2^2 + C_3^2 + \ldots} \qquad (135.3)$$

for narrowing as much as possible the limits $\pm l$ within which the error of the value $x$ should be contained with probability $P$. These limits are then obtained from the equation

$$t = lg\sqrt{C_1^2 + C_2^2 + C_3^2 + \ldots} \qquad (135.4)$$

where $g$ is the modulus of convergence for the common law of probabilities of the errors $\delta_1$, $\delta_2$, ...

When the correction $x$ is defined by formula (135.2) or by the Cotes rule the limits $\pm l$ are connected with probability $P$ by the equation

$$t = lg\frac{C_1 + C_2 + C_3 + \ldots}{\sqrt{m}}$$

where $m$ denotes the number of observations. And by the virtue of the principle which we repeatedly cited (§§ 73, 77) the coefficient of $l$ in that formula is always smaller than in formula (135.4). For an approximate value of the unknown modulus $g$ in that formula we can take the number

$$\gamma = \frac{\sqrt{m(C_1^2 + C_2^2 + \ldots)}}{\sqrt{2[(C_1^2 + C_2^2 + \ldots)(\Delta_1^2 + \Delta_1^2 + \ldots) - (C_1\Delta_1 + C_2\Delta_2 + \ldots)^2]}}$$

given by the observations themselves[6] and the product

$$\gamma\sqrt{C_1^2 + C_2^2 + \ldots}$$

expresses the weight of the correction $x$.

If the law of probabilities of the error and therefore the modulus of convergence vary from one observation to another, but the curve of probabilities nevertheless remains symmetric and formulas (135.3) and (135.4) will be replaced by

$$x = \frac{g_1^2 C_1 \Delta_1 + g_2^2 C_2 \Delta_2 + g_3^2 C_3 \Delta_3 + \ldots}{g_1^2 C_1^2 + g_2^2 C_2^2 + g_3^2 C_3^2 + \ldots},$$

$$t = l\sqrt{g_1^2 C_1^2 + g_2^2 C_2^2 + g_3^2 C_3^2 + \ldots}.$$

The value of $x$ provided by the formula (135.3) results from the condition that the sum of the squares of the errors $\delta_1, \delta_2, \ldots$ or the sum

$$(C_1 x - \Delta_1)^2 + (C_2 x - \Delta_2)^2 + \ldots$$

is minimal so that that formula is called the *rule of least squares of errors*[7].

**136.** The same series of observations can often be applied for determining corrections to many elements at the same time. For example, observations of the places of a planet should serve for simultaneously correcting the six elements of its elliptic motion or the masses of perturbative planets if perturbations are considered when studying that motion[9]. Equation (135.1) will be replaced by

$$\delta_1 = C_1 x + C'_1 x' + C''_1 x'' + \ldots - \Delta_1.$$

Here, $x, x', x'', \ldots$ denote the corrections sought, $C_1, C'_1, C''_1, \ldots$, numbers provided by the theory and $\Delta_1$ is a number derived from the observations themselves. The Cotes rule is not here applicable. How to combine equations of that form for obtaining the same number of resulting equations as there are unknowns $x, x', x'', \ldots$ remained absolutely indeterminate until Legendre proposed the rule of least squares for eliminating that indeterminacy. Actually, the minimal value of $\delta_1^2 + \delta_2^2 + \ldots$ can always be derived whatever is the number of the unknowns $x, x', x'', \ldots$ which enter in each function $\delta_1, \delta_2, \ldots$

**137.** We should not forget that all these rules are based on the hypothesis that we already have the values of the elements which should be corrected with a good approximation. It is then possible to obtain a very high probability $P$ that the error of the corrected value is contained within very narrow limits. This, however, does not prevent later observations from discovering that the corrected value is still very inexact if its reputedly very close value applied for calculating the correction $x$ and the probability $P$ was, on the contrary, considerably corrupted.

That circumstance occurred not so long ago when determining an important element of the Solar system[10]. By measuring the elongation of Jupiter's satellites the astronomer Pound, a contemporary of Newton, found that the mass of Jupiter was equal to $1/1067 = 0.00093721$ of the mass of the Sun. By applying a method which I sketched above to solve 126 initial conditions complied by Bouvard for the motion of Jupiter in longitude and 129 equations for the motion of Saturn, Laplace decreased that value by 0.00000294 fixing the mass of Jupiter at 0.00093427. He then found probability *P* equal to [found odds of] $10^6$:1 that the relative error of the corrected value was less than 1/100 in either direction, or that the mass of Jupiter was contained in the interval 0.00092493 − 0.00094261.

Nevertheless, it was later discovered that the perturbations in the motion of the small planets and of the Encke comet[11] occasioned by the action of Jupiter demanded an attribution of a more considerable mass to Jupiter. Finally, Airy discussed the Pounds observations anew and discovered an error so that the mass of Jupiter became 0.00095357. This value is in accord with the result of calculations of the perturbations of the small planets and is now admitted by the astronomers although it very appreciably exceeds the limits assigned by Laplace.

The defect of his conclusion could have been occasioned by wrong or incomplete terms in the expansion of the approximate formulas providing planetary perturbations or by errors in treating the initial equations or even by an error in the hypothesis according to which the law of probabilities of the errors of the applied observations remained without change and the positive and negative errors of equal magnitude were equally probable. However, this defect can also be due to the prior premise that the Pounds measures only demanded very small corrections.

**138.** This is the place to return to the remark of § 130 that it is absurd to pretend that a determination of a continual magnitude with an arbitrary precision is possible by indefinitely multiplying its observations or measurements. If, for example, an observer wished to measure a distance to within 20 decimals[12], it will be at once evident that at least the ten last digits were quite arbitrary and have no relation to the veritable numerical expression of that distance. After a large number of such expressions is calculated, it will be discovered either that each observer had irregularly and randomly chosen the last 10 digits so that the mean of the numerical values of the digits of the same order was the mean of 0, 1, 2, …, 8, 9, that is, 9/2; or, that the same causes led all the observers to prefer some digits to other and then that mean could have much differed from 9/2 but, being quite independent from the pertinent veritable digits, will not be less random.

It would be a grave mistake to infer from our example that, under conditions of actual experimentation, it is practically possible to measure a magnitude to within 10 decimals. On the contrary, we introduced an exaggerated hypothesis showing as clearly as possible the need to admit a limit to precision whatever is the number of particular measures or observations.

A singular fact! There are no continuous magnitudes yielding to empirical determination with an arbitrary approximation with the exception of constant chances of the arrival of phenomena resulting from the concurrence of constant and fortuitous causes. If, for example, in a given climate the mean chance of a male birth does not vary in time or with a sudden change in social habits, we imagine that it can be determined with arbitrary precision by an indefinite accumulation of observations. Indeed, we really understand that a perfectly exact enumeration is possible (? - O.S.) but at the same time we deny that a continuous magnitude can be thus determined by measuring it with our instruments and applying our organs of sense.

**139.** Errors in measuring, strictly speaking, and the indeterminacy necessarily affecting the *reading* of the measure should be clearly distinguished. Thus, when measuring an angle, independently from errors of sighting, centring, etc, which can vary the measure much greater than the value of the smallest division of the limb, there is the indeterminacy attached to reading each measure. It occurs since the observer either neglects[13] or arbitrarily *estimates* the fraction smaller than that value. Since this is an important point for sensibly interpreting all the results of experimental sciences and since it did not gain its merited attention, I will be excused for entering here into some minute details.

To each kind of magnitudes and manner of adopted measurement or observation there corresponds a limit after which the fractions of unity become indiscernible. Imagine such a magnitude as a straight line [as a segment] one of whose ends is fixed, and the other mobile. That line is divided into equal parts, $w$ [in length], counted from the fixed end, and so short that the observer can not discern their subdivisions. And if the point on the variable end of the line is situated more closely to the $n$-th division rather than to any other, the observer either assumes that the value of the measured magnitude is $nw$ thus neglecting the difference smaller than $w/2$ in either direction or wishes to estimate this difference, but his vision fails him and the *estimation* is purely arbitrary and fortuitous; he has no intrinsic reason to choose one value rather than another. The mean of a large number of such estimations is $\pm w/4$ for positive and negative differences.

Suppose that OA is the veritable value of the measured magnitude. The errors of measurement in its proper sense can extend to points $m$ and $n$ to the left and right from A. If there is no indeterminacy in the *reading*, the mean of a large number of measures will be equal to the distance from O to the centre of gravity of the bar $mn$ (§ 67) if $fx$ or the probability of the error $x$ of measurement expresses the density of the bar at point $x$. That mean coincides with the magnitude OA if the density varies symmetrically with respect to point A. And this is what becomes with the mean because of the indeterminacy of reading.

Suppose that marks 1, 2, 3, … separate into two equal parts the consecutive divisions of the bar from $m$ to $n$. We should imagine the mass of the bar between 1 and 2 concentrated not in its centre of gravity but in point $a$ separating the interval between these marks into two equal parts. The same consideration concerns marks 2 and 3 and the point $b$ etc. The distance of O from that system of material points represents the mean of an infinite number of observations and in general differs from OA or the measured magnitude.

What we called reading can be replaced by another method of estimation by the human organs of sense and affected by similar indeterminacies. Thus, a keen sense of hearing allows people having a good ear for music to estimate with a good approximation the tonic syllables and therefore to measure the duration of the vibratory motions producing the sound. However, this faculty has its limits so that we arrive at such small fractions of the tone that even the best exercised ear can not estimate them anymore or the estimation becomes quite arbitrary.

We have supposed that an arbitrary estimation of differences smaller than $w/2$ repeated a large number of times by the same observer should result in $\pm w/4$ for the mean of positive or negative differences. Had it been otherwise with the same observer invariably tending to *estimate* wrongly in either direction[14], the means would have been different and irregularly changing from one observer to another. It is therefore necessary to accumulate many series of measures taken by different observers or made under different systems of experimentation for providing a fixed mean, no less prone, however, to differ from the veritable value by magnitudes of the order of the indeterminacy of reading.

In general, whether the operation leading to the measurement of a continuous magnitude is delicate or crude, complicated or simple, we should recognize that independently from the causes of random or constant errors occurring because of the defects of the method, of the construction of the instruments, perturbative influence of the environment, distractions of the observer and brief or permanent disorder of his organs, − there independently exists an indeterminacy inherent in the reading or its substitute. It is therefore impossible to exceed a certain degree of precision whatever is the number of partial measures combined for determining the mean.

**140.** If desirable, the errors due to the indeterminacy of reading can be combined with those caused by the defects in the construction of the instruments and by the limits at which the hand of their constructor absolutely cease to be guided by his senses (§ 43) and only obeys blind and fortuitous causes.

If it is only necessary to establish differences between magnitudes, instruments can aid the senses of man unboundedly perfecting them. Thus, optical instruments which magnify ever more allow a classification of distances or volumes by the order of their magnitude, a feat impossible for the naked eye or less powerful instruments. Balances and thermometers of ever greater sensitivity establish inequalities of weights or temperatures imperceptible for less delicate instruments. However, when it is required to measure those differences, to compare them with some magnitude or with a unit of magnitude, these very sensitive instruments, excellent as *indicators*, do not help anymore because of the error of comparison.

Imagine indeed a thermometer capable of indicating variations of temperature to $0.01°$ whose range is very restricted because of its sensibility, for example to the interval of ca. [40, 45°]. The divisions of that thermometer are prone to irregularity and in addition, even if they are perfectly regular, before measuring a temperature it is necessary to establish the point of the scale corresponding to point zero on the standard thermometer and the ratio of the former's divisions to that of the latter. The precision of all these operations which we understand as comparison depends on the precision of the ordinary thermometer rather than on the degree of sensibility of the auxiliary thermometer.

When measuring a long distance by successively laying out a unit of length, or a large volume by scooping out its contents etc, the error inevitably affecting the comparison of the applied instrument as a unit of measure, affects all the partial measures in the same way. Compensation does not at all destroy the partial errors in the final result. Thus, when measuring a distance of about a kilometre with an error of the applied metre being 0.1 *mm*, the resulting error will be a decimetre.

In general, the measure of magnitudes with an arbitrary unit or such whose zero point should be fixed by preliminary experience is less precise than the measure of other magnitudes since the errors of the measurement itself are complicated by the errors of comparison. For this reason, other things being equal, distances are measured less precisely than angles and densities, more precisely than temperatures.

**141.** In the actual state of physical sciences, there are hardly any measures the errors of whose measurements are not complicated by errors affecting many other elements, whether the observer accepts their previously known numerical expressions or determines them himself. It is evident that the degree of precision of a final result in general depends on that of the least precise element and that calculation itself, however precise, can not render greater precision to the final numerical expression than inherent in the initial numerical data. It is therefore illusory to divide or extract a root up to the seventh decimal when having the initial data exact to within the fourth.

It should be recognized that the possible degree of precision is attained when a new series of measuring a magnitude *a* maintains intact the decimals of its previous numerical expression up to *n* inclusive but imprints quite irregular variations on the next digits. The ($n$ + 1)-st digit can be, for example, 6 in the first series of trials, 2 and 9 and … in the second, the third, … series, and their mean for a large number of series will approximately be 9/2. It is therefore illusory to preserve digits of the ($n$ + 1)-st and next orders. On the contrary, the previous digits should be thought free from all the errors occurring because of the indeterminacy of reading and similar causes. They are *certain* if only all observations are free from the causes of a kind which an exact theory renders impossible and which affects the decimals of those previous orders.

If the digits of the ($n$ + 1)-st order vary in new observations in such a manner that some of them occur oftener than the others, or if their mean appreciably differs from 9/2, they are still uncertain but should be preserved hoping that after sufficiently multiplying the measures they will be exactly determined. Advances in measurement and experimentation draw back the limits of possible precision. Digits of the ($n$ + 1)-st order impossible to be determined by previous procedures can become derivable and usefully combined with those established by new methods provided that they were not assigned purely arbitrarily.

**142.** If a magnitude is composed of two parts, one of them computed and the other one measured, the study of its degree of precision can only concern the latter part. Thus, the mean tropical year is composed of 365 and a fraction of mean solar days. There is evidently no incertitude in the whole number of days and the measuring of that fraction means evaluating those continuous magnitudes which only admit of a limited precision.

It is quite understandable that if the calculated parts are not, or not thought to be strictly equal to each other, the study of the degree of precision should concern, as usually, the entire magnitude.

**143.** Digits which we call decimals can denote natural numbers or decimal multiples of unity. Suffice it to replace the decimal point or change the unit with which a magnitude *a* is compared. If *a* is a distance, the unit can be a millimetre, a metre or a kilometre. Its choice is not absolutely arbitrary since it is absurd to evaluate in millimetres a geodetic distance which can only be determined to within a few metres or measure in metres an astronomical distance in which remains an incertitude of many kilometres. Even the radius of the Earth's orbit becomes too small as a unit if required to express in numbers the distance from the Sun to the nearest stars.

In each case the limit imposed on the degree of precision determines the limits for the choice of the unit with which the pertinent magnitude is to be compared and the number of digits with which it is to be exactly determined. That indicator of precision varies from one epoch to another and even in the same epoch from one branch of natural philosophy to another[15]. It evidently is not the same in astronomy and chemistry, in optics and acoustics or electricity. We do not go further in that delicate discussion and will only provide a few examples suitable for orienting the readers' minds.

**144.** *First example*. After collecting in Paris a large number of particular pendulum observations made in different places according to the rules of an exact theory and taking into account all their necessary corrections to render them quite comparable, we will obtain the length of a simple seconds pendulum. Here are the results. [Cournot provides a table of 10 values obtained by various authors from J. Picard to Bessel.]

The first measure by Picard in the 17$^{th}$ century barely differs from the mean of the results achieved by such able observers who had not neglected any improvements prompted by the advances in physical sciences and calculations. We conclude that it is impossible to vouch either for hundredths of a millimetre in those measurements or for obtaining them with more than 4 exact digits.

*Second example*. In his memoir on the determination of the density of the Earth, Cavendish (1798) reported 29 results in which the density of water was taken as unity. [Cournot provides a table of those results.]

The mean is 5.48, the sum of the squares of the deviations [from the mean] is 1.1967, the weight of the [final] result is 18.745. Even money can be bet on the deviation of the mean 5.48 from the absolute mean to be contained within interval ± 0.026, and the existence of a six times larger deviation is extremely improbable.

In 1837, F. Reich [1799 – 1882] in Freiberg repeated the Cavendish experiment and obtained [Cournot provided a table of means of 14 groups consisting of 2 – 6 observations.] The general mean of the 57 observations was 5.44 in good agreement with Cavendish. The difference is of the order easily explained by anomalies of chance. However, at least the third digit should be considered uncertain.

At the same time, being encouraged by the Royal Society, F. Baily made more numerous experiments varying the manner of suspending the oscillating needle and applying balls of differing materials and diameters. The summary of his 2004 experiments is shown below. [Cournot provides a table showing the materials and diameters of the balls, the number of experiments made and the appropriate means, all this for each of the three manners of suspending the balls.]

The general mean was 5.67. However, the divergence of the results for the materials of differing densities and manners of suspending the balls leaves no doubt in the existence of causes of constant errors. It is only possible to conclude from a very large number of experiments similar to those of Cavendish that the mean density of the Earth little differs from 5.5. This, however, does not yet prove that the veritable value of that mean density actually very little differs from 5.5. Indeed, the ingenious Cavendish method is connected with a cause of a constant error which affects all the similar experiments, and Maskelyne, by his very precise observations of the deviation of a plumb line made in 1774 at the foot of mount Shehallien[16] in Scotland is known to have derived 4.5 for that same density.

*Third example*. Consider finally the series of experiments made recently by J. B. A. Dumas for determining with a higher degree of precision the composition of water and for verifying the theoretical law of Dr. Prout according to whom the *chemical equivalents* (§ 136, Note 9) are exact multiples of those equivalents for hydrogen. Therefore, if the weight of the hydrogen entering into a given weight of water is represented as 1, according to Prout the weight of oxygen then combined with the hydrogen should be represented by the natural number 8. Or, if the weight of oxygen is 100, the corresponding weight of the hydrogen will be represented by 12.5.

However, a series of 19 experiments made by Dumas, after introducing all the corrections, resulted in [Cournot provided a table of the results of those experiments ranging from 12.472 to 12.562 indicating in each case the applied dehydrating acid, sulphuric or phosphoric]. The mean is 12.515, the sum of the squares of the deviations is 0.0173, the weight of the [final] result, 102.145. When substituting that weight instead of $\gamma\sqrt{m}$ in the equation

$$t = l\gamma\sqrt{m}$$

and assuming that $l = 0.015$, we get $t = 1.532$ and the corresponding value of $P$ will be 0.969. Thus, admitting that it is possible to apply here an approximate formula (see below), and supposing that the numbers provided above are free from the influence of all causes of constant error, 32 can be bet against 1 on the error of the mean 12.515 to be smaller than 0.015 in either direction. Reciprocally, admitting, according to a prior viewpoint, and for satisfying the Prout law the value 12.500, 32 can be bet against 1 that the Dumas numbers, however carefully that able chemist had carried out his research, are still affected by a constant error tending to increase his values.

In the Table above, 9 experiments out of 19 provided numbers smaller than 12.500 with their mean being 12.486, which is almost 12.480 as resulted from previous experiments by Berzelius and Dulong. Until the last experiments of Dumas, they were adopted by all those who did not share the theoretical ideas of Prout as well as by the school of British chemists. In other words, the law of probabilities of errors that apparently results from that Table, tends to group the errors in the vicinity of the extreme values which is contrary to our idea about that law in the strict sense and which does not originate from the limit imposed by the nature of things on the possible precision. We simply conclude that it is impossible to vouch for the fourth digit in that analysis or rather in the synthesis of water by the actual method of experimenting[17]. The mean of 10 experiments with sulphuric acid as the dehydrating substance was 12.520, and the mean of the 9 other experiments in which phosphoric acid was applied, was 12.511. By the formulas of § 128 we are authorized to regard the difference as purely fortuitous rather then indicating any variation in the system of constant causes in those series.

**Notes**

**1.** The ordinates *fx* of a bell-shaped curve of probabilities can be represented by a function such as $K\exp(-k^2x^2)$, cf. formula (33.2). The coefficient $K$ measures the maximal ordinate. […] If the function *fx* is really of that form, the probability $P$ corresponding to limit $l$ of the deviations can be calculated by formula (69.1) without having a large number of observations $m$. That formula becomes exact rather than approximate for any values of $m$. This case necessarily occurs when each of the particular values applied for calculating the definitive mean is by itself a mean of a large number of measures taken under similar circumstances. A. A. C.

**2.** That rule is due to Laplace, 1811, see his *Théorie* (1812/1886, pp. 345 – 348) and is included in Poisson (1837, § 108). The mean is taken as the value of $g$ rather than of $a$ as stated by Cournot. [B. B.] I am not convinced. O. S.

**3.** By the principles of the differential calculus the coefficients $C_1$, $C_2$, … are the numerical values of the differential coefficients

*df(a, b, c, …)/da, df(a, b, c, …)/db,* …

with $a = \alpha$, $b = \beta$, …

We can also determine these coefficients approximately and empirically of sorts, without knowing the rules of the differential calculus. To this end, we take $a = \alpha$, $b = \beta$, … in formula (133.1) and calculate η. Then, without changing the values of *b, c,* … we assume that $a = \alpha + \delta_1$ where $\delta_1$ is an arbitrary and very small fraction of α, for example a hundred or a thousand times smaller. We thus get $\eta_1$ very little differing from η, and $(\eta_1 - \eta)$ will be very near to $C_1\delta_1$ in equation (133.2). Therefore, the value of the quotient $(\eta_1 - \eta)/\delta_1$ will be very near to $C_1$. In the same way we successively determine $C_2$, $C_3$, … A. A. C.

**4.** Cotes published his proposal in 1722, see Gowing (1983). Without any justification he advised to regard the weighted arithmetic mean, which he compared with the centre of gravity of the system of given points, as the most *probable* estimator of the constant sought. Cournot's description follows Laplace (1814/1995, p. 121).

**5.** Laplace first provided that formula in 1811, see his *Théorie* (1812/1886, p. 325). [B. B.]

**6.** See Laplace (1812/1886, p. 327). [B. B.]

**7.** Legendre was the first to propose the rule of least squares, but only as an empirical procedure for introducing more symmetry into the calculations. Then Gauss proved that that rule ought to be considered as the most *advantageous* in virtue of the principles of the theory of chances if the law of probabilities of the errors is $K\exp(-k^2x^2)$, see Note 1.

Finally Laplace demonstrated the same for any law of probabilities under the conditions: **1.** The law of probabilities is the same for all the observations and the same for positive and negative errors. **2.** The number of observations reaches the order of magnitude allowing the application of approximate formulas. But we ought to admit that, having a very large number of observations, this rule becomes barely practicable because of long calculations. This circumstance essentially restricts the practical value of Laplace's theory. In this case the benefit of the rule depends mainly on the fact that the form of the function expressing the law of probabilities for observations of astronomical precision should little deviate from the law initially assigned by Gauss[8]. A. A. C.

**8.** Cournot did not know about Gauss's decisive contribution (1823), had not properly described Laplace's stillborn theory of errors and had not noticed that the treatment of a large number of observations is difficult in any case. Bru noted that Cournot had ignored Gauss in his other contributions and that in the first half of the 19th century Laplace's work had been considered as the last word of the pertinent sciences.

Cournot apparently thought that the number of observations was usually very large (§ 130) and that errors can be (always) judged by their relative values (§ 130). He did not explain how to compare the precision of measured distances and angles (§ 140), did not mention errors caused by external conditions of observation and his discussion of the errors of reading (§ 139) is wishful thinking. He discussed the measurement of the height of a tower but did not refer to Cotes (Gowing 1983), and he could have mentioned Bessel (1839) who had found out the most suitable points for supporting a measuring bar so that its weight will least corrupt its length. My criticism concerns geodesy, the main pertinent field both at that time and for many later decades. Also see Note 3 in Chapter 6.

**9.** Until now, the theory discussed in this chapter has been only applied in astronomy. When admitting the Prout theory (see below) in full, would not it be proper to apply it for determining the atomic weights or chemical equivalents by subjecting to simultaneous calculations the corrections, to which those weights are liable according to a large number of analyses of various substances. The effects of different causes of errors will probably be compensated since the weight of the same chemical radical is determined by analysing many compounds in which it is included in various combinations. On the other hand, certain causes of errors can act constantly in the same sense when the analysis of the same compound is repeated. A. A. C.

**10.** See Poisson (1837, § 113).

**11.** Encke calculated the path of that comet discovered in 1818 by Pons. The Airy memoir mentioned a few lines below was published in the *Mem. Roy. Astron. Soc.* for 1837. [B. B.]

**12.** Cournot possibly had in mind that in 1799, the standard metre was fixed as being 443.295936 lines. [B. B.] The line was approximately equal to 1/12 of the inch. O. S.

**13.** Bru refers to Saigey's paper of 1842 and provides information about that scholar. He many times mentions Saigey (1832) as the source of the measurements discussed by Cournot in the sequel.

**14.** Cournot could have referred to the personal equation discovered in astronomical observations by Bessel in 1823.

**15.** The research of the estimation of precision of measurement is a subject with which my friend and comrade at the *Ecole Normale* Saigey is particularly occupied. I have borrowed from him the most essential in the preceding remarks and the two first examples in § 144. A. A. C.

**16.** Bru indicates that various authors had applied different spelling of that mountain.

**17.** *The Enc. of Chemical Technology*, vol. 21, 1970, stated that that constant was 12.59. [B. B.]

## Bibliography

**Bessel F. W.** (1839), Einfluß der Schwere auf die Figur eines … Stabes. *Abhandlungen*, Bd. 3. Leipzig, 1876, pp. 275 – 282.

**Cavendish H.** (1798), To determine the density of the Earth. *Phil. Trans. Roy. Soc.*, vol. 88, pp. 469 – 526. Also in *Phil. Trans. Roy. Soc. Abridged*, vol. 18, 1809 and in author's *Scient. Papers*, vol. 2. Cambridge, 1921, pp. 249 – 286.

**Gauss C. F.** (1823, Latin), in *Theory of Combinations of Observations Least Subject to Error*. Philadelphia, 1995. Latin & English. Translated by G. W. Stewart.

**Gowing R.** (1983), *Roger Cotes − Natural Philosopher*. Cambridge.

**Laplace P. S.** (1812), *Théorie analytique des probabilités*. *Œuvr. Compl.*, t. 7. Paris, 1886.

--- (1814, French), *Philosophical Essay on Probabilities*. New York, 1995. Translated by A. I. Dale.

**Poisson S.-D.** (1837), *Recherches sur la probabilité des jugements etc*. Paris, 2003. English translation: www.sheynin.de   downloadable file 53.

**Saigey J.-F.** (1832), *Petit Physique du Globe*. Paris, 1842.

# Chapter 12. Application to Problems in Natural Philosophy[1]

**145.** For a long time the calculus of chances had only been applied to games of chance and therefore to purely speculative problems of no practical interest[2], and to facts in social statistics whose causes escape any mathematical investigation because of their complication so that about them we only have the materials of experience. Little was done to adapt that calculus to problems of natural philosophy, of, so to say, a mixed nature. It was possible to hope that in that field observational data can be confronted with relations provided by the theory. If there is a branch of natural philosophy to which that kind of research (? - O.S.) can be applied with hope for success, it is surely astronomy.

That science, distinguished among all the others because of its simplicity and grandeur of its studied phenomena, should therefore offer the most remarkable examples of rapid separation of regular causes counter to anomalies of chance. The immensity of distances separating celestial bodies and the relative smallness of their dimensions maintain the simplicity of their motion, introduce admirable geometric regularity and render those bodies more independent from each other, freer to arrange themselves according to the mathematical laws of combinations and causes influencing the initial conditions of their motion. In a word, just as the observational astronomy is a model of sciences of observation, theoretical astronomy is the model of scientific theories, and statistics of the heavenly bodies (if this association of words is possible[3]) will someday serve as a model for other statistics.

The example of a randomly thrown globe with a regular or irregular structure (§§ 71 and 81) is actually only a method of representing by a visual image a problem that should be frequently reproduced in astronomy. It consists of discovering whether a number of points was distributed over a sphere under the influence of regular or irregular causes. These points can be real as when discussing the distribution of the fixed stars of various magnitudes, of double stars[4], nebulae etc. They can also exist only geometrically when we study the intersection of the celestial sphere with a number of straight lines or radii vectors originating in a common point, the centre of the sphere, or the *poles* of a number of planes passing through that centre, i. e., the points in which the perpendiculars erected from the centres of these planes cut the spherical surface. It is evident that to each direction of the plane corresponds a particular polar point and that if the planes are uniformly distributed in all directions the polar points will be uniformly placed on the spherical surface.

The three angular magnitudes used in astronomy for determining the planetary and cometary orbits are the inclination of the orbital plane to the plane of the ecliptic, the longitudes of the ascending node and of the perihelion[5]. The first magnitude is simply the angular distance of the pole of the pertinent orbit to the pole of the ecliptic; the second magnitude differs exactly by 90° from the longitude of the orbital pole. If the *sense* of the motion of a heavenly body is not taken into account, the orbital pole can be taken as one of the two opposite points in which the perpendicular to the orbital plane cuts the sphere. In the contrary case one of those points should be fixed by a convention similar to those which geometers apply in mechanics and which provide clarity to the exposition of theorems.

Thus, we can agree to choose the pole at the north of the ecliptic if the heavenly body has *direct* motion, a motion from west to east, and at the south of the ecliptic when the body has *retrograde* motion.

**146.** The solar system, such as we know it, consists of 11 primary planets. Here are their three elements defined above. [Cournot provides a table[6] for the main planets up to and including Uranus and four minor planets.] These values which slowly change in time are shown for 1 Jan. 1801 and 1 Jan 1820 for the main and the minor planets respectively. All planets have direct motion and except Pallas, which appears to form a separate group with the three other minor planets, have very small inclinations of their orbits to the plane of the ecliptic. Theory proves that these inclinations, varying in time due to the mutual attraction of the planets, are invariably very small.

However, the theory does not tell us either why all the planets have direct motion or why the mutual inclinations of their orbital planes were initially very small. These are very remarkable features of the constitution of the solar system, and it is natural to require of the calculus of chances whether they can be attributed to fortuitous causes or causes which acted quite independently on each planet separately.

**147.** Let us first examine the circumstance of all the planets having direct motion. Denote by $p$ the chance of such a motion for each planet. The posterior probability that $p > 1/2$ is, by calculations based on the Bayes rule (§§ 88 and 92), $(2^{12} - 1)/2^{12} = 4095/4096$. In other words, 4095 can be bet against 1 on some cause favouring the appearance of direct motions more than that of retrograde motion. When only attributing to chance $p$ two possible values, 1 and 1/2, that is, when only admitting two hypotheses, that the motion will necessarily be direct, or that the directions of motion are indifferent, the relative probability of the first hypothesis will be, once more according to that rule, $2^{11}/(2^{11} + 1) = 2048/2049$.

However, to understand how fragile is the basis of these calculations suffice it to have a look at the table above. We see that for all the planets except the Earth the longitude of the ascendant node is contained between 0 and 180° so that the Bayes rule provides the probability $(2^{11} - 1)/2^{11} = 2047/2048$; 2047 can be bet against 1 on the causes being not independent since the planets favoured the concentration of their ascending nodes in that half of the ecliptic where the longitudes were less than 180°.

At the same time, that concentration is quite certainly an effect of chance. A very small displacement of the orbital planes of the order caused by planetary perturbations is sufficient for that fortuitous accumulation to disappear someday. Even an accumulation of a much more considerable number of planets than actually observed and an enormously high probability then following from the Bayes rule would have provided no meaningful consequence as is clearly shown by the construction served for representing the geometrical conditions of our problem.

Suppose indeed that after 11 random throws of a globe we determined 11 points of its contact with the floor accumulated in a certain region of the globe only a few degrees distanced from each other[7]. We will conclude with a very high probability that that accumulation was not random, but, on the contrary, followed from the structure of the sphere or manner of the throws. However, when constructing at random an arc of a great circle passing through one of those points, there will be nothing singular in that all the other 10 points are situated on the same side of that arc, nothing that can not be easily explained by a freak of chance. Indeed, for that to happen it is sufficient that, 1) the chosen point was the vertex of an angle *a* of a spherical convex polygon with the other 11 points situated within it or on its sides; 2) the arc of the great circle randomly passing through *a* did not intersect that polygon.

In our astronomical problem, the point *a* represents the pole of the ecliptic, i. e., of the orbit of the Earth; the other points are the poles of the orbits of those other planets; the arc of the great circle passes through the pole *a* in a plane perpendicular to the ecliptic and the equinoctial line. There is nothing extraordinary about the Earth being one of the planets whose orbital poles are the vertices. And it is not more singular that the equinoctial line is situated in a manner that leads the arc of the great circle to be beyond the vertex of the angle *a*. Even were there much more planets and the orbital poles and the equinoctial line had not been displaced in time, the considered fact reduced to its veritable sense by the discussion above does not at all justify an intervention of a special cause, and all the calculations of probabilities on which its existence is founded are illusory.

On the contrary, we recognize that the existence of 5 or 6 planets with their orbital poles accumulated within a very small spherical polygon would not have reasonably been attributed to chance.

Supposing that the chances of direct and retrograde motions are the same, and the chances of the longitude of the ascending node either larger or smaller than 180° are also the same, there will doubtless be the same prior probability both for having direct motion 11 times in succession and for that longitude to be smaller than 180°. However, it will not follow that these two occurred opposing events lead to the same posterior probability of their having unequal chances. At least that posterior probability only has a subjective value for someone who had to consider these two events on a par because of his complete ignorance of the other features distinguishing them. Such a probability can not lead to any objectively significant consequence.

**148.** The sum of the inclinations of the planetary orbits with respect to the ecliptic or the sum of the distances of the pole of the ecliptic to the poles of the other orbits is equal to 82°14′19″.1. If all the values of the polar distances contained between 0 and 180° are equally probable[8], the probability that that sum is less than 90° will be (Note 5 in Chapter 6) an extraordinarily small fraction $1/2^{10}10! = 1/3,715,891,200$. However, if the purely fortuitous causes determine the directions of the orbital planes, the probability that the value of each polar distance will increase proportionally to the sine of that value and the values near to 0 or 180° will be much less probable than those near 90°. Therefore, the probability that the sum is less that 90° or the mean less than 9° will be considerably lower than the preceding fraction. Its exact determination would have been tiresome and useless.

Calculation according to our inexact hypothesis, much less favourable for an equal probability of all the values, will suffice for proving that the accumulation of the orbital poles around the poles of the ecliptic can not be regarded as random. Instead of relating the orbital planes of other planets to the plane of the ecliptic, it is natural to relate the latter, just as the other planetary orbits, to the plane of the solar equator. Then, after considering that all the motions of the rotation of the Sun and the planets and of the translations and rotations of the satellites are direct[9], and in most cases occur in planes little inclined to the solar equator, we will not doubt that an initial cause tended to bring together the orbital planes and the plane of the solar equator and to imprint on all those bodies translations and rotations directed in the same sense as the rotation of the great mass dominating the system.

**149.** It is natural to inquire whether the same or a similar cause acted on the comets[10], or, on the contrary, whether these bodies differ from the planets in everything. Their volumes are enormous and masses imperceptible, they move along hyperbolas or excessively eccentric ellipses, some in the direct, and some in the retrograde sense. In addition, their orbital planes are indifferently directed to all regions of the space. When calculating the mean of the inclinations of the cometary orbits to the ecliptic, we find that its value is[11] ca. 50° thus exceeding 45°. And Laplace, owing to a very singular oversight, concluded that the comets, *far from following the tendency of the bodies of the planetary system to move in planes little inclined to the ecliptic, seem to have a contrary tendency*[12].

For purely empirically becoming convinced in the defect of that conclusion, suffice it to calculate the inclinations of the cometary orbits with respect to one or two planes perpendicular to the ecliptic. It will occur that the [their] means exceed 60° and it becomes sensible that the orbital planes are not uniformly distributed in all regions of the space and that they seem, on the contrary, to indicate a marked tendency to approach the plane of the ecliptic. This becomes theoretically evident also by considering the orbital poles. Had the orbital planes been indifferently directed to all the regions of the space, the poles would have been uniformly distributed over the regions of the celestial sphere and then the mean distance of these poles to the pole of the ecliptic or the mean inclination of the planes to the ecliptic would converge to 57°17′44″.8 (§ 71), i. e., to a value much exceeding 50°.

**150.** For establishing the laws to which the distribution of the cometary orbits in space can be subjected and for applying the calculus of chances for determining the particulars presented by that distribution, the system of coordinates used by astronomers (§ 145) is of little use. On the one hand, we have remarked (§ 71) that the mean of a number of longitudes converges less rapidly to a fixed number and again less rapidly than the mean of a number of latitudes or of the corresponding polar distances gets rid of the anomalies of chance. On the other hand, we know that the best method for discovering the laws governing the combinations consists in combining the symmetrical elements whereas the system of longitudes and latitudes or longitudes and polar distances do not satisfy that condition of symmetry. Therefore, we[13] have applied the distances of the orbital poles from three points symmetrically situated on the heliocentric celestial sphere: the north pole of the ecliptic, the vernal equinox or the first point of *Aries*, and the summer solstice or the first point of *Cancer*. These distances, which we denote respectively by $\theta$, $\theta'$ and $\theta''$, measure the angles that a straight line perpendicular to the orbital plane makes up with these mutually perpendicular straight lines drawn from the centre of the sphere to the three indicated fixed points. They also measure the inclinations of the orbital plane to the ecliptic and to the two other mutually perpendicular planes, both of them perpendicular to the ecliptic.

Since the sense of the cometary motion is not taken into account, it is possible to choose as the orbital pole any of the two points in which the straight line perpendicular to the orbital plane cuts the sphere, and assume that $\theta$, $\theta'$ and $\theta''$ take values from 0 to 90° without attaching signs to them. In the contrary case (§ 145) the pole should be attributed either to the north, or to the south hemisphere depending on whether the comet has direct or retrograde motion. However, we will still assume that $\theta$, $\theta'$ and $\theta''$ take values from 0 to 90°, but distinguish them by their signs and regard as positive the angles between the straight lines drawn from the centre of the sphere to the orbital poles and the radii vectors of the north pole of the ecliptic and the first point of *Aries* and *Cancer*, and as negative the angles between those straight lines and the extensions of the same radii vectors to the opposite regions of the celestial sphere. Finally, we denote by $t$, $t'$, and $t''$ the analogues of $\theta$, $\theta'$ and $\theta''$ for the perihelia instead of the orbital poles and apply to them the same system of notation[14].

We base our calculations on the catalogue of comets published in 1823 by Olbers[15] in the *Astronomische Abhandlungen* with a supplement extending to 1825 and followed up by the catalogue of Santini (1830) extended to comet No. 137. Two more comets observed in 1830 and 1832 were added, so now 139 orbits are known. On the other hand, we thought it necessary to exclude as very uncertain the Chinese, Arabic and European observations made before the 16th century and only take into account 125 orbits. The dates of their first observed appearance together with the elements attributed by calculation to that date were only preserved in the catalogue for comets with a constant period of returning. These are the comets of Halley of 1607, of Encke of 1786 and Biela of 1772. This manner of applying the elements of the periodic comets seems to be most free from arbitrariness.

**151.** Suppose that by means of the catalogue mentioned above we calculated and arranged in a table in the chronological order of the appearance of the comets the pertinent angles $\theta$, $\theta'$ and $\theta''$ and $t$, $t'$, and $t''$ and that we calculated the [their] means for 10, 20, … of the first comets until exhausting all 125 without taking into account the signs of each angular value. We will have the following table. [Cournot provided a table of those 6 magnitudes for 10, 20, …, 110, 120, 125 orbits.]

When inspecting this table we are first of all surprised at the small intervals within which oscillate those consecutive means after the 30 first ones. From the beginning of the 18th century when there were only 30 comets with somewhat known elements, those values have remained almost the same up to our epoch. According to this table, it is impossible to deny the existence of constant causes, real or apparently real, having to do with the conditions of that phenomenon or of observations, which maintain the mean of θ smaller than the means of θ′ and θ″, and those of t″ smaller than the means of t and t′.

If no constant cause would have opposed the uniform distribution of the orbital planes and of their major axes, the absolute mean (? - O.S.) will be 57°14′44″.8 and the modulus (? - O.S.) becomes 2.9518 if a quarter of a circumference is chosen as unity (§ 71). We will then have in advance the probabilities 0.99991, 0.938, 0.982, 0.939, 0.892 and 0.9986 that the means of the 125 trials of θ, θ′ and θ″ and t, t′, and t″ will not randomly deviate from the absolute means in either direction as much as they did in the last line of the Table. Therefore, and as always bearing in mind the explications repeatedly provided here of the essence of posterior probabilities, we can bet 9,999 against 1 on the probability that until today, judging by the mean value of angle θ, the constant causes acting on separate observations tended to bring together the plane of the ecliptic and the planes of the observed cometary orbits. Just the same, when taking into account the mean value of angles t″, we can bet 715 against 1 on the probability that constant causes tended to bring together the major axes and the solsticial line.

**152.** At present, it should be remarked that the laws of probabilities of the angles θ, θ′, … are not independent from each other. Indeed, on the one hand the sum of the numerical values of the angles θ, θ′ and θ″ ought to remain within the limits 180° and 164°13′ (thrice the angle whose tangent is $\sqrt{2}$) and it is the same for the sum of the angles t, t′, and t″. On the other hand, the system of the angles θ, θ′ and θ″ acts on the system of angles t, t′, and t″ and vice versa. Denote for the sake of brevity by Θ the point of the celestial sphere corresponding to the system of the angles θ and by T the point corresponding to the system of the angles t. If point Θ is given, the point T should be situated on the great circle of the sphere with pole Θ; and when T is given, the orbital plane can only rotate about the diameter passing through T so that the point Θ should also be situated on the great circle of the sphere with point T as its pole.

Deliberating about the essence of this mutual dependence, we are bound to think that it can be most suitably revealed[16] by decomposing the series of the values of each angle into two partial series formed respectively by the values exceeding and smaller than 60°. According to the hypothesis of a uniform distribution these two series should consist of the same number of angles. We see however, on the contrary, that after that decomposition the results can be expressed by the following notation convenient by its brevity.

θ, 48:77; θ′, 65:60; θ″, 69:56 and *t*, 77:48; *t*′, 66:59; *t*″, 44:81

It means that, for example, the series of angles θ contains 48 exceeding 60° and 77 smaller than that value. Now we clearly see the influence that the manner of distribution of the values of θ exerts on the distribution of the values of *t* and we find a similar dependence in comparing θ′ with *t*′, and θ″ with *t*″. The drawing together of the numbers in the last case is less surprising and we will soon reveal the cause of this particular circumstance[17].

**153.** After establishing the existence of constant causes which influence the chronological series of observations of the cometary orbits we should wish to go further and study by certain decompositions of the series the nature of the influencing causes and whether their influence is the same for various angular regions and different senses of motion.

To this end we can at first decompose the series in two others. One of them is formed by all comets whose passage of their perihelion was observed in the half-yearly winter period from 22 September to 22 March, and the other, from 22 March to 22 September. That first separation is only relative to the situation of the observer, and should be regarded as accomplished by chance and therefore having no influence on the laws of probabilities of the elements, at least if these laws are not themselves subordinated to causes only depending on the circumstances of observation.

**154.** Out of 125 comets of our list, 71 belong to the winter series and 54, to the summer series. That inequality could have been expected. Thus, Arago remarked that during the summer months *the long day properly understood and the twilight can not fail to conceal from us a certain number of those acts*[18]. Suppose that the existence of this cause had not been manifested in advance; then, after deriving the probability of the ratio 71/54 by the ordinary formulas, we will obtain 0.924, a number of an order which is not usually thought as decisive in natural philosophy.

However, that probability can be much heightened (§ 116) when studying the ratio of the numbers in both series of 10, 20, … orbits taken in chronological order of appearance and beginning from the first 10 ones. This is indicated in the following table. [Cournot provides a table showing both the total number of cometary orbits and their number observed in winter and summer: 40 (24 and 16), 50, … 110, 120 and 125 (71 and 54).]

**155.** When taking into account the signs of the angles (§ 150), the distribution of the orbits will be [Cournot provides a table showing the distribution of the 125 orbits among positive and negative values of the set of the six angles θ, θ′, θ″ and *t*, *t*′, *t*″ separately for winter and summer.]

The most apparent result is that in winter the number of positive values of *t*″ exceeds the number of its negative values whereas the opposite takes place in summer. This phenomenon is certainly due to the optical circumstances of observation. Actually, in winter the radius vector passing from the Sun to the Earth forms an acute angle with the line, again issuing from the Sun, to the first point of *Cancer*; in summer, however, that angle becomes obtuse. And when considering a plane perpendicular to the ecliptic and the solstician line and passing through the centre of the Sun, it becomes clear that a terrestrial observer has more chances to discover comets at their perihelia situated on the same side of that plane as the Earth than those being on the other side of that plane, in the hemisphere opposed to the heliocentric celestial sphere[19].

**156.** Bearing in mind these preliminary remarks, after constructing a table of series for each half-yearly period similar to those of § 151 for the general series, we obtain [Cournot provides a table showing the values of θ, θ′, θ″ and $t$, $t'$, $t''$ for 10, 20, …, 50, 60, 71 orbits and the winter period, and for 10, 20, …, 40, 50, 54 orbits and the summer period.]

The same remarks concerning the small extent of the oscillations of the means and of their rapid convergence as in § 151 are possible. Angles θ and θ″ are apparently the only ones for which we can conclude from the final means of the Table that their chances in those periods are unequal. According to the usual formulas whose sense we have explained, the probabilities of these inequalities are 0.927 for θ and 0.946 for θ″.

A very simple consideration suggests the idea that there exists an optical cause which should modify the means of the two series of angles θ. Actually, the effect of the solar light which conceals from the European observer more comets in summer than in winter is caused by the plane itself of the ecliptic. Consequently, mostly affected are the comets less deviating from the ecliptic or moving along orbits little inclined to that plane. The natural means for checking that conclusion consists in comparing the number of the orbits in both series whose inclinations are contained within certain limits. In addition, important consequences that can be connected with this result compel us to consider it doubtless when comparing the orbits, beginning after the first 30 of them, as indicated in the following table. [Cournot provides a table showing the number (30, 40, …, 110, 120, 125) of orbits having inclinations 0 – 40, 40 – 60, and 60 – 90°, separately for winter and summer.]

The table makes it evident that the difference in the action of the considered optical cause in winter and summer almost exclusively tells on the comets whose orbits are inclined from 0 to 40° to the ecliptic and ceases to be noticed when inclinations exceed 60°. We should conclude that, without the influence of the solar light which acts in the same sense both in winter and summer, although less intensively in winter, the accumulation of the cometary orbits in the zodiacal regions would have been much more appreciable for observing. Actually, it seems impossible to assign mean inclinations obtained after completely eliminating that optical influence. We can not even assert that there is no other optical cause producing the accumulation in the zodiacal regions. However, at least the opinion that that accumulation is real and due to some cosmological cause becomes very highly probable.

**157.** The optical influence is principally manifested by the inequality of its action in winter and summer and modifies the distribution of the angles θ apparently without appreciably acting on the angles *t*. It should therefore corrupt the relations which would have been naturally established between the two laws of distribution, and these relations are more exactly seen in the winter rather than in the summer or in the total series. Actually, when continuing to apply the notation of § 152, we find for the winter series this remarkable result:

θ, 24:47; θ′, 36:35; θ″, 46:25 and *t*, 46:25; *t*′, 36:35, *t*″, 27:44.

It can be shown in quite a symmetrical form:

θ, *m:n*; θ′, *p:p*; θ″, *n:m* and *t*, *n:m*; *t*′, *p:p*, *t*″, *m:n*.      (157.1)

Here, $p = (m + n)/2$. However, the surprise possibly caused by the extreme simplicity of this statistical law derived from such a small number of elements, will be still strengthened when decomposing the winter series in two others depending on whether the angles θ are positive or negative. When compiling as many such parts as there are angular elements and when finally compiling a last partition depending on perihelion distances being smaller or larger than 3/4 of the semimajor axis of the terrestrial orbit, we will obtain the following result. [Cournot provides a table showing the number of orbits for each of the six angular magnitudes taken separately with either sign and broken down into those exceeding 60° and smaller than that.]

All these results so well agree with formula (157.1), the deviations are so small in spite of the small number of considered orbits and of the combinations taking place when having *seven* different partitions of the same series, that it is very difficult to attribute that coincidence to chance[20].

**158.** The particular position of a European observer provides occasion for studying the possible differences between two series consisting of comets with a north (+ *t*) and south (− *t*) perihelia. Indeed, in Europe, the retained elevation of the north pole of the ecliptic above the horizon can not fail to conceal from us a certain number of comets which, while near their perihelia, do not leave the southernmost regions of the heliocentric sphere.

These considerations perfectly agree with the experience as we can show by a table classifying the angles *t* in both series similar to what was done in § 156. For the sake of brevity we restrict our description by providing the final result. The first number refers to the series of comets with a north perihelion and the second, to the series of those with a south perihelion.

For *t* = 0 − 40, 40 − 60, 60 − 90 and 0 − 90° the respective series are 21:6, 10:11, 37: 40 and 68:57.

It is easily seen that the comets whose angle *t* is negative (? - O.S.) and less than 40°, or, in other words, whose perihelia have a south heliocentric latitude larger than 50°, can only be seen in Europe under very rare circumstances. Only one such comet among the 60 (most of them barely, or not at all visible by the naked eye) was observed after 1780.

That circumstance could not have failed to influence the inclination to the ecliptic because of the relation between the angles θ and *t*. If we compare the means of the elements for north and south perihelia both in the total, and in each half-yearly series we will obtain the following results. [Cournot provides a table showing for each of those three series the mean values of each of the six angular magnitudes, separately for the two different kinds of perihelia.] The differences are of the same sense for series of the half-yearly periods except for the angle *t* whose deviations can very probably be attributed to anomalous causes[21].

**159.** And so, the two optical influences which are evidently connected with the local situation of the European observer, act in contrary senses. The result of our discussions is that the influence tending to diminish the mean inclination is insufficient for explaining, at least nowadays, the difference observed between that mean and the absolute mean under the hypothesis of uniform distribution. However, the main problem remains unsolved and we are still justified in requiring whether the purely optical causes connected with the position of the terrestrial orbit in the celestial space do not occasion the observed difference. For Lambert, that point did not seem probable but he had not provided any arguments supporting his opinion which is however strongly corroborated by the remark in § 156.

On the one hand, it seems difficult to take into account, precisely and in advance, the chances of visibility; on the other hand, we should recognize that the number of observations (sufficient for establishing the apparent laws of distribution in the general series) is not enough for providing a definitive solution of this interesting problem[22].

**160.** It has been remarked long ago that the numbers of comets with direct (+ θ) and retrograde (− θ) motions are appreciably equal whenever the series was terminated. Out of the 125 comets we have 65 and 60 respectively. There is therefore room for believing that the chances of visibility are not less for the former than for the latter. On the other hand, we have 69 comets whose direction of motion renders the angle θ′ positive (§ 150) and 56 causing that angle to be negative. We are therefore authorized to suppose that the chances of visibility are greater for the comets of the series (+ θ′) than for those of (− θ′) or that at least these chances are the same.

However, suppose that the smallness of the mean inclination occurs because a certain number of comets among those whose orbits are more inclined to the ecliptic escape the conditions of visibility, then the mean θ will not be less for the comets with direct, than for those with retrograde motion, and it will be larger for the series (+θ′) than for those of (− θ′). But this is exactly contrary to the observed: the means of θ for series of (+ θ) and (+ θ′) are noticeably smaller than for the series (− θ) and (− θ′), and their deviations are sufficiently large and sufficiently stable for rendering probable the inequalities between the laws of distribution depending on whether we consider one of the partial series or a [partial] series with a contrary sign as can be judged by the following table.

[Cournot provides a table showing the number of orbits 1) 30, 40, …, 60, 65 and 2) 30, 40, 50, 56, 60, 69 and the respective mean values of θ for series (+ θ) and (− θ) and (+ θ′) and (− θ′).]

Suppose that someday that result becomes doubtless[23]. It will then be natural to conclude that the inequality of the distributions between the series does not occur because of an optical influence so that real or cosmological causes modified the law of distribution depending on the sense of motion, and that that law for the general series is definitively affected by the influence of real causes to which the accumulation of the orbits in the zodiacal regions should be indeed attributed[24].

**161.** If the general series is decomposed in two others, one of them consisting of comets with perihelia distances smaller than 0.75 (with the semimajor axis of the terrestrial orbit chosen as unity), and the other, with larger perihelia distances, these partial series will include an almost equal number of comets, 65 and 60. However, for recognizing that this result is only temporary, suffice it to remark that for the first 60 appearances until 1772 that ratio was 40:20 whereas for the 65 later appearances it was 25:40. A more thorough study of the sky by more powerful instruments evidently discovered a much larger number of comets among those which previously had escaped the conditions of visibility because of their large perihelion distances. When taking the means of the two partial series thus formed we obtain the following result.

[Cournot provides a table showing the six angular magnitudes for both series.] We should conclude, at least for the time being, and contrary to Lambert's opinion, that the laws of distribution apparently do not essentially change with the perihelion distances[25].

**162.** We have dealt in detail with the problems of the distribution of cometary orbits because it seems to have provided a remarkable type of analysis and a discussion of a statistical fact, and because it clearly shows (contrary to the received prejudice) that in this type of discussion it is not always necessary to have a very large number of observations. Large numbers are really needed for the stability of the means[26] and only when the ties of solidarity between individual observations irregularly change from time to time and place to place.

We will only add a few words about the application of the calculus of chances to phenomena produced by various physical agents on the surface of our globe and on the fluids covering it. Take for example the variations of the atmospheric pressure indicated by the barometer[27]. Many accidental and irregular causes incessantly vary, especially in our climate, the barometric height, so that, when being content with measuring this height at certain hours of the day for a small number of weeks or even months, we can only arrive at discordant results and be absolutely unable to single out the laws of *diurnal* barometric variations, or the small but constant influence of the hour of the day on the barometric height.

Suppose, on the contrary, that we collect a large number of observations and take the means of the heights observed at two different hours of the day, for example, at 9 o'clock in the morning and at 3 o'clock in the afternoon. Suppose also that the days of afternoon observations are chosen fortuitously and independently from the days of the morning observations. The difference of the observations will indicate the existence of constant causes by whose virtue the atmospheric pressure depends on the hour of the day. The anomalies due to the intervention of causes more energetic but irregular and independent from the hour of observation will disappear. Depending on the magnitude of the difference and the number of values taken to form each mean, calculations will indicate the probability with which we can decide about the existence of constant causes to which that difference was due.

However, when proceeding in that manner, we often find it necessary to have a very large number of observations for arriving at appreciably fixed means and to be able by formulas of the theory of chances to become assured with sufficient security of the existence of the sought constant causes.

Now suppose that instead of registering the absolute values of the barometric heights as described above we note the excess of the morning heights for each day. The evening height is not in general independent from the height in the morning since the action of anomalous causes which heighten or lower the morning barometric column usually does not cease at 3 o'clock in the afternoon. It thus becomes evident that the difference of the two consecutive heights is largely free from the influence of causes irregularly affecting each absolute height. Consequently, the mean difference between consecutive heights more rapidly and more surely indicates the laws of diurnal variation or the influence of the hour on the barometric height.

Such procedures have actually established the laws of the diurnal variations of the atmospheric pressure in different seasons and climates. And, had it been required, we would have investigated in the same way the diurnal variations of the magnetic declination or of any other similar phenomenon. The same principle which we apply to consecutive observations at a given place is usable for treating simultaneous observations made in different places sufficiently near to each other for being influenced by the same perturbative causes. Thus, thermal inequalities in two neighbouring places are more rapidly and more surely established by taking the means of the differences of temperatures observed at the same time in both those places rather than when determining them by the means of independent observations of the temperatures[28].

**163.** The theory of the tides, so important in itself and because of its connection with great astronomical phenomena, offers an opportunity of one more excellent application of isolating the effects of regular causes from those due to perturbative forces[29]. When a sequence of observed heights of the sea is formed in a port, at first the results are usually very irregular and discordant because of the action of winds, tempests, currents and many accidental circumstances. They corrupt the regular action of the Moon and the Sun which tend periodically to rise and abate the waters of the Ocean.

However, when taking the excess of the height of the high tide over that at the preceding and following low tide, and comparing its mean values during syzigies and quadratures near the equinoxes and then near the solstices and, finally, comparing the series of observations made when the Sun or the Moon are near their apogees and perigees, we will reveal the proper influence of each of those attracting celestial bodies depending on their distance and remoteness from the equatorial plane.

We thus obtain, as Laplace[30] had shown, results sufficiently precise for determining the value of the mass of the Moon agreeing with what was established by other purely astronomical phenomena. Besides that, it was attempted to find out whether barometric observations also indicate the existence of periodic movements of the atmosphere caused by the lunar attractive force. Laplace compared the diurnal barometric variations observed in Paris for the days of syzigies and quadratures. He studied a series of 298 days of either and his calculations led to the *lunar atmospheric flux* equal to 0.0176 *mm*[31]. This magnitude is too small for definitively concluding that the lunar attractive force appreciably influences the atmospheric pressure, at least in our climate.

On the other hand, by issuing from a long sequence of noon observations at Viviers, Flaugergues concluded that the variation of the diurnal barometric heights amounted to 1.48 *mm* and seemed to be very regularly connected with the lunar phases and therefore, because of the chosen hour of observation, with the distances of the Moon from the meridian, see the following table[32]. [Cournot provides a table showing the number of observations (246 – 248) and the respective barometric heights during various lunar phases.] By interpolation we find that the minimal value of the barometric heights at Viviers (754.78 *mm*) occurred at $9^h 18^{min}$ before the Moon passes the meridian, and the maximal value (756.26 *mm*) at $6^h 12^{min}$ after that passage.

## Notes

**1.** Except for some details, this chapter repeats Cournot (1834). [B. B.]

**2.** Games of chance had been in the social order of the day. Leibniz (1704/1996, p. 506) advocated the creation of *a new type of logic* and therefore recommended to study all kinds of games. In 1713 Nikolas Bernoulli invented the celebrated Petersburg game and thus started a extremely useful discussion and De Moivre (1718) described many games which compelled him to develop the theory of recurrent sequences and to introduce the definitions of fundamental notions (probability, independence). All this profited probability and mathematics in general.

**3.** Heavenly bodies are too different and their statistics was never developed, but Cournot failed to mention William Herschel's statistical investigation of the starry heaven as well as the study of the proper motions of stars which began in 1837 (F. Argelander).

**4.** Michell (1767) was the first to inquire whether the existence of double stars can be explained by a random distribution of stars over the celestial sphere. Many later authors discussed this problem (Sheynin 1984a, § 5).

**5.** The longitude of the ascending node is the angle between two straight lines passing from the centre of the Sun to the vernal equinox and to that node. Recall that the astronomers understand the longitude of the perihelion as the longitude of the ascending node augmented by the angular distance from the perihelion to the node; that distance is measured in the orbital plane. A. A. C.

**6.** Cournot borrowed that Table from John Herschel (1834). [B. B.]

**7.** D'Alembert (1768) did not agree; his main criticism was however directed against Daniel Bernoulli (1734). [B. B.]

**8.** This is the hypothesis of Daniel Bernoulli (1734) maintained by Laplace and contested by Cournot (1834, p. 505) and here in § 149. [B. B.] Poisson (1837, § 110) also rejected it. O. S.

**9.** Several satellites are now known to have a retrograde motion (Blazko 1947, p. 343).

**10.** This is what Daniel Bernoulli did as well as Lambert and Laplace. [B. B.]

**11.** Previously Cournot (1834) wrote 45° and Poisson (1837, § 111) provided 48°55′ for 138 comets. [B. B.]

**12.** See Laplace [1812/1886, p. 263]. Only respect for the memory of Laplace prevented Poisson to recognize the palpable error of that great geometer. Contrary to what Poisson (1837, § 111) apparently had in mind, the problem was not to find whether all the inclinations to the ecliptic were equally probable. No, as Laplace expressly indicated, it was to determine whether the orbital planes were indifferently directed to all regions of the space, or, on the contrary, have the tendency to approach the plane of the ecliptic. A. A. C.

Cournot (1834, p. 505) noted that Laplace made that mistake in his *Mécanique Céleste* and that Poisson (1837, § 111) criticized Cournot stating that he confused *the assumption that all the points of the celestial sphere can with the same probability be the poles of the cometary orbits with the hypothesis of equal probability of all possible inclinations of the comets*. Cournot's answer is not convincing. On the inclinations of the comets see for example Richter (1963, pp. 12 – 16). [B. B.]

**13.** See my *Additions* to Herschel (1834). A. A. C.

**14.** Denote the longitude of the ascending node by $\lambda$; by $l$, the longitude of the perihelion in the sense indicated in Note 5; by $w$, the projection of the angle $(l - \lambda)$ on the ecliptic. Then $(\lambda + w)$ will be the angle between two circles of latitude, one of them passing through the perihelion of the comet, and the other, through the vernal equinox. We will then have

$\cos\theta' = \sin\theta \sin\lambda, \quad \cos\theta'' = -\sin\theta \cos\lambda,$
$\tan w = \cos\theta \tan(l - \lambda), \quad \cos t = \sin\theta \sin(l - \lambda),$
$\cos t' = \sin t \cos(\lambda + w), \quad \cos t'' = \sin t \sin(\lambda + w).$ A. A. C.

**15.** Olbers regularly published his pertinent calculations but did not compile any catalogue. [B. B.]

**16.** This is a rudimentary test. Galton, in 1877, provided the first contingency table. [B. B.]

**17.** Cournot returned to that circumstance in § 157. [B. B.]

**18.** See *Annuaire* for 1832, 2nd edition, p. 357. Lambert (p. 209 of the French edition of 1801 of his *Lettres*) annotated by d'Utenhoven made the same remark even earlier. A. A. C. I have only found an *Annuaire historique ou histoire politique et littéraire* but this reference is only a conjecture. O. S.

**19.** Cournot (1834, p. 520) also stated that a *confrontation* of statistical data with the appropriate prior theory could be applied for appraising future statistical results. [B. B.]

**20.** Cournot (1834, pp. 524 – 525) also stated that the existence of such a simple law could be partly due to the relations similar to the discussed not necessarily being continuous. [B. B.]

**21.** Cournot (1834, pp. 526 – 527) also stated that with probability 0.948, as followed from the total series and augmented by both partial series providing results of the same sense, the law of distribution of the angles $\Theta$ is not the same for comets with north and south perihelia. [B. B.]

**22.** Cournot (1834, p. 527) also stated that in a hundred years the number of the observed comets will double and his results will be seen as a good example of combining the calculus of chances with statistical analysis. In 1843, he apparently lost his confidence and optimism. [B. B.]

**23.** It is thought that a sample of a few hundred comets out of many millions is too small and too heterogeneous. However, almost all comets with a short period and small inclination have a direct motion (Richter 1963, p. 16). [B. B.]

**24.** Cournot (1834, p. 529) remarked that it was time to find out the nature of those causes. [B. B.]

**25.** Cf. Lambert's *Lettres* [1761?], pp. 222 – 223 and 67 – 69. [B. B.]

**26.** Cournot (1834, p. 503) remarked that for obtaining stable means much more trials were needed if made under differing circumstances by different people and globes which is the reason why social statistics should be based on such a large number of observations. [B. B.] Cournot had unjustifiably restricted the need to have a large number of observations. O. S.

**27.** See Laplace (1812/1886, pp. 355 – 358, 1814/1995, pp. 54 – 55) and his *Oeuvr. Compl.*, t. 13, 1904, pp. 342 – 358, Stigler (1975) and the article *Laplace* in the *Dict. Scient. Biogr*. [B. B.]

**28.** Lamont (ca. 1839, p. 263) stated without proof that a year of simultaneous observations at different places was tantamount to 30 years of ordinary observations. Then, after 30 years of experience, he (1867, p. 245) declared that that difference method will enable meteorology to become a mathematical discipline. Also see Sheynin (1984b, pp. 71 – 72).

Chetverikov, the translator of Cournot into Russian, referred to many later authors beginning with Cave-Browne-Cave (1904) and added that the mentioned new practice became the embryo of the *variate difference method*.

**29.** This was the main goal of De Moivre (1718) as he stated in a *Dedication* of that book to Newton. That *Dedication* was reprinted in 1756 (p. 329).

**30.** See *Méc. Cél*, livre 4, and *Oeuvr. Compl.*, t. 12, 1898, p. 480. [B. B.]

**31.** Cournot provided a reference, but see Note 32.

**32.** *Bibliothèque universelle*, Dec. 1827, p. 264, and Apr. 1829, p. 265. A. A. C. It was the astronomer H. Flaugergues (1755 – 1835) who observed in Viviers, Ardèche.

## Bibliography

**Bernoulli Daniel** (1734), Recherches physiques et astronomiques etc. *Werke*, Bd. 3. Basel, 1987, pp. 303 – 326.

**Blazko S. N.** (1947), *Kurs Obshchei Astronomii* (Course in General Astronomy). Moscow – Leningrad.

**Cave-Browne-Cave F. E.** (1904), On the influence of the time factor on the correlation between the barometric heights of stations more than 100 miles apart. *Proc. Roy. Soc.*, vol. A74, pp. 403 – 413.

**Cournot A. A.** (1834), Commentary to Herschel (1834).

**D'Alembert J. le Rond** (1768), Doutes et questions sur le calcul des probabilités. In author's *Mélanges de litterature, d'histoire et de philosophie*, t. 5, pp. 239 – 264. Amsterdam.

**De Moivre A.** (1718), *Doctrine of Chances*. London, 1756, 3rd edition, reprinted: New York, 1967

**Herschel J.** (1834), *Traité d'astronomie*. Translated by Cournot from *Outlines of Astronomy*. Second edition, London, 1849. Cournot also added a commentary (*Addition*).

**Lambert J. H.** (1761 German), *Cosmological Letters*. Edinburgh – New York, 1976.

**Lamont J.** (ca. 1839), Nachricht über die meteorologische Bestimmung des Königreiches Bayern. In author's *Jahrb. Kgl. Sternwarte bei München* für 1839, pp. 256 – 264 and 247 – 249 (paging wrong).

--- (1867), Über die Bedeutung arithmetische Mittelwerthe in der Meteorologie. *Z. Öster. Ges. Met.*, Bd. 2, No. 11, pp. 241 – 247.

**Laplace P. S.** (1812), *Théorie analytique des probabilités*. *Œuvr. Compl.*, t. 7. Paris, 1886.

--- (1814 French), *Philosophical Essay on Probabilities*. New York, 1995. Translated by A. I. Dale.

**Leibniz G. W.** (1704), *Neue Abhandlungen über menschlichen Verstand*. Hamburg, 1996.

**Michell J.** (1767), An inquiry into the probable parallax and magnitude of the fixed stars. *Phil. Trans. Roy. Soc. Abridged*, vol. 12, 1809, pp. 423 – 438.

**Poisson S.-D.** (1837), *Recherches sur la probabilité des jugements* etc. Paris, 2003. English translation: www.sheynin.de downloadable file 53.

**Ramond de Selas L. F.** (1811), *Mémoires sur la formule barométrique etc*. Clermont-Ferrand.

**Richter N. B.** (1963), *The Nature of Comets*.

**Santini G.** (1830), *Elementi di astronomia*. Padova.

**Sheynin O.** (1984a), On the history of the statistical method in astronomy. *Arch. Hist. Ex. Sci.*, vol. 29, pp. 151 – 199.

--- (1984b), On the history of the statistical method in meteorology. Ibidem, vol. 31, pp. 53 – 95.

**Stigler S. M.** (1975), Napoleonic statistics: the work of Laplace. *Biometrika*, vol. 62, pp. 503 – 517.

# Chapter 13. Application to Problems about the Elements of Population and the Duration of Life

**164.** An analysis of the work of statisticians on the elements of population and on mortality of humankind can be the sole aim of a considerable treatise. However, the plan of our book mostly envisages formulation of principles and does not allow us to enter into details. We restrict our deliberations by indicating a summary of some results since the interest in these problems was one of the main causes of the progress of the pertinent theory[1].

## 13.1. Births of Both Sexes

**165.** In our ignorance of the physiological conditions determining the birth of one or the other sex, a most curious problem of zoology[2] would have been the establishment by observation the chances of a male and female birth for each species of animals. In advance, the probability of a strict equality of these chances is infinitely low because mathematical rigour never occurs in the complex phenomena of natural forces. On the other hand, the admirable harmony dominating the works of nature is quite sufficient for believing that we will discover the ratios between the habits[3] of each species and the role which fell on them, and the values of the chances of the procreation of both sexes provided by observation.

Until now, that research had only been studied by very imperfect essays; with respect to domestic animals, it apparently is not really difficult. For other species it will be necessary to give up such research had it been absolutely necessary to treat a very large number of observations. However, owing to the less complicated causes dominating here the random phenomenon whose chances we are looking for, even with medium numbers interesting and sufficiently probable results can be likely achieved.

In the beginning of the 18th century[4], when scientific curiosity led to the study of public registers of vital statistics, it was noticed that more boys had been born than girls. Nowadays, no other fact is statistically established firmer than that. According to the official register of the movement of population in France and published by the *Annuaire* of the Bureau of Longitudes[5], the following table shows [Cournot provides a table showing the annual number of male and female births in France and the corresponding ratios for 1817 – 1840].

The mean ratio 1.0631 results from more than 23 *mln* births; it differs almost by 0.009 from the maximal value in 1817 and by 0.0075 from its minimal value in 1830. For determining whether that ratio varies in time appreciably and progressively, Charles Dupin [(1842) − B. B.] calculated it for five-year periods from 1801 to 1840 inclusive, both for France and England as shown in the table below. [Cournot provides a table showing these data for France and, only for 1801 – 1830, England & Wales. The ratio for the latter had invariably been smaller (mean value 1.0442).] However, the table is corrupted by noticeable anomalies so that a progressive diminution of the ratio[6] can not be considered as sufficiently established.

**166.** For finding out whether climate influences the studied ratio, the Editors[7] of the *Annuaire* separately considered two groups of departments, 8 in the north and 15 in the south. [Cournot provides a table showing that data totally from 1817 to 1839. The ratios were 1.0635 and 1.0618 respectively.] That table leads to the conclusion that the chances of birth of either sex are appreciably free from variations due to the influence of latitude [Quetelet (1836, t. 1, p. 42) − B. B.]. The same follows, at least for the time being, from the discussion of the statistical documents of the principal European states. [Cournot provides a table[8] showing the data for France, Belgium (for 1915 – 1939!), Holland, Great Britain, Sweden, Russia, Portugal and several German states. Great Britain and Sweden had the smallest ratio (ca. 1.045) and Russia (too large for having a single ratio), the largest (1.089).]

Observations in Egypt [Laplace (1814/1995, p. 39) − B. B.] and the Cape of Good Hope seem, on the contrary, to indicate an appreciable variation directly or obliquely caused by the influence of the climate. According to what follows below, it can not at all be doubted that the climate influences the chances of the birth of both sexes only when it is sufficiently strong for profoundly modifying the morals and habits.

It is also probable in advance that for such different races as the white and the black, the chances ought to be subjected to appreciable variation [ought to be appreciably different] independent from any moral or climatic influences[9].

Instead of grouping the departments according to their latitudes, Ch. Dupin got the fortunate idea [Quetelet (1836, t. 1, p. 78) − B. B.] to form two groups of 24 maritime and 62 interior departments[10]. Again discussing five-year periods, he constructed the following table. [Cournot provides a table showing the period 1801 – 1840. The mean ratios were 1.0574 and 1.0677 for the former and latter groups respectively.] In spite of some objections made by Demonferrand (1842), it seems that these results should be seriously considered owing to their concordance.

**167.** The influence of morals and social habits on the chances of birth of both sexes became doubtless[11] after babies born in and out of wedlock as well as those born in towns and rural areas began to be distinguished. From 1817 to 1839 illegitimate births in France amounted to 814,524 boys and 781,238 girls, ratio 1.0426. A similar result follows from the discussion of statistical documents of the main European states.

The influence of living in towns[12], and especially in cities, is not stronger contestable. The Belgian ratio 1.0654 (§ 166) rises to 1.0670 for those living in villages and decreases to 1.0607 for the urban population. In Paris, the number of male and female births during 1817 – 1840 inclusive was 340,817 and 329,142, ratio 1.0355. That ratio varies as shown in the table[13] [Cournot provides a table showing the annual births of either sex during that period. The totals coincide with the numbers provided above, so that the table is compiled for Paris (which is not expressly stated). The corresponding ratios vary from 1.0163 (in 1838) to 1.0572 (1829). Cournot remarks that the numbers of girls born in 1822 and 1823 are the same, but that this equality is not a mistake; it is confirmed by their being broken down in Fourier's *Recherches statistiques* … The volume of that source is not provided.]

The number of illegitimate births is proportionally much larger in Paris than for France in its entirety. This, however, is not at all the only cause of the difference between the means for Paris and the whole kingdom. The number of illegitimate babies born in Paris during the period under consideration was 117,605 boys and 114,031 girls, ratio 1.0319. It becomes 1.0377 when only considering those born in wedlock. Therefore, the causes diminishing in Paris the excess of male births act upon all the babies and even more strongly on those latter and similar results are observed in the main European cities.

Finally, if believing some observers, all causes tending to weaken physical force [Quetelet (1836, t. 1, p. 49) − B. B.], tend also to decrease the preponderance of male births. On the contrary, according to Sadler and Hofacker [Quetelet (1836, t. 1, pp. 51 – 53) − B. B.] the cause of the excess of male births only consists in that usually the father is older than the mother. The chance of the birth of either sex does not depend on their absolute ages, but only on that difference. Regrettably, the numbers on which these statisticians had supported their conclusion were insufficient[14] for solving such an important problem.

**168.** It should have been natural to inquire whether the parents' preference for boys is the cause or one of the causes of the excess of male births. *As a consequence of that preference*, asks Prévost (1829) from Geneva, *will not the parents prevent after male births the increase of their family?* […] *The parents have a son; when various causes hinder the increase of their family, they will be perhaps less worried if their main expectation was fulfilled which would not have happened had they no male infants* [Quetelet (1836, t. 1, pp. 48 – 49) − B. B.].

Long before Prévost formulated that argument, Laplace had refuted it. Actually, it is the same as made by a gambler who claims to have changed the chances of the game and the mean result of a large number of trials by adopting a system of quitting the game after gaining once or in many sets (§ 62)[15]. However, the problem can be studied from another point of view, not noticed yet which leads to an explanation of one of the causes of the excess of male births.

Suppose that the parents believe that setting the life of a girl causes more serious demands than settling a boy. Then there will be families whose multiplication stops after the birth of one or many girls although this will not happen after the birth of the same number of boys.

It is however probable that the chances of a male conception vary from one couple to another and that they are higher for a couple which had given birth to one or many boys than for another couple which had the same number of girls. Then those marriages will have more chances to extend their fertility whose first children were boys so that the superiority of the chance of a male birth is thus augmented.

The same consequence can be reached by quite different considerations. Nowadays it is well known that the mortality at earliest ages is appreciably higher for boys than for girls [Quetelet (1836, t. 1, p. 156) − B. B.]. It is even less doubtless that the babyhood mortality suppresses one of the obstacles limiting the number of births. The mean fertility of marriages is therefore higher for those marriages which had given birth to boys; that is, for those for which the mean chance of procreating a boy is higher.

If these two conclusions are justified, the causes rendering the setting up of boys in the world more difficult or increasing the cost of educating them, as well as those that decrease the infant mortality tend to diminish the excess of male births. However, to all appearances that excess is maintained in virtue of purely physiological causes inherent in the constitution of mankind. Above all, it would be necessary to determine separately the chances of the sex of the firstborn and to see whether they appreciably differ from the general mean. Many statisticians[16] actually thought that for those firstborn the excess of male births was apparently less, but that result was contested and we regard it as one of those which will be most interesting to resolve definitively.

**169.** Conjectures were also made about the causes of the decrease of the chance of illegitimate male birth [see Villermé (1832) − B. B.]. And, first of all, if this chance actually diminishes among the firstborn because of the mentioned causes, it should therefore decrease when passing from legitimate births to illegitimate in which the proportional number of the firstborn is certainly much larger. If the excess of the father's age over the mother's, the severity of morals and hard physical work (if it only maintains and develops the body's vigour) are the causes of the preponderance of the male births, they should evidently tend to increase that excess in legitimate births.

In addition, the mortality of the newborns and the proportion of stillbirths are considerably higher for the male sex[17] and it is certain that the ratio of male/female conceptions very appreciably exceeds the ratio of the respective births. And if during uterine life the foetus of the male sex is less viable or has less chances for resisting the causes of destruction; and if on the other hand (which can not regrettably be doubted) more numerous of these threaten the embryo of an illegitimate conception [C. Bernoulli (1838) − B. B.], the combined action of both these circumstances should tend to diminish the proportion of illegitimate male births.

It is also possible, as can be suspected, that false declarations or utterances conceal from the registers of vital statistics the veritable proportion of the sexes among those born out of wedlock. It is not unlikely that among the abandoned babies registered as illegitimate although born in wedlock boys are proportionally less numerous than in general. It is more difficult to admit, as some authors had proposed [C. Bernoulli (1838, pp. 62 – 63) − B. B.], that in countries with a proper (? - O.S.) registration of vital statistics illegitimate births of a considerable number of babies, especially boys, is concealed.

**170.** The mean ratio of the male/female births for France in entirety for 1817 – 1840 inclusive (§ 165) was 1.0631 and the mean value of the probability of a male birth was 0.51529. The *weights* of this result (§ 107) is expressed by the number

$$\frac{23{,}215{,}333^{3/2}}{\sqrt{2 \cdot 11{,}962{,}811 \cdot 11{,}252{,}522}} = 6817.$$

Suppose that the causes influencing in the same way all the births or a group of births are not undergoing irregular or progressive variations, so that large numbers only compensate the influence of causes varying irregularly from one birth to another. Then we can bet 45,000 against 1 on the error due to anomalies of chance of that number 0.51529 not to exceed 0.00044 in either direction, pending that the initial documents were exact.

After excluding the last year of that period, the numbers of male and female births will be 11,473,437 and 10,789,578, ratio 1.0636, and the probability of a male birth, 0.51541 with weight 6676. In particular, for 1840 the ratio of the male births to all of them is 0.51388.

Insert

$m = 22{,}263{,}015$, $n = 11{,}473{,}437$, $m' = 952{,}318$,
$n' = 489{,}374$, $\delta = 0.51541 - 0.51388 = 0.00153$

in formulas of § 108, then $P = 0.1834$ and $\Pi = 0.5917$ for the probability that the deviation should not be attributed to chance. It is seen that that probability is very low and has no possible meaning although for the year 1840 that ratio had attained a value only exceeding it in 1830 and 1828.

The total number of births for 23 years from 1817 to 1839 inclusive for France in entirety amounted to 22,263,015 and the mean annual number of births was 967,957. 1/86 of it, or 11,255, we consider as the mean number of these births in a department with a mean population. Insert $m = 11{,}255$, $p = 0.51541$ and $l = 0.01541$ in formula (33.1), then $t = 2.313$. The corresponding $P = 0.99893$ well enough measures the probability that the yearly ratio of male births to all births for such a department is contained within interval $0.51541 \pm 0.01541$ or 0.5 [0.50000] and 0.53082.

A half of $(1 - P)$ or 0.00053 well enough measures the probability that that ratio is smaller than 0.5 or that the number of female births will fortuitously exceed the number of male births. That singular fact should therefore only occur about once in 2000 cases, but, according to the yearly tables of the movement of population compiled by the prefectures, it happened 37 times in 23 years, − in 86·23 = 1978 trials. Independently from the remark made in § 166, we ought to conclude that the mean chance of male births experiences very considerable variations from one department to another and year to year.

Nevertheless, many of those registers presenting this anomaly are among those suspected by Demonferrand [(1838, p. 33) − B. B.] in other respects. For being better assured in their accuracy, it should be interesting to follow the succession of those deviations (? - O.S.) and to see whether they indicate perturbations in the causes dominating all the births at once or whether they can not unlikely be attributed to causes whose variations are fortuitous and independent from one birth to another.

### 13.2. Laws of Mortality and Population

**171.** Were it possible to guard a living being against all accidental causes of destruction whether leading to sudden and violent death or calling forth diseases with death following after more or less time, we will observe the *natural* duration of life determined by intrinsic conditions of his organization. That duration undoubtedly differs from one individual of a species to another but a rather restricted number of observations is sufficient for obtaining a considerably stable mean, a measure of the longevity of the species. We can more or less approach those conditions and therefore approximately solve one of the most interesting problems of zoology, of comparing various species of animals with respect to their longevity and discover, if possible, the law according to which it depends on the variety of organization and functioning of the species and the action of the environment. In this direction we only have sketchy information.

Actually, there are only a few or none at all individuals dying naturally because of failed vitality. All of us are incessantly exposed to causes of destruction against which we struggle more or less successfully according to our power. An atmospheric variation not influencing a young man or leading to his passing indisposition causes the death of an old man. In this sense, as Bichat [(1799 – 1800/1822, p. 2) − B. B.] put it, life is resistance to death. Therefore, even without distinguishing the causes of death, variations of mortality due to age, sex and other conditions provide precious indications about variations of vitality.

It is certain that mortality can decrease either because the action of destructive forces becomes less intensive or the resistance of the vital forces strengthens. And it is also doubtless that abundant energy of those forces during certain periods of life can multiply the dangers and indirectly contribute to the increase in mortality.

When considering mankind in particular, the knowledge of the chances of mortality is not only highly important for a physician, an administrator, an economist; it is of most lively interest for each of us. It can prevent us in ordinary life from exaggerating fears and hopes, can facilitate our submission to the severe laws of nature.

**172.** At the mid-17th century the celebrated Jean de Wit, a statesman and geometer, had been studying the probabilities of human life for calculating [the price of] annuities. Priority in such problems naturally belonged to a nation that passed ahead of all others in banking and credit operations [in particular]. As testified by Leibniz, Hudde, another Dutch geometer, initiated like Jean de Wit in management, also wrote on the same subject [see Haas (1956)].

However, the first mortality table [life table] by the astronomer Halley compiled from the registers of the town of Breslau, appeared in 1693 [in 1694][18]. That subject had always preferentially attached the attention of statisticians and nowadays the number of the published pertinent tables is considerable. Nevertheless, their compilation is so difficult that they greatly diverge, and much time will undoubtedly pass before we arrive at quite satisfactory results[19].

Suppose we select at random a large number of newborns, 10,000 for example, and follow their lives until death. It will then be possible to show the number of survivors to each age, then to compile a mortality table after which it will be easy to obtain a table of probabilities of the duration of human life. That table will be corrupted by fortuitous errors affecting each table of probabilities compiled from a restricted number of observations even when the conditions of randomness remain invariable during all the trials. Moreover, the table will be affected by sudden and irregular variations happening in the causes of mortality during the same time. Thus, a visitation of an epidemic in the twentieth year of observation [of study] will perturb the table at years 19 and 20. Furthermore, if during the time necessary for the compilation of the table the causes of mortality could have experienced slow and progressive variations, all the numbers in the table will be affected; its part describing old age will not anymore correspond with its other parts.

Organizations or associations uniting large numbers of people pursuing common interests can have their own registers and compile mortality tables for their members which will be perturbed by the same causes. Moreover, their initial dates are usually uncertain since members do not enter unions at the same age, they say nothing about the law of mortality for the first years of life, and they can only be applied to people of the same stratum.

**173.** The law of population of a given state, or of distribution of the entire population by ages, is evidently connected with the law of mortality. Population arrives at a *stationary number* when the yearly number of births equals the mean yearly number of deaths and there is no *external movement*, or emigration and immigration compensate one another. However, that double condition is not sufficient for the population to reach the *stationary state*. We can suppose, for example, that the number of births decreases but that at the same time the causes of death become less energetic and the number of deaths also diminishes. The mean duration of life therefore lengthens whereas the number of the population does not change. The law of population can vary even if there is no variation in the numbers of births and deaths. Mortality will, so to say, transfer elsewhere and will not affect the various ages the same way as before.

If the law of population of a state can be thought stationary, and emigration and immigration compensate one another at each age, a census discovering that law will at the same time and with the same precision indicate the law of mortality. Suppose for example that having 10,000 yearly births we find on 1 Jan. 1841 6000 people of age 20 – 21 years, 6000 remaining out of the 10,000 born in 1820. Because of the stationary state of the population, 6000 people will remain on 1 Jan. 1862 out of the 10,000 which should be born in 1841.

Suppose that the same census indicated 5900 people aged 21 – 22 years; we will conclude that out of the 10,000, 100 die during their 22$^{nd}$ year and that out of the 6000 of 21 years 100 die before reaching age 22. The ratio 100/6000 is the *yearly danger*, or the probability for individuals reaching a certain age to die the next year. Its variation can be regarded as the indicator of the variation of the action of causes of mortality at various ages.

The preceding suffices for understanding the compilation of mortality tables by a census [providing the distribution of the population] by ages under the hypothesis of a stationary law of population. We also ought to remark that for the earliest ages, when mortality varies very rapidly, it is proper to compile these tables for each month rather than year.

**174.** Consider a table indicating the number of survivors at a given age. Suppose that out of 10,000 babies registered as having been born at the same time, 6000 reach the age of 21 and 3000 live at 65 years. We conclude that 34 years is the median value of the time left for those aged 21 years. This is what the authors usually call *probable life* (§ 68). It is ordinarily indicated in mortality tables for each age. The probable life at the moment of birth, or the median duration of life, is the age at which the number of those born at the same time is halved.

When following separately 10,000 babies born during the same year the sum of their ages at death divided by 10,000 is the mean duration of life, or the *mean life* (§ 67). And if we follow in the same manner the 6000 people reaching 21 years, the ratio of the sum of their ages at death divided by 6000 will be their *resting mean life*.

Mortality tables indicate the number of babies dying before reaching their first year, their number dying during their second year, etc. Suppose that all deaths during a year happened at the same time, for example, in the mid-year. When calculating the mean, we will obtain then an approximate value of the mean life. That value will be more exact if the table provides deaths by months, and in any case its error can be made negligible as compared with those occurring because of the imperfection of the data.

Under the hypothesis of a stationary law of population the mean life, in its proper sense, reckoned from the moment of birth is equal to the quotient of the number expressed by the total population divided by the number of yearly births. Suppose as above that there are 10,000 yearly births and 100 of them die during their $22^{nd}$ year. Suppose also for the sake of simplification that all births and deaths happen during the same day of the year. Out of the total population on 1 Jan. 1841 those who ought to die at age 22 consist of a 100 born in 1840 and dying in 1862; of a hundred born in 1839 and dying in 1861; …; and finally of a 100 born in 1819 and dying in 1841. The total number of those people is equal to 100·22 or 22·100.

Consequently, we will again estimate the total population by summing the products of each age by the number of babies born in the [chosen] year and destined to die at that age. However, that sum divided by the number of yearly births is precisely the mean life reckoned from birth. Inversely, the mean life is expressed by the ratio of the total population divided by the yearly births[20].

**175.** Practical difficulties of a general census[22] and especially of censuses by age are great. Interior and exterior movements of the population as well as municipal and private interests prevent a high precision of the total numbers whereas an exact distribution by ages should be considered impossible. However, a frequent repetition of that operation carried out in France each five years, perfection of the administrative system, and means of checking provided by the registers of vital statistics and tables of military recruitment inspire hope that the results of the five-year censuses will someday serve as a solid foundation for statistical investigations.

Instead of distribution by ages it is possible to apply the yearly acts of deaths. By the terms of the French law, they should indicate the age of the deceased so that the numbers of dying at each age can be immediately found. After comparing them with the numbers of yearly births we can compile a mortality table, once more under the hypothesis of a stationary law of population. In turn, it can serve for compiling a table of population such as provided by censuses by age. That method, much more practicable than the previous, has been actually applied from the very beginning[23]. Its inexactitude is caused by inexactitudes in the registers themselves, especially concerning ages, and by the chances of error inherent in compiling the acts of deaths.

**176.** Without, for the moment, taking into account these sources of error the determination of the law of mortality and all its connected elements will still be affected because of 1) Anomalies properly fortuitous and resulting from the restricted number of trials of the same randomness. Their influence can be indefinitely decreased by applying ever larger numbers. 2) Anomalies occasioned by sudden changes in the conditions of randomness or causes of mortality.

For example, in 1832 an epidemic [of cholera − B. B.] led to devastation and the census of 1842 will indicate its traces in all strata of population older than 10 years. It will register fewer individuals than there ought to have been as compared with those of a younger age group. The census will indicate a much higher mortality in the group situated between younger and older ages. On the contrary, the registers of the deaths in 1842 will indicate a very low mortality in the old ages which consists of a fewer number of individuals and therefore has fewer deaths than it should have had ordinarily. This remark is all the more applicable to the lacunas left in the virile population by long wars such as those that racked Europe for 24 years after the French revolution.

**177.** Emigration and immigration without being compensated for each age follow an appreciably constant law, since the law of population remains stationary [Fourier (1821b, § 67, p. 51) − B. B.]. However, for deriving in this case the law of mortality either by taking censuses by ages or studying the registers of deaths, we ought to know the laws which emigration and immigration are obeying.

For example, the population of cities consists of people of every occupation coming to find a job or seeking pleasure whereas many babies are sent far away to wet nurses and die there. Because of these various causes the acts of deaths of a city, when compared with the numbers of yearly births, can not provide a fair idea about the chances of death without taking into account the exterior movement. Inversely, the law of the exterior movement supposed constant can be derived by comparing those acts with a census by ages [Fourier (1821a) − B. B.].

**178.** Finally, all nations influenced by our European civilization are still far from that stationary state to which we have supposed the laws of population and mortality are led. On the other hand, the law of mortality being influenced by changing progressive and secular causes indirectly varies the law of population whereas the conditions of the development of the population are changing directly, for example because of ploughing up[24] or introduction of new cultures, independently from the variation of the causes of mortality.

In all European states the population is nowadays increasing, although in different nations the rapidity of this process is very unequal. If we wish to compile a mortality table from the yearly acts of deaths, it is absolutely necessary to calculate and take into account that secular increase[25]. It is very easy to derive that correction if the growth of the population only results from the increase in births without modifying the chances of death at different ages. However, this is what can not be admitted. For obtaining a precise correction we should know exactly the required law. Nevertheless, when allowing for the incertitude inherent in other initial materials, it is possible to make that correction by trial and error with a sufficient precision.

**179.** We have mentioned the main causes of error and incertitude. Because of all of them mortality tables compiled until now greatly diverge. They are usually distributed in two categories, of slow and rapid extinction. Companies whose speculations include payment of annuities or life insurances base their calculations on either of these depending on their interests. In France, the table of Deparcieux compiled in 1746 had been applied for a long time and is still applied as a table of the first type. It was based on lists of the tontines of 1689 and 1696 containing about 9000 deaths.

On the contrary, the table of Duvillard published in 1806 was compiled, as he stated, from a list containing about 100,000 deaths occurring before the revolution. Today, it certainly indicates a far too rapid mortality. In 1835, Bienaymé compared the tables of recruitment with official documents about the movement of the population and left no doubt about that conclusion.

Quetelet's calculation for Belgium [(1836, t. 1, pp. 161 – 164) − B. B.] also lead to an appreciably slower law of mortality than Duvillard's, and finally Demonferrand (1838) exhaustively analysed official materials and derived a still much slower mortality for France in its entirety. Nevertheless, it considerably differed from one department to another. Here is a summary table comparing the results of Duvillard and Demonferrand and thus providing an idea of their divergence and of its limits. [Cournot provided a table showing the mean life and the yearly danger at ages 0(1)5, 10(5)105 years according to both authors, separately for each sex.]

**180.** The divergence of the tables mostly concerns the first years of life. During that period mortality is so rapid that, according to Quetelet [(1836, t. 1, p. 167) − B. B.], the number of babies is reduced by 1/10 by the end of the fist month and by 1/4 by the end of the first year. During the first month 4 times more babies die than during the second month and almost as many as during two years after the first year although then also mortality is still very high.

The maximal value of the resting mean life occurs between ages 5 and 6 years [Quetelet (1836, t. 1, p. 168) − B. B.] and the minimal yearly danger, at ages 12 – 14 years, just before reaching puberty. After age 50 the discordance between the tables becomes much less pronounced and we can regard the law of mortality as very well known for the period extending from then to very advanced ages. Then the incertitude of the tables again becomes very large because of the small number of cases of anomalous longevity.

**181.** The law of mortality is not at all the same for both sexes [Quetelet (1836, t. 1, p. 155) − B. B.]. Even before birth, as we have remarked in §§ 168 and 169, the chances of death preferentially affect babies of the male sex. This becomes doubtless by the acts of deaths for the stillborn showing 13 or 14 boys for 10 girls [Quetelet & Smit (1832) − B. B.]. According to Quetelet's studies, during the first 10 months of life mortality of boys continues to exceed remarkably the mortality of girls although relatively the former gradually lowers [Quetelet (1836, t. 1, p. 158) − B. B.]. At the age of 2 years these mortalities almost equalize.

By the age of puberty mortality of women exceeds the mortality of men. During ages 20 – 25 years, the period of the most lively passions, mortality of men becomes once more predominant after which mortality of women again takes the lead until age 50, or during all the time of fecundity. Statistics does not at all confirm the usual prejudice about the existence of a maximal value of feminine mortality at the critical period when that fecundity ceases. The yearly danger increases with age for women from the age of puberty to death but for men it passes its maximal value at 24 years and its minimal value at 30.

These results of Demonferrand's calculations [(1838, p. 44) − B. B.] agree with those of Quetelet and are among the most interesting which statistics was able to obtain. They should be thoroughly verified when it will not be feared that the law of population is still affected by the lacunas left in the virile generation by wars.

According to Quetelet [(1836, t. 1, p. 141) − B. B.], after age 50 the yearly danger is approximately the same for both sexes. Once more contrary to common prejudice, the Demonferrand tables show that it continues to be a bit greater for women.

Mortality is higher in towns than in rural areas. It considerably varies depending on local and climatic influences, on morals and occupations. We will however deviate from our plan by entering here into more detail.

### Notes

**1.** See especially Quetelet (1836) and the recently published treatise of C. Bernoulli (1840 – 1841), a professor in Basel. A. A. C. *Elements of population* apparently mean *basic principles governing* it. O. S.

**2.** Bru refers to several sources on the sex ratio at birth of domestic animals and mentions Darwin, see Sheynin (1980, pp. 346 – 347).

**3.** With respect to mankind, if replacing habits by morals or manners, this was the idea of Montesquieu. [B. B.]

**4.** Rather in the mid-17[th] century (Graunt). [B. B.]

**5.** Official registers had been regularly published since 1817 in the *Annuaire* of the Bureau of Longitudes. See Levasseur (1889). [B. B.]

**6.** On the changes of that ratio in time see Worms (1912) and Schwartz (1975, pp. 95 – 111). [B. B.]

**7.** The Editors were Poisson and the astronomer Mathieu. [B. B.]

**8.** Cournot borrowed that table from C. Bernoulli (1840, pp. 139 – 140). See also Quetelet (1836, t. 1, p. 43). [B. B.]

**9.** Schwartz (1975) expressed an opposite opinion. [B. B.]

**10.** The *C. r. Acad. Sci. Paris* mentions 52 which is undoubtedly wrong. A. A. C.

**11.** See Poisson (1830) and Quetelet (1836, t. 1, pp. 46 – 51). [B. B.]

**12.** Noted by Graunt in 1662. See also Laplace (1814/1995, pp. 40 – 41) who remarked that in Paris the ratio male/female for abandoned babies was smaller than in general. [B. B.]

**13.** Cournot borrowed his table from the *Annuaires* of the Bureau of Longitudes reproduced by Fourier in his *Recherches statistiques …* [B. B.]

**14.** Bru refers to Quetelet (1836, t. 1, pp. 56 and 57) and to both mentioned authors.

**15.** I do not agree with Cournot and I do not think that Laplace (see just above) or anyone else could have refuted Prévost.

**16.** Worms (1912, Chapter 10) cites Bertillon the father. [B. B.]

**17.** See especially Quetelet (1836, t. 1, p. 156) and Demonferrand (1838, pp. 16 – 17). [B. B.]

**18.** Graunt is known to have compiled the first (very imperfect) life table.

**19.** In 1834, the Paris Academy of Sciences set up a commission (Poisson, Mathieu, Dupin) with the aim of indicating methods for obtaining more exact mortality tables. [B. B.] But what did it decide? O. S.

**20.** Denote by $fxdx$ the probability for a person aged $x$ to die during time $dx$ and by $Fx$ the probability that a newborn baby attains age $x$ [and dies at once]. Then

$$dFx = -Fx fx dx, \quad Fx = \exp[-\int_0^x fx dx]$$

with $Fx = 1$ at $x = 0$. If $fx = f_1 x + f_2 x + \ldots$ where $f_1 x, f_2 x, \ldots$ correspond to causes of death which act independently from each other and possess their own proper laws. Then $Fx = F_1 x \cdot F_2 x \cdot \ldots$ where $F_1 x$ would be the initial function $Fx$ had there only been one cause of mortality to which corresponds $f_1 x$. If that cause is suppressed, the altered function $Fx$ denoted by $(Fx)$ will be

$$(Fx) = Fx/F_1 x.$$

The probability for a newborn baby to die at age $x$ is $-dFx$ so that his mean duration of life is

$$M = -\int_0^\infty x dFx = -\int_0^\infty Fx dx$$

where we applied integration by parts after noting that the product $xFx$ should disappear at both limits of integration. At the same time the mean duration of life at age $x$ is

$$\frac{1}{Fx}\int_x^\infty Fx dx.$$

Its median duration at birth is the root $\xi$ of equation $Fx = 1/2$ and the same duration at age $x$ is obtained from the equation $F(x + \xi) = Fx/2$.

Suppose that the law of population is stationary and disregard exterior movement, then that law will be expressed by $Fx$ so that

$$Fx dx \div \int_0^\infty Fx dx$$

is the probability that a man chosen by chance will be aged $x$ and the number of people of ages between $x_1$ and $x_2$ is proportional to the integral

$$\int_{x_1}^{x_2} Fx dx$$

and if a year is taken as a unit of time, and the number of yearly births, *N*, is multiplied by that integral it will express the total population

$$P = N \int_0^\infty Fxdx.$$

The mean life under the assumed hypothesis is the quotient of that total population divided by the number of [yearly] births. The mean age of the population is

$$\int_0^\infty xFxdx \div \int_0^\infty Fxdx$$

whereas the median age ξ is obtained from the equation

$$\int_0^\xi Fxdx = \frac{1}{2} \int_0^\infty Fxdx.$$

Denote by $D_0, D_1, D_2, \ldots$ the numbers of yearly deaths for individuals aged at death less than a year, from 1 to 2 years, … Then

$$D_0 = N[1 - \int_0^1 Fxdx], D_1 = N[\int_0^1 Fxdx - \int_1^2 Fxdx], D_2 = N[\int_1^2 Fxdx - \int_2^3 Fxdx], \ldots$$
$$N = D_0 + D_1 + D_2 + \ldots$$

conforming to the hypothesis of stationary population.

We can consider men and women separately. Then, still by the same hypothesis,

$$M' = \int_0^\infty F'xdx, \quad M'' = \int_0^\infty F''xdx, \quad P' = N'M', \quad P'' = N''M''$$

where one and two strokes are applied for those sexes respectively. Now, $N' > N''$ but $M''$ possibly exceeds $M'$. Adopting the numbers calculated by Demonferrand [(1838, pp. 46 – 47) − B. B.], although for a non-stationary population, we get

$M'$ = 38 years 4 months, $M''$ = 40 y. 10 m., $M''/M'$ = 1.0652.

According to Rickman [1831 − B. B.], in England

$M'$ = 32 years, $M''$ = 34 years, $M''/M'$ = 1.0625.

The values of the ratio $M''/M'$ are quite near to the likely values of $N'/N''$ so that there is some room to believe that for a stationary population the ratio $P'/P''$ little differs from 1. In any case, it differs less than $N'/N''$ does as though the excess of male births was the final cause for approximately compensating the shorter mean life or the higher chances of death of men. It is certainly unreasonable to admit a strict compensation and to suppose that the functions $F'x, F''x$ are adjusted so as the integrals $M', M''$ will be precisely in the inverse ratios of the numbers $N', N''$ since being determined by a system of efficient causes very different from those influencing the integrals $M', M''$. Whatever ideas we have about the final aim of nature, it certainly does not proceed with such mathematical exactitude.

Quite recently Pouillet [1842 − B. B.] proposed a law for a population exposed to the influence of perturbative causes and gradually approaching a stationary state. He expressed it by the proportion

$(D'/P')/(D''/P'') \approx N'/N''$

so that in the limit $P' = P''$ when $D'$ and $D''$ become respectively equal to $N'$ and $N''$. The objections to this law[21] proved that it was not strictly mathematical and that it was not necessary to calculate six decimals because, even after all the anomalies of chance were compensated, the limiting equality $P' = P''$ could evidently be only considered as an approximation.

Nevertheless, his indicated relation understood as a result of approximate compensation is still not less remarkable exactly because it could have denoted a tendency of natural causes to maintain an almost strict equality $P' = P''$ in a stationary population. A. A. C.

**21.** Objections were made the same year by Mathieu, Dupin and Demonferrand (*C. r. Acad. Sci. Paris*, t. 15, pp. 1021 – 1028, 1028 – 1036 and 1097 – 1106). [B. B.]

**22.** Censuses in France had [have?] been carries out in 1801, 1806, 1821 and each 10 years beginning with 1831. Plausible results were first obtained in 1841. [B. B.]

**23.** Beginning with Halley (1694). [B. B.]

**24.** The French term was *défrichements*. Bru noted that Laplace (1814, p. CV) had mentioned it, but Dale (Laplace 1814/1995, p. 85) had translated it as *easily tilled*.

**25.** In this connection Lacroix cited Euler. [B. B.]

## Bibliography

**Bernoulli C.** (1838), Sur la différence dans la population sexuelle des naissances légitimes et illégitimes. *Annales d'hygiène*, t. 19, pp. 60 – 65.

--- (1840 – 1841), *Handbuch der Populationistik*. Ulm.

**Bichat X.** (1799/1800), *Recherches physiologiques sur la vie et le mort*. Paris, 1822.

**Bienaymé I. J.** (1835), Sur les erreurs présumées des documents à l'aide desquels on a calculé en France les tables de population. *C. r. Acad. Sci. Paris*, t. 1, p. 364.

**Demonferrand F.** (1838), Essai sur les lois de la population et de la mortalité en France. *J. Ecole Polyt.*, t. 16, No. 26.

**Deparcieux A.** (1746), *Essai sur les probabilités de la durée de la vie humane*. Paris. Supplement, 1760.

**Dupin C.** (1842), Observations et calculs sur les variations du rapport entre le nombre des naissances du sexe masculine et du sexe féminin. *C. r. Acad. Sci. Paris*, t. 15.

**Duvillard E. E.** (1806), *Analyse et tableaux d'influence de la petite vérole sur la mortalité …* Paris.

**Fourier J. B. J.** (1821a), Notions générale sur la population. In author's (1821b, pp. IX − LXXIII).

--- (1821b), *Recherches statistiques sur la ville de Paris …* t. 1. Paris.

**Haas K.** (1956), Die mathematischen Arbeiten von J. H. Hudde, Bürgermeister von Amsterdam. *Centaurus*, t. 4, pp. 235 – 284.

**Halley E.** (1694), *An Estimate of the Degree of Mortality of Mankind*. Baltimore, 1942.

**Laplace P. S.** (1814 French), *Philosophical Essay on Probabilities*. New York, 1995. Translated by A. I. Dale.

**Levasseur E.** (1889), *La population française*, t. 1.

**Poisson S.-D.** (1830), Sur la proportion des naissances des filles et des garçons. *Mém. Acad. Sci. Paris*, t. 9, pp. 239 – 308.

**Pouillet** (1842), Sur les lois générales de la population. *C. r. Acad. Sci. Paris*, t. 15, pp. 861 – 874.

**Prévost** (1829), *Bibl. universelle*, Oct., p. 140.

**Quetelet A.** (1836), *Sur l'homme*, tt. 1 − 2. Bruxelles.

**Quetelet A., Smits E.** (1832), *Recherches sur la reproduction et la mortalité de l'homme …* Bruxelles.

**Rickman J.** (1831), *Preface to the Abstract of the Population of Great Britain*.

**Schwartz D.** (1975), *La fécondation*.

**Sheynin O.** (1980), On the history of the statistical method in biology. *Arch. Hist. Ex. Sci.*, vol. 22, pp. 323 – 371.

**Villermé** (1832), *Annales d'hygiène*, t. 8, pp. 445 – 459.

**Worms R.** (1912) *La sexualité dans les naissances françaises*.

## Chapter 14. On Insurance

**182.** A contract for *insurance* had been barely known to ancient lawyers[1], but, with the development of commercial institutions and spirit of association, it tends to become one of the most ordinary acts of life. The theory of insurance considered in its generality is very tightly connected with the mathematical doctrine of chances and randomness so that we devote to it a separate chapter.

Suppose, for treating at first the simplest case, that a large number, $m$, of individuals insure a thing A prone to perish because of a chance event. Let $a$ be its sale price had there been no chance of its loss; $p$, the chance of its loss during a given time, a year for example, or what can be called for short the *yearly risk*. Then there will be probability $P$ that the total yearly value of accidents or the annual sum of indemnities payable by the insurance office is contained within the limits

$$mpa \pm ta\sqrt{2mp(1-p)}.$$

Here, $P$ and $t$ are connected as indicated in § 33.

If the insurance is based on *fixed premiums* and $w$ denotes the premium rate then $wa$ is the yearly payment made by the insurant[2]. Evidently, $w$ is larger than $p$ not only because of management expenses and the interest on the circulating and reserve capital, but for providing the insurer the benefit due him owing to his industry. The yearly profit (*boni*) of the insurance office or the sum from which the management expenses, the interest on the capital and the profit (benefices) of the insurer should with probability $P$ oscillate between the limits

$$ma(w-p) \pm ta\sqrt{2mp(1-p)}. \qquad (182.1)$$

In general, when considering the insurance of *things* [of property] the risk $p$ is a very small fraction[3]. Suppose for example that $p = 0.001$, $w = 0.0015$ and $m = 10,000$. Then there will be probability 0.571 that the number of *accidents* is not either larger than 12 or smaller than 8 so that the *boni* is contained within $3a$ and $7a$ or $(5 \pm 2)a$. And rather often (about 48 times in 1000) the insurance office will suffer a loss. With $w = 0.002$, twice exceeding $p$, a deficit becomes extremely unlikely[4].

If $m = 100,000$, and $w$ remains equal to 0.0015 and still exceeds the risk $p$ by 0.0005, there will be probability 1/2 that the *boni* is contained within the limits $(50 \pm 6.742)a$. We will be authorized to consider the possibility of a loss as almost physically impossible. Admitting that the insurer wishing to expand his business or get ahead of his competitors lowers the premium rate by 0.00125, his mean benefit will become 5 times larger and in addition he could regard himself sufficiently guaranteed against losses.

**183.** It seems at first sight that if the individual or the insurance office continues operations for many years in succession with sufficient capital, the loss during one year can be compensated by the previous or future profits. Take, for example, a series of 10 years. The insurer will enjoy the same security and can lower the premium rate just as if he is in business only for one year with ten times more insurants. Here, however, we make a grave mistake [Laplace (1812/1886, p. 432) − B. B.] as Bienaymé (1839) justly noted. It occurs because of the compound interest[5]. Suppose that the office liquidates business after a decade. If an extraordinary loss had occurred during the first year, the compound interest on the advanced capital will considerably burden the office. If, on the contrary, that first year was extraordinarily favourable, the compound interest on the profit will considerably increase the definitive profit of the shareholders.

Therefore, during the same period many offices with the same number of cases and prone to the same chances can either enjoy large profits or suffer large losses, even exhaust their reserve capital and, if the number of yearly cases is not sufficient for appropriately narrowing random oscillations, find themselves compelled to liquidate their business before the fixed time. These enterprises will actually become speculations in chances.

Being placed in such conditions, they will only be able to secure themselves against the chances of ruin by exaggerating the rate. Even if they have sufficient capital for guaranteeing them from a forced liquidation, an indefinite extension of the duration of their work will not indefinitely narrow the limits of their losses or secure a mean fortuitous gain. The action of the compound interest tends to amplify the deviations almost inversely to their decrease owing to the increase in the number of cases with time.

**184.** Some authors had attempted to explain the insurant's advantages of a contract for insurance, in spite of its being a source of the insurer's profit, by the distinction between the *mathematical* and *moral expectation* (§ 52)[6]. However, in our opinion these explanations, the names of their proponents notwithstanding, are vague and arbitrary and there are no real reasons for resorting to them. Security provided by the contract for insurance to the insurant is undoubtedly a boon, but moral and inappreciable. It can not be entered in the balance-sheet of a particular insurant, nor does it increase directly his fortune or the social wealth, that is, the sum of the fortunes separately possessed by the members of the society.

Nevertheless, the institute of insurance indirectly provides an appreciable increase in both fortunes. Thus, a building which can be destroyed by a chance event is less commercially valued than another one providing the same profit, utilized under the same conditions but not running the same risk of perishing. If $p$ denotes the yearly chance of perishing; $a$, the value of the second building; $r$, its profit rate, then the depreciation of the first one will in general be much larger than $pa/r$[7]. One or another buyer will undoubtedly be disposed according to his temper and manner to depreciate either more or less considerably, but these are individual estimations as all of them are on isolated markets.

In the course of commerce they will, or tend to be reduced to a normal and mean rate. Only experience, that is, the quoting can find out the extent of the depreciation of a building because of the chance of its perishing. Depreciation is a function, which we have elsewhere [1838] called *empirical*, of the elements [arguments] *p* and *a,* whose determination is nowadays impossible because the rapid and progressive variations in the system of values require numerous observations and an appreciably stationary economy. As we said above, it is only permissible to state that such depreciation considerably exceeds *pa/r*.

Actually, as we saw, the more does the number of the insurants increase, the more can the insurer lower the premium rate just as competition also compels. The yearly sum demanded from the insurants will ever less exceed *pa* which corresponds to a depreciation of the capital ever less exceeding *pa/r*. Conversely, the more the number of insurants decreases, the stronger is the insurer compelled to increase that rate so as to secure sufficiently against chances of loss. And the increased incertitude of his speculations reduces competition which compelled him to be content with a smaller rate.

Therefore, the depreciation of buildings exceeds *pa/r*. In an imagined case of one single insurant the premium rate and the corresponding depreciation attain their maximal values. Otherwise it should be admitted (which is impossible) that that rate, after increasing with a graduate decrease of the number of insurants, will decrease once more with a further decrease of that number. However, when a buyer of a building subject to chance of an accident can not insure it, he is compelled to act disorderly both as an insurer and insurant. He is an insurer with one single insurant for whom the premium rate attains its maximal value. The corresponding depreciation of the building with its interest representing that fictitious premium should exceed its value calculated as an insurance premium in the most unfavourable case of insurance[8].

We certainly do not at all wish to say that each person acquiring a building prone to a chance of its loss will undoubtedly take into account that depreciation or even that no one will pay for it [in full] as though that chance did nor exist. We are discussing the laws of economics which dominate general and mean results compensating random deviations rather than accidental causes which determine conditions on some market.

From the moment when some office offers full insurance of a building prone to risk for a yearly premium *w* little differing from *p*, there is no more reasons for the depreciation to exceed *wa/r* which is the capital ensuring interest calculated according to the rate for such buildings located in the neighbourhood and sufficient for the payment of the yearly premium. The fortune of the owner of the endangered building and therefore the social wealth considered as the sum of the wealth of individuals will increase by the whole difference between *wa/r* and the depreciation of the building had the owner no means for getting rid of the risk.

That increase of the social wealth, which recent advances will certainly make ever more evident, is the result of the institution of insurance rather than of a contract for insurance. Indeed, it little matters for fixing the commercial value of the building whether it is actually insured. Suffice it that its owner can insure it at will. This is a point where our theory is essentially distinguished from those based on the alleged measure of moral expectation.

Apart from this effect of the institution of insurance there is another resulting from the contract itself and tending to prevent the diminution of the social capital. If the insured [uninsured] building is burned down, the owner's capital is lost, and if he builds it anew, he does not at all spend his savings but withdraws capital from other productive investments. On the contrary, if the building is insured, the capital is returned to its owner. It is obtained from the premiums of other insurants, deducted from the income of many, is a saving which, strictly speaking, could have been made but which ordinarily will not be done without a contract for insurance.

**185.** We have explained how the institution of insurance tends to increase the sell price of endangered funds and therefore the social wealth. By itself, insurance does not create actual values, nor does it oppose their loss. We can even fear that a slackened supervision will cause oftener destruction. However, we should not confuse an increase of wealth owing to a rise in values [in cost] and to increase in material production. The sum of wealth or of the values capable of entering commerce can considerably change from one epoch or country to another without the variations in material production being of the same order or even when it varies in the opposite direction. The science known as political economy ought to develop the consequences of this fundamental distinction[9].

When the contract for insurance deals with risks of fabrication and commerce it becomes or should become the cause of an increase in material production and consumption and at the same time of an increase in the wealth or the value since it encourages industrial and commercial activity (§ 54 [§ 53 − B. B.]). Because of the great risk a prudent speculator will refuse to undertake a venture capable of providing him a large profit, supporting other activities, improving the well-being of consumers and making everything profitable for the entire society. But he will not hesitate to start, if he wishes, such a venture after ensuring a guarantee against the risk by paying a premium.

**186.** In addition to the institution of *insurance by payment of premiums* there appear, both in theory and practice, institutions of *mutual insurance*. If a large number $m$ of individuals, each possessing a thing A (§ 182), unite for jointly covering the losses which some of them experience by chance during a year, there will be probability $P$ that the contribution of each of them will be contained within the limits

$$pa \pm ta\sqrt{2p(1-p)/m}.$$

In addition, each ought to participate in covering the expenses of the association. With an increase of membership the fortuitous variations of their contributions will become contained in narrower limits and [together with those contributions] tend to change into a fixed premium.

On the other hand, with that increase it will be more difficult for each interested member to supervise effectively the agents of the association. In general, it could be said that for very large numbers $m$, when competition narrowed the profits of the insurer and contained them within fair limits, the great economic principle of division of responsibilities and work[10] tends to prefer the system of insurance with fixed premiums to mutual assurance.

On the contrary[11], when the number of the insurants is very small for the competition to restrict the insurer's profit, or even when the competition exists, the extremely restricted scope of operations compels the insurer to increase the mean profit for sufficiently guaranteeing himself against the chances of ruin (§ 183), − in these circumstances the system of mutual insurance should be preferred.

**187.** In ordinary insurance it is opportune to encounter such simple problems which we first of all treat for conveniently describing general principles. Usually insurance deals with unequally valued properties not exposed to the same chances of destruction. It also often happens that the insured property can be lost completely or partly depending on the gravity of the accident. Finally, complete independence of the causes of an accident for each insured property does not at all exist frequently. For example, the same fire can destroy at once a large number of houses insured against such accidents by the same office.

When taking into account the inequality of risks for different insured properties the insurance offices usually range them in classes and sharply change the premiums from one class to another. Maintaining here the law of continuity is known to be impossible and the errors inherent in the hypothesis of a fictitious discontinuity merge with many others from which calculations can not be freed.

Denote by $m_i$ the number of the insured buildings whose value is $a_i$ and the yearly risk, $p_i$. There will be probability $P$ that the sum of the indemnities paid yearly by the insurance office will be contained within the limits [Laplace (1812/1886, pp. 431 − 432) − B. B.]

$$m_1 p_1 a_1 + m_2 p_2 a_2 + ... \pm$$
$$t\sqrt{2[m_1 p_1 (1-p_1) a_1^2 + m_2 p_2 (1-p_2) a_2^2 + ...]}. \qquad (187.1)$$

The amplitude of the fortuitous oscillations is smaller[12] when the values $a_i$ are equal with their sum remaining constant and in addition when all the buildings are running the same risk equal to the mean of $p_i$. To put it otherwise, the amplitude of fortuitous oscillations increases with an unequal distribution of either the total insured value or of what can be called the common fund of risks (§ 82).

**188.** If the insurance is mutual, and the contributions of each member of the association is proportional to the numbers $a$ and $p$, that is, to the insured value and the extent of risk, a member belonging to class $i$ will have probability $P$ that his yearly contribution is contained within the limits

$$p_i a_i \pm \frac{t p_i a_i \sqrt{2[m_1 p_1(1-p_1)a_1^2 + m_2 p_2(1-p_2)a_2^2 + ...]}}{m_1 p_1 a_1 + m_2 p_2 a_2 + ...}.$$

His mean contribution will always be $p_i a_i$ with any number of members, any risks and insured values. However, the fortuitous variation of his contribution will increase if new members deposit too large or too risky properties as compared with previously insured values. From this viewpoint, when following the usual indicated rule, admissions to the association can be more detrimental than advantageous for the existing members.

To fix this idea, suppose that there are 2000 members each of them having deposited value $a_1$ and 1000 to be admitted and to deposit value $a_2$. For the sake of simplicity suppose also that all the risks are the same. The existing members will benefit by accepting newcomers if $a_2 < 4a_1$. Otherwise admittance will contradict the aim of insurance by increasing the future yearly irregularities in the contributions.

I owe an ingenious method of overcoming that difficulty (? - O.S.) to Bienaymé's friendly communications. The supposed yearly loss $\mu$ is the sum of the mean loss

$$M = m_1 p_1 a_1 + m_2 p_2 a_2 + \ldots$$

and, with probability $P$, the random deviation $(\mu - M)$ is contained within the limits

$$\pm t \sqrt{2[m_1 p_1(1-p_1)a_1^2 + m_2 p_2(1-p_2)a_2^2 + ...]}.$$

With the total loss being $\mu$ its part burdening members of class $i$ is $\mu_i$. Multiplying each $\mu_i$ by its corresponding probability we will obtain the mean[13]

$$m_i p_i a_i + \frac{[\mu - (m_1 p_1 a_1 + m_2 p_2 a_2 + ...)]m_i p_i(1-p_i)a_i^2}{m_1 p_1(1-p_1)a_1^2 + m_2 p_2(1-p_2)a_2^2 + ...} =$$

$$M \frac{m_i p_i a_i}{m_1 p_1 a_1 + m_2 p_{12} a_2 + ...} +$$

$$(\mu - M) \frac{m_i p_i(1-p_i)a_i^2}{m_1 p_1(1-p_1)a_1^2 + m_2 p_2(1-p_2)a_2^2 + ...}. \qquad (188.1)$$

Now, the contribution of members of the class $i$ to the mean loss $M$ is proportional to the factor multiplied above by that loss and is therefore reduced to $m_i p_i a_i$ whereas their contribution to the random deviation $(\mu - M)$ is proportional to the factor multiplied above by that deviation.

When following that rule, two or more classes of members will always benefit by joining up since such a union or an accession of new members always tends to decrease the amplitude of random variations in each contribution. Suppose for the sake of simplicity that there are only two classes of insurants. The amplitude of those variations for those of the first class is proportional to

$$\sqrt{2m_1 p_1(1-p_1)a_1^2}.$$

After joining up the amplitude of the variations in the total deviation, with *t* and *P* remaining constant, will be proportional to

$$\sqrt{2[m_1 p_1(1-p_1)a_1^2 + m_2 p_2(1-p_2)a_2^2 + ...]}.$$

By the Bienaymé's rule the part falling on the first class is proportional to the factor multiplied by (μ − M) in formula (188.1). Necessarily

$$\sqrt{2m_1 p_1(1-p_1)a_1^2} >$$
$$\frac{m_1 p_1(1-p_1)a_1^2 \sqrt{2[m_1 p_1(1-p_1)a_1^2 + m_2 p_2(1-p_2)a_2^2]}}{m_1 p_1(1-p_1)a_1^2 + m_2 p_2(1-p_2)a_2^2}$$

since this inequality, after cancelling out the common factors, is reduced to[14]

$$1 > \sqrt{\frac{m_1 p_1(1-p_1)a_1^2}{m_1 p_1(1-p_1)a_1^2 + m_2 p_2(1-p_2)a_2^2}}.$$

The Bienaymé rule is very simple, but practice all by itself would have never suggested it and nothing indicates that it is known to those who had until now treated insurance. It should become the basis of contracts for mutual insurance and be a fundamental condition for uniting several insurance offices, pension or mutual aid funds. It is not seen whether that rule can be precisely applied to contracts with fixed premiums.

However, in such cases the insurer, wishing to be better guaranteed against the chances of loss, can compile the premium rate from two parts, one proportional to *pa*, and the other, to $p(1 − p)a^2$ so that with *p* always less than 1/2, the sum paid by the insurants will increase greater than proportionally to the risk and the insured value. Still, the usual exaggeration of the premium shows that insurance offices find more advantageous to offer a discount on high values for increasing the extent of their business.

**189.** For showing by a very simple example how can the unequal burden of accidents be allowed for in calculations, suppose that all the insured properties, *m* in number, are of equal value and run the same risks. Let *p* be the probability of a total destruction of a thing, or the probability that the insurer loses capital equal to the value *a* of the insured thing; *p′*, the probability of loss *a′* < *a*, etc. Then there will be probability *P* that the sum of indemnities paid by the insurance office is contained within limits

$$m(pa + p'a' + p''a'' + ...) \pm$$

$$t\sqrt{2m[pp'(a-a')^2 + pp''(a-a'')^2 + ... + p'p''(a'-a'')^2 + ...]}.$$

Actually, the range of the risks varies from one property to another. All things being equal, a large property should have more chances than a small one to avoid total destruction. If the insurance is mutual, the contribution of each member for covering the mean loss should be proportional to the value of the function ($pa + p'a' + $ …) for the insured property. The other part of the contribution covering deviations is calculated by a rule similar to that indicated above.

**190.** In insurance, all the probabilities of risk entering like elements [like arguments] in the preceding formulas can not at all be assigned in advance. However, it is possible to determine them by experience[15] and the more exactly the larger is the number of compiled facts and the better they are classified. The registers of insurance offices provide valuable pertinent documents. However, there is a certain circumstance affecting the formulas estimating error and impossible to be taken into account theoretically. It is the lack of complete independence between the chances of accidents for different insured properties. A fire destroying a whole town will derange all the commercial calculations based on the mathematical law of compensation in a succession of random independent events.

The chance of such accidents can not be reasonably ignored so that an insurance office dealing with fixed premiums will more or less resemble a gambler's lot. What we say about accidents caused by fire is all the more applicable to those falling upon harvest and possibly extending over a vast territory. If an office insuring against such chances does not considerably multiply its cases and cover a much greater territory than possibly dominated by solidary accidents and in addition does not have sufficient capital for bearing the burden of disastrous years, it can not provide its insurants a complete guarantee or sufficiently avert chances of ruin. In any case it should considerably exaggerate the premiums.

In such circumstances the system of mutual insurance shows all its advantages. The helpful effects of a social contract of which mutual insurance is only a particular case, are based on human nature itself and persist whichever is the burden of chances owing to which people unite, and whether the harmful causes against which they join their efforts and resources are solidary or independent[16].

**191.** The contract for insurance can take most various forms. The bank that countersigns commercial promissory notes and thus ensuring trust guarantees the necessary payment in due time. A fraction of the ordinary interest on the capital can be considered as an insurance premium against the debtor's insolvency. Economists disapproving of laws restricting the interest rate based their arguments on that very principle. In case of loans connected with *great risks* the probable profit of the creditor should at least be equivalent to the insurance premium for the sum of the risk plus the interest on the same capital had it not been risked or ran the usual commercial risk.

By extracting yearly deductions from the salary of an employé a pension fund assures him of an annuity after reaching a certain age or completing a determined period of working. Yearly, monthly or weekly dues paid by a worker [see Hubard (1852) and Debouteville (1844) − B. B.] to a mutual aid fund assures him of an allowance in case of illness or an aid for his widow and children in case of his death etc. A large number of establishments created for protecting the interests of the public or individuals base their activities on general principles of insurance applied to probabilities of human life. We restrict our account to indicating their main operations.

**[1]** Determination of annuities on one or more lives, reversible partly or entirely, and ensured by a capital paid at the conclusion of the contract.

**[2]** Assurance of a life pension or a retirement pension paid out after a certain age or a determined period in exchange for yearly deductions or for an annuity paid out until pensionable age.

**[3]** Assurance of a life pension for the widow in case of a preceding death of the husband or after a determined or undetermined period ensured by investing a capital or by yearly payments made by him during his lifetime.

**[4]** Assurance of a capital for the widow or children at the death of husband or father ensured by husband or father who pay annuities during their lifetime or pay a lump capital.

**[5]** [Managing] tontines, or many individuals uniting their funds with the interest on the common fund to be shared by the survivors. The capital itself is given back to the last survivor or shared by several of them according to various rules (§ 52)[17].

For solving all the problems connected with, or similar to those operations account should be taken not only of chances of death but, because of the compound interest, also of various values of the same sum payable at different moments. Denote by $A$ the sum actually paid for establishing a life annuity $a$, by $r$, the interest rate; and by $p_1$, $p_2$, …,
the probabilities that the annuitant will live 1, 2, … years. Then [Lacroix (1816/1833, § 125, p. 222) − B. B.]

$$A = a[\frac{p_1}{1+r} + \frac{p_2}{(1+r)^2} + ...]. \qquad (191.1)$$

For the sake of simplicity we have neglected here the proportional possible arrears down to the payer at the death of the annuitant.

It is well known that a company speculating in annuities should demand the payment of a capital much exceeding the capital $A$ determined by the preceding formula. It is also known that the probabilities $p_1$, $p_2$, … should be calculated by mortality tables representing the law of mortality of annuitants rather than by tables of mean mortality. For being better guarded against the chances of loss companies usually calculate $A$ by tables of slow mortality (§ 179) considered as limiting the veritable numbers.

Suppose on the contrary that it is required to determine the yearly life payment $b$ to a company for giving up capital $A$ to the payer's heirs. This problem is solved very simply if imagining as Laplace did that the payer of the yearly dues pays [instead?] capital *A* to be reimbursed at his death. Or, rather, he actually pays out the value stipulated by the company for an annuity *a* on his own life which exceeds the dues *b* by the interest on capital *A* because this interest should compensate a part of the annuity equivalent to capital *A*. Therefore, $b = a - rA$ and is always determined as a function of *A* by formula (191.1). For pertinent numerical calculations the company determines the values of $p_1, p_2, \ldots$ by tables of rapid mortality.

More details can be found in special contributions [Baily (1836) − B. B.]. Here, it was sufficient to describe our subject in the most general way.

**Notes**

**1.** On the history of insurance in France see Richard (1956). [B. B.]

**2.** The first French company to introduce fixed premium appeared in 1819, see Hamon (1896). [B. B.]

**3.** Lacroix (1816/1833, p. 241) described a manuscript in which the risk of fire for stone buildings with slated roofs was estimated as 1/20,000. According to a note that he kindly sent me, he paid an agency for mutual fire insurance for a house in Paris with evaluated value of 300,000 *francs* [Cournot mentioned the payments made in 1836 – 1842 amounting to 44 *fr* 95 *c* in 1836 and decreasing, largely monotonically, to 29 *fr* 90 *c*.] or about 12 *c* for 1000 *francs* of the property which agrees with the rate stated by Lacroix [not at all!]. As an extreme case, British offices estimate the risk of losing a whaler as 1/100. A. A. C.

**4.** When the product *mp* is not very large, either it is necessary to take into account the second term of formula (33.2) or, for achieving a more precise result, to apply the formula

$$P = e^{-pm}[1 + \frac{pm}{1} + \frac{p^2 m^2}{2!} + \ldots + \frac{p^n m^n}{n!}]$$

provided by Poisson (1837, § 81) in which *P* is the probability that there will not be more than *n* accidents for *m* insured properties.

For better judging the degree of approximation ensured by that formula we compare that probability for $p = 0.01$ and $m = 200$ with rigorous formulas (Lacroix 1816/1833, p. 244). In the following table (*p*) denotes the probability of *n* accidents, and *P* is the sum of the numbers (*p*) or the probability that the number of accidents will not exceed *n*. [Cournot provides a table for $n = 0(1)11$. The values of (*p*) and *P* are given both approximately and exactly to 6 decimals, and Bru, having checked them against Pearson (1914, Table 51), states that all of them are correct.]

The differences are of the order quite permissible to neglect. We are all the more authorized to apply that formula for larger values of *m* or smaller values of *p* than they usually are. For avoiding excessive complications we suppose that it is still possible to apply formula (33.2) reduced to its first term. A. A. C.

**5.** Lacroix had remarked on the effect of compound interest whereas Cournot followed Bienaymé (1839) whose note is not easy to understand at once. The same is true with respect to Cournot's pertinent paragraph, see below.

Apply the notation of § 182 and denote by *r* the interest. Each year the profit of the insurer oscillates with probability *P* within the limits (182.1). In the *n*-th year it will oscillate with the same probability within those limits divided by $(1 + r)^n$ so that the total profit for *n* first years will be contained within

$$ma(w-p)[\frac{1}{1+r} + \ldots + \frac{1}{(1+r)^n}] \pm ta\sqrt{2mp(1-p)[\frac{1}{(1+r)^2} + \ldots + \frac{1}{(1+r)^{2n}}]}.$$

For a large *n*, when neglecting $r^2$, we have

$$ma\frac{w-p}{r} \pm at\sqrt{\frac{mp(1-p)}{r}}.$$

The limits do not narrow and the enterprise becomes a *speculation in chances*. [B. B.]

**6.** See especially Fourier (1819). A. A. C. In addition, see Laplace (1812/1886, pp. 447 and 454). [B. B.]

**7.** Indeed, *pa/r* is the sum of the infinite progression $pa/(1+r) + pa/(1+r)^2 + \ldots$ See also Cournot's formula (191.1). [B. B.]

**8.** Someone will perhaps object to the administration of a lottery, operating under conditions of monopoly, fixing the price of its tickets in excess of their mathematical expectation (§ 55). [Some excess is necessary.] However, such an objection is not solidly based. Lotteries are directed towards a very exceptionable [and numerous − O. S.] kind of people bewildered by ignorance or passion. For them, games of chance became an unnatural necessity and they are compelled to satisfy it by coming to the administration that enjoys a monopoly on chances.

On the contrary, the man raising the price of a building he is selling, in general does not at all experience any need to gamble. Far from studying the chances inherent in realty, he is subjected to them against his will. If he accidentally does have a passion for games of chance he will rather speculate in chances luring him by the prospect of unlikely profit. A. A. C.

**9.** See especially my contribution (1838, Chapters 1 and 12). A. A. C.

**10.** That principle is usually attributed to A. Smith, but Diderot (*Encyclopédie*, t. 1, 1751, article *Art*) forestalled him. [B. B.]

**11.** J.-N. Nicollet (*Enc. moderne*, 1828 – 1830, article *Assurance*) expressed an opposite opinion. [B. B.]

**12.** This paragraph is unclear. In the literal sense the mean loss in case of constant insured values and risks is $\sum m_i p_i \sum m_i a_i / m^2$ which generally differs from $\sum m_i p_i a_i$. On the other hand, in this case the amplitude of the fortuitous oscillations can be either larger or smaller depending on the values of $p_i$ and $a_i$. For example, if $a_i = a$ and $p_i$ are variable, we have by formula (77.2), when denoting $p = \sum m_i p_i / m$,

$mp(1-p)a^2 \geq \sum m_i p_i (1-p_i) a_i^2$

and the amplitude of the oscillations is smaller than in the simple case considered in § 182.

On the contrary, if $p_i = p$ and $a_i$ are variable, denoting $p$ as above, we have by formula (77.2) the same inequality of the opposite sense. This occurs just because the function $x(1-x)y^2$ changes its concavity.

Cournot refers to § 82. It is possible that he hints at the two-stage pattern there, but

$$\sum_{i<j} \frac{m_i m_j}{m}(p_i a_i - p_j a_j)$$

should be added in the square brackets of formula (187.1) which does not further clear the situation.

The same confusion is also present in the works of Poisson, Bienaymé and almost all statisticians of the 19[th] century. [B. B.]

**13.** As usual, this is an asymptotic equality by which the mean $\mu_i$ is linear in the total loss μ:

$E(\mu_i/\mu) = E\mu_i + [(\mu - E\mu)/\text{var}\mu]\text{var}\mu_i,$

$\mu = \sum \mu_i$ and $\mu_i$ for normal laws are independent. This remarkable Bienaymé formula apparently was not directly provided by Laplace or Gauss and seems to have little interested the organizers of mutual insurance. [B. B.]

**14.** There were no common factors and the transformation was not explained. O. S.

**15.** This was understood by Condorcet. Lacroix thought that the lack of realistic evaluations of risks was one of the most serious hindrances to the development of insurance. [B. B.]

**16.** Note Laplace's statement (1814/1995, p. 89): *One may look upon a free people as a large association whose members mutually protect their property by proportionally supporting the costs of the protection.* O. S.

**17.** Why should the capital been given back?

## Bibliography

**Baily F.** (1836 French), *Doctrine of Life Annuities and Assurances Analytically Investigated*. 1866. A later edition of the original English source.

**Bienaymé I. J.** (1839), Effets de l'intérêt composé. *Soc. Philomat. Paris Extraits*, ser. 5, pp. 60 – 65. *L'Institut*, 286, t. 7, pp. 208 – 209.

**Cournot A. A.** (1838 French), *Researches into the Mathematical Principles of the Theory of Wealth*. London, 1897.

**Debouteville** (1844), *Recherches sur les sociétés de secours mutuels et de prévoyance*. Paris.

**Fourier J. B. J.** (1819), Sur la théorie analytique des assurances. *Annales de chimie et de physique*, t. 10, pp. 177 – 189.

**Hamon G.** (1896), *Histoire générale de l'assurance*. Paris.

**Heyde C. C., Seneta E.** (1977), *I. J. Bienaymé*. New York.

**Hubard G.** (1852), *De l'organisation des sociétés de prévoyance et de secours mutuels*. Paris.

**Lacroix S.-F.** (1816), *Traité élémentaire du calcul des probabilités*. Paris, 1822, 1833, 1864.

**Laplace P. S.** (1812), *Théorie analytique des probabilités*. *Œuvr. Compl.*, t. 7. Paris, 1886.

--- (1814 French), *Philosophical Essay on Probabilities*. New York, 1995. Translated by A. I. Dale.

**Pearson K., Editor** (1914), *Tables for Statisticians and Biometricians*. Cambridge, 1924, 1930.

**Poisson S. –D.** (1837), *Recherches sur probabilité des jugements* … Paris, 2003. English translation: www.sheynin.de   downloadable file  53.

**Richard P.-J.** (1956), *Histoire des institutions d'assurance en France*. Paris.

# Chapter 15. Theory of Probability of Judgements. Applications to Judicial Statistics of Civil cases[1]

**192.** It is obvious that the conditions of majority imposed on the decisions of judiciary personnel or a deliberating assembly should have a relation with the mathematical theory of chances. An accused[2], knowing nothing about his judges, about their favourable or unfavourable disposition towards him, who is not informed about the procedures following investigation and pleadings, or about the manner in which the judges communicate with each other and collect their votes, − even he will not indifferently regard whether he is tried by a tribunal of three judges who condemn by a majority of two votes against one or by six judges who may only condemn by a majority of four votes against two.

It follows that the number of the judges and the established necessary majority are all by themselves arithmetical conditions independent from the judges' qualities and personal disposition, always influencing a long series of decisions. They ought to prevail in the long run over variable circumstances concerning the composition of the tribunal in each particular case. A purely arithmetical problem is therefore present in the foundation of each law regulating the votes of tribunals. It essentially belongs to the theory of chances. However, calculations necessarily depend on certain observational material, on judiciary statistics which summarizes and coordinates sufficiently numerous facts so that the anomalies of chance will not appreciably influence the mean results.

A large country such as France[3] governed by a strictly uniform legislation and a centralized administration finds itself in most favourable circumstances for collecting judiciary statistics. It was in France that the administration of justice initiated in 1825[4] the publication of the *Comptes rendus*. Someday it will become a source of a large number of documents precious for perfecting the legislation and studying the society from the viewpoint of morals and civic duties and responsibilities[5].

**193.** When treating probabilities of judgements we always bear in mind the application of that theory (? - O.S.) to civil and criminal cases, but it is not out of order to consider the word *judgement* in its widest sense both in ordinary and philosophical language and study in the most general manner the consequences following from associating the ideas of chance and judgement.

Such an investigation is interesting all by itself and will prepare us to understand better the special theory of judgements of the tribunals. For fixing the ideas by an example suppose that a man living in the countryside is always turning his attention to the state of the sky so as to forecast the next day weather by looking at the sunset. If he registers his forecasts or judgements, and in a large number $m$ of them the event had been verified in $n$ times, the fraction $n/m = v$ will express the probability that a new forecast by the same observer will also be verified.

In other words, if he had not become either less or more able to forecast, we will note, when continuing to register his judgements, that the probability $n_1/m_1$ of verified forecasts is appreciably equal to $n/m$ for sufficiently large numbers $m_1$ and $n_1$. Suppose now that two people, A and B, independently make the same observations resulting in probabilities $v_1$ and $v_2$. If the causes (for example, depending on their physical and moral dispositions, on their health, degree of their attention) influencing the verity or error of their judgement are completely independent, then, evidently,

[1] The probability that A and B agree either when forecasting correctly or mistakenly, is

$$p = v_1 v_2 + (1 - v_1)(1 - v_2) = 1 - (v_1 + v_2) + 2 v_1 v_2. \qquad (193.1)$$

[2] The probability of the contrary case is

$$q = v_1(1 - v_2) + v_2(1 - v_1) = (v_1 + v_2) - 2 v_1 v_2 = 1 - p.$$

[3] The probability that the forecasts in the first case were correct

$$V_1 = v_1 v_2 / [v_1 v_2 + (1 - v_1)(1 - v_2)].$$

[4] The probability that A's forecasts in the second case were correct

$$V_2 = v_1(1 - v_2) / [v_1(1 - v_2) + v_2(1 - v_1)].$$

These conclusions[6] should be understood in an objective and absolute sense. They signify that when actually registering the forecasts of both A and B and comparing them with the future event for a large number $N$ of simultaneous observations, we will approximately get for the first case

$$pN = [v_1 v_2 + (1 - v_1)(1 - v_2)]N$$

and $V_1 N$ will be the number of verified forecasts, etc. The magnitudes $v_1$ and $v_2$ are determined by a previous series of observations as explained above (? - O.S.).

**194.** In our imagined example the verity or error of judgement of each observer can be confirmed by a faultless *criterion*, by observing the event. However, in many other cases similar criteria do not exist, their existence would have even contradicted the essence of the matter. For example, when a physician prescribes a treatment for his sick patient, no faultless criterion of the verity or error of his judgement can be elicited from the event. The patient can die although the prescribed treatment was really the best possible, or can recover in spite of its being wrong.

Suppose that two physicians are called to consult either jointly or separately on a long series of pathological cases and that it is impossible to determine directly the numbers $v_1$ and $v_2$ expressing the probabilities of their correct forecast or judgement. However, the register of the consultations will tell us how many times the physicians agreed, or not. If that series is sufficiently long, we will have an appreciably exact number $p$ entering equation (193.1) and therefore the equation between the numerical values of $v_1$ and $v_2$ which are impossible to assign by direct observations.

We should not forget that the existence of this latter equation is founded on the hypothesis that the anomalous causes tending to dispose the judgements of [physicians] A and B towards verity or error are independent from each other. We will first study the consequences of that hypothesis tacitly admitted by all those [Condorcet, Laplace, Lacroix, Poisson − B. B.] who had previously treated the probabilities of judgements.

**195.** We return to our first example about meteorological forecasts and suppose that there is a register of a series of them made by A, B and C. We retain the notation $v_1$ and $v_2$ and additionally introduce $v_3$. It can happen that all three of them agree; that A, B, and C decided contrary to B and C, to A and C; and to A and B. Denote the four respective probabilities by $p, a, b$, and $c$. Then

$$p = 1 - (v_1 + v_2 + v_3) + v_1 v_2 + v_1 v_3 + v_2 v_3$$
$$a = v_1(1 - v_2 - v_3) + v_2 v_3$$
$$b = v_2(1 - v_1 - v_3) + v_1 v_3 \qquad (195.1a, b, c, d)$$
$$c = v_3(1 - v_1 - v_2) + v_1 v_2$$

Had the numbers $v_1, v_2, v_3, p, a, b, c$ been determined by direct observation of a long series of trials, they would appreciably satisfy equations (195.1). A deviation too large to be attributed to anomalies of chance[7] will prove that the admitted hypothesis about the independence of the causes of error of each observer does not conform to reality.

On the contrary, when there does not exist a proper criterion (? - O.S.) for directly determining the numbers $v_1, v_2, v_3$, they can be found indirectly by means of the values of $p, a, b, c$ given by observation and any three equations (195.1) since the four magnitudes are connected by an equation

$$p + a + b + c = 1.$$

Because of symmetry it is convenient to have

$$v_1 = 1/2 + z_1, v_2 = 1/2 + z_2, v_3 = 1/2 + z_3,$$
$$a - 1/4 = \alpha, b - 1/4 = \beta, c - 1/4 = \gamma$$

so that

$$\alpha = z_2 z_3 - z_1 z_2 - z_1 z_3, \beta = z_1 z_3 - z_1 z_2 - z_2 z_3,$$
$$\gamma = z_1 z_2 - z_1 z_3 - z_2 z_3, \qquad (195.2a)$$

$$z_2 z_3 = -(\beta + \gamma)/2, z_1 z_3 = -(\alpha + \gamma)/2,$$
$$z_1 z_2 = -(\alpha + \beta)/2, \qquad (195.2b)$$

$$v_1 = \frac{1}{2} \pm \sqrt{\frac{(a + b - 1/2)(a + c - 1/2)}{1 - 2(b + c)}},$$

$$v_2 = \frac{1}{2} \pm \sqrt{\frac{(a+b-1/2)(b+c-1/2)}{1-2(a+c)}},$$

$$v_3 = \frac{1}{2} \pm \sqrt{\frac{(a+c-1/2)(b+c-1/2)}{1-2(a+b)}}.$$

For $v_1$, $v_2$, $v_3$ to be real, it is necessary that

$$a + b - 1/2, \; a + c - 1/2, \; b + c - 1/2 \qquad (195.3)$$

are negative, or that two of them are positive and one, negative. And for $v_1$, $v_2$, $v_3$ to be contained within the interval 0, 1, it is necessary, as can be easily shown, that all the magnitudes (195.3) should be less than 1/2 in absolute value, a sufficient condition for which is $a + b < 1$, $a + c < 1$, $b + c < 1$.

If these various conditions are not satisfied by the values of $a, b, c$ as given by observations, the mentioned hypothesis should be rejected. The magnitudes $v_1$, $v_2$, $v_3$ are double-valued since the radicals are; however, equations (195.2) do not allow to combine them indifferently. The remark about the signs of magnitudes (195.2) requires that all the three products

$$z_1 z_2, \; z_1 z_3, \; z_2 z_3 \qquad (195.4)$$

should be positive, or two negative and one positive.

Suppose they are all positive, then $z_1$, $z_2$, $z_3$ are of the same sign, positive or negative. If another hypothesis about the signs of magnitudes (195.3) or (195.4) is admitted, then either two hypotheses about the signs of $z_1$, $z_2$, $z_3$, or two systems of values of the unknowns $v_1$, $v_2$, $v_3$ will correspond to each of these magnitudes.

**196.** This analysis is naturally applicable to judgements in tribunals composed of three judges as they are in most French tribunals of the first instance. If the legal secretary may register the votes of each judge, his list summarizing a very large number of cases will provide the numbers $a, b, c$ and the values of $v_1$, $v_2$, $v_3$ can then be derived by the preceding formulas.

It is impossible to determine these values directly since the verity or propriety of a tribunal's judgement can only be checked by another tribunal, also exposed to error however supposedly enlightened can its members be. Calculations indeed provide two systems of values for the numbers $v_1$, $v_2$, $v_3$ but in most cases one of them is at once inadmissible which eliminates any ambiguity. For example, according to those systems each of the magnitudes $v_1$, $v_2$, $v_3$ is either larger or smaller than 1/2, but it is repulsive to admit that each of the tribunal's three judges is oftener mistaken than correct. It will mean taking seriously Rabelais' jibe at the judge who decided his cases by throwing dice. The first system is therefore the only admissible. And calculations indirectly provide the values of $v_1$, $v_2$, $v_3$ as surely as by direct observation provided we had a faultless criterion for checking such judgements.

We should note that all these consequences rest on the hypothesis, whose righteousness is discussed below, on the independence of causes individually predisposing an error of each judge. This means that such an error is indifferently combined with an erroneous or proper decision of another judge. However, even before calculations the values of *a, b, c* often directly prove[8] that that hypothesis is inadmissible since $v_1, v_2, v_3$ become imaginary or negative, or exceed unity which can at least limit from above their veritable values, see below.

If it is possible to consider in advance that these three magnitudes are equal one to another or that the chances of error of each judge are the same, then, after rejecting the values of $v < 1/2$, equation (195.1a) will provide

$$v_1 = \frac{1}{2} + \frac{1}{2}\sqrt{\frac{4p-1}{3}}. \qquad (196.1)$$

By that hypothesis it will suffice, as Laplace had remarked, to know the rate *p* of the unanimously decided cases. It will not be either difficult or inconvenient to determine that rate which should exceed 1/4 for $v_1$ not to become imaginary, not to signify that at least one of the admitted hypotheses (on the independence of the causes of error of each judge or on the equality of the chances of error, again of each judge) was inadmissible. However, although the second hypothesis is certainly arbitrary and in general wrong, it is easy to see that $v_1$ in formula (196.1) appreciably expresses the mean of the veritable values of the three magnitudes $v_1, v_2, v_3$, at least when the differences between them is not relatively very large. Suppose for example that they are 0.6, 0.7 and 0.8. The mean is 0.7 and the corresponding value of *p* is 0.36. Substituting that value in (196.1) we find that $v_1 = 0.692$, only by 1/125 less than the veritable mean.

It is doubtless interesting to find out for each tribunal composed of permanent [irremovable] judges a value so close to the mean probability of verity or error of each of them. It is therefore desirable that the administration takes steps ensuring the knowledge for each such tribunal the rate of unanimous decisions reached during a decade. Purely formal decisions (conciliatory, by default etc) should certainly be excluded.

**197.** In addition, we should remark that formula (196.1) and all those supposing that the chances of error for all the voters are the same become applicable if the council or the tribunal is composed not of permanent judges, but of those chosen randomly from a long list. In such formulas, the letter *v* denotes the mean of the true values of that magnitude for each person entered there. So, if that list contains $n_1$ people having $v_1$ as the value of *v*, $n_2$ people having $v_2$, … then *v* will be [the generalized arithmetic mean of $v_1, v_2, …$].

Actually, it is possible to imagine that the first judge chosen by chance deposits his vote in urn A, the second, in urn B, … Then these urns substitute judges A, B, … of a permanent tribunal. But then, that generalized mean will evidently express the probability of the propriety of the votes left in urn A as also in urns B, C, … if only the number of people entered in their list is sufficiently large so that the removal of those chosen will not appreciably alter the value of the mean $v$.

**198.** The probability that a tribunal of three judges pronounces its judgement unanimously and correctly is $v_1v_2v_3$; that it judges still properly but not unanimously is

$$(1 - v_1)v_2v_3 + (1 - v_2)v_1v_3 + (1 - v_3)v_1v_2.$$

Therefore, the probability of a correct judgement is

$$V = v_1v_2 + v_1v_3 + v_2v_3 - 2\, v_1v_2v_3. \qquad (198.1)$$

In other words, the tribunal is a moral person for whom $V$ represents that which for each of the judges A, B, C we denoted by letters $v_1$, $v_2$, $v_3$. There should always be such relations between $v_1$, $v_2$, $v_3$ that the value of $V$ in (198.1) will exceed each of them. Indeed, if $V$ is, for example, less than $v_1$, it will be unreasonable to join judges B and C to A since this will only diminish the chance of a *proper verdict*. And, representing equation (198.1) in the form

$$V = v_1(v_2 + v_3) - (2v_1 - 1)v_2v_3,$$

we see that necessarily $V < v_1$ if $v_2 + v_3 < 1$, $v_1 > 1/2$. However, in the most probable case in which each of the three numbers $v_1$, $v_2$, $v_3$ are larger than $1/2$, $V$ necessarily exceeds the largest of them.

**199.** For indicating at least the general course of calculation, consider in addition the case in which the tribunal consists of four judges A, B, C, D having chances $v_1$, $v_2$, $v_3$, $v_4$ of deciding properly. Assume that their votes in a long series of judgements are established and denote by $a$ the rate of judge A opposing all the others and by $b$, $c$ and $d$ the similar rates for the other judges. Then we will have 4 equations for determining the 4 unknowns $v_1$, $v_2$, $v_3$, $v_4$

$$a = (1 - v_1)v_2v_3v_4 + v_1(1 - v_2)(1 - v_3)(1 - v_4)$$
$$b = (1 - v_2)v_1v_3v_4 + v_2(1 - v_1)(1 - v_3)(1 - v_4)$$
$$c = (1 - v_3)v_1v_2v_4 + v_3(1 - v_1)(1 - v_2)(1 - v_4)$$
$$d = (1 - v_4)v_1v_2v_3 + v_4(1 - v_1)(1 - v_2)(1 - v_3)$$

Let

$$v_1 = 1/2 + z_1,\ v_2 = 1/2 + z_2,\ v_3 = 1/2 + z_3,\ v_4 = 1/2 + z_4,$$
$$2a - 1/4 = \alpha,\ 2b - 1/4 = \beta,\ 2c - 1/4 = \gamma,\ 2d - 1/4 = \delta$$

so that

$$\alpha + \beta = 2(z_3z_4 - z_1z_2) - 8z_1z_2z_3z_4, \qquad (199.1a)$$
$$\gamma + \delta = 2(z_1z_2 - z_3z_4) - 8z_1z_2z_3z_4 \qquad (199.1b)$$

Therefore

$$8z_1z_2 = (\gamma + \delta) - (\alpha + \beta) \pm \sqrt{[(\alpha+\beta)-(\gamma+\delta)]^2 - 4(\alpha+\beta+\gamma+\delta)},$$

$$8z_3z_4 = (\alpha + \beta) - (\gamma + \delta) \pm \sqrt{[(\alpha+\beta)-(\gamma+\delta)]^2 - 4(\alpha+\beta+\gamma+\delta)}$$

and the products $8z_1z_3$, $8z_2z_4$, $8z_1z_4$, $8z_2z_3$ are obtained similarly. Each of them is double-valued because the radicals are. However, equations (199.1) can only be satisfied when choosing the same radical with the same sign for $z_1z_2$ and $z_3z_4$. At the same time there will only be two systems of values, positive or negative, for each of the groups ($z_1z_3$ and $z_2z_4$) and ($z_1z_4$ and $z_2z_3$). Then,

$$z_1 = \pm\sqrt{\frac{z_1z_2z_1z_3}{z_2z_3}} = \pm\sqrt{\frac{z_1z_2z_1z_4}{z_2z_4}} = \pm\sqrt{\frac{z_1z_3z_1z_4}{z_3z_4}}$$

with the other three unknowns $z_2$, $z_3$, $z_4$ being expressed similarly.

Considering the cases in which A and B are of the same opinion and oppose C and D; A and C oppose B and D; A and D oppose B and C; and when all four are of the same opinion, we will have 4 other equations from which the values of $v_1$, $v_2$, $v_3$, $v_4$, can be determined. If they do not agree with those derived from the values of *a, b, c, d*, we will reject the hypothesis about the independence of the causes of error of each judge.

**200.** For avoiding an equal separation of votes and the necessity of attributing a preponderant vote to one of the judges, or of inviting other judges for removing the impasse, tribunals are usually composed of an odd number of judges. Denote by $V_m$ the probability of a proper judgement for a tribunal composed of $(2m + 1)$ judges having the same chance *v* of faultlessness. Then[9]

$$V_m = v^{2m+1} + \frac{2m+1}{1}v^{2m}(1-v) + \ldots$$
$$+ \frac{(2m+1)2m(2m-1)\ldots(m+2)}{m!}v^{m+1}(1-v)^m.$$

Had this probability been known, the value of *v* would be found by solving that equation of $(2m + 1)$-st degree. The hypothesis on which it rests, the equality of the values of *v* for each judge, is in general certainly inadmissible. However, that formula is nevertheless applicable, as explained in § 197, if the tribunal is composed of judges chosen by chance from a long list.

For determining *v* in a similar case when $V_m$ is not known, the simplest procedure is to determine by experience the rate *q* of majority judgements. We will have

$$q = \frac{(2m+1)2m(2m-1)\ldots(m+2)}{m!}[v^{m+1}(1-v)^m + v^m(1-v)^{m+1}] =$$

$$\frac{(2m+1)2m(2m-1)\ldots(m+2)}{m!}v^m(1-v)^m,$$

$$v = \frac{1}{2} \pm \sqrt{\frac{1}{4} - \sqrt[m]{\frac{qm!}{(2m+1)2m(2m-1)\ldots(m+2)}}}.$$

If the tribunal is composed of an even number $2m$ of voters, and if $q$ is the rate of an equal separation of votes, then

$$v = \frac{1}{2} \pm \sqrt{\frac{1}{4} - \sqrt[m]{\frac{qm!}{2m(2m-1)\ldots(m+1)}}}.$$

The expansion of the binomial (§ 5) is symmetric and it is easy to see[10] that the probability of a proper decision by a majority of $i$ votes is

$$v^i / [\, v^i + (1-v)^i\,]$$

so that it is the same as though the decision was reached unanimously by a tribunal only composed of $i$ judges for whom the chance of proper voting is $v$. In other words the probability of a proper vote depends not on the absolute number of votes but on the difference between the positive and negative votes.

Imagine therefore two tribunals the second of which is composed of a lesser number of judges and admit as previously that the probability $v$ of a proper individual decision is the same for both. If each tribunal has judged a very large number of cases from which $N_1$ and $N_2$ cases respectively were decided by a majority of $i$ votes, then, in general, $N_2 < N_1$, but the ratios of proper to mistaken judgements will be appreciably the same.

**201.** If a large number of the same cases are successively tried in many tribunals we can calculate the probability of a proper decision for each court just as that probability for the different judges of a tribunal by the concordance and discordance of the votes in a long series of cases. It seems that the institution of appeals and the publication of judiciary statistics in such a country as France should lead to the determination of $V$ and therefore of the mean value of $v$. However, many important remarks ought to be made.

First of all, a process, and especially a civil process, often presents complicated problems and is transformed during different phases of the procedure. The factual or legal essence of a case submitted to the judges of an appeal court can considerably differ from the subject decided in the court of first instance and the appellant can win his case without the judgement in the appeal court being a reversal in the proper sense of that of the court of first instance.

It is different with regard to cassation complaints since the plaintiff can only refer to the legal aspects of the contested sentence. However, on the other hand, if the cassation indicates that the appeal court judged mistakenly (or at least contrary to the doctrine of the cassation court) concerning a point of the complex problem submitted to it, a rejection of a complaint, as is known to all those for whom the principles of our French legal system is not alien, does not indicate either that the appeal court decided properly or that the cassation court had in essence approved the former's sentence.

If wishing to eliminate the first consideration, against which in any case neither the judiciary statistics nor the combinatorial analysis can do anything, it should be remembered that the appeals only concern a minority of cases decided in courts of first instance. The method under discussion can in any case only determine the magnitude $V$ for the cases appealed against the decisions of the tribunals of first instance. In reality, had the appeals only been determined by the plaintiff's whim, that magnitude would be the same whether the cases were appealed or not. It is the same if the decision to appeal or not also depends on the pecuniary importance of the process. Indeed, it is natural to suppose that a process of little pecuniary importance presents the same difficulties as one of great importance of that kind and that conscientious magistrates take the same care to resolve them according to the principles of equity and law.

However, we should also admit that, after losing his case, the plaintiff often acquiesces in believing that his case is weak and we find that there are considerably more properly decided processes in courts of first instance than among those sent to judges of the courts of second instance.

The personnel of tribunals are renewed in time, legislation varies, the law consolidates in some respects and new controversial problems appear so that the magnitudes $V$ and $v$ should vary in time. For covering only one period during which they remain appreciably invariable but nevertheless collecting a sufficient number of decisions we should not restrict our study to a small number of tribunals of first instance or appeal courts. We ought to, for example, apply the yearly data of the public administration for France in its entirety.

This is tantamount to supposing that there is only one French court of first instance and one appeal court. […] If tribunal $i$ whose chance of a proper decision is $V_i$ considers $m_i$ yearly appeals, the required magnitude $V$ for a fictitious court of first instance will be

$$V = (m_1 V_1 + m_2 V_2 + \ldots)/(m_1 + m_1 + \ldots).$$

For a fictitious appeal court similarly $V$, $m_i$ and $V_i$ are replaced respectively by $V'$, $m_i'$ and $V_i'$.

**202.** According to the law of 16 Aug. 1790, the district tribunals had at the same time been reciprocal appeal courts. The constitution of year III maintained the same system but left only one tribunal for each department and magnitudes $V$ and $V'$ became equal to each other. Denote by $q$ the rate of reversed appeals, then

$$q = 2V - 2V^2, \quad V = \frac{1}{2} \pm \sqrt{\frac{1}{4} - \frac{q}{2}}.$$

Nothing would have been easier than to determine the mean $V$ for that period had the judiciary statistics then, during the time of civil disorder and considerable perturbations even in the ordinary courts, been organized.

A more complicated system although similar to a certain point is still governing in France with regard to appealing the decisions of police courts. In those departments where there are no Royal courts, decisions of such kind reached at courts of first instance by circuit tribunals are appealed to the tribunal of five judges of the main city of the department. Decisions reached by courts of first instance by the tribunal of three judges of the main city are appealed, depending on the distances, either to the tribunal of a neighbouring main city or to the Royal court which indifferently receives appeals against the tribunals of the police courts of the same department. The very complication of that system provides the necessary materials for determining magnitudes $v$ and $V$ with a very good approximation, but we will not enter here into details[11].

**203.** In civil matters, denoting by $V$ the mean value of the chance of a proper decision for the kingdom's tribunals of the first instance and for the appealed cases; by $V'$; the same mean chance for the Royal courts; and by $q$, the rate of the reversed appeals, we will have

$$q = V + V' - 2VV'. \qquad (203.1)$$

This equation is not, however, sufficient for determining $V$ and $V'$ separately and the same uncertainty occurs with respect to the appeals against sentences reached by justices of the peace and sent to the district tribunals.

That uncertainty can only be eliminated by a hypothesis. At first, following Poisson, we suppose that the mean chance $v$ is the same for judges of the courts of first instance and appeal courts. We will also admit that all the decisions at those former courts are reached by 3 judges and by 7 at the latter courts. These are actually the minimal numbers stipulated by law, and they are rarely exceeded. This double assumption leads to

$$V = v^3 + 2v^2(1 - v), \qquad (203.2a)$$
$$V' = v^7 + 7v^6(1 - v) + 21v^5(1 - v)^2 + 35 v^4(1 - v)^3. \qquad (203.2b)$$

By issuing from the *Comptes généraux* of the administration of French civil justice, from the beginning of the judicial year 1830/31 to the end of the civil year 1840, we compiled the following table. [Cournot provides a table showing years, number of judgements, both confirmed and reversed; reversed wholly or partly, total numbers 85,161 and 27,141; and the rate $q$, mean value 0.3187.]

In this table, the appeals against the decisions of commercial tribunals and the judgements of the civil tribunals of first instance composed of permanent magistrates are combined. Their separation testifies about the remarkable fact that the value of $q$ is appreciably the same in both. It is as though the advantages of the civil judges caused by the permanence of their jobs and their professional studies are almost exactly compensated by a fair appreciation achieved by the knowledge of commercial deals performed by eminent merchants temporarily empowered to solve disputes occurring in such cases[12].

Substitute 0.3187 instead of *q* in equation (203.1), then this equation and equations (203.2) will provide[13]

$v = 0.686, V = 0.766, V' = 0.855.$

The mean probability of a proper decision of a Royal court will be 0.950 in cases of confirmed appeals and only 0.642 in the opposite cases[14].

**204.** However, the hypothesis on which these values are based is evidently very unfavourable for the Royal courts in that their superiority over tribunals of first instance only depends on the larger number of their judges. The hierarchal constitution of the judiciary corps must nevertheless concentrate more experience and enlightenment in superior tribunals and other statistical documents wholly confirm this opinion and determine the limits within which the values of *v* and *V* ought to be contained.

These documents consist of the lists of reversed complaints [concerning various types of courts] as indicated in the following table.

[Cournot provides a table showing the number of reversed complaints about the decisions of the Royal courts and other superior tribunals and, separately, of civil and commercial tribunals for the years 1830/31 − 1840.]

It follows that the ratio of the cassations to the number of complaints about the decisions of the Royal courts concerning civil and commercial cases is 0.202; and the ratio of complaints about the judgements of the tribunals of first instance or commercial tribunals on which appeals are not allowed is 0.467. The first ratio can be considered quite exact; on the contrary, the value of the second can only be accepted for the time being.

Denote as previously by *V* and *V'* the mean values of the chances of proper decisions by the tribunals of first instance and Royal courts; by *V''*, the similar chance for the cassation court; by *q''* and *q'*, the ratios whose numerical values were given above. Then

$q' = V + V'' - 2VV'', q'' = V' + V'' - 2V'V''$ (204.1a, b)

and after eliminating *V''*

$V(1-2q'') - V'(1-2q') = q' - q''.$

Suffice it to combine this with equation (203.1) for separately determining *V* and *V'* independently from the hypothesis of § 203 if only admitting that these values are the same for the series of cases appealed to the Royal courts and cassations. However, calculations provide imaginary values for *V* and *V'* and it suffices to be a bit familiar with the principles of our judiciary organization for presuming in advance that that hypothesis is inadmissible: the delicate problems most often causing requests for cassation should be expounded to the tribunals of first instance and the Royal courts whereas in the mean the chances of error are greater for the latter than in cases which are only appealed.

**205.** After successively substituting $V'' = 1$ and $V'' = V$ in equation (204.1a) we obtain two values of $V$, one of them certainly smaller, and the other certainly larger than its true value. Indeed, the cassation court can itself be mistaken and it sometimes reforms its own practice. On the other hand, it is absurd to suppose that $V'' < V$. We conclude, after assuming that $q'' = 0.467$, that

$$0.533 < V < 0.630. \qquad (205.1)$$

The same reasoning applied to equation (204.1b) leads to

$$0.798 < V' < 0.866. \qquad (205.2)$$

However, in virtue of that equation $V''$ increases when $V'$ decreases so that $V'' > 0.886$ and we obtain a superior limit for $V$ when inserting $V'' = 0.886$ in equation (204.1a). Indeed, there is no reason to suppose that in a series of complaints against the decisions of the tribunals of first instance $V''$ can decrease lower than the limit not passed in a series of complaints against the decisions of Royal courts.

Inequalities (205.1) should therefore be replaced by

$$0.533 < V < 0.543$$

to which correspond inequalities

$$0.520 < v < 0.528.$$

We see that this assumption narrows the unknown values of $V$ and $v$ relative to a series of appealed cases. The residual uncertainty is less than that caused by the incertitude of the statistical materials themselves.

Inequalities

$$0.649 < v' < 0.710$$

correspond to inequalities (205.2). Here, $v'$ denotes for the judges of the Royal courts the same as $v$ denotes for the judges of the courts of first instance, and it can be thought that the true values (? - O.S.) are nearer to the inferior rather than to the superior limits.

For a series of appealed cases the value of the ratio $V'/V$ should be less than for those which led to requests for retrials. On the one hand, they do not in general present such difficulties, and the more difficulties there are, the more should be felt the enlightenment of the judges of the appeal courts. On the other hand, for the interested side the causes for resorting to requests for cassation, when cases are unimportant, as those ordinarily finally decided by inferior tribunals are, should seem very serious.

In a series of requests, when taking the inferior limits of $V$ and $V'$, which can not differ much, we will appreciably have $V'' = 3V/2$. Together with equation (203.1) we will have

$V = 0.668$, $v = 0.614$

and $V'$ will be appreciably equal to unity. And so, for a series of appealed cases we will certainly assume that

$0.668 < V < 0.766$, $0.614 < v < 0.686$.

We should also be on our guard against confounding, just as I remarked in § 201, the values of $V$ and $v$ for the appealed cases or those under cassation, with cases which are in general tried by courts of first instance.

### Notes
**1.** This chapter and the next one are essentially a reprint of Cournot (1838). [B. B.] *Explanation*. Courts of appeal, appeals or appellate courts review decisions reached by lower courts. In turn, their own decisions (as the decisions of lower courts) can be reviewed by courts of cassation which only verify the correctness of the interpretation of the law. O. S.
**2.** Since mentioning an accused, Cournot began discussing criminal cases which was contrary to the title of the chapter. O. S.
**3.** Cournot (1838) referred here to the previous work of Condorcet and Poisson. At the time of Condorcet and Laplace, as he stated, there was no judiciary statistics and they had to introduce arbitrary hypotheses. [B. B.]
**4.** Actually, the first volume of the *Comptes généraux* for the judiciary year 1825 appeared in 1827. [B. B.]
**5.** Cournot (1838) stated that he had postponed its publication to read Poisson's forthcoming book (1837) about which that illustrious author had informed him in *quelques commucations*. [B. B.]
**6.** The formulas above are due to Condorcet; then they appeared in the works of Lacroix and Poisson. [B. B.]
**7.** The expression *anomalies of chance* is due to Laplace (1814/1995, p. 43). [B. B.]
**8.** Boole (1854, Chapter 18, No. 4 and Chapter 21, No. 5) criticized that test. [B. B.]
**9.** That formula is due to Condorcet. [B. B.]
**10.** Same comment. [B. B.]
**11.** Concerning the other developments which can not be included either here or in the next chapter, see my paper (1838). A. A. C.
**12.** In our opinion, that point is one of the most curious occurring in the *Comptes généraux* and it is desirable to justify it in more detail. To achieve this, we (1838) distinguished more clearly than in the official reports the appeals against two jurisdictions, as shown in the table.
[Cournot provides a table showing the appeals against decisions of civil and commercial tribunals, see previous table, for 1830/31 − 1834.]
In 1840, for the first time, the appeals were thoroughly distinguished not only by those two jurisdictions, but in addition against civil tribunals with respect to civil cases proper and commercial judgements in towns where there were no commercial tribunals. The results are shown in the following table […]. A. A. C.
**13.** Poisson (1837, § 151) derived $v = 0.6832$ and $V = 0.7626$. [B. B.]
**14.** Poisson (Ibidem) derived, respectively, 0.9479 and 0.6409. [B. B.]

### Bibliography
**Boole G.** (1854), *Laws of Thought*. New York, 2003.

**Cournot A. A.** (1838), Sur les applications du calcul des chances à la statistique judiciaire. *J. math. pures et appl.*, sér. 1, t. 3, pp. 257 – 334.

**Laplace P. S.** (1814 French), *Philosophical Essay on Probabilities*. New York, 1995. Translated by A. I. Dale.

**Poisson S.-D.** (1837), *Recherches sur la probabilité des jugements* … Paris, 2003. English translation: www.sheynin.de   downloadable file 53.

## Chapter 16. The Theory of Probability of Judgements Continued. Applications to Judicial Statistics of Criminal Cases. Probability of Testimonies

### [16.1. Applications to Judicial Statistics of Criminal Cases]

**206.** In the preceding chapter, we issued from a hypothesis that the causes of error of each judge are independent so that the cases in which judge A decided properly or mistakenly were indifferently combined with those concerning judges B, C, … It was as though each face of a die indifferently combined with each face of another die. This is only true with regard to causes of errors which occur because of circumstances dominating each judge individually, such as the state of his physical and moral health, the degree of the encouragement of his attention, the habits of his mind, individual prejudices, etc.

However, there exist other causes of error which at the same time influence all those who become acquainted with the case so that the error of judge A more easily or oftener combines with the error of judges B, C, … than with the contrary event. I return to my initial example (§ 193) and suppose that two observers of the same degree of perspicacity and experience simultaneously register their meteorological predictions. Denote by $v$ the rate of successful predictions for each of them and by $p$, the ratio of the number of observations leading to coinciding predictions to the total number of observations. In virtue of equation (193.1) we will have

$$p = 1 - 2v + 2v^2, \quad v = \frac{1}{2} \pm \sqrt{2p - 1}. \qquad (206.1)$$

The register of predictions after being compared with the register of subsequent observations will determine the numbers $v$ and $p$ so that, if the causes of error of those observers are independent, the relation between them will be expressed by equation (206.1).

Suppose now that we separate the series of predictions in many categories by months or by distinguishing predictions of fine and rainy weather, etc. Each category should be quite numerous[1] for a sufficiently exact determination of the ratios $v$ and $p$. Denote by $p_1, p_2, …, p_n$ and $v_1, v_2, …, v_n$ the values of those ratios for the different categories and by $k_i$ the ratio of the number of observations in category $i$ to the total number of them, so that

$$k_1 + k_2 + … + k_n = 1. \qquad (206.2)$$

The veritable value of $v$ will be

$$k_1 v_1 + k_2 v_2 + … + k_n v_n.$$

Assume now that in each category the causes of the observers' errors are independent, then strictly for each of them

$$v = \frac{1}{2} + \frac{1}{2}[k_1\sqrt{2p_1 - 1} + k_2\sqrt{2p_2 - 1} + … + k_n\sqrt{2p_n - 1}] \qquad (206.3)$$

so that the influence of causes inclining both observers to err simultaneously is eliminated.

For simplifying the reasoning we suppose that each number $v_1$, $v_2$, …, $v_n$ exceeds 1/2 and that all radicals are taken with a positive sign. A contrary hypothesis will be examined later. If categories were not separated or if this was impossible to achieve, and if we admit that the causes of error for the general series were independent, then equation (206.1) will lead to

$$v = \frac{1}{2} + \frac{1}{2}\sqrt{2(k_1 p_1 + k_2 p_2 + \ldots + k_n p_n) - 1}$$

and it is easy to prove (§ 77) that this approximate value of $v$ always exceeds its true value (206.3).

**207.** Suppose that a direct determination of $v$ is possible. Before any classification of judgements by categories we should warn against an error of the hypothesis about the independence of the causes of error for each judge. The value of $v$ derived from equation (206.1) will exceed that provided by direct determination and the difference between these values can still increase if the value of $v$ for certain categories can fall lower than 1/2 with its being nevertheless larger than 1/2 for the general series.

On the contrary, in the ordinary case in which the value of the number $v$ can only be known indirectly by means of equation (206.1) or similar to it, nothing will warn the calculator about the mistake of his hypothesis if the general series of judgements is not sufficiently numerous[2] for allowing to subdivide it by issuing from statistical documents into partial series quite numerous for the ratios to be appreciably permanent.

When such a classification is possible, it occurs in general that the ratio $p$ varies from one category to another and becomes, in succession, $p_1$, $p_2$, …, $p_n$. And then the value of $v$ can be calculated by equation (206.3). That second value, always smaller than the first one, will still exceed the true value if for each category or for some of them the hypothesis about the independence of the causes of errors is not yet appreciably true. When statistics will be enriched by a larger number of observations the number of categories could be multiplied for obtaining a value of $v$ smaller than the preceding and closer to the true value.

**208.** For simplifying calculations suppose that

$$\frac{1}{2}\sqrt{2p-1} = z, \quad \frac{1}{2}\sqrt{2p_1-1} = z_1, \quad \frac{1}{2}\sqrt{2p_2-1} = z_2, \ldots$$

$$k_1 z_1 + k_2 z_2 + \ldots = \varsigma$$

so that the values of $v$ derived from equations (206.1) and (206.3) will be respectively

$$v = 1/2 + z, \quad v = 1/2 + \varsigma.$$

Then

$$z^2 - \varsigma^2 = k_1 k_2 (z_1 - z_2)^2 + k_1 k_3 (z_1 - z_3)^2 + \ldots + k_2 k_3 (z_2 - z_3)^2 + \ldots$$

Nevertheless, we can prove (Note 5 to Chapter 6) that if the numbers $k_1, k_2 \ldots$ are to remain positive and satisfy equation (206.2), and if, on the other hand, $a$ and $b$ ($a < b$) denote the limits of $z_1, z_2, \ldots$ the value of the right side of the preceding equation will be contained between 0 and $(b - a)^2/4$.

In addition, the numbers $p_1, p_2, \ldots$ are necessarily contained between 1/2 and 1 so that $a = 0$, $b = 1/2$ and

$$z^2 - \varsigma^2 < 1/16, \varsigma > \sqrt{z^2 - 1/16}.$$

If, for example, the general series applied without distinguishing categories provides $v = 0.9$ or $z = 0.4$, the veritable value of $v$ (obtained were it possible to multiply the categories until no other causes of error are left except the influence of those varying from one judge to another) will be necessarily contained between 0.9 and

$$1/2 + \sqrt{0.16 - 1/16} = 0.81225.$$

This formula relative to the inferior limit of $v$ will become illusory and useless when $z < 1/4$ or only a bit exceeds that fraction. In such a case we can only wait for the improvement of the exact statistical value of the ratio $v$. When this value remains stationary even if the increased number of observations permits to multiply the number of categories, we will be assured that the limit is attained, that in each category the causes of error can be thought to act irregularly and variably for each judge or observer.

**209.** Consider now a tribunal of three judges for each of whom we are justified to attribute the same value $v$. The judges are permanent and equally enlightened or randomly selected from a general list for each case with $v$ (§ 197) then being the mean of the appropriate values for each individual included in that list. Denote by $p$ the rate of unanimous decisions determined from a long succession of observations, then

$$v = \frac{1}{2} \pm \frac{1}{2}\sqrt{\frac{4p-1}{3}}. \tag{196.1}$$

The form of this expression [of this function] is the same as of (206.1) and we can apply here all the previous considerations about the successive lowering of the value of $v$ with the multiplication of the categories and about the limits of this decrease.

The probability or the rate of a proper decision is

$$V = 3v^2 - 2v^3$$

if all decisions can be included in a single series. Suppose now that that series is separated in two categories with $v$ taking values $v_1$ and $v_2$ so that the rate of proper judgement will be

$$k_1 V_1 + k_2 V_2 = 3(k_1 v^2_1 + k_2 v^2_2) - 2(k_1 v^3_1 + k_2 v^3_2)$$

and we ought to prove that it is smaller than $V$,

$$k_1 V_1 + k_2 V_2 < V, \qquad (209.1)$$

at least when $v, v_1, v_2 > 1/2$.

In this case the ratio $p$ for the general series is the mean of $p_1$ and $p_2$ and $v$ is contained between $v_1$ and $v_2$ so that we can suppose that

$$v_2 < v < v_1. \qquad (209.2)$$

On the other hand, as proved above,

$$k_1 v_1 + k_2 v_2 < v, \; k_1 < (v - v_2)/(v_1 - v_2). \qquad (209.3)$$

If the inequality (209.1) is not satisfied, and if, on the contrary,

$$k_1 V_1 + k_2 V_2 > V, \; k_1(V_1 - V_2) > (V - V_2),$$

we can conclude that $k_1 > (V - V_2)/(V_1 - V_2)$ and, because of inequality (209.3)

$$(v - v_1)(V_1 - V_2) > (v_1 - v_2)(V - V_2) \qquad (209.4)$$

since $V$ increases with $v$ and inequalities (209.2) lead to

$$(V_1 - V_2) > 0, \; (V - V_2) > 0.$$

Without changing the sense of the inequalities we can multiply or divide [them] by $(V_1 - V_2)$ and $(V - V_2)$. After substituting the values of $V$, $V_1$ and $V_2$ inequality (209.4) is reduced to

$$(v_1 - v)(v_1 - v_2)(v - v_2)[3 - 2(v + v_1 + v_2)] > 0.$$

However, $v, v_1, v_2 > 1/2$ as assumed above, and the factor in square brackets is negative whereas all three binomials are positive in virtue of (209.2). Therefore, the preceding inequality can not be valid and finally inequality (209.1) is verified if the infinitely low probable case in which both its sides coincide is excluded.

It follows that when categories are multiplied, the value of the chance of proper judgement of a tribunal provided by calculation lowers just as does the value assigned for the mean chance of proper decision for each judge.

**210.** Until now, we supposed that the numbers $k_1, k_2, \ldots$ are known and directly provided by statistical documents. They can be called the coefficients of the categories expressing the probabilities that a judgement chosen randomly from the general series belongs to categories 1, 2, … However, we can also suppose that these numbers are unknown and propose to determine them by a sufficient number of elements chosen from those directly provided by observation. On solving such a problem depends the application of the theory of chances to judiciary statistics in criminal cases.

Actually, a series of accused brought before a criminal tribunal is naturally separated in two categories, guilty and innocent with $k_1$ and $k_2$ denoting the probabilities that a person randomly chosen from the general list of the accused belongs to the respective category. Numbers $k_1$ and $k_2$ whose sum is unity are not, and can not be given directly because we never have a pertinent faultless criterion.

It is only quite likely that, taking into account our morality and our judiciary institutions, $k_1$ considerably exceeds $k_2$. The accused only appear before a tribunal after preliminary investigation which excludes those against whom there are no really serious charges.

The chance $v$ of a proper decision by each of the tribunal's judges is $k_1 v_1 + k_2 v_2$ where $v_1$ and $v_2$ are the values of that chance for the series of guilty and innocent accused. We can believe that in general these numbers are not equal at all, or, in other words, that the ratio of the condemned guilty accused to the acquitted guilty is not equal to the ratio of the acquitted innocent accused to the condemned innocent. In any case, an equality of $v_1$ and $v_2$ can only be established by experience. In general we have three unknowns, $k_1$, $v_1$ and $v_2$ to be determined by observations.

**211.** Once more we consider a tribunal of 3 judges and we are justified to attribute to each of them the same values of $v_1$ and $v_2$. For the time being assume also that the judge who condemns or acquits an accused asserts by that same act that the accused is guilty or innocent respectively. Denote by $c_1$ the ratio of the number of unanimously condemned accused to the whole number of the accused, by $c_2$ the same ratio for those condemned by a majority verdict, and finally by $a$ the ratio of the unanimously acquitted accused to the total number of the accused. Then we will have three equations in three unknowns, $k_1$, $v_1$ and $v_2$:

$$k_1 v_1^3 + (1 - k_1)(1 - v_2)^3 = c_1 \qquad (211.1a)$$
$$3 k_1 v_1^2 (1 - v_1) + 3(1 - k_1) v_2 (1 - v_2)^2 = c_2 \qquad (211.1b)$$
$$k_1 (1 - v_1)^3 + (1 - k_1) v_2^3 = a \qquad (211.1c)$$

**212.** It is reasonably indicated that a judge who acquits an accused ordinarily does not at all mean that the accused is innocent, but that in his eyes the indicators of guilt are not sufficient for determining a conviction. Inversely, a judge condemning an accused does not at all affirm with absolute certitude that the accused is guilty, but only then there exist such indications, such a strong presumption of guilt that he can not acquit the accused against whom such indications and strong presumptions are levelled without paralyzing the action of justice and compromising public security.

Those criticisms lead to the numbers $k_1$, $v_1$ and $v_2$ determined by the previous equations being relative not at all to guilty and innocent accused, but to two other categories, *convictable* and *not convictable* or *absolvable* accused. The first can, strictly speaking, include innocent accused, and the second, quite likely, many really convictable. On the other hand, if the mind clearly and at once comprehends the absolute distinction between the guilty and innocent accused, how nice it would also be to form easily a precise idea about a categorical separation of the accused into convictable and absolvable. This is the most delicate point of the theory, delicate as it is, and we ought to turn extreme attention to it.

On this point we ask permission to return once more to the fictitious example. If someone engaged in meteorological predictions predicts fine weather for tomorrow, he does not surely affirm his statement in an absolute manner but only believes that the chances of such weather are very good, good enough, let us say, for not hesitating to undertake a voyage or climb a mountain.

Just the same, a surgeon who believes that an amputation of an injured limb is necessary does not state that another choice is absolutely impossible. He only affirms that otherwise the chances of death are in his opinion sufficiently high for sacrificing the affected limb. The same remark is applicable to most of human judgements and there is nothing special about judgements in criminal cases.

**213.** And so, for returning to our subject which demanded that digression let us imagine that the accused are separated into a sufficiently large number of categories containing the guilty and the innocent so that in each the causes of error act fortuitously and independently on each judge. Suppose that in each category the value of the ratio $v$ can not decrease lower than 1/2 either for the guilty or for the innocent. Then equations (211.1) applied to each category will determine the numbers $k_1$, $v_1$ and $v_2$ according to the distinction just mentioned. Calculation will provide double values of $v_1$ and $v_2$, $1/2 \pm z_1$ and $1/2 \pm z_2$, and the positive sign should always be chosen.

According to the theory elucidated above, the same equations when applied to the general series of the accused will certainly only provide approximate values of the required ratios but that approximation will always concern the classification of the accused into guilty and innocent. On the contrary, we should admit that for numerous categories of the accused the chance $v$ of a vote conforming to reality drops lower than 1/2 and even indefinitely approaches zero.

There are doubtless many guilty accused who will almost certainly be acquitted either because of the weakness of the legal charges levelled against them or due to various causes (such as the excessive harshness of the penal law) which predispose most judges to indulgence[3]. We can not at all refuse to admit that a very small number of innocent accused will almost certainly be convicted because, owing to a fatal coincidence of circumstances, charges were heavy and compelled even the most enlightened and impartial judges to convict them. As a consequence, there are accused for whom, when classifying them as innocent and guilty and applying equations (211.1), it is necessary to choose for $v_1$ and even for $v_2$ their calculated values which are less then 1/2.

Therefore, those equations should not be applied to the general series of accused even for obtaining a first approximation or at least, when having a good reason for doubt, perhaps choose the lesser calculated value $v_1$ as being nearer to its true value. There is only one way to get over this difficulty and include the second case into the first one, viz, to consider absolvable those guilty accused for whom the chance of conviction is smaller than 1/2 and at the same time regard as convictable those innocent accused (actually and happily there being a very small number of them) for whom that chance is larger than 1/2.

Then, when changing the initial sense of the letters $v_1$ and $v_2$ and assuming that $v_1$ is the chance of convicting convictable accused and $v_2$, the chance of absolving absolvable accused, for no category of the accused these numbers by their very definition will be smaller than 1/2. And if equations (211.1) are applied to the general series of the accused as a first approximation, it will be necessary to choose those of the calculated values for $v_1$ and $v_2$ which exceed 1/2.

**214.** These explanations are useful in that they provide a precise mathematical definition of the sense attached to the words *convictable* and *absolvable*. They allow seeing clearly how the pertinent classification of the accused is connected with the state of the enlightenment and moral disposition of the stratum of citizens from whose midst jurymen and judges are chosen. It follows that if the judges belong to another stratum, or even to the same stratum as the accused, but are influenced otherwise, the accused can pass over from being convictable to absolvable or vice versa.

And so, the rate of conviction in Belgium amounting to 0.83 when crimes had been heard by magistrates fell to 0.60 after the French system of jurisprudence was re-established there. It follows, as Poisson had remarked, that the proportion of convictable accused (in our sense) sharply declined because of that change although the form of preliminary investigation remained as it was previously, so that the proportion of the really guilty accused did not considerably vary. Actually, jurymen are more inclined towards indulgence than permanent magistrates, and there are many categories of the guilty accused for whom the chance of conviction was higher than 1/2 under the previous system but became lower than that value when the vote was transferred to the jurymen. Those accused are considered convictable when formulas (211.1) or similar are applied with the vote carried out by permanent magistrates but they pass over to the absolvable when the same formulas are applied with the vote being granted to the jurymen.

This theory also allows us to foresee in what sense the results of calculation will be modified depending on the essence of the variations of criminal legislation or other circumstances influencing the votes of the jurymen. All that tends to heighten their enlightenment should increase the values of $v_1$ and $v_2$. Therefore, all things being equal, we find lesser $v_1$ and $v_2$ for jurymen voting without communication between themselves than for those who deliberate jointly and can more clearly grasp the situation[4]. On the contrary, a softening of penal legislation leading to a greater number of deserved convictions and therefore to a more efficient suppression of certain offences should be regarded as an incontestable improvement but it can lower the values of $v_1$ and $v_2$ as related to the classification of the accused into convictable and absolvable.

There had existed a category of guilty accused almost sure to be acquitted for which $v_2$ had a very large value, but the cause of the pertinent constant error which almost certainly determined the acquittal, was excluded and the fate of those accused became influenced by the causes of error acting irregularly and independently on each juryman. For those accused the chance $v_2$ decreased; if it fell below 1/2, they passed over to the convictable for whom $v_1$ could nevertheless have a value a bit higher than 1/2. These mean values of $v_1$ and $v_2$ for the general series of the accused can thus lower after the softening of the penal legislation although the number of proper judgements increases.

In general, ignorance is a cause of error acting irregularly and variably from one juryman to another. All that tends to heighten their enlightenment inclines to diminish the risky part in their verdicts and to increase the proportion of unanimous verdicts or those carried by a strong majority, and therefore to increase the calculated values of $v_1$ and $v_2$. On the contrary, the elimination of the causes of error following from dominant prejudices and natural dispositions of the human heart while increasing the number of proper judgements can also increase the role of chance, diminish the proportion of verdicts carried unanimously or by a strong majority and therefore decrease the calculated values of $v_1$ and $v_2$.

**215.** If the judicial statistics provides for our correctional tribunals composed in general of three judges the values of the elements $c_1$, $c_2$ and $a$, equations (211.1) can be applied for determining the ratios $k_1$, $v_1$ and $v_2$. Those elements are not however given and, according to our laws, can not be given. On the contrary, the judicial statistics provides all necessary documents concerning the appeals against correctional police and those ratios can be determined for the accused brought before the two degrees of jurisdiction. We can not enter into all the necessary details about this curious application and refer to our memoir (1838) already cited in § 202.

**216.** The most interesting application of the theory of probability of judgements is that whose subject is the decisions pronounced by our jurymen in criminal cases. A tradition that goes back to the Middle Ages established the jury panel of 12 members, just like in England. Otherwise, however, the jury system in those two countries rests on very different foundations.

In French legislation, the established majority for pronouncing a convictive verdict had varied many times. According to the present law, a simple majority of 7 votes against 5 is sufficient. Suppose that $N$ is the total number of the accused; $N_1$ of them, convictable by our definition, and $N_2$, absolvable. Then, $C_1$ of the accused are convicted by a majority stronger than 7 votes, and $C_2$, by a simple majority; $A$ are acquitted because the votes were equally divided; $V_1$ and $V_2$ are the probabilities of convicting and acquitting verdicts for those who are convictable or absolvable respectively, and the meaning of $v_1$ and $v_2$ remains without change. Let also

$N_1/N = k_1$, $N_2/N = k_2$, $C_1/N = c_1$, $C_2/N = c_2$, $A/N = a$.

Then $k_1 + k_2 = 1$,

$$k_1[v_1^{12} + 12v_1^{11}(1-v_1) + 66v_1^{10}(1-v_1)^2 + 220v_1^9(1-v_1)^3 + 495v_1^8(1-v_1)^4] + k_2[(1-v_2)^{12} + 12(1-v_2)^{11}v_2 + 66(1-v_2)^{10}v_2^2 + 220(1-v_2)^9 v_2^3 + 495(1-v_2)^8 v_2^4] = c_1,$$

$$792[k_1 v_1(1-v_1)^6 + k_2 v_2^5(1-v_2)^7] = c_2,$$
$$924[k_1 v_1^6(1-v_1)^8 + k_2 v_2^6(1-v_2)^8] = a.$$

If the numbers $a$, $c_1$, $c_2$ are provided by statistics, these equations suffice for determining $k_1$, $k_2$, $v_1$ and $v_2$ and therefore $V_1$ and $V_2$. Statistical documents provide at once the number $c_1 + c_2$ as well as $c_2$ at least according to the actual legislation which obliged the jury panels, as did one of the previous legislations, to indicate whether a conviction was carried by a simple majority. However, the legislation is always opposed to an indication of the majority in cases of acquittal.

Therefore, neither the ratio $a$, nor any other analogue can be provided by juridical statistics. For numerical determinations it is thus necessary to reduce the number of the unknowns just as Poisson (1837) did it since he tacitly supposed that $v_1 = v_2$. Nevertheless, our analysis above concerning the appeals against the correctional police well agrees with considerations according to which we should suppose in advance that $v_2$ and $V_2$ exceed $v_1$ and $V_1$ respectively, and that $v_2$ and $V_2$ very little differ from unity.

The causes leading to this result for permanent judges such as those who pronounce judgement in correctional police courts should stronger influence jurymen. Still, the indication of such a result by direct observation would have been so interesting, that we should ardently desire the adoption of a measure which, without revealing the separation of votes in each particular acquittal, will provide the lacking element of criminal statistics for a long series of cases. Thus, for example, for each convicted or acquitted accused the foreman of the jury panel can be obliged to deposit in a sealed box white and black tickets according to the number of acquitting and convicting votes, with those tickets to be counted yearly out of interest for the judicial statistics but without violating the secret of voting in each case. It is not difficult to show that for our goal a registration of the results of such counts will be tantamount to the knowledge of the element $a$.

Not knowing its value[5] but having every reason to believe that the value of $v_2$ is contained between that of $v_1$ and unity we can only formulate two hypotheses, $v_2 = 1$ and $v_2 = v_1$. The true values of the unknowns $k_1$, $v_1$ and $V_1$ are contained between those which correspond to the extreme hypothesis above.

According to the hypothesis $v_2 = 1$, $v_1$ will be provided by an equation of the fifth degree

$$v_1^5 + 12v_1^4(1-v_1) + 66v_1^3(1-v_1)^2 + 220v_1^2(1-v_1)^3 +$$
$$495v_1(1-v_1)^4 - 792v\frac{c_1}{c_2}(1-v_1)^5 = 0 \qquad (216.1)$$

and we will then have

$$V_1 = v_1^7[v_1^5 + 12v_1^4(1-v_1) + 66v_1^3(1-v_1)^2 + 220v_1^2(1-v_1)^3 +$$
$$495v_1(1-v_1)^4 + 792(1-v_1)^5, \qquad (216.2)$$

$$k_1 = (c_1 + c_2)/V_1. \qquad (216.3)$$

If, on the contrary, $v_2 = v_1 = v$,

$$c_2 = 792v^5(1-v)^5[k_1(2v-1) + (1-v)^2], \qquad (216.4)$$

$$c_1 + c_2 = k_1[1 - 924v^6(1-v)^6 - (2k_1 - 1)(1-v)^7[(1-v)^5 +$$
$$12(1-v)^4v + 66(1-v)^3v^2 + 220(1-v)^2v^3 +$$
$$495(1-v)v^4 + 792v^5]. \qquad (216.5)$$

Without an appreciable error we may neglect the negative term multiplied by $(1-v)^7$, and then the elimination of $k_1$ from equations (216.4) and (216.5) will be very simple and the root of the final equation in $v$ could be obtained by trial and error the more easily since we know in advance that it can not much differ from that of equation (216.1).

The value of $V_1$ is always provided by equation (216.2) in which $v$ can be substituted instead of $v_1$. Although we assumed that $v_2 = v_1$, the value of $V_2$ is not the same as that of $V_1$ since acquittals do not require the same majority as convictions. We have

$$V_2 = V_1 + 924\ v^6(1-v)^6.$$

Out of $N$ accused the number of the acquitted although convictable will be

$$P = k_1(1 - V_1)N \qquad (216.6)$$

and of the convicted although absolvable

$$Q = (1 - k_1)(1 - V_2)N \qquad (216.7)$$

which disappears if $V_2 = 0$.

**217.** In our memoir (1838) we have applied these formulas to the criminal statistics for the six years 1825 – 1830 during which the legislation admitted convicting verdicts by simple majority, although only if the majority of the 5 magistrates then forming the assize courts confirmed the majority of the jurymen. During that period we had $N = 42,300$ with the number of the convicted amounting to 25,777. That number however can not be adopted[6] for $C_1 + C_2$ since it did not include the accused in whose favour the majority of the court differed from the majority of the jurymen. On the other hand, statistics did not then directly indicate the number $C_2$ and it is only possible to calculate it indirectly by adopting a hypothesis which leaves incertitude in the result. We concluded from our calculations that $c_1 + c_2 = 0.621$, $c_2 = 0.071$ and $c_1 = 0.550$.

During four years 1832 – 1835 the law of 4 March 1831 stipulated a majority stronger than 7 votes for conviction and the new penal code permitted the jury panel to soften the penalty because of mitigating circumstances. We had then N = 28,702, $C_1 = 11,116$ so that $c_1 = 0.596$. We ought to conclude that the possibility of declaring the presence of mitigating circumstances granted to the jury panels and other extenuations introduced in the penal legislation increased the ratio $c_1$ by about 0.046.

The law of 9 Sept. 1835 introduced secret vote, or rather allowed it for the jurymen. It stipulated a simple majority for conviction but left the possibility for the majority of the court to annul such verdicts when returned by a simple majority and compelled the jury panels to mention this circumstance [simple majority] of the verdict. The *Comptes généraux* for the years 1836 – 1840 of that latest phase of the criminal legislation directly provide the numbers $C_2$, see the following table.

[Cournot provided a table showing the numbers of the accused separately for crimes against the person and against property; of the convicted for these crimes, separately and without distinguishing the categories mentioned, all this yearly for 1836 – 1840 as well as the total figures.]

We conclude that for those five years, for the total series of the accused and without distinguishing the category of crime, $c_1 + c_2 = 0.645$, $c_2 = 0.026$ and $c_1 = 0.619$. The value of the number $c_1$ increased from the second phase [of the legislation] to the third just as it did from the first phase to the second. The value of $c_2$ became much less than its calculated (hypothetically though) value for the period before 1831. This somewhat justifies the opinion that the jurymen during the previous legislation, when being perplexed, often *agreed* to return their decision by a *simple majority* to compel the court to pronounce the final judgement and thus to remove their responsibility of announcing the [final] verdict[7].

**218.** Assuming $c_1 = 0.619$, $c_2 = 0.026$ in equation (216.1) we will have by one of the extreme hypotheses ($v_2 = 1$)

$v_1 = 0.816$, so that $V_1 = 0.987$, $k_1 = 0.653$.              (218.1)

From equation (216.6) we will then have during that five-year period $P = 335$ acquitted but convictable accused with the total number of the acquitted being 13,984.

According to the other extreme hypothesis ($v_1 = v_2 = v$) $v = 0.817$ which only insignificantly differs from the previous value of $v_1$ when taking into account the degree of approximation connected with such determinations. The same remark is applicable to the magnitude $k_1$

=0.652. We have $V_2 = 0.997$ so that equation (216.7) provides $Q = 41$ convicted absolvable accused with the total number of those convictions being 25,440. We should not forget our definitions of *convictable* and *absolvable*, and we ought to be especially on our guard over confounding absolvable and innocent accused.

The *Comptes généraux* tell us that during those five years the total number of the accused convicted by a simple majority amounted to 1023 and that the assize courts used the right conferred on them by the law of 1835 in favour of 20 of them. Out of those absolved by the majority of 3 magistrates 12 were finally acquitted by other jury panels and 8 were convicted anew in spite of the influence exerted on the new jurymen by the decision of the assize court and although magistrates just as jurymen are more inclined to indulge when much time passes between crime and judgement.

**219.** The preceding results concern the general series of the accused without distinguishing their categories and assuming an essentially faulty hypothesis according to which all the causes of error act fortuitously and independently on each juryman. We (§ 207) saw how the improvement of the judicial statistics allowing to subdivide the general series of judgements into ever more categories at the same time provided means for determining the role of those causes of error which influence all the judges at once.

The *Comptes généraux* first of all separates the accused into two main categories depending on the separation of crimes into those *against the person* and *against property*. Each of these is further subdivided in many others, and still many other divisions can be established depending on the sex, age, degree of education, existence of repeated offences, essence of punishment etc. For the sake of brevity we are only discussing the first two categories. The table [of § 217] provides for them respectively

$c_1 + c_2 = 0.556$, $c_2 = 0.032$, $c_1 = 0.524$;
$c_1 + c_2 = 0.679$, $c_2 = 0.024$, $c_1 = 0.655$.

We will accompany letters $v, k, …$ with one or two strokes depending on the category of crime. First of all, we have by the hypothesis that $1 = v'_2 + v''_2$

$v'_1 = 0.796$, $v''_1 = 0.821$ so that
$V'_1 = 0.979$, $k'_1 = 0.568$; $V''_1 = 0.989$, $k''_1 = 0.682$

We should remark that the large difference between the values of $c_1$ and $c_2$ in each category affects the element $k_1$ much more than $v_1$. For the general series the values of $v_1$ and $k_1$ calculated as above and formulas

$v_1 = k'v'_1 + k''v''_1$, $k_1 = k'k'_1 + k''k''_1$

where $k'$ and $k''$ are the rates of the accused for crimes in the two main categories ($k' + k'' = 1$), provide

$k' = 0.2731$, $k'' = 0.7269$ so that $v_1 = 0.814$, $k_1 = 0.651$.

When comparing these values with the system (218.1), we see that the difference is very small and it is possible to regard the mean values of the elements $v_1$ and $k_1$ for the general series as determined with sufficient approximation in accordance with the essence of the data without previously multiplying the number of categories. We also have $P' = 127$, $P'' = 215$, so that $P' + P'' = 342$ which little differs from $P = 335$ provided by the general series.

As remarked above, we can without an appreciable error or with an error of the order supposed for the uncertainty of the data assume the previous values of $v'_1$, $v_1''$, $k'_1$, $k_1''$ as $v'$, $v''$, $k'_1$, $k_1''$ in accord with the other extreme hypothesis, $v'_1 = v'_2 = v'$, $v''_1 = v''_2 = v''$. Then

$V'_1 = 0.996$, $V'''_1 = 0.998$, and $Q'' = 18$, $Q' = 19$, $Q' + Q'' = 37$

instead of $Q = 41$ determined for the general series. This means no more than one accused in a thousand is convicted although absolvable; convicted, although according to our definition, the chance of a convicting vote for them fell below 1/2.

Concerning that small number of the accused it is legitimate and comforting to believe that the majority is guilty, but it is impossible to evaluate even approximately the probability of their real guilt. On the other hand, that category of convicted although absolvable accused does not necessarily include all the accused who can be convicted although innocent. It is regrettably possible that for some innocent accused the chance of a convicting vote exceeds 1/2 and is even very close to unity. Calculation applied to judicial statistics has no means for revealing this possibility and assign a chance to it.

**220.** Concluding what we have to say on this subject, we believe it useful to add some explanations to those given above about the meaning of the letters $v_1$ and $v_2$ and about the sense of the fundamental distinction established between the convictable and absolvable accused.

For simplifying the discussion, we will at first only consider accused of the same category at whose trial all the causes of error act fortuitously and variably from one judge to another. We will also admit that with regard to those accused it is only possible to distinguish one category of citizens chosen, or possibly chosen to perform the duties of a juryman. The ratio of the number of convicting and acquitting votes will be the same whether for the same randomly chosen juryman trying successively a very large number of accused or for a very large number of jurymen asked about the same accused.

In either case that ratio is $v_1/v_2$ with $v_1$ and $v_2$ being the chances of a convicting and acquitting vote for the mentioned categories of accused and jurymen. Therefore, since we understand and ought to understand convictable accused as such for whom $v_1$ and therefore $v_2$ exceed ½, those convictable accused will certainly be convicted at least by a simple majority if the pleadings are conducted before a very large number of jurymen for each of whom the chances $v_1$ as well as $v_2$ have the same values. And it is not difficult to see that this conclusion also persists when it is not anymore permissible to admit the just mentioned condition for all the citizens from whom jurymen are chosen by chance.

Actually, $v_1$ and $v_2$ therefore denote means

$$k^{(1)}v_1^{(1)} + k^{(2)}v_1^{(2)} + \ldots, \; k^{(1)}v_2^{(1)} + k^{(2)}v_2^{(2)} + \ldots$$

where $v_1^{(1)}, v_1^{(2)}, \ldots v_2^{(1)}, v_2^{(2)}, \ldots$ are the values of $v_1, v_2, \ldots$ for each category of jurymen and $k^{(1)}, k^{(2)}, \ldots$ express for each category the rate of citizens composing those categories. However, the same mean values also express the probabilities that a juryman chosen by chance from the general list will convict or acquit an accused and when the former exceeds 1/2 (i. e., when the accused is convictable in the sense of the definition) we are sure that conviction will follow at least by simple majority if a very large number of jurymen randomly chosen from the general list can be called to the pleadings.

According to this manner of defining the magnitudes $v_1$ and $v_2$ and their analogues, the problems treated in this chapter assume a purely arithmetical sense easily understood even by people remote from mathematical analysis. Eliminated are delicate considerations following from the use of the words *truth* and *error*, when applied to judgements such as those pronounced by tribunals, for which there is no general criterion of verity. Nevertheless[8], we believe to be duty-bound to prefer a method which connects the conveniently modified theory of judgements of tribunals with the theory of chances of verity or error in judgements considered generally. That method facilitates the comparison of our analysis with that of the authors treating the same subject, and more clearly shows, as we think, the imperfection of the previous theories and can direct them towards an improvement.

**221.** The contempt for the calculus of judicial chances felt by certain lawyers is unfounded and I believe that the legislators' viewpoint about the organization of tribunals is in essence the same as that of geometers. The former are only interested in the mean and general results of the system they establish; the latter know that their formulas are only useful when applied to large numbers without being applicable to a particular case. If desiring to confirm authentically their assumptions, legislators can only study judicial statistics; and without statistics the formulas of the geometers remain sterile or at least only some general propositions rather than numerical results can be derived from them.

The former know or ought to know that judicial institutions can never prevent fatal blunders when all the manifestations of a crime point to an innocent accused; that in civil cases those institutions do not impede errors due to dominant prejudices; that their sole aim is to guarantee a judgement conforming to the opinion of a majority of impartial and enlightened contemporaries; to offer even in criminal cases a sufficient assurance that a convicting decision will be approved by a great majority; and to restrain the anomalous influence of chance on the fate of the accused.

All the events which the legislator can not perceive by his means just the same can not be submitted to calculation by a geometer, but what is understandable to one of them is assessable to the other by studying statistical documents.

## 16.2. On the Probabilities of Testimonies

**222.** Our long discussion of the probabilities of judgements allows us to restrict the present subject to a few considerations.

Unreliability of a conscientious witness is akin to that of an honest judge. The former's testimony, like the latter's vote can only be wrong due to an error of judgement. There is therefore room for applying our theory of the probability of judgements to probabilities of testimonies provided that the witnesses are not suspected of ill will.

Suppose that a large number of times the same person A is involved as a witness and that by some means we are able to distinguish quite certainly true and mistaken testimonies. Denote by $m_1$ the total number of testimonies, $n_1$ of them true. The fraction $n_1/m_1 = v_1$ expresses the chance [the probability] of truth of A's testimony. In other words, if A will testify again under similar circumstances the ratio $n_2/m_2$ of his testimonies admitted as being true will not appreciably differ from $v_1$ if only the numbers $m_2$ and $n_2$ are sufficiently large as were $m_1$ and $n_1$.

Suppose that $v_2$ is analogues to $v_1$ for witness B and assume as in § 193 that the causes influencing A's verity or error are completely independent from those influencing B. Then

[1] The probability of an agreement between A and B is

$$p = 1 - (v_1 + v_2) + 2 v_1 v_2. \qquad (222.1)$$

[2] The probability of their disagreement is

$$q = v_1 + v_2 - 2 v_1 v_2 = 1 - p.$$

[3] The probability of the testimony's verity in case [1] is

$$V = \frac{v_1 v_2}{v_1 v_2 + (1-v_1)(1-v_2)}$$

etc. If the numbers $v_1$ and $v_2$ are not known in advance, but the number $p$ is quite precisely determined by experience we will at least know that those numbers ought to satisfy equation (220.1).

**223.** Suppose that the chance of verity of a third witness C is $v_3$ and that all three are simultaneously testifying about a large number of cases. Each of them can successively find himself in opposition to the others or all of them can agree. Denote by $a, b, c, d$ the probabilities of these four combinations, then

$$\begin{aligned}a &= v_1(1 - v_2 - v_3) + v_2 v_3 \\ b &= v_2(1 - v_1 - v_3) + v_1 v_3 \\ c &= v_3(1 - v_1 - v_2) + v_1 v_2\end{aligned} \qquad (195.1b, c, d)$$

and $p = 1 - (a + b + c)$. Assume that observations provided the numbers $a, b$ and $c$, then the values of the chances $v_1, v_2$ and $v_3$ can be derived from equations (195.1b, c, d). Indeed, it is impossible to determine these values directly since a criterion for discerning correct and erroneous testimonies is lacking. However, we will not discuss the consequences of this remark; we considered them for the probabilities of judgements proper whereas no statistics of testimonies is practically possible and no means exist for deriving numbers from formulas.

If C is opposed to A and B, the probability of his error is

$$\frac{v_1 v_2 (1-v_3)}{v_1 v_2 (1-v_3) + (1-v_1)(1-v_2) v_3}.$$

It is reduced to $v_1$ if $v_2 = v_3$. The testimonies of B and C are contradictory and of equal value and therefore neutralize each other and the probability of the verity of A's testimony remains without change as though he is the only witness.

For understanding this proposition in its veritable sense we should suppose that the number of trials was very large and the number of those in which A agreed with B and disagreed with C was registered. The rate of A's true testimonies in that partial series will not appreciably differ from the value $v_1$ derived from the complete series.

If the three witnesses agree, the probability of the truth of their testimonies becomes

$$V = \frac{v_1 v_2 v_3}{v_1 v_2 v_3 + (1-v_1)(1-v_2)(1-v_3)}.$$

In general, if all $n$ witnesses [!] agree, that probability will be

$$V = \frac{v_1 v_2 v_3 \ldots v_n}{v_1 v_2 v_3 \ldots v_n + (1-v_1)(1-v_2)(1-v_3)\ldots(1-v_n)}.$$

Suppose that all fractions $v_1, v_2, v_3, \ldots v_n$, exceed 1/2, then as $n$ indefinitely increases the value of $V$ will in general approach unity. We can however suppose that $v_n$ decreases with an increasing $n$ and indefinitely tends to 1/2 according to such a law that $V$ converges to a value differing from unity.

**224.** Suppose now that we separate the witnesses depending on the essence of the attested fact. We will certainly see that the chance of verity of the testimony of A, whom we believe to be invariably honest, is not the same in all selected categories. Experience would show this, had it been possible to check verity by a criterion in a numerous enough series of testimonies. Lacking such an experience, our knowledge of the laws of human nature will sufficiently well indicate this circumstance. Fondness for the marvellous, force of prejudices, exaltation of emotions by sectarian and party views, − all that which brings into play sympathies and antipathies of the human heart influences witnesses, most often unaccountably, and leads them to illusions, bewilder them, expose them to err involuntarily.

All previous authors, when appreciating testimonies, understood the need to allow for the essence of the attested fact. However, they reasoned as though the probability of the fact in itself should be combined with the probability of the witness's verity supposed to be invariable for the same person whatever the essence of that fact. Actually, it is this latter element which varies with each category of facts. In addition, the impossibility of assigning the law of that variation from one category of facts to another evidently renders impracticable any numerical applications.

**225.** According to the theory of probability of judgements we did not explicitly allow for the possibility of judges' malfeasance and did not distinguish between interior decisions suggested to them by their enlightenment and their votes. In essence, this distinction is not necessary and we can consider malfeasance as an error of vote in an exterior decision. Even if there is no malfeasance in its proper sense the judge's or juryman's vote in criminal cases can oppose his interior judgement. Thus, a feeling of pity determined him to acquit an accused although internally he believed that person to be guilty.

We can just as well consider the chance of a witness deliberately deceiving us, perhaps being corrupted by a bribe, interlaced with other chances influencing his verity or error. On certain occasions and especially when a solemn oath is lacking, a witness can lie, just as a judge without there being any malfeasance in its proper sense because he believes, reasonably or not, to have a good reason for concealing the truth. Previous authors[9] distinguished the witness's chances of error and lie without proposing a similar difference for judges, undoubtedly because we ought to fear a witness's lie much more than a judge's formal malfeasance. However, this distinction can be admitted in theory but it is barely useful in practice since it is absolutely impossible to determine either chance quantitatively and separately.

We will not say anything more about our subject and we guard ourselves against the desire to apply the calculus of probability to facts supposed to be known to us through a chain of witnesses or by *tradition*[10]. The values of the elements included in such calculations are not at all assignable and in addition the very combinations of these elements rests on arbitrary hypotheses establishing fictitious independence between actually solidary facts whose solidarity prevents any legitimate application of the theory of chances.

**Notes**

**1.** In this connection Cournot (1838) earlier mentioned Poisson's law of large numbers. [B. B.] It is generally known that, being under Bienaymé's influence, Cournot here ignored that law. O. S.

**2.** Same remark. [B. B.]

**3.** Montesquieu noted that an inevitable harsh punishment often compelled judges to acquit the accused. [B. B.]

**4.** This seems to be one-sided and in any case unjustified.

**5.** Gelfand & Solomon (1974) applied American statistics and evaluated the parameter $v_1$ without supplementary hypotheses. They concluded that $v_1$ and $v_2$ little differed from each other. [B. B.]

**6.** Poisson (1837, § 135) did just that and concluded that $c_1 + c_2 = 0.6094$. [B. B.]

**7.** Earlier Cournot (1838) pronounced a contrary opinion. [B. B.]

**8.** Earlier Cournot (1838, p. 333) stated that he had studied Poisson's book (1837) *with all attention of which I was capable* and largely followed him. [B. B.] See also Note 5 to Chapter 15. O. S.

**9.** Bru referred to Laplace (1812/1886, pp. 455 – 458) and noted that Lacroix had thought that Laplace's pertinent formula was practically useless.

**10.** This is what Craig in 1699 and Laplace had done. [B. B.]

On Craig see Stigler (1986). Poisson (1837, §§ 39 and 40) should also be mentioned. O. S.

## Bibliography

**Cournot A. A.** (1838), Sur les applications du calcul des chances à la statistique judiciaire. *J. math. pures et appl.*, t. 3, pp. 257 – 334.

**Gelfand A. E., Solomon H.** (1974), Modelling jury verdicts in the American legal system. *J. Amer. Stat. Assoc.*, vol. 69, pp. 32 – 37.

**Laplace P. S.** (1812), *Théorie analytique des probabilités*. *Œuvr. Compl.*, t. 7. Paris, 1886.

**Poisson S.-D.** (1837), *Recherches sur la probabilité des jugements* … Paris, 2003, 2012. English translation in www.sheynin.de downloadable file 53.

**Stigler S. M.** (1986), John Craig and the probability of history. In author's *Statistics on the Table*. Cambridge (Mass.) − London, 1999, pp. 252 − 273.

# Chapter 17. On the Probability of Our Knowledge and on Judgements Based on Philosophical Probabilities.

## Summary.

### [17.1. On the Probability of Our Knowledge and on Judgements Based on Philosophical Probabilities]

**226.** All our faculties by which we acquire knowledge are or seem to be subjected to error. Senses are illusory, attention dulls, memory is capricious, and the faculty of calculating or reasoning escapes us many times in succession. Thus, we justly do not believe ourselves and regard that verity is only established after being checked and admitted by a large number of competent judges situated in various circumstances.

During each period of the history of philosophy sceptics boasted about this maxim of common sense for denying the possibility of distinguishing the true and the false. Other philosophers concluded that our knowledge, although never completely certain, can acquire a probability ever nearer certainty, still others regarded a unanimous agreement as the sole and solid foundation of our knowledge. Philosophical criticism is beyond the scope of this book, but we ought to say a few words about its fundamental issues in so far as they are connected with the theory of chances and probabilities whose principles and all its important applications we desire to indicate.

Let us admit that each faculty by which we acquire knowledge can be likened to a fallible judge or witness. A superior intelligence with an unlimited scope of understanding which penetrates for example the mysterious skill of the memory will be able to assign the chance of verity or error attached to the action of each function, to the application of any of our faculties for each individual under certain determined circumstances. We will perhaps recognize that for some people under certain circumstances the chances of error disappear. Indeed, nothing authorizes us to affirm absolutely that no intellectual operation, even the simplest of them, is free from a chance of error.

An intelligence lacking that capacity but possessing an infallible criterion can therefore experimentally determine the chances of error inherent in exercising each of our faculties had it been possible to carry out a sufficiently numerous series of trials under suitably determined psychological conditions. However, we will never acquire a posterior absolute certainty that, under certain conditions, the chance of error disappears or that the possibility of an error is exactly zero.

Even when that intelligence does not have such a criterion of verity, observations can lead it to a numerical determination of the unknown chances of error, or at the very least of the mean of the values which these chances can take when passing from one person to another, from one category (? - O.S.) to another if only we assume that the chance of verity is invariably higher than the chance of error. This restriction is necessary if we agree that human faculties are normally destined for, and result in leading us to verity. Mistaken perception or judgement is an anomaly caused by an accidental disturbance of faculties and functions. We thus return to the mathematical theory of judgments or testimonies which was the subject of the [two] preceding chapters.

**227.** However, we should not be mistaken. That theory is of little interest here even if we have the necessary knowledge for assigning in advance the numerical values of the chances which are entering [the formula] as elements, or if supposing that we are able to determine these values experimentally. Indeed, it is important to weigh in each case in particular the force of the reasoning leading us to believe, to reject or to abstain from agreeing (? - O.S.), and the mathematical theory as expounded until now most often only provides deceptive indications.

Suppose for example that it is perfectly established by experience that each of the two people, A and B, is subject to err only once in twenty numerical calculations of a well determined kind, such as a solution of a right triangle. It does not follow that since B attentively checked A's calculation and found it correct, the probability of its error is exactly $(1/20)^2 = 1/400$. Actually, by the very fact that B proposed to check an already obtained result, we may think that he will be more attentive and better guarded against the chances of error.

Even when B had not known A's result and did not wish to check it, it will be extraordinary that from all possible calculative mistakes the same one remains unnoticed by both or that B overlooks another mistake affecting the final result in the same way[1]. Therefore, if the two calculators exactly agree, the probability of the correctness of their common result as derived by applying these notions (? - O.S.) about combinations and chances, can much exceed 399/400. The calculation of that probability is a complicated problem whose solution depends on the kind of the numerical calculations which provided the common result, on the number of the retained digits etc.

If, on the contrary, the faults of the calculation depended on some mistake in the common method applied by both A and B, or some error in the tables applied by them, the probability of the same error in the coinciding result can exceed 1/400. Suppose now that the result obtained by those two calculators satisfies some simple law suggested by the theory, already verified by similar cases and expected to be confirmed once more. Everyone will agree to regard extremely unlikely or even impossible that an accidental calculative error provides exactly that, which brings the result in compliance with the theoretical law. No one will ever doubt the correctness of the obtained result and will never inquire whether those two calculators are subject to err once in 20 or a 100 times.

**228.** We considered an example of a numerical calculation, i. e., of the most mechanical intellectual operation of sorts, but a similar reasoning is evidently applicable to all the actions of the mind directed to obtaining knowledge. Nevertheless, the evaluation of the chances of error whether in advance or after the calculation seems to present difficulties the less surmountable the more complex are the operations or the involved areas of our intellectual organization are more concealed.

Even greatest geometers fall into error and propositions admitted in pure mathematics are later abandoned as wrong or inexact. Nevertheless, it would be extraordinary and therefore very improbable that so many geometers for more than 20 centuries have been mistakenly thinking that the Euclidean demonstration of the Pythagorean proposition was irreproachable. Indeed, when considering that that theorem had been proved in many ways, and that it conforms to the entire system of perfectly connected propositions, we become totally convinced in that the demonstration conforms to the laws regulating human thinking and that this theorem belongs to the rank of truths subsisting independently from the faculties revealing them to us and from the laws governing them.

Similar remarks are applicable to historical testimonies. We firmly believe in the [former] existence of that person who was called Augustus not only because many historians mention him and agree about the main circumstances of his life, but also because Augustus is not an isolated figure but renders meaning to many contemporaneous and subsequent events which would have remained groundless and unconnected when such an important link is removed from the historical chain.

If some unusual minds doubt the Pythagorean proposition or the [former] existence of Augustus, it will not at all shake our belief; we will not hesitate to decide that some of their intellectual faculties are disordered, that they overstepped the normal conditions necessary for performing their destination.

**229.** Therefore, neither the repetition of the same judgements[2] nor a unanimous or nearly unanimous agreement is the sole foundation of our belief in certain truths. It mainly rests on a perception of a rational order of linking the truths and on a conviction that the causes of error are anomalous, irregular and subjective, and can not engender such a regular and objective coordination.

This is indeed the principle of philosophical criticism[3] of our knowledge. Our senses, and in general the performance of all the faculties by which our knowledge is extended or perfected, are guided and controlled by a superior and regulating faculty called *reasoning*. And the human mind is able to rouse us to inquire about the reasonable in various things, is the faculty of perceiving the chain of causes and effects, of principles and consequences.

And an aberration in the sensibility of some individuals being in certain anomalous physiological states or even in those normally and periodically reproduced in our sleep, in spite of the objections of ancient sceptics are certainly unable to shake our belief in the testimony of the mind. The notions about the external objects when we are awake and our senses function normally, perfectly agree one with another. Impressions of various nature rendered by our different senses are well enough connected, systematized and coordinated.

Memory identifies the notions provided by our senses beginning from the obscure period of our babyhood when their training is completed in spite of the variability of painful or agreeable feelings which, during the different periods of life, accompany for each of us the perception of those same external objects. The same identification of those objects unites in each of us our faculties and clearly manifests itself in our continuous encounters with other human beings whereas no regular connection exists between our dreaming today and tomorrow or between dreams of different people.

And finally, although little do we know about the principle of sensibility and of our psychological functioning, we know enough for discerning that the perturbations of sensibility during sleep or because of other circumstances of the organic life result from suspending or obliging certain faculties, from an injury of certain organs. *Exceptio firmat regulam* [Exception proves the rule].

Sometimes our senses expose us to illusions, but they can be called normal because they are universally adopted and do not result from accidental disturbances of their functioning. Such are optical illusions causing the sky to look like a flattened vault whereas the Moon seems to be much larger at the horizon than near the zenith. Numerous explanations have been proposed of these and of many other illusions. However, even if they had remained inexplicable, the concurrence [of the indications] of our other senses and the intervention of our mind would have hastened to rectify the errors of judgement which can at first accompany them.

If one faculty apparently contradicts another one, our mind, without being embarrassed, will decide between them. It discerns which of them is preeminent and does not hesitate to perceive the phenomena in such a manner that solely suits a systematic and regulated coordination and solely satisfies the supreme laws of the mind.

**230.** It is also the mind understood as a faculty judging all the others which the philosopher asks whether the notions provided to us by the system of inferior faculties about the external objects are only true as human verities suited to our condition, to the laws of our proper nature, or, on the contrary, whether these faculties are granted to humans for attaining in a certain measure an effective knowledge of what are things externally and objectively.

We see some celestial objects traversing the terrestrial atmosphere which deflects their luminous rays and alters for us [their] relative position. Owing to the same disturbing cause called astronomical refraction the stars do not seem to rotate uniformly in perfect circumferences about the axis of the world. However, even had not the experience of the physicists instructed us in the laws of refraction, it would have been sufficient for us to remark that the anomalies of the diurnal motion of the stars change with the observer's horizon [with their zenith distances] for concluding without hesitation that these anomalies are only apparent, depend on the conditions of observation and have no objective reality.

Assume now as Fr. Bacon did the possibility of such a structure of the human eye that the relative positions of the stars corrected for the refraction are still wrong and that therefore the laws of the diurnal motion in all their imposing simplicity are only illusory. Conclude that perhaps all the edifice of the astronomical sciences resting on those laws is baseless? These consequences are contrary to the mind. Indeed, how could have happened such a faulty structure of the human eye that, without being able to disturb the external phenomena's order and regularity it introduces (? - O.S.) the lacking order, regularity and simplicity?

Just the same, we are firmly convinced that observation does not fail us at all, that, after taking into account the deflection caused by the intermediate atmosphere and some other disturbances originated by the Earth's motion, the stars are shown to us in their veritable optical positions. Had we enough place, we would have indicated that similar conclusions certainly justify even for the most meticulous mind the received objective reality of the fundamental notions of space and time. We would thus dispel the systems of the modern sceptical school which only desires to understand those motions as laws proper to the human mind, as forms of our thoughts, having no external reality.

**231.** That discussion of our intellectual faculties and their engendered ideas which seems to us as an essential object of philosophical speculations, unlike geometrical theorems, is not at all proceeding by demonstrations, and, unlike formal reasoning, does not lead to conclusions from premises. The existence of bodies, the objective reality of space and time can not be proved and the same should be said about the most certain physical laws, such as the law of gravitation. Indeed, what can prevent a misconstrued mind demanding in such cases geometrical demonstrations to attribute to chance the invariable agreement between the Newtonian hypothesis and observations of phenomena?

Therefore, independently from *apodictic* proofs, or formal demonstrations, there exist *philosophical* or *rational* certitudes resulting from judgements of the mind by appreciating various suppositions or hypotheses, admitting some of those introducing order and rational chains in the system of our knowledge and rejecting those incompatible with that rational order of the world pursued by human intelligence as avidly as possible.

Thus, certain natural and instinctive beliefs are legitimized in the eyes of the mind whereas others are rejected as prejudices or illusions of our senses. In the final analysis all our knowledge is based on that philosophical certitude since all demonstrated truths issue from primary truths accepted but indemonstrable.

**232.** The derived secondary certitude acquired by logical demonstration is fixed and absolute, does not admit any nuances or degrees, but the judgements of the mind on which we insist are the foundation of the certitude of the received verities. Under certain conditions they produce unshakeable conviction, but in many cases they seem only to lead to probabilities lowering by imperceptible nuances and they act differently on various minds.

For example, in the present state of science some physical theories are thought to be more probable than others since they apparently better satisfy the rational chain of the observed facts, since they are simpler or discover more remarkable similarities. However, these similarities and inductions do not strike all minds, even most enlightened and impartial of them, with the same force. Intellect recognizes certain probabilities which nevertheless are insufficient for establishing complete conviction. They change with the progress of science. Some contested theories become unanimously accepted, some are sooner or later rejected which proves that their probabilities include elements varying from one mind to another.

In other cases, we are condemned only to possess probabilities and such is the problem of inhabitancy of [other] planets. We are astonished by the similarity between those planets and our Earth, and disgusted to admit that, according to the design of nature, a tiny globe lost in the midst of immense spaces is the only one on whose surface wonders of organisation (? - O.S.) and life have developed. However, we can not at all expect that the progress of science will throw new light on things which nature apparently designed to leave beyond any means of observation. Relatively near us a globe whose dimensions are comparable to terrestrial seems to be situated in such physical conditions under which any life organized similarly to that peopling our Earth is impossible. As the mind becomes more astonished by similarities and dissimilarities, it more or less firmly adheres to the philosophical opinion about most worlds.

**233.** That subjective probability is variable and sometimes either excludes doubt or only glimmers, and we wish to call it *philosophical*. Indeed, it is caused by that superior faculty by which we judge the nature of things. Should it be thought to be in essence the same probability with which we have been dealing until now and connected with the notions of chances and randomness, or, as we explained many times, should it be linked with the concept of independence of combined causes? The identity is somewhat obscure and, since it complies with the rules of philosophical criticism, it suffices for us to remain on this side of possible reductions rather than to risk a confusion of really distinct principles. Let us see whether it will not be possible to further the analysis and formally isolate this distinction.

**234.** For better fixing the ideas, we consider a fictitious and very simple example. Suppose that a variable magnitude can take values expressed by numbers from 1 to 10,000 and that 4 of its observations or measures formed a geometric progression. We will be really inclined to think that that result was not fortuitous, that it could not have been caused by an operation comparable to 4 drawings by chance from an urn containing 10,000 tickets numbered from 1 to 10,000; that, on the contrary, it indicates the existence of a regular law in the variation of the measured magnitude and in the order of the succession of its measures.

Instead of a progression of the kind we call geometric, the 4 numbers obtained by observation could have offered some other arithmetical law. For example, they can form 4 terms of a progression with equal differences [arithmetic progression], or a series of square, cubic, triangular, pyramidal etc. numbers. The number of laws of this kind is unbounded and, as indicated by the theory of interpolation, we can even always find one which mathematically connects the 4 obtained numbers.

However, if the mathematical law to which we have to turn for connecting those 4 numbers, is expressed in an ever more complicated way, it becomes ever more probable that their succession is due to chance or a combination of independent causes. On the contrary, if the law appears very simple, and even if sufficiently numerous observations do not strictly satisfy it, we will not hesitate to admit its existence and to attribute the pertinent deviations to errors of observation or some disturbing causes too late indicated by the theory.

But what exactly means the simplicity of a law? How to compare and rank in this sense the infinitely varying laws capable to be imagined by the mind and presented mathematically? This problem can seem insoluble in itself, and so it is for us owing to the imperfection of our knowledge. And even if soluble, with the ranking established, it will not at all result in a numerical estimation of the probability of the existence of a law which a restricted number of observations satisfy strictly or approximately.

**235.** The planetary system offers a very remarkable example which almost exactly returns us to the abstract hypothesis of § 234 and can further elucidate our subject[4]. It was noted long ago that, after ranging the planets (except Mercury) by their distances from the Sun, the intervals between their orbits or the differences between the mean distances of two consecutive planets from the Sun almost follow a geometric progression with ratio 2. Thus, assuming the interval between the orbits of Venus and the Earth as unity, the following intervals will be expressed by numbers 2, 4, 8, …:

Venus − Earth, 1; Earth − Mars, 2;
Mars − Vesta etc − 4; Vesta etc − Jupiter, 8;
Jupiter − Saturn, 16; Saturn − Uranus − 32

This progression which the mean distances in the astronomical sense or the semimajor axes only approximately satisfy, is strictly obeyed by the *limits of the eccentricities*[5]. Indeed, assign for each planet a value of the radius vector between the perihelion and the aphelion distances so that the series satisfies the progression of those double intervals. It is apparently very difficult to attribute to chance such a simple ratio and not to see here a law of the construction of the planetary system although the theory does not indicate the causes which governed its formation.

Before the four telescopic planets [Vesta etc] were discovered and even before Uranus became known and added a new term to the series thus essentially corroborating the probability of the discussed law, eminent minds[6], surprised by the gap in the series of intervals between Mars and Jupiter, had suspected the existence of an intermediate planet. The mean distances from the Sun of the four telescopic planets little differ one from another; more than one indicator led to regard them as the debris from a destroyed planet and the previous conjectures were thus verified. It became much more difficult to regard the progression of the double intervals as a fortuitous coincidence.

Nevertheless, Mercury is an exception because the interval between its orbit and that of Venus is approximately equal to the interval between Venus and the Earth, not twice smaller as required by the presumed law. To solve this anomaly Bode, a German astronomer of the 18th century, imagined another law. Express the mean distance of the Earth from the Sun by 10, then those distances for Mercury will be approximately 4, for Venus, 4 + 3 = 7, and for planet number *i beginning with Venus*,

$$4 + 3 \cdot 2^{i-1}: \qquad (235.1)$$

Mercury, 4; Venus, 7; Earth, 10; Mars, 16;
Vesta etc., 28; Jupiter, 52; Saturn, 100; Uranus, 196

Because of the constant 4, the formula (235.1) is not as simple as a geometric progression whose terms are free from a constant. In addition, the anomaly presented by Mercury is not completely eliminated: we are unable to derive Mercury's distance from the Sun by attributing a suitable value to $i$. In addition, Mercury is an exception to the system of the 7 main planets both by the eccentricity of its orbit almost equal to those of the orbits of Juno and Pallas and by the notable distance of the pole of its orbit from the region of the sky in which the poles of the 6 of the planetary orbits are now situated (§§ 145 and 156).

We will now inquire whether it is possible to assign a numerical value to the probability of the Bode law presented in one or another form, whether it is possible to take numerically into account its confirmation by discovery of new planets or its refutation by Mercury's anomaly. All this discussion evidently leads to probabilities which natural philosophy [philosophy of science?] can not neglect. Still, by their essence they can not lead us to complete certainty and a possibility of numerically expressing them can only be an illusory wish.

**236.** When replacing purely arithmetical notions by geometric considerations the remarks of § 234 acquire a new power. Suppose that 10 points determined on a plane surface by so many observations are situated on a circumference. We will not hesitate to admit that that coincidence has nothing fortuitous but indicates a law. If these points very little deviate from a circumference, some of them in one direction and some, in another, we will attribute the deviations to errors of observations or disturbing causes of an inferior order rather than abandon the law.

We will be even more surprised and hesitated still less to attribute the result to a regular cause if the pertinent circle occupies certain remarkable positions, if, for example, its centre coincides with the centre of the figure on the plane (? - O.S.) on which all the points should be situated. Instead of a circumference the points can be situated on [the perimeter] of an ellipse, on a parabola, on an infinity of various curves susceptible of a mathematical definition. And the theory informs us that it is always possible to find a curve among those qualified as mathematical passing through all the observed points whatever is their number even when their individual situation depends on fortuitous and independent causes.

The probability that the observed points were situated under the influence of regular causes depends however on the simplicity attributed to the curve that connects them exactly or approximately. Any classification of lines in this sense is incontestably artificial whether it depends on the degree or the number of terms of their equations or of the number of the included parameters. From a certain point of view a parabola can be regarded as being a simpler curve than a circle (? - O.S.). A curve having a transcendental equation, a spiral for example, can in a sense be regarded as simpler and more proper to express a law of nature governing the production of certain phenomena than an infinite number of curves having algebraic equations [algebraic curves]. However, the feeling of simplicity of an observed curve as opposed to that of an infinite multitude of possible curves leads to a judgement of the probability (? - O.S.) which is not at all expressed in numbers[7] as resulting from an enumeration of favourable and unfavourable cases for the production of an event when cases are equally possible or at least when we have no reason to prefer one case to another.

**237.** Poisson (1837, § 42) proposed to assign a probability that a *remarkable* event was due to a special regular cause rather than to a combination of chains. [Cournot inserts a long passage in which Poisson provided examples of remarkable extractions of 30 balls of two colours from an urn and of remarkable arrangements of printed letters.]

By the usual rules admitted in the theory of posterior probability and issuing from this statement and assuming that we know how many events are *remarkable* and how many, *not remarkable,* Poisson determined the probability that the occurrence of a remarkable event is not at all due to chance. The defect of this reasoning consists in supposing, first, that we can demarcate remarkable and unremarkable events[8], and, second, that that events believed to be remarkable are supposed remarkable to the same degree and placed on the same level. Who can tell, after exhausting all possible combinations of the order of the occurrence of those 30 balls, when a combination ceases to be *remarkable*? And how about the combinations of printed letters: if a traveller recognizes a few words spoken by a savage tribe, will it be for us as remarkable as another combination offering us usual words in our language? Will we consider equally probable that neither was caused by fortuitous causes?

**238.** In addition, the number of combinations in Poisson's examples, whether remarkable or not, is restricted whereas we saw above that in most cases a judgement about a similar probability was based on the simplicity presented by an observed law as compared with an infinite number of laws regarded equally possible if only the law indicated by observations had no intrinsic raison d'être but was a result of a fortuitous combination of independent causes.

Here, the number of remarkable laws, just as that of unremarkable laws (supposing that they can be separated with laws of each category placed on the same level) is unbounded and indefinite and we can not imagine how the ratio of those two numbers converge to a finite and assignable limit when both indefinitely increase. Any possible application of the notions of mathematical probability to judgements about the discussed probability is therefore corrupted and illusory.

**239.** Geometers have nothing to do with probabilities resisting the application of calculation, but we should never decide that in the eyes of the philosopher those probabilities ought to be considered useless. As we have indicated, all the criticism of human knowledge beyond the narrow path of logical deductions is based on probabilities of that nature. Sometimes they puzzle all minds by determining or justifying the irresistible conviction called *common sense*, in other cases they are only appreciated by trained intelligence. When investigating new truths, the geometer himself is most often only guided by probabilities of that kind which allow him to feel the searched truth before becoming able to demonstrate it and to offer it in that form to all minds capable of comprehending a succession of rigorous reasoning.

### [17.2. Summary]

**240.** We will summarize in a few words the main points of the doctrine which we attempted to establish in this essay.

**240/1.** The notion of *randomness* is conveyed by a coincidence of independent causes producing a determined event. Combinations of such diverse causes equally contributing to its occurrence should be understood as its *chances*.

**240/2.** If among infinitely many chances only one is possible to produce the event, that event is *physically impossible*. The notion of physical impossibility[9] is neither a fiction of the mind, nor an idea only valuable relative to the state of our imperfect knowledge. It should be included as an essential element in the explanation of natural phenomena whose laws do not depend on that knowledge.

**240/3.** When considering a large number of trials of the same randomness the rate of the occurrences of an event becomes appreciably equal to the rate of its favourable chances, or to that, called *mathematical probability* of the event. If trials can be repeated infinitely many times, it will be physically impossible for those rates to differ from each other by a finite magnitude.

In this sense mathematical probability can be considered as measuring the *possibility* of an event or the facility of its production. Also in that sense it expresses a ratio existing beyond the mind perceiving it, a law to which the phenomena are subjected and whose existence does not depend on the extension or restriction of our knowledge about the circumstances of the production of phenomena.

**240/4.** If, in the state of imperfect knowledge we have no reason to suppose that some combination will arrive rather than another one, although in reality they can be connected with so many events possibly having unequal mathematical probabilities or possibilities. And if we understand the *probability* of an event as the ratio of the number of favourable combinations to the number of all equally ranked combinations, that probability, lacking anything better, can also be applied for fixing the conditions of a bet or of some random business deal, but it will not anymore express a ratio really and objectively existing between things.

It will acquire a purely subjective essence possibly varying from one individual to another according to the measure of their knowledge. If only we desire to avoid confusion and error whether in elucidating theory or in applying it, nothing will be more important than a thorough separation of the double meaning of the term *probability* sometimes understood in an objective and sometimes in a subjective sense.

**240/5.** Mathematical probability taken objectively is understood as measuring the possibility of events produced by coinciding independent causes. In general, when discussing natural events, physical and moral, it can only be determined by experience. If the number of trials of the same randomness increases to infinity, it will be determined exactly with a certitude comparable to that of an event whose contrary is physically impossible. When the number of the trials is only very large, the probability will only be approximate but we will still be authorized to regard as extremely unlikely a considerable difference between its real value and the value derived from observations. In other words, when equating the observed and the real values we will be very rarely considerably mistaken.

**240/6.** If the number of trials is small the usual formulas for evaluating probabilities by experience will become illusory, they will only indicate subjective probabilities proper for regulating the conditions of a bet but inapplicable to the production of natural phenomena.

**240/7.** We should not however conclude from the previous remark that for providing with sufficient exactitude and sufficient likelihood the real values of probabilities of events the number of trials always ought to be very large. Such likelihood is nevertheless not the same as an objective probability. We can not assign a chance of being in the right when pronouncing that the real value is contained within determined limits. In other words, we can not assign the rate of mistaken judgements made in similar circumstances.

**240/8.** Independently from mathematical probability understood in the two senses discussed above there are probabilities not reducible to enumeration of chances. For us, they substantiate many judgements including those most important, they are largely based on our idea about the simplicity of the laws of nature, of order and rational connections between phenomena and they can be called *philosophical probabilities*[10].

All reasonable people have an obscure feeling about these probabilities. When it becomes distinct or is applied to delicate subjects, it only belongs to cultivated intelligence or it can even constitute an attribute of a genius. It provides a foundation for a system of philosophical criticism vaguely felt in the most ancient schools which suppressed or reconciled scepticism and dogmatism (? - O.S.), but under the threat of strange corruption it [certainly not the feeling but philosophical probability] can not be included in the field of the applications of mathematical probability.

### Notes

**1.** A few lines below Cournot refuted himself.
**2.** Poisson (1837, §§ 63 – 64) was of the same opinion. [B. B.]
**3.** Cournot many times applies the word *criticism* likely meaning *discussion*.
**4.** Cournot is discussing the Titius − Bode law which remained topical at least until recently and about which contradictory opinions have been formulated, see Nieto (1972). Gauss thought that the observed regularity was only coincidental and claimed that he was the first who noted Mercury's anomaly (Sheynin 1984, pp. 153 – 155).
**5.** The perihelion and aphelion distances of a planet from the Sun (a few lines below) are the distances when it intersects the appropriate ends of its major axis, and the *limits of eccentricities* (a most unusual term) are apparently connected with these distances. Bru noted that Cournot had used the same expression previously.

**6.** Notably Kepler whom Cournot mentioned later. [B. B.]

**7.** On induction and simplicity of laws see Laplace (1814/1995, pp. 112 – 113), also Boole and Mill. [B. B.]

**8.** Astrologers call remarkable arrangements of the Sun, the Moon and the planets *aspects*, and here again the same question is essential: which arrangements are remarkable? Being an astrologer, Kepler (1601/1997, § 38, p. 97) *added* three aspects to the five recognized by the ancient astrologers.

**9.** Physical impossibility is contrary to moral certainty which was introduced much earlier, see *Note by Translator*. See also §§ 43 and 233.

**10.** In § 233 Cournot noted that philosophical probabilities are subjective.

## Bibliography

**Kepler J.** (1601 Latin), On the most certain foundation of astrology. *Proc. Amer. Phil. Soc.*, vol. 123, 1997, pp. 85 – 116.

**Laplace P. S.** (1814 French), *Philosophical Essay on Probabilities*. New York, 1995. Translated by A. I. Dale.

**Nieto M. M.** (1972), *The Titius − Bode Law*. Oxford.

**Poisson S.-D.** (1837), *Recherches sur la probabilité des jugements* … Paris, 2003. English translation www.sheynin.de   downloadable file 53.

**Sheynin O.** (1984), On the history of the statistical method in astronomy. *Arch. Hist. Ex. Sci.*, vol. 29, pp. 151 – 199.

i (? - O.S.) symbols mean here and later, without additional notice, translator's concern that the original text allows some ambiguity or unclearness.